\documentclass[11pt]{article}
 \usepackage[latin1]{inputenc}
  \usepackage[T1]{fontenc}
  \usepackage{amsmath,amssymb,color}
 
  \usepackage[english,frenchb]{babel}
  \usepackage[all]{xy}
 \newtheorem{thm}{Th\'eor\`eme}
 
 \newtheorem{cor}[thm]{Corollaire}
   
 \newtheorem{lem}[thm]{Lemme}
  
 \newtheorem{prop}[thm]{Proposition}
   
 \newtheorem{rmq}[thm]{Remarque}

\title{Stabilisation de la formule des traces tordue X: stabilisation spectrale}
\date{}
\author{Colette Moeglin et Jean-Loup Waldspurger\\ CNRS, Institut de math\'ematiques de Jussieu}
\begin{document}
\maketitle{}
\section{Introduction}
Le but de cette partie est de finir la stabilisation   de   la formule des traces tordues.

La m\'ethode est tr\`es voisine de celle de \cite{partie3}. Le principe m\^eme de cette m\'ethode due \`a Langlands et Arthur, est de mener de front des r\'eductions pour l'expression g\'eom\'etrique et l'expression spectrale de cette stabilisation,
$$
I^{\tilde{G}}(\omega,f)-\sum_{{\bf G}'}i(\tilde{G},{\bf G}')SI^{{\bf G}'}(f^{\bf G'}),\eqno(1)
$$
o\`u les notations sont expliqu\'ees dans le texte mais sont \`a peu pr\`es standard, en fait, pour \^etre correct, il faut fixer un ensemble fini de places $V$ suffisamment grand (la condition pr\'ecise est que $V$  contient l'ensemble $V_{ram}$ d\'efini en \cite{stabilisationVI} 1.1) et ne consid\'erer que les fonctions $f$ qui hors de $V$ sont les fonctions caract\'eristiques d'un espace compact hyperspecial; ce n'est qu'avec ce choix de $V$ que les distributions (vues uniquement pour les places dans $V$) sont invariantes pour le c\^ot\'e gauche et stables pour le c\^ot\'e droit. Dans ce cas seules les donn\'ees endoscopiques elliptiques non ramifi\'ees hors de $V$ interviennent dans le membre de droite et il n'y a qu'un nombre fini de telles donn\'ees ensdoscopiques. 

Tr\`es sch\'ematiquement, la r\'eduction spectrale montre que (1) est une distribution qui s'exprime avec des caract\`eres de repr\'esentations (ce serait m\^eme une somme discr\`ete de traces de repr\'esentations  si le sous-groupe  du  centre de $G$ invariant sous $\tilde{G}$ \'etait un groupe alg\'ebrique compact); la r\'eduction g\'eom\'etrique avec d'autres r\'esultats montre, elle, qu'en fait (1) n'est pas du tout une distribution discr\`ete si elle n'est pas nulle. Cette incompatibilit\'e prouve  la nullit\'e cherch\'ee.

Ce qui est int\'eressant est la traduction de (1) en une \'egalit\'e de transfert spectral qui ne fait intervenir que la partie discr\`ete du c\^ot\'e spectral de la formule des traces; ce r\'esultat fait partie de la d\'emonstration et se raffine de la fa\c{c}on suivante. En toute place archim\'edienne, on fixe un caract\`ere infinit\'esimal et on note le produit de ces caractères infinitésimaux, $\nu$. On fixe aussi $V$ un ensemble fini de places suffisamment grand (pr\'ecis\'ement $V$ doit contenir $V_{ram}$ comme ci-dessus) et pour toute place $v$ hors de $V$, on fixe un caract\`ere de l'alg\`ebre de Hecke sph\'erique de $G(F_{v})$, $c_{v}$. On note $c^V$ le produit de ces caract\`eres pour toutes les places $v$ non dans $V$. On note alors $\pi^{\tilde{G}}_{\nu}[c^V]$ la somme des $\omega$ repr\'esentations de $\tilde{G}$ intervenant dans la partie discr\`ete du c\^ot\'e spectral de la formule des traces (on ne voit \'evidemment que la trace tordue de ces repr\'esentations), qui ont $\nu$ comme caract\`ere infinit\'esimal et $c^V$ comme action de l'alg\`ebre de Hecke sph\'erique hors de $V$. Alors on consid\`ere l'ensemble des fonctions $f$ qui en une place $v$ finie ont une composante cuspidale nulle  (au sens de la d\'ecompostion de Paley-Wiener) et pour ces fonctions on  montre le transfert:
$$
tr\, \pi^{\tilde{G}}_{\nu}[c^V](f)=\sum_{{\bf G}'}i(\tilde{G},{\bf G}')tr \pi_{\nu,st}^{G'}[c^V](f^{\bf G'}), \eqno(2)
$$
o\`u  ${\bf G}'$ parcourt l'ensemble des donn\'ees endoscopiques elliptiques et non ramifi\'ees hors de $V$ de $\tilde{G},\omega$ et o\`u $$\pi_{\nu,st}^{G'}[c^V]=\sum_{\nu',c^{V,G'}}\pi_{\nu',st}^{G'}[c^{V,G'}]$$
la somme portant sur les caract\`eres infinit\'esimaux $\nu'$ de $G'$ se transf\'erant en $\nu$ par la fonctorialit\'e entre alg\`ebre de Lie et sur les caract\`eres des alg\`ebres de Hecke hors de $V$ pour $G'$ se tranf\'erant en $c^V$ pour la fonctorialit\'e non ramifi\'ee de Langlands (il y a une inclusion de $L$-groupe qui fixe cette fonctorialit\'e et les facteurs de transfert de fa\c{c}on compatible) et  $\pi_{\nu',st}^{G'}[c^{V,G'}]$ est l'analogue stable de $\pi_{\nu}^{\tilde{G}}[c^V]$. En d'autres termes on montre que le terme de gauche de (2) moins le transfert du terme de droite est une somme discr\`ete de $\omega$-repr\'esentations elliptiques en toute place de $V$. 
  Comme on a toujours le droit d'ajouter \`a $V$ un ensemble fini de places $V_{0}$ et de n'appliquer la formule qu'aux fonctions non ramifi\'ees en les places $V_{0}$, on conclut facilement si   les fonctions non ramifi\'ees en une place  de $V_{0}$ ont  leur composante elliptique nulle (dans la r\'ealisation de Paley-Wiener de cet espace); ceci est exactement \'equivalent au fait qu'il n'y ait pas de donn\'ee endoscopique elliptique qui soit un tore non ramifi\'e et n'est donc pas toujours vrai.  On se tire de cet ennui en faisant agir le groupe adjoint de $G$ ou plut\^ot dans le cas tordu le groupe $(G/Z(G)^\theta)({\mathbb A}_{F})$. Ce groupe agit du c\^ot\'e gauche via son action sur les fonctions et on montre que cette action se d\'ecompose suivant un nombre fini de caract\`eres automorphes non ramifi\'e hors de $V$ (du moins si $V$ est suffisamment grand). Il agit aussi du c\^ot\'e droit en fait sur les facteurs de transfert; en \cite{stabilisationI} 2.7, il est montr\'e que chaque donn\'ee endoscopique elliptique de $\tilde{G}, \omega$ donne un caract\`ere de $(G/Z(G)^{\theta})({\mathbb A}_{F})$ dont la restriction \`a l'image de $G({\mathbb A}_{F})$ est le caract\`ere $\omega$; on note $\omega_{\sharp}^{{\bf G}'}$ ce caract\`ere. Et l'\'egalit\'e (2) se raffine en une \'egalit\'e, pour tout caract\`ere automorphe $\chi_{\sharp}$ de $(G/Z(G)^\theta)({\mathbb A}_{F})$:
$$
tr\, \pi^{\tilde{G}}_{\nu, \chi_{\sharp}}[c^V](f)=\sum_{{\bf G}';\omega_{\sharp}^{{\bf G}'}=\chi_{\sharp} }i(\tilde{G},{\bf G}')tr \pi_{\nu,st}^{G'}[c^V](f^{\bf G'}).\eqno(3)
$$Et c'est cette \'egalit\'e plus fine que l'on d\'emontre; on v\'erifie qu'elle est vraie pour toute fonction $f=f_{V}1_{\tilde{K}^V}$ si l'une des composantes $f_{v}$ a sa composante elliptique nulle; c'est le raffinement de l'\'egalit\'e (2). 
Mais pour $\chi_{\sharp}$ fix\'e, il existe une place $v$ dans $V$ o\`u ce caract\`ere est trivial. Le c\^ot\'e gauche appliqu\'e \`a $f$ valant $1_{\tilde{K}_{v}}$ prend la m\^eme valeur qu'en la fonction o\`u on remplace simplement la composante $1_{\tilde{K}_{v}}$ par la fonction
$$
\gamma\in \tilde{G}(F_{v}) \mapsto \sum_{g_{\sharp}}1_{\tilde{K}_{v}}(ad(g_{\sharp})\gamma) \chi^{-1}_{\sharp}(g_{\sharp}),
$$
o\`u $g_{\sharp}$ parcourt un ensemble de repr\'esentants de $(G/Z(G)^\theta)(F_{v})/G(F_{v})$ (l'image dans $I(\tilde{G}(F_{v}))$ de cette fonction ne d\'epend pas de l'ensemble de repr\'esentants); le c\^ot\'e droit ne change pas quand on remplace $f_{v}$ par cette nouvelle fonction. On gagne car   la nouvelle fonction a sa composante elliptique nulle. en la place $v$. On conna\^{\i}t donc l'\'egalit\'e pour cette fonction ce qui permet de conclure.

Les r\'eductions g\'eom\'etriques sont les plus difficiles et ont \'et\'e men\'ees dans les articles pr\'ec\'edents. Elles ram\`enent le probl\`eme \`a deux assertions: d'une part la stabilisation des int\'egrales pond\'er\'ees locales invariantes pour des \'el\'ements semi-simples r\'eguliers et d'autre part une identification des coefficients pour les int\'egrales orbitales ordinaires mais correspondant aux \'el\'ements exceptionnels au sens de \cite{stabilisationIII} 6.3 (c'est \cite{stabilisationVII} 3.3, 3.4, 3.5 qui donne les r\'eductions). 

La r\'eduction spectrale est effectu\'ee ici et est bien plus simple (cf \ref{redstabspecglob}); on montre que par r\'ecurrence on sait stabiliser la partie continue de cette forme spectrale de la formule des traces et donc la formule (1) \'ecrite s'\'ecrit en termes de caract\`eres de repr\'esentations. 

Admettons pour le moment
le point cl\'e, c'est-\`a-dire la preuve de la stabilisation locale des int\'egrales pond\'er\'ees invariantes. Alors les r\'eductions g\'eom\'etriques montrent que (1) est nulle si en une place, $v$, que l'on peut choisir arbitrairement, la fonction $f$ est nulle pr\`es des \'el\'ements exceptionnels de $\tilde{G}(F_{v})$. Quand on fixe $f$ aux autres places et que l'on fait varier $f_{v}$ disons parmi les fonctions non ramifi\'ees, alors (1) d\'efinit une distribution en $f_{v}$ qui est une combinaison lin\'eaire des int\'egrales orbitales en les composantes en $v$ des \'el\'ements exceptionnels (il n'y en a qu'un nombre fini); or les int\'egrales orbitales sont des transform\'ees de Fourier de caract\`eres de repr\'esentations elliptiques de sous-espaces de Levi de $\tilde{G}$ (\cite{arthurcrelle} g\'en\'eralis\'e en \cite{appendice}) et on peut supposer que ces Levi sont propres car on a suppos\'e $f_{v}$ non ramifi\'e. En revenant \`a l'\'ecriture spectrale de (1) on montre une incompatibilit\'e entre ces deux expressions qui forcent leur nullit\'e. Quand on fixe $f$, on peut toujours trouver $v$ tel que $f_{v}$ soit non ramifi\'e et c'est ce qui fournit suffisamment de fonctions pour lesquelles on sait stabiliser la formule des traces et permet de conclure comme expliqu\'e plus haut.

Il faut donc montrer la stabilisation locale g\'eom\'etrique; la premi\`ere partie de l'article  relie cette stabilisation g\'eom\'etrique  \`a la stabilisation de la formule des traces locales; la stabilisation g\'eom\'etrique entra\^{\i}ne la stabilisation de la formule des traces locales, la r\'eciproque est loin d'\^etre claire mais c'est bien cela que l'on d\'emontre (cf. \ref{redlocgeo} et \ref{stabgeohyp}). La stabilisation de la formule des traces locales a elle aussi une \'ecriture g\'eom\'etrique et une \'ecriture spectrale. L'\'ecriture spectrale est simple, c'est essentiellement une combinaison lin\'eaire de traces de repr\'esentations discr\`etes. L\`a aussi il faut montrer que l'\'ecriture g\'eom\'etrique s'apparente \`a une distribution continue au sens que c'est une int\'egrale pour des espaces de Levi propre de repr\'esentations discr\`etes de ces Levi. On aura encore une incompatibilit\'e entre les deux formes de la distribution qui assure sa nullit\'e. Malheureusement la d\'emonstration est compliqu\'ee par le fait que  les hypoth\`eses cl\'es de \ref{stabgeohyp} sont d\'emontr\'ees par voie globale en \ref{preuvegeoloc} et n\'ecessitent donc elles aussi de jouer avec la r\'eduction du c\^ot\'e spectral. De plus dans le cas tordu, il est encore plus difficile de globaliser une situation locale que dans le cas non tordu; on utilise pour cela les travaux de Kottwitz et Rogawski (\cite{kottwitzrogawski}) qu'il faut compl\'eter (cf \ref{globalisation}).

\tableofcontents
\section{Notations g\'en\'erales\label{notations}}
 Lle corps de base, toujours not\'e $F$ sera parfois local (p-adique ou archim\'edien) et parfois un corps de nombres. Si $F={\mathbb R}$ ou si $F$ est un corps de nombres, on doit travailler avec des $K$-espaces, cf. \cite{stabilisationI} 1.11 et \cite{stabilisationVI} 1.16. Pour notre propos, cela ne change rien. Aussi, pour simplifier, on n\'egligera dans la notation ce passage aux $K$-espaces que l'on notera comme des espaces -appel\'es aussi bitorseurs- connexes. On fixe $\tilde{G}$ un bitorseur sous un groupe alg\'ebrique $G$; 
  On reprend la notion $\tilde{M}$ pour les espaces de Levi de $\tilde{G}$ et on note $\mathcal{L}(\tilde{M})$ les espaces de Levi de $\tilde{G}$ qui contiennent $\tilde{M}$; on renvoie \`a \cite{stabilisationI} pour ces notations/d\'efinitions. Pour $\tilde{M}$ un espace de Levi de $\tilde{G}$, on note $M$ le sous-groupe de Levi de $G$ sous-jacent et on pose: $$W(\tilde{M}):= Norm_{G(F)}(\tilde{M})/M(F).$$ Tout objet invariant attach\'e \`a $\tilde{M}$ sera invariant par ce groupe donc dans les formules faisant intervenir  une somme sur les espaces de Levi $\tilde{M}$ pris \`a conjugaison pr\`es de tels objets, cette somme sera tr\`es naturellement quotient\'ee par $\vert W(\tilde{M})\vert^{-1}$. Mais, il est en fait plus simple de sommer sur les Levi semi-standard, ce qui a l'avantage de donner des formules qui fonctionnent aussi dans certains cas o\`u les objets attach\'es \`a $\tilde{M}$ ne sont pas invariants: c'est-\`a-dire on fixe un espace de Levi $\tilde{M}_{0}$ minimal de $\tilde{G}$ et un Levi semi-standard est un Levi qui contient $\tilde{M}_{0}$. Pour un tel Levi,  on note $\tilde{W}_{M}:=Norm_{M(F)}(\tilde{M}_{0})/M_{0}(F)$ et $\tilde{W}_{G}:= Norm_{G(F)}(\tilde{M}_{0})/M_{0}(F)$, c'est le groupe que l'on avait auparavant not\'e $W(\tilde{M}_{0})$. Dans les formules qui font intervenir des sommes de Levi semi-standard, le coefficient qui revient quasiment en permanence est:
$$
w(\tilde{M}):=\vert \tilde{W}_{M}\vert/\vert \tilde{W}_{G}\vert.
$$
Pour tout espace de Levi $\tilde M$, on note $A_M$ le tore d\'eploy\'e maximal du centre de $M$ et $A_{\tilde{M}}$ le tore d\`eploy\'e maximal de $A_M$ inclus  dans le centralisateur de $\tilde{M}$.
Pour certaines formules, on a une somme altern\'ee sur les espaces de Levi semi-standard, le signe est alors $(-1)^{a_{\tilde{M},\tilde{G}}}$ o\`u
$$
a_{\tilde{M},\tilde{G}}:= rang\, ({{A}}_{\tilde{M}}/{{A}}_{\tilde{G}}
).
$$
\section{Stabilisation de la formule des traces locales tordues\label{stablocgeo}}
Dans ce paragraphe, le corps $F$ est local.   
On   va montrer que la stabilisation de la formule des traces locales tordues est \'equivalente \`a la stabilisation des int\'egrales orbitales pond\'er\'ees tordues invariantes pour les \'el\'ements semi-simples r\'eguliers, modulo des hypoth\`eses de r\'ecurrence tout \`a fait naturelles, pr\'ecis\'ement sur cette stabilisation des int\'egrales orbitales pour des groupes ''plus petits''.

\subsection{Le c\^ot\'e g\'eom\'etrique de la formule des traces locales}
\subsubsection{Rappel du c\^ot\'e g\'eom\'etrique de la formule des traces locales et de sa variante endoscopique\label{rappelgeolocal}}

En \cite{ftlt} 6.6, la formule des traces locale invariante est d\'efinie. Le c\^ot\'e g\'eom\'etrique est, pour toute paire de fonctions $f_{1},f_{2}\in I(\tilde{G})$:
$$
I_{geo}^{\tilde{G}}(f_{1},f_{2})=\sum_{\tilde{M}}(-1)^{a_{\tilde{M},\tilde{G}}}w(\tilde{M})I_{\tilde{M}}^{\tilde{G}}(\omega,f_{1},f_{2}),
$$
o\`u $w(\tilde{M})$ et $a_{\tilde{M},\tilde{G}}$ sont d\'efinis dans le paragraphe \ref{notations} et o\`u 
$$
I_{\tilde{M}}^{\tilde{G}}(\omega,f_{1},f_{2})=\int_{\tilde{M}_{ell}/\sim}i'(\gamma)I_{\tilde{M}}^{\tilde{G}}(\gamma,\omega,f_{1},f_{2})\, d\gamma,
$$
o\`u on renvoie \`a loc.cite pour la description de la mesure sur $\tilde{M}_{ell}/\sim$, o\`u $$i'(\gamma)= mes(Cent_{M}^0(\gamma,F)/A_{\tilde{M}}(F)) \vert Cent_{M}(\gamma,F)/Cent^0_{M}(\gamma,F)\vert^{-1},$$
et $
I^{\tilde{G}}_{\tilde{M}}(\gamma,\omega,f_{1},f_{2})=$
$$\sum_{\tilde{L}_{i}\in {\mathcal L}(\tilde{M}); i=1,2}d_{\tilde{M}}^{\tilde{G}}(\tilde{L}_{1},\tilde{L}_{2})
\overline{I^{\tilde{L}_{1}}_{\tilde{M}}(\gamma,\omega,f_{1,\tilde{L}_{1}})}I^{\tilde{L}_{2}}_{\tilde{M}}(\gamma,\omega,f_{2,\tilde{L}_{2}}), \eqno(S)_{1}
$$avec $d_{\tilde{M}}^{\tilde{G}}(\tilde{L}_{1},\tilde{L}_{2})$ vaut z\'ero si l'application naturelle ${\mathcal{A}}_{\tilde{M}}^{\tilde{L}_{1}}\oplus {\mathcal{A}}^{\tilde{L}_{2}}_{\tilde{M}}$ dans ${\mathcal{A}}_{\tilde{M}}^{\tilde{G}}$ n'est pas bijective et vaut le jacobien de cette application sinon, ce qui refl\`ete le rapport des mesures sur ces espaces indispensables pour construire des int\'egrales orbitales pond\'er\'ees.

\

Pour le c\^ot\'e endoscopique,
on pr\'ecise ici quelques notations: on est dans une situation locale. Pour ${\bf G}'$ une donn\'ee endoscopique elliptique du bitorseur $\tilde{G}$ et du caract\`ere $\omega$, on pose $i(\tilde{G},{\bf G}'):=$
$$\vert det_{{\mathcal{A}}_{G}/{\mathcal{A}}_{\tilde{G}}} (1-{\theta})\vert^{-1}\vert \pi_{0} (Aut({\bf {G}'})/\hat{G}')\vert^{-1}
$$
$$\vert \pi_{0}(Z(\hat{G})^{\Gamma_{F}})\vert\vert \pi_{0}(Z(\hat{G}')^{\Gamma_{F}})\vert^{-1} \vert \pi_{0}(Z(\hat{G})^{\Gamma_{F},0}\cap Z(\hat{G}'))\vert.
$$C'est la notation $c(\tilde{G},{\bf G}')$ de \cite{stabilisationI} 4.17 (ce qui suit la formule (3)).

Il faut aussi d\'efinir la variante stable de la formule des traces locales pour un tel ${\bf G}'$. En \cite{stabilisationI} 4.17 ( formule (3)), cette formule stable est \'ecrite pour la partie elliptique de ${\bf G}'$. Avec un choix de mesure compl\`etement explicit\'e, c'est l'int\'egrale sur les points elliptiques de ${\bf G}'$, des int\'egrales orbitales stables. Pour ${\bf M}'$ un Levi de ${\bf G}'$, on g\'en\'eralise la d\'efinition de fa\c{c}on imm\'ediate en rempla\c{c}ant l'int\'egrale sur les \'el\'ements elliptiques de ${\bf G}'$ par l'int\'egrale sur les \'el\'ements elliptiques de ${\bf M}'$ et en rempla\c{c}ant les int\'egrales orbitales stables par les int\'egrales orbitales pond\'er\'ees stables, c'est-\`a-dire
$$SI_{{\bf M}',geo}^{{\bf G}'}(f'_{1},f'_{2})=$$
$$ \int_{\tilde{M'}(F)_{ell}/st-conj}k_{M'}(\delta)^{-1}mes(A_{{\bf M'}}(F)\backslash M'_{\delta}(F))SI^{{\bf G}'}_{{\bf M}'}(\delta,f'_{1},f'_{2})
$$o\`u $k_{M'}(\delta)$ est le nombre de classes de conjugaison par $M'(F)$ contenues dans la classe de conjugaison stable de $\delta$ (la d\'efinition doit \^etre adapt\'ee dans le cas $F={\mathbb R}$, cf. loc. cit.) et o\`u $SI^{\bf{G}'}_{\bf{M}'}(\delta,f'_{1},f'_{2})$ se calcule par la formule de scindage
$$
SI^{\bf{G}'}_{\bf{M}'}(\delta,f'_{1},f'_{2})=$$
$$\sum_{{\bf{L'}}_{i}\in {\mathcal L}({\bf M}'); i=1,2}e_{{\bf M}'}^{{\bf G'}}({\bf L}'_{1},{\bf L}'_{2})
\overline{SI^{{\bf G}'}_{{\bf M}'}(\delta,f'_{1,{\bf L}'_{1}})}SI^{{\bf G}'}_{{\bf M}'}(\delta,f'_{2,{\bf L}'_{2}})\eqno(S)_{2}
$$
o\`u les int\'egrales orbitales pond\'er\'ees stables sont celles de \cite{stabilisationII} 1.10 (8) (c'est la d\'efinition standard due \`a Arthur g\'en\'eralis\'ee au cas tordu) et o\`u $e_{{\bf M'}}^{{\bf G}'}({\bf L}'_{1},{\bf L}'_{2})$ est d\'efini en \cite{stabilisationII} 1.14, ce sont les constantes ''universelles'' qui interviennent dans des formules de scindage pour des distributions stables et que l'on rappelle m\^eme si on va utiliser (S)$_{3}$ ci-dessous au lieu de (S)$_{2}$; 
$$e_{{\bf M'}}^{{\bf G}'}({\bf L}'_{1},{\bf L}'_{2})=d_{{\bf M'}}^{{\bf G}'}({\bf L}'_{1},{\bf L}'_{2}) \vert (Z(\hat{L}'_{1})^{\Gamma}\cap Z(\hat{L}'_{2})^{\Gamma})/Z(\hat{G}')^{\Gamma}\vert^{-1}.$$ 
On a aussi $d_{\tilde{M}}^{\tilde{G}}(L_{1},L_{2})=d_{M'}^{G'}(L'_{1},L'_{2})$ si  ${\bf M}'$ est une donn\'ee endoscopique elliptique de $\tilde{M}$ et  si ${\bf L}'_{i}$ pour $i=1,2$ est une donn\'ee endoscopique elliptique de $\tilde{L}_{i}$.

Pour toute paire de fonctions $f_{i}$ pour $i=1,2$ sur $\tilde{G}$, et pour ${\bf G}',{\bf M}'$ comme ci-dessus, on pose:
$$I_{{\bf M'},geo}^{\tilde{G},{\mathcal{E} }}(f_{1},f_{2})=0,\eqno(1)$$
si ${\bf M}'$ n'est pas une donn\'ee endoscopique elliptique d'un espace de Levi de $\tilde{G}$ et $$
I_{{\bf M'},geo}^{\tilde{G},{\mathcal{E} }}(f_{1},f_{2})=\sum_{{\bf G}'} i_{\tilde{M}}(\tilde{G},{\bf G}')SI_{{\bf M}',geo}^{{\bf G}'}(f_{1}^{\bf G'},f_{2}^{{\bf G}'}),\eqno(2)
$$
o\`u la somme porte sur les donn\'ees endoscopiques elliptiques de $\tilde{G}$ contenant ${\bf M}'$ et o\`u $i'_{\tilde{M}}(\tilde{G},{\bf G}')=$
$$
j(\tilde{G})^{-1}j(\tilde{M}) \vert Z(\hat{G}')^{\Gamma_{F}}/Z(\hat{G})^{\Gamma_{F}}\cap Z(\hat{G}')^{\Gamma_{F}}\vert^{-1} \vert Z(\hat{M}')^{\Gamma_{F}}/Z(\hat{M})^{\Gamma_{F}}\cap Z(\hat{M}')^{\Gamma_{F}}\vert,
$$
avec $j(\tilde{G})=\vert det_{{\mathcal{A}}_{G}/{\mathcal{A}}_{\tilde{G}}} (1-{\theta})\vert$ et $j(\tilde{M})$ est son analogue pour $\tilde{M}$.
On remarque que (2) ne d\'epend que de $\hat{M}$ et pas de l'espace $\tilde{M}$ et est donc d\'efini m\^eme si ${\bf M}'$ n'est pas relevant pour un espace de Levi de $\tilde{G}$ (cf. \cite{stabilisationVI} 6.6). Mais   une preuve analogue \`a celle de la proposition de \cite{stabilisationVI} 6.6 montre  que (2) est alors nul, ce qui est compatible avec (1).

On donne tout de suite la formule de scindage suivante, pour ${\bf M}'$ un donn\'ee endoscopique elliptique de $\tilde{M}$ et pour $\delta$ une classe de conjugaison stable d'\'el\'ements elliptiques de $M'$ 
$$I_{\tilde{M}}^{\tilde{G},{\mathcal E}}(\delta,f_{1},f_{2})=$$
$$
\sum_{\tilde{L}_{i}\in {\mathcal L}(\tilde{M}); i=1,2}d_{\tilde{M}}^{\tilde{G}}(\tilde{L}_{1},\tilde{L}_{2})
\overline{I^{\tilde{L}_{1},{\mathcal E}}_{\tilde{M}}(\delta,f_{1,\tilde{L}_{1}})}I^{\tilde{G},{\mathcal L}_{2},\mathcal E}_{\tilde{M}}(\delta,f_{2,\tilde{L}_{2}}). \eqno(S)_{3}
$$
Ceci est similaire \`a \cite{stabilisationVI} 4.5 proposition (i), o\`u $V$ n'a que deux places.
\subsubsection{S\'eparation suivant les espaces de Levi\label{separationsuivantlevi}}
\begin{prop} Pour toute paire de fonctions $f_{i}$ pour $i=1,2$ sur $\tilde{G}$, on a l'\'egalit\'e
$$I^{\tilde{G}}_{geo}(\omega,f_{1},f_{2})-\sum_{{\bf G}'}i(\tilde{G},{\bf G}') SI^{{\bf G}'}_{geo}(f_{1}^{{\bf G}'},f_{2}^{{\bf G}'})=$$
$$\sum_{\tilde{M}} w(\tilde{M})(-1)^{a_{\tilde{M},\tilde{G}}}
\biggl(I_{\tilde{M},geo}^{\tilde{G}}(\omega,f_{1},f_{2})-\sum_{{\bf M}'}i(\tilde{M},{\bf M}')I_{{\bf M'},geo}^{\tilde{G},{\mathcal{E} }}(f_{1},f_{2})\biggr)\eqno(3),
$$
o\`u $\tilde{M}$ parcourt l'ensemble des espaces de Levi de $\tilde{G}$ semi-standard et o\`u ${\bf M}'$ parcourt l'ensemble des donn\'ees endoscopiques elliptiques de $\tilde{M}$ relevantes.
\end{prop}
Il y a deux \'etapes dans cette proposition; la premi\`ere consiste \`a v\'erifier que les sous-groupes de Levi des donn\'ees endoscopiques elliptiques de $\tilde{G}$ qui ne sont pas relevants pour un sous-espace de Levi de $\tilde{G}$ ne contribuent pas, c'est ce que l'on a expliqu\'e avant l'\'enonc\'e. Et la deuxi\`eme partie est purement combinatoire: on utilise \cite{stabilisationVI} 6.3 pour \'echanger la somme sur les donn\'ees endoscopiques elliptiques pour $\tilde{G}$ et celles sur d'abord les espaces de Levi $\tilde{M}$ puis les donn\'ees endoscopiques elliptiques de $\tilde{M}$; les constantes $i(\tilde{G},{\bf G}')$ n'y \'etaient pas les m\^emes mais le quotient $i(\tilde{G},{\bf G}')/i(\tilde{M},{\bf M}')$ est bien le m\^eme (les $Ker^1$ qui sont les objets globaux ne jouent pas de r\^ole dans cette combinatoire).
\subsection{Stabilisation du c\^ot\'e g\'eom\'etrique de la formule des traces locales et stabilisation des int\'egrales orbitales pond\'er\'ees}
\subsubsection{}
On note $d(\theta)$ le d\'eterminant d\'efini pr\'ecis\'ement dans ce qui pr\'ec\`ede  le th\'eor\`eme de \cite{stabilisationI} 2.4 et qui si $\theta$ stabilise un \'epinglage est la valeur absolue du  d\'eterminant de $1-\theta$ dans $\mathfrak{t}/\mathfrak{t}^{\theta}$ (o\`u $T$ est le tore de l'\'epinglage). Cette constante va redispara\^{\i}tre aussi vite qu'elle est intervenue.

Soit $\tilde{M}$ un espace de Levi de $\tilde{G}$ et $\gamma\in \tilde{M}$; on a d\'ej\`a d\'efini $i'(\gamma)$.
On note ${\mathcal X}(\gamma)$ l'ensemble des classes de conjugaison sous $M(F)$ \`a l'int\'erieur de la classe de conjugaison stable de $\gamma$ et on voit cet ensemble comme un ensemble fini de cardinal $k(\gamma)$. On note ${\mathcal X}^{\mathcal E}(\gamma)$ l'ensemble des couples $({\bf M}',\delta')$ o\`u ${\bf M}'$ est une donn\'ee endoscopique elliptique relevante de $\tilde{M}$ et o\`u $\delta'$ est une classe de conjugaison stable de $M'$ qui se transf\`ere en la classe de conjugaison stable de $\gamma$; on prend ces couples \`a conjugaison pr\`es, c'est-\`a-dire que le groupe des automorphismes stabilisant $\bf{M}'$ agit sur la classe de $\delta'$ et c'est le quotient $Out({\bf M}')$ qui agit et il agit librement. Pour $\gamma$ un \'el\'ement elliptique de $\tilde{M}$, on pose:
$$
x(\gamma):= \vert Cent_{M}(\gamma,F)/Cent^0_{M}(\gamma,F)\vert k(\gamma)^{-1}.
$$
\begin{lem}Le terme (3) dans la proposition de \ref{separationsuivantlevi} est la somme sur les espaces de Levi $\tilde{M}$ semi-standard du produit de $w(\tilde{M})$ par l'int\'egrale sur les classes de conjugaison d'\'el\'ements elliptiques de $\tilde{M}$, dont un repr\'esentant  est not\'e $\gamma$ de la fonction
$$
i'(\gamma)\biggl(
I_{\tilde{M}}^{\tilde{G}}(\gamma,f_{1},f_{2})-d(\theta)^{-1}x(\gamma)^2\sum_{({\bf M}',\delta')\in {\mathcal X}^{\mathcal E}(\gamma)}I_{M'}^{\tilde{G},{\mathcal E}}(\delta',f_{1},f_{2})\biggr).
$$
\end{lem}
Le deuxi\`eme membre de l'\'egalit\'e (3) de la section \ref{separationsuivantlevi} est une int\'egrale sur les classes de conjugaison des  \'el\'ements elliptiques dans $\tilde{M}$ avec une mesure (et des coefficients) \'ecrits en \cite{stabilisationI} 4.17, o\`u il faut g\'en\'eraliser de $\tilde{G}$ \`a tous ses espaces de Levi; comme expliqu\'e en loc. cite les choix de mesures sont coh\'erents avec ceux de \cite{ftlt}. Le terme index\'e par $\tilde{M}$ est donc, pour $I_{\tilde{M},geo}^{\tilde{G}}(\omega,f_{1},f_{2})$ une int\'egrale sur les classes de conjugaison elliptiques de $\tilde{M}$, repr\'esent\'ees par un \'el\'ement not\'e ${\gamma}$ de l'int\'egrale orbitale pond\'er\'ee $I_{\tilde{M}}^{\tilde{G}}({\gamma},\omega,f_{1},f_{2})$ affect\'ee du coefficient $i'(\gamma)$.

Le deuxi\`eme terme intervenant est une somme sur ${\bf M}'$, les groupes endoscopiques elliptiques relevants de $\tilde{M}$, pris \`a isomorphisme pr\`es, de l'int\'egrale sur les classes de conjugaison stable d'\'el\'ements elliptiques de $M'$ avec un repr\'esentant not\'e $\delta'$,   chaque terme \'etant  affect\'e du produit du coefficient g\'en\'eral $i(\tilde{M},{\bf M}')$ avec $k(\delta')^{-1}mes(A_{{M}'}(F)\backslash M'_{\delta'}(F))$ o\`u $k(\delta')$ est le nombre de classes de conjugaison sous $M'(F)$ \`a l'int\'erieur de la classe de conjugaison stable de $\delta'$ (dans le cas o\`u le corps de base est archim\'edien, c'est la classe dans le $K$-groupe qui intervient \'evidemment). On reprend les calculs de \cite{stabilisationI} preuve de la proposition 4.16; on fixe $\gamma$ un \'el\'ement elliptique de $\tilde{M}$, ou plut\^ot sa classe de conjugaison sous $M(F)$.  On peut r\'ecrire le deuxi\`eme terme en faisant une somme sur les classes de conjugaison stable d'\'el\'ements elliptiques de $\tilde{M}$ comme 
$$
\sum_{({\bf M}',\delta)\in {\mathcal X}^{\mathcal E}(\gamma)}i(\tilde{M},{\bf M}')\vert Out({\bf M}')\vert k(\delta)^{-1}mes(A_{M'}(F)\backslash M'_{\delta}(F))I_{M'}^{\tilde{G},{\mathcal E}}(\delta,f_{1},f_{2}),\eqno(1)
$$
ce qui ressemble \`a \cite{stabilisationI} 4.17 (4). On peut simplifier les coefficients car il est d\'emontr\'e en \cite{stabilisationI} 4.17 (5), pour $({\bf M}',\delta) \in {\mathcal X}^{\mathcal E}(\gamma)$, l'\'egalit\'e:
$
i(\tilde{M},{\bf M}')=$
$$d(\theta)^{-1}\vert Out({\bf M}')\vert^{-1}k(\delta) mes (A_{M'}(F)\backslash M'_{\delta}(F))^{-1} k(\gamma)^{-1} $$
$$mes (A_{\tilde{M}}(F)\backslash Cent_{M}({\gamma},F)).
$$
 Ainsi (1) devient
$$
\sum_{({\bf M}',\delta)\in {\mathcal X}^{\mathcal E}(\gamma)} d(\theta)^{-1}k(\gamma) x(\gamma)^{-1}i'(\gamma)I_{M'}^{\tilde{G},{\mathcal E}}(\delta,f_{1},f_{2}).
$$
Cela montre l'\'enonc\'e puisque dans l'\'enonc\'e on somme sur les classes de conjugaison et non comme ci-dessus sur les classes de conjugaison stable (d'o\`u la disparition du $k(\gamma)$)
\subsubsection{}
On peut encore simplifier l'\'enonc\'e du lemme pr\'ec\'edent; on fixe $\tilde{M}$ et $\gamma$ une classe de conjugaison d'\'el\'ement elliptique de $\tilde{M}$; on reprend les notations pr\'ec\'edentes en particulier $x(\gamma)$. On pose
$$
I_{\tilde{M}}^{\tilde{G},{\mathcal E}}(\gamma,f):= d(\theta)^{-1/2}x(\gamma)\sum_{({\bf M}',\delta)\in {\mathcal X}^{\mathcal E}(\gamma)} \Delta(\delta,\gamma)^{-1} I_{M'}^{\tilde{G},{\mathcal E}}(\delta,f).
$$
On g\'en\'eralise la d\'efinition \`a $\tilde{G}$ remplac\'e par un de ses sous-espaces de Levi contenant $\tilde{M}$. Pour $i=1,2$ fixons des espaces de Levi $\tilde{L}_{i}$ contenant $\tilde{M}$. Alors on a:
\begin{rmq} 
$
d(\theta)^{-1}x(\gamma)\sum_{({\bf M}',\delta)\in {\mathcal X}^{{\mathcal E}}(\gamma)} \overline{I_{M'}^{\tilde{L}_{1}}(\delta,f_{1,L_{1}})}{I_{M'}^{\tilde{L}_{2}}(\delta,f_{2,L_{2}})}$
$$=\sum_{\gamma'}\overline{I_{\tilde{M}}^{\tilde{L_{1}},{\mathcal E}}(\gamma',f_{1,L_{1}})}{I_{\tilde{M}}^{\tilde{L_{2}},{\mathcal E}}(\gamma',f_{2,L_{2}})}
,
$$
o\`u la somme porte sur les classes de conjugaison $\gamma'$ stablement conjugu\'ees de $\gamma$.
\end{rmq}
C'est simplement la formule d'inversion des facteurs de transfert.

\begin{cor}Avec les notations pr\'ec\'edentes et $\gamma'$ parcourant le m\^eme ensemble que ci-dessus, 
$$\biggl(\sum_{\gamma'}I_{\tilde{M}}^{\tilde{G}}(\gamma,\omega,f_{1},f_{2})-\sum_{({\bf M}',\delta)\in {\mathcal X}^{\mathcal E}(\gamma)} d(\theta)^{-1}x(\gamma)^2I_{M'}^{\tilde{G},{\mathcal E}}(\delta,f_{1},f_{2})\biggr)=
$$
$$
\sum_{\tilde{L}_{1},\tilde{L}_{2}}d_{M}^{G}(L_{1},L_{2})$$
$$\sum_{\gamma'}
\biggl(\overline{I_{\tilde{M}}^{\tilde{L}_{1}}(\gamma,f_{1,L_{1}})}I_{\tilde{M}}^{\tilde{L}_{2}}(\gamma,f_{2,L_{2}})-\overline{I_{\tilde{M}}^{\tilde{L}_{1},{\mathcal E}}(\gamma,f_{1,L_{1}})}
I_{\tilde{M}}^{\tilde{L}_{2},{\mathcal E}}(\gamma,f_{2,L_{2}})\biggr),
$$o\`u la somme porte sur les couples d'espaces de Levi de $\tilde{G}$ contenant $\tilde{M}$.
\end{cor}
Chacun des termes de la premi\`ere \'egalit\'e v\'erifie une formule de scindage (cf \ref{rappelgeolocal} $(S)_{1}$ et $(S)_{3}$). Pour $I_{\tilde{M}}^{\tilde{G}}$, on scinde en sommant sur les couples d'espaces de Levi contenant $\tilde{M}$ comme dans l'\'enonc\'e. Pour $I_{M'}^{\tilde{G},{\mathcal E}}$, il s'agit des couples de Levi $L'_{1},L'_{2}$ contenant $M'$ et inclus dans un groupe endoscopique elliptique $G'$ de $\tilde{G}$. Dans ce dernier cas, on note \'evidemment $L_{1},L_{2}$ les espaces de Levi de $\tilde{G}$; ainsi, pour $i=1,2$, $L'_{i}$ et ${\bf M}'$ d\'efinissent une donn\'ee endoscopique elliptique pour $\tilde{L}_{i}$. Et il r\'esulte alors de la d\'efinition m\^eme donn\'ee en loc.cite que $d_{M}^{G}(L_{1},L_{2})=d_{M'}^{G'}(L'_{1},L'_{2})$.  Le corollaire r\'esulte alors de la remarque pr\'ec\'edente.

\subsubsection{R\'eduction pour la stabilisation g\'eom\'etrique \label{redlocgeo}}
\bf Hypoth\`ese de r\'ecurrence locale g\'eom\'etrique: \rm 
on suppose que l'on conna\^{\i}t l'\'egalit\'e $I_{\tilde{M}}^{\tilde{L}}(\gamma,\omega,f)=I_{\tilde{M}}^{\tilde{L},{\mathcal E}}(\gamma,f)$ pour tout espace de Levi propre $\tilde{M}$ de $\tilde{G}$ et pour tout espace de Levi propre $\tilde{L}$ de $\tilde{G}$ contenant $\tilde{M}$.

On remarque  que l'on a de toute fa\c{c}on  cette \'egalit\'e si $\tilde{M}=\tilde{L}=\tilde{G}$ ce qui permet de commencer la r\'ecurrence quand $\tilde{L}$ n'a pas de Levi propre.
\begin{prop} Avec les hypoth\`eses faites, pour toute paire de fonctions $f_{i}$ pour $i=1,2$ sur $\tilde{G}$, on a l'\'egalit\'e
$$I^{\tilde{G}}_{geo}(\omega,f_{1},f_{2})-\sum_{{\bf G}'}i(\tilde{G},{\bf G}') SI^{{\bf G}'}_{geo}(f_{1}^{{\bf G}'},f_{2}^{{\bf G}'})=$$
$$
\sum_{\tilde{M}} w(\tilde{M})^{-1} (-1)^{a_{\tilde{M},\tilde{G}}}\int_{\tilde{M}_{ell}/\sim}i'(\gamma)
$$
$$\overline{
\biggl(I_{\tilde{M}}^{\tilde{G}}(\gamma,\omega,f_{1})-I_{\tilde{M}}^{\tilde{G},{\mathcal E}}(\gamma,f_{1})\biggr)}{I^{\tilde{M}}(\gamma,\omega,f_{2,\tilde{M}})
}+$$
$$\overline{
I^{\tilde{M}}(\gamma,\omega,f_{1,\tilde{M}})}
{\biggl(I_{\tilde{M}}^{\tilde{G}}(\gamma,\omega,f_{2})-I_{\tilde{M}}^{\tilde{G},{\mathcal E}}(\gamma,f_{2})\biggr)}
$$
\end{prop}
C'est un corollaire des paragraphes pr\'ec\'edents.

\subsection{Le c\^ot\'e spectral de la formule des traces locales et sa stabilisation}
\subsubsection{Rappel des notations}
On note $Rat(G(F))$ le groupe des caract\`eres rationnels $Hom_{alg}(G(F),F^*)$; c'est un ${\mathbb Z}$-module et on pose ${\mathcal{A}}^*_{G,{\mathbb C}}:= Rat(G(F))\otimes_{{\mathbb Z}}{\mathbb C}$. On a ainsi d\'efini un espace vectoriel avec une structure r\'eelle; $\tilde{G}(F)$ op\`ere sur cet espace vectoriel de fa\c{c}on semi-simple et on note ${\mathcal{A}}^*_{\tilde{G},{\mathbb C}}$ le sous-espace vectoriel des \'el\'ements invariants pour cette action; ce sous-espace vectoriel a aussi une structure r\'eelle et il existe une application injective, $\tilde{G}(F)$ invariante de ${\mathcal  A}^*_{\tilde{G},{\mathbb C}}$ dans ${\mathcal{A}}^*_{{G},{\mathbb C}}$. Les \'el\'ements de ${\mathcal{A}}^*_{G,{\mathbb C}}$ s'identifient \`a des caract\`eres continus \`a valeurs dans ${\mathbb C}^*$  de $G(F)$ par l'application $$\forall g\in G(F), (\chi\otimes c)(g)=\vert \chi(g)\vert^{c}.
$$
On notera $i{\mathcal{A}}^*_{\tilde{G}}$ la partie imaginaire de ${\mathcal{A}}^*_{\tilde{G},{\mathbb C}}$, c'est-\`a-dire $Rat(G(F))\otimes_{{\mathbb Z}} i{\mathbb R}$.

Par les d\'efinitions ci-dessus, on a une application (c'est l'application usuelle) $H_{G}$ de $G(F)$ dans l'espace vectoriel dual de ${\mathcal{A}}^*_{G}$ et qui, par projection donne une application $H_{\tilde{G}}$ dans l'espace vectoriel dual de ${\mathcal{A}}^*_{\tilde{G}}$; on note ${\mathcal{A}}_{G}$ et ${\mathcal{A}}_{\tilde{G}}$ ces espaces vectoriels duaux.  On note ${\mathcal L}$  l'ensemble des \'el\'ements de ${\mathcal{A}}^*_{\tilde{G},{\mathbb C}}$ qui envoient dans $2\pi i{\mathbb Z}$  l'image de $H_{\tilde{G}}$ dans ${\mathcal{A}}_{\tilde{G}}$; ${\mathcal L}$ est trivial si $F$ est un corps archim\'edien et est un r\'eseau de $i{\mathcal{A}}^*_{\tilde{G}}$ si $F$ est $p$-adique. Dans la suite, on note $i{\mathcal{A}}^*_{\tilde{G},F}$ le quotient de $i{\mathcal{A}}^*_{\tilde{G}}$ par ce r\'eseau (trivial si $F$ est un corps archim\'edien).

On veut faire op\'erer ${\mathcal{A}}^*_{\tilde{G},{\mathbb C}}$ sur toute $\omega$-repr\'esentation de $\tilde{G}$; on ne peut pas le faire canoniquement mais d\`es que l'on fixe $\tilde{g}_{0}\in \tilde{G}(F)$, on pose pour tout \'el\'ement $g\tilde{g}_{0}\in \tilde{G}(F)$, $H_{\tilde{G}}(g\tilde{g}_{0})=H_{\tilde{G}}(g)$ et en composant avec $H_{\tilde{G}}$, tout \'el\'ement de ${\mathcal{A}}^*_{\tilde{G},{\mathbb C}}$ s'identifie en une fonction sur $\tilde{G}(F)$. La mutiplication d'une $\omega$-repr\'esentation par une telle fonction est encore une $\omega$-repr\'esentation.

On a ainsi d\'efini une action de $i{\mathcal{A}}^*_{\tilde{G},F}$ sur l'ensemble des $\omega$-repr\'esentations.

On g\'en\'eralise toutes ces notations aux sous-espaces de Levi, $\tilde{M}$ de $\tilde{G}$ en rempla\c{c}ant l'indice $G$ par un indice $M$.

La partie spectrale de la formule des traces locales est \'ecrite en  \cite{ftlt} 3.25, 3.26 et 6.6, pour toute fonction $f_{1},f_{2}\in I(\tilde{G})$

$$
I_{spec}^{\tilde{G}}(\omega,f_{1},f_{2})=\sum_{\tau} 	\int_{i{\mathcal{A}}^*_{\tilde{G},F}} d\lambda \iota(\tau) \vert Stab_{W^G\times i{\mathcal{A}}_{\tilde{G},F}^*}(\tau)\vert^{-1}\overline{\tau_{\lambda}(f_{1})}\tau_{\lambda}(f_{2}),
$$
o\`u $\tau$ parcourt un syst\`eme de repr\'esentants de l'action de $i{\mathcal{A}}^*_{\tilde{G}}$ agissant par tensorisation sur l'ensemble des $\omega$-repr\'esentations discr\`etes de $\tilde{G}$; $\iota(\tau)$ est d\'efini en \cite{ftlt}  \`a la fin de 2.11  et le stabilisateur  \`a la fin de \cite{ftlt} 2.9 et leurs valeurs explicites n'ont gu\`ere d'importance pour nous ici.
\subsubsection{Une remarque sur les mesures}
Le c\^ot\'e spectral de la formule des traces locale pour $\tilde{G}$ fait intervenir des int\'egrales sur les espaces $i{\mathcal{A}}^*_{\tilde{G},F}$ de $\omega$-repr\'esentations de $\tilde{G}$.  Tous les termes ont un coefficient qui tient compte du stabilisateur de la repr\'esentation sous $i{\mathcal{A}}^*_{\tilde{G},F}$ ce groupe agissant par tensorisation. 

Dans le cas des corps $p$-adiques, on peut donc modifier l'espace d'int\'egration en int\'egrant sur $i{\mathcal{A}}^*_{\tilde{M}}/{\mathcal L}'$ o\`u ${\mathcal L}'$ est un sous-r\'eseau inclus dans le r\'eseau ${\mathcal L}$ d\'ecrit ci-dessus. Cela modifiera le coefficient calculant le stabilisateur et ne modifiera pas l'int\'egrale.

C'est une remarque \`a faire car on va comparer le c\^ot\'e spectral de la formule des traces pour $\tilde{G}$ \`a celle de ses donn\'ees endoscopiques elliptiques; notant ${\bf G}'$ une telle donn\'ee. L'\'ellipticit\'e assure que ${\mathcal{A}}_{\tilde{G}}={\mathcal{A}}_{G'}$ mais n'assure \'evidemment pas que $i{\mathcal{A}}_{\tilde{G},F}$ soit \'egal \`a son analogue pour $G'$, $i{\mathcal{A}}_{{G}',F}$ si le corps de base est p-adique.Et la remarque pr\'ec\'edente montre que l'on peut quand m\^eme int\'egrer sur le m\^eme espace pour $\tilde{G}$ et pour ses donn\'ees endoscopiques elliptiques \`a condition de prendre les stabilisateurs dans l'espace sur lequel on int\'egre.

\subsubsection{D\'efinition du c\^ot\'e spectral stable de la formule des traces locales, pr\'eliminaires \label{defspeclocstabpreliminaire}}
Ici on suppose que $\tilde{G}$ est un bi-torseur sous un groupe $G$ quasi-d\'eploy\'e, qu'il est \`a torsion int\'erieure et que  $\omega$ est trivial. Il faut d\'efinir $SI^{\tilde{G}}_{spec}$. Par tensorisation $i{\mathcal{A}}^*_{\tilde{G},F}$ op\`ere dans l'ensemble  des caract\`eres unitaires de $A_{\tilde{G}}(F)$ et on fixe un ensemble, ${\mathcal X}$, de repr\'esentants pour ces orbites. Pour chacun de ces caract\`eres, $\chi$, on fixe  une base des repr\'esentations temp\'er\'ees et stables de $\tilde{G}$ se transformant sous $A_{\tilde{G}}(F)$ par ce caract\`ere $\chi$ et on note ${\mathcal B}_{\chi}$ cette base. On impose en plus \`a cette base d'avoir la propri\'et\'e d'orthogonalit\'e suivante: on consid\`ere le produit scalaire elliptique d\'efini pour les repr\'esentations elliptiques. On commence d'abord par choisir une base des repr\'esentations elliptiques stables orthogonale pour ce produit scalaire, ceci est possible gr\^ace \`a \cite{stabilisationI} 4.17 qui montre la compatibilit\'e de ce produit scalaire avec le transfert des fonctions cuspidales. Ensuite on utilise le fait que  la d\'ecomposition d'une repr\'esentation temp\'er\'ee stable en induites de repr\'esentations elliptiques (modulo le centre)   ne fait intervenir que des induites de repr\'esentations elliptiques stables.  En tensorisant par des \'el\'ements de $i{\mathcal{A}}^*_{\tilde{G},F}$ on obtient alors une base pour les repr\'esentations elliptiques stables sans hypoth\`ese sur le caract\`ere de $A_{\tilde{G}}(F)$. On fait la m\^eme construction pour les espaces de Levi de $\tilde{G}$ et en induisant on a ainsi une base des repr\'esentations temp\'er\'ees stables gr\^ace aux propri\'et\'es des d\'ecompositions des repr\'esentations temp\'er\'ees stables en combinaison lin\'eaires d'induites de repr\'esentations elliptiques n\'ecessairement stables (cf. le paragraphe 3 de \cite{prolongementdescaracteres} et \cite{stabilisationIV} 3.2).

\bf Remarque: \rm
Fixons un ensemble fini de couples $(\tilde{M},\chi_{M})$ form\'es d'un espace de Levi de $\tilde{G}$ et d'un \'el\'ement de $i{\mathcal{A}}^*_{\tilde{M},F}$. Fixons aussi un $K$-type. Alors il n'existe qu'un nombre fini de repr\'esentations temp\'er\'ees irr\'eductibles admettant un vecteur invariant sous ce $K$-type et qui soient un sous-module d'une combinaison lin\'eaire d'induites de repr\'esentations elliptiques \`a partir de donn\'ees $\tilde{M},\sigma$ o\`u en notant $\chi_{\sigma}$ le caract\`ere de $\sigma$ restreint \`a ${A}_{\tilde{M}}(F)$, les couples $(\tilde{M},\chi_{\sigma}
)$ soient dans l'ensemble fix\'e.

\

Cette propri\'et\'e de finitude est une g\'en\'eralisation simple d'un r\'esultat analogue pour les s\'eries discr\`etes d\^u \`a Harish-Chandra: en effet le r\'esultat d'Harish-Chandra se g\'en\'eralise aux repr\'esentations elliptiques (cf. \cite{appendice} paragraphe 6) et donc aux sous-modules des  induites comme dans l'\'enonc\'e.

\

On fixe  $\chi$ et on consid\`ere un  ensemble de nombres complexes $b_{\phi,\phi'}$ index\'es par les couples $\phi,\phi'$ d'\'el\'ements de ${\mathcal B}_{\chi}$; on dit qu'un tel syst\`eme est centralement fini si il  existe un nombre fini de couples ${\tilde{M},\chi_{M}}$ tel que $b(\phi,\phi')=0$ sauf \'eventuellement si $\phi$ et $\phi'$ ont une composante sous-quotient d'une induite d'une repr\'esentation elliptique de l'un des $\tilde{M}$ avec comme restriction \`a $A_{\tilde{M}}(F)$ le caract\`ere $\chi_{M}$. 

Il r\'esulte de la remarque que, pour un tel ensemble de nombres complexes, et pour toute paire de fonctions  $f_{1},f_{2}$ sur $\tilde{G}$, la somme
$$
\sum_{\chi}\sum_{\phi,\phi'\in {\mathcal B}_{\chi}} b(\phi,\phi')tr\, \phi(f_{1})\overline{tr\, \phi'(f_{2})}\eqno(1)
$$
est finie. Elle d\'efinit alors une distribution stable.

\subsubsection{Quelques propri\'et\'es de finitude\label{finitudecarcent}}
Le lecteur est en droit de se demander pourquoi ces d\'efinitions techniques. C'est \`a cause de probl\`eme de d\'efinition du c\^ot\'e spectral de la formule des traces locales, on ne veut que des sommes finies et c'est justifi\'e par le lemme suivant:
\begin{lem} Soit $\chi$ un \'el\'ement de $i{\mathcal{A}}^*_{\tilde{G},F}$. Alors il existe un ensemble fini de couples $(\tilde{M},\chi_{M})$ comme ci-dessus tel que toute repr\'esentation discr\`ete de $\tilde{G}$, dont le caract\`ere se restreint en le caract\`ere $\chi$ de $A_{\tilde{G}}(F)$, soit sous-module d'une combinaison lin\'eaire d'induites de repr\'esentations elliptiques des sous-groupes de Levi $\tilde{M}$ en n'autorisant que l'un des caract\`eres $\chi_{M}$ en restriction \`a $A_{\tilde{M}}(F)$.
\end{lem}
La d\'emonstration repose sur la remarque suivante:
\begin{rmq}
Soit $\pi$ une repr\'esentation discr\`ete associ\'ee \`a un triplet discret comme dans \cite{ftlt} 2.11. Alors il existe un unique couple $(\tilde{M},\sigma_{M})$ (unique \`a conjugaison pr\`es sous $G$) form\'e d'un espace de Levi $\tilde{M}$ de $\tilde{G}$ et d'une repr\'esentation elliptique $\sigma_{M}$ de $\tilde{M}$ tel que $\pi$ soit l'induite de $\sigma_{M}$.
\end{rmq}
On prend $\tilde{M}$ minimal contenant le triplet d\'efinissant $\pi$. Comme expliqu\'e en \cite{ftlt} 2.11, ce triplet pour $\tilde{M}$ correspond \`a une repr\'esentation elliptique $\sigma_{M}$ de $\tilde{M}$ et \cite{ftlt} 2.12 montre alors que $\pi$ est l'induite de $\sigma_{M}$. Et cette r\'ef\'erence montre aussi l'unicit\'e \`a conjugaison pr\`es.

La remarque entra\^{\i}ne le lemme. En effet si $\pi$ comme dans la remarque a pour restriction \`a $A_{\tilde{G}}(F)$ le caract\`ere $\chi$, $\sigma_{M}$ a la m\^eme propri\'et\'e. La repr\'esentation $\sigma_{M}$ a pour restriction \`a $A_{\tilde{M}}(F)$ la restriction du caract\`ere central d'une repr\'esentation d'un sous-groupe de Levi de $G$ invariant sous un \'el\'ement r\'egulier du groupe de Weyl de ce sous-groupe de Levi. Par d\'efinition de la r\'egularit\'e cela ne laisse qu'un nombre fini de possibilit\'es.

\subsubsection{D\'efinition du c\^ot\'e spectral stable de la formule des traces locales}
On revient \`a la situation de \ref{defspeclocstabpreliminaire}
On suppose d\'efini par r\'ecurrence le c\^ot\'e spectral stable pour les groupes endoscopiques elliptiques de $\tilde{G}$ propre, not\'e $SI^{G'}_{spec}$, de la forme:
$$
SI^{G'}_{spec}(f'_{1},f'_{2})=\sum_{\chi}\sum_{\phi,\phi'\in {\mathcal B}^{\bf G'}_{\chi}}\int_{i{\mathcal{A}}^*_{{\bf G}',F}}d\lambda b^{{\bf G}'}(\phi,\phi') \overline{tr\, \phi_{\lambda}(f'_{1})}{tr\, \phi_{\lambda}'(f'_{2})}
$$ avec un ensemble de nombre complexes $b(\phi,\phi')$ centralement fini.

Il faut remarquer que le sous-groupe de $i{\mathcal{A}}^*_{{\bf G}',F}$ form\'e des caract\`eres dont la restriction \`a ${A}_{{\bf G}'}(F)$ est trivial op\`ere sur l'ensemble des repr\'esentations stables ayant un caract\`ere sur ${ A}_{{\bf G}'}(F)$ fix\'e. Comme on int\'egre sur $i{\mathcal{A}}^*_{{\bf G}',F}$, on demande aussi \`a $\sum_{\phi,\phi'}b^{\bf G'}(\phi,\phi')\phi\otimes \overline{\phi'}$ d'\^etre invariant sous cette action.
\begin{prop}
Il existe un syst\`eme de coefficients centralement fini $b(\phi,\phi')$ associ\'e au triplet $\chi,\phi,\phi'$ comme ci-dessus tel que en posant $SI_{spec}^{\tilde{G}}$ la distribution obtenue en int\'egrant (1) sur $i{\mathcal{A}}^*_{\tilde{G},F}$, on ait pour toute paire de fonctions $f_{1},f_{2}$ sur $\tilde{G}$
$$
I^{\tilde{G}}_{spec}(f_{1},f_{2})-\sum_{{\bf G}'\neq \tilde{G}}SI^{G'}_{spec}(f_{1}^{\bf G'},f_{2}^{\bf G'})=SI^{\tilde{G}}_{spec}(\omega,f_{1},f_{2}).
$$
\end{prop}
Ce que dit la proposition est que le c\^ot\'e gauche est stable et s'\'ecrit sous la forme (1) pour des bons choix. 
On montre d'abord que le c\^ot\'e gauche de l'\'egalit\'e est une int\'egrale sur $i{\mathcal{A}}^*_{\tilde{G},F}$ d'une somme comme dans (1) mais o\`u $\phi$ et $\phi'$ parcourt une base de l'ensemble des repr\'esentations temp\'er\'ees et pas seulement des repr\'esentations stables. Pour cela, on transf\`ere terme \`a terme les \'el\'ements de la somme sur ${\bf G}'$ (cf. le th\'eor\`eme du paragraphe 3 de \cite{prolongementdescaracteres} et \cite{stabilisationIV} 3.2 dans le cas archim\'edien); un caract\`ere $\chi'$ pour $A_{G'}(F)$ se transf\`ere en un caract\`ere $\chi_{1}$ de $A_{\tilde{G}}(F)$ par l'application naturelle du deuxi\`eme groupe dans le premier. On tensorise par un \'el\'ement de $i{\mathcal{A}}^*_{\tilde{G},F}$ pour le ramener en l'un des repr\'esentants fix\'es, $\chi$. On fait la somme sur toutes les donn\'ees et tous les caract\`eres et on obtient l'assertion (sans la stabilit\'e). On montre maintenant que les sommes que l'on int\'egre  pour $\chi$ fix\'e sont des distributions stables; c'est ici qu'il faut commencer par faire agir le sous-groupe de $i{{\mathcal{A}}^*_{\tilde{G},F}}$ qui agit trivialement. Ce sous-groupe est fini et apr\`es cette op\'eration l'int\'egrale porte sur le quotient de $i{{\mathcal{A}}^*_{\tilde{G},F}}$ par ce sous-groupe.

Dans le cas de torsion int\'erieure que nous sommes en train de consid\'erer, il a \'et\'e montr\'e en \cite{stabilisationV} 1.13, que les int\'egrales orbitales pond\'er\'ees \'equivariantes pour les \'el\'ements fortement r\'eguliers sont stabilisables c'est-\`a-dire que $I_{\tilde{M}}^{\tilde{G}}(\gamma,\omega,f)=I_{\tilde{M}}^{\tilde{G},{\mathcal E}}(\gamma,\omega,f)$ en de tels points $\gamma$. Ainsi la partie g\'eom\'etrique de la formule des traces locales est stabilisable ou encore  le membre de gauche de (3) dans \ref{separationsuivantlevi} est nul (cf. \ref{redlocgeo}).

Ainsi on a une distribution  du c\^ot\'e spectral qui est stable c'est-\`a-dire nulle si on l'applique \`a  un couple de fonctions $f_{1},f_{2}$ o\`u l'une des fonctions annule toutes les int\'egrales orbitales semi-simples stables.

Fixons $f_{1}\in I(\tilde{G})$ d'image nulle dans $SI(\tilde{G})$; $f_{2}$ varie librement donc l'application $\pi\mapsto tr\, \pi(f_{2})$ d\'ecrit toutes les fonctions de Paley-Wiener sur l'ensemble des repr\'esentations temp\'er\'ees de $\tilde{G}$. Cela prouve que si $\phi'$ parcourt une base des repr\'esentation temp\'er\'ees (ayant les propri\'et\'es d'orthogonalit\'e que nous avons fix\'ees), le coefficient de $\phi'$ qui est une fonction de $f_{1}$ est nulle quand on l'\'evalue en $f_{1}$. Ainsi on peut faire parcourir \`a $\phi$ une base des repr\'esentations stables (cf. \cite{prolongementdescaracteres} paragraphe  2 et \cite{stabilisationIV} 2.8 dans le cas archim\'edien) modulo l'action de $i{{\mathcal{A}}^*_{\tilde{G},F}}$; comme l'application $I(\tilde{G})\rightarrow SI(\tilde{G})$ est surjective, on \'echange les r\^oles de $f_{1}$ et $f_{2}$ pour obtenir le r\'esultat annonc\'e.

\subsubsection{Description plus fine du c\^ot\'e spectral de la formule des traces locale\label{locspecfin}}
On peut am\'eliorer la proposition pr\'ec\'edente en d\'ecoupant suivant les espaces de Levi de $\tilde{G}$ et leurs repr\'esentations elliptiques. On fixe toujours un ensemble de repr\'esentants des caract\`eres unitaires de $A_{\tilde{G}}(F)$ modulo l'action de $i{\mathcal{A}}^*_{\tilde{G},F}$. Pour $\chi$ un tel repr\'esentant, on consid\`ere l'ensemble fini de couples $(\tilde{M},\chi_{M})$ satisfaisant le lemme de \ref{finitudecarcent}. Pour $(\tilde{M},\chi_{M})$ dans cet ensemble on fixe Quand on a fix\'e un tel caract\`ere et si $\tilde{G},\omega$ est \`a torsion int\'erieure comme dans le paragraphe pr\'ec\'edent, on fixe aussi une base orthogonale (pour le produit scalaire elliptique) de l'ensemble des repr\'esentations elliptiques stables de $\tilde{M}$ sur lesquelles $A_{\tilde{M}}(F)$ agit par ce caract\`ere; on note ${\mathcal B}(\chi_{{M}})$ cette base.

On utilise la notation ambigu\"e suivante: soit $\tilde{L}$ un espace de Levi de $\tilde{G}$ et soit $\tau$ une repr\'esentation elliptique de $\tilde{L}$, c'est \`a dire correspondant \`a un triplet elliptique de \cite{ftlt} 2.11 qui avec les notations de loc.cite est de la forme $(M,\sigma,\tilde{w})$; ici $M$ est un sous-groupe de Levi de $L$ et n'a rien \`a voir avec les $\tilde{M}$ de ce paragraphe. On dit que l'induite de $\tau$ de $\tilde{L}$ \`a $\tilde{G}$ est irr\'eductible si le $\tilde{R}$-groupe de $\sigma$ calcul\'e dans $\tilde{L}$ est aussi celui calcul\'e dans $\tilde{G}$, c'est-\`a-dire que le morphisme de  $ W^{\tilde{L}}(\sigma)$ dans $W_{0}^G(\sigma)\backslash W^{\tilde{G}}(\sigma)$ est un isomorphisme; ce morphisme est injectif car $\tau$ a \'et\'e suppos\'e elliptique et c'est la surjectivit\'e qui est la condition. La repr\'esentation induite de $\tau$ est une repr\'esentation discr\`ete de $\tilde{G}$ pour laquelle $\iota({ind\, \tau})$ est bien d\'efini (\cite{ftlt} 2.11). 
\begin{prop}
(i) On ne fait pas d'hypoth\`ese sur $\tilde{G},\omega$.
$$
I^{\tilde{G}}_{spec}(f_{1},f_{2})=\int_{i{\mathcal{A}}^*_{\tilde{G},F}} d\lambda
$$
$$\sum_{\tilde{M},\chi_{\tilde{M}},\tau} \vert Stab_{W^G\times i{\mathcal{A}}_{\tilde{G}}^*}(ind\, \tau)\vert^{-1}\iota(ind\, \tau) \overline{tr\, \tau_{\lambda}(f_{1,\tilde{M}})} tr\, \tau_{\lambda}(f_{2,\tilde{M}})
$$
o\`u $\tilde{M}$ parcourt les classes de conjugaison d'espaces de Levi et o\`u   les $\tau$ parcourt l'ensemble des repr\'esentations elliptiques de $\tilde{M}$ dont l'induite  \`a $\tilde{G}$ est irr\'eductible (modulo conjugaison sous le normalisateur de $\tilde{M}$ dans $G$) de caract\`ere central $\chi_{{M}}$.

(ii) Ici on suppose que $\tilde{G},\omega$ est \`a torsion int\'erieure avec $\omega=1$ et $G$ quasi-d\'eploy\'e. Il existe un syst\`eme de coefficients $b(\tilde{M}, \chi_{\tilde{M}},\phi,\phi')$ centralement fini index\'e par les  donn\'ees pr\'ec\'edant l'\'enonc\'e (les espaces de Levi sont pris \`a conjugaison pr\`es), identiquement nul pour presque tout $\chi_{{M}}$ tel que 
pour toute paire de fonctions $f_{1},f_{2}$, 
$$
SI^{\tilde{G}}_{spec}(f_{1},f_{2})=\int_{i{\mathcal{A}}^*_{\tilde{G},F}}\sum_{\tilde{M},\chi_{{M}},\phi,\phi'}b(\tilde{M}, \chi_{{M}},\phi,\phi') \overline{tr\, \phi_{\lambda}(f_{1,\tilde{M}})} {tr\, \phi_{\lambda}'(f_{2,\tilde{M}})}\, d\lambda,
$$o\`u l'indice $\lambda$ est la tensorisation par $\lambda$.
\end{prop}
(i)
Il faut d'abord remarquer que dans $I_{spec}^{\tilde{G}}$ seules des repr\'esentations discr\`etes de $\tilde{G}$ interviennent. De plus, on a bien une somme d'int\'egrales sur $i{\mathcal{A}}^*_{\tilde{G},F}$ de produit $\overline{tr\, \tau_{\lambda}(f_{1})}tr\, \tau_{\lambda}(f_{2})$, o\`u comme on vient de le dire $\tau$ est discr\`ete. Or une repr\'esentation discr\`ete est induite d'une repr\'esentation elliptique comme on l'a vu dans \ref{finitudecarcent} pour un unique espace de Levi de $\tilde{G}$  (unique \`a conjugaison  pr\`es). Cela donne (i)

(ii)
On raisonne comme dans la proposition pr\'ec\'edente: on applique le r\'esultat par r\'ecurrence aux donn\'ees endoscopiques elliptiques propres. On remarque qu'un sous-groupe  de Levi, ${\bf M}'$ d'une donn\'ee endoscopique elliptique ${\bf G}'$ qui n'est pas relevant pour $\tilde{G}$ a une contribution nulle car il s'applique \`a $f^{\bf G'}_{M'}=0$. On n'a donc qu'\`a transf\'erer des repr\'esentations elliptiques stables de sous-groupes endoscopiques elliptiques d'espaces de Levi de $\tilde{G}$. Un tel transfert est une combinaison lin\'eaire de repr\'esentations elliptiques de ce sous-espace de Levi. Cela montre une formule comme dans l'\'enonc\'e mais sans savoir que les $\phi$ et $\phi'$ sont stables.

Pour avoir cette stabilit\'e, on utilise encore \cite{stabilisationIV} 2.8 dans le cas archim\'edien et  le paragraphe 2 de \cite{prolongementdescaracteres}  dans le cas $p$-adique.
\subsubsection{Remarque sur la stabilisation locale spectrale\label{stabilisationspecgrossiere}\label{redlocspec}}
Ici on revient \`a un espace $\tilde{G}$ g\'en\'eral, donc on ne conna\^{\i}t plus la stabilisation des int\'egrales orbitales pond\'er\'ees et on n'a ni la stabilisation g\'eom\'etrique ni la stabilisation spectrale. Comme dans le cas \`a torsion int\'erieure, on fixe un ensemble de repr\'esentants dans le quotient du groupe des caract\`eres unitaires de $A_{\tilde{G}}(F)$ modulo l'action par tensorisation de $i{\mathcal{A}}^*_{\tilde{G},F}$ et pour $\chi$ un tel repr\'esentant, on consid\`ere les couple ${\tilde{M},\chi_{M}}$ form\'es d'un espace de Levi de $\tilde{G}$ et d'un caract\`ere unitaire de $A_{\tilde{M}}(F)$ prolongeant $\chi$.  On fixe une base ${\mathcal B}(\tilde{M},\chi_{M})$ du groupe de Grothendieck des repr\'esentations elliptiques de $\tilde{M}$ se transformant via $\chi_{M}$ sous $A_{\tilde{M}}(F)$. 
\begin{lem}Pour tout $\chi$ et pour tout ${\tilde{M},\chi_{M}}$ et $\tau,\tau' \in {\mathcal B}(M,\chi_{M})$, il existe un syst\`eme de coefficients $a(\tilde{M},\chi_{M},\tau,\tau')$ nul pour presque tout $\chi_{M}$,   tel que pour toute paire de fonctions $f_{1},f_{2}$ sur $\tilde{G}$, on ait:
$$
I^{\tilde{G}}_{spec}(\omega,f_{1},f_{2})-\sum_{{\bf G}'}i(\tilde{G},{\bf G}')SI^{\bf G'}_{spec}(f_{1}^{\bf G'},f_{2}^{\bf G'})=\eqno(*)
$$
$$
\int_{i{\mathcal{A}}^*_{\tilde{G},F}}d\lambda \sum_{\tilde{M},\chi_{M}}\sum_{\tau,\tau'\in {\mathcal B}^{\tilde{M}}(\chi_{M})}a(\tilde{M},\chi_{M},\tau,\tau')\overline{ tr\, \tau_{\lambda}(f_{1,\tilde{M}})}{tr\, \tau_{\lambda}' (f_{2,\tilde{M}})}.\eqno(1)
$$
\end{lem}
Cela se d\'emontre comme dans le cas de torsion int\'erieure,  et ici on s'attend \`a ce que les coefficients soient tous nuls et on le d\'emontrera.

\
\begin{rmq} Dans l'\'enonc\'e pr\'ec\'edent, la somme ne porte que sur les $\tilde{M}$ sous-espace de Levi propre de $\tilde{G}$.
\end{rmq}
En effet, supposons que $f_{1}$ et $f_{2}$ soient cuspidales. Alors dans (1) ne reste que le terme pour $\tilde{M}=\tilde{G}$. On \'ecrit (*) \`a l'aide de l'\'ecriture g\'eom\'etrique de \ref{redlocgeo}. Comme on suppose que $f_{1}$ et $f_{2}$ sont cuspidales les termes constants de $f_{1}$ et $f_{2}$ pour les espaces de Levi propres de $\tilde{G}$ ont toutes leurs int\'egrales orbitales  nulles. Avec \ref{redlocgeo}, on sait que (*) est nul. Donc (1) est nul pour tout $f_{2}$ cuspidale.  Via la trace tordue, les fonctions cuspidales s\'eparent les \'el\'ements d'une base ${\mathcal B}(\tilde{G},\chi)$ pour $\chi$ fix\'e.  Par inversion de Fourier, on conclut alors \`a la nullit\'e des $a(\tilde{G},\chi,\tau,\tau')$ quand $\tau$ et $\tau'$ parcourt une base des repr\'esentations elliptiques de $\tilde{G}$ de caract\`ere $\chi$ en restriction \`a $A_{\tilde{G}}(F)$ comme dans l'\'enonc\'e du lemme.

\subsection{Elimination de certaines conditions\label{eliminationsingularite}}

On va aussi en d\'eduire le corollaire ci-dessous. On a besoin de la terminologie suivante.
Une fonction $f\in I(\tilde{G})$ se d\'ecompose suivant le th\'eor\`eme de Paley-Wiener en une composante cuspidale et une composante dite non elliptique dont la trace est nulle sur toute les $\omega$ repr\'esentations elliptiques de $\tilde{G}$. 
On a aussi besoin de la notion d'\'el\'ements isol\'es ou plus exactement exceptionnels. Ils sont d\'efinis en \cite{stabilisationVIII} 4.4: ces \'el\'ements n'apparaissent que pour certains choix de $\tilde{G},\omega$ (donn\'es en \cite{stabilisationIII} 6.3) et sont tels que leur partie semi-simple stabilise une paire de Borel \'epingl\'ee; appelons un tel \'el\'ement un \'el\'ement isol\'e. C'est une g\'en\'eralisation de la notion d'\'el\'ement unipotents modulo le centre du cas non tordu mais ils ne g\^enent que pour certains $\tilde{G},\omega$ et on ne les consid\`ere que  dans ces cas. On a besoin de les consid\'erer s\'epar\'ement \`a cause des hypoth\`eses de la proposition de \cite{stabilisationVIII} 4.4. C'est pour cela qu'il vaut mieux les appeler exceptionnels. Ceci dit pour nous la seule  chose qui nous importe est qu'ils n'apparaissent pas sur un corps local archim\'edien.

On fixe un espace de Levi $\tilde{M}$ de $\tilde{G}$ et une fonction $f\in I(\tilde{G})$. On dit que $f$ est $\tilde{M}$ cuspidale si les termes constants de $f$ sont nuls pour tout espace de Levi qui ne contient pas, \`a conjugaison pr\`es, $\tilde{M}$ (c'est la g\'en\'eralisation imm\'ediate de la d\'efinition donn\'ee en \cite{partie3} avant le lemme 2.3)
\begin{cor} On fixe $\tilde{M}$ un espace de Levi de $\tilde{G}$ et on fixe $f_{2}\in I(\tilde{G})$ en supposant que $f_{2}$ est $\tilde{M}$-cuspidale.

(i) $$I^{\tilde{G}}_{spec}(\omega,f_{1},f_{2})-\sum_{{\bf G}'}i(\tilde{G},{\bf G}')SI^{\bf G'}_{spec}(f_{1}^{\bf G'},f_{2}^{\bf G'})=0\eqno(1)$$ pour toute fonction $f_{1}\in I(\tilde{G})$ si et seulement si cela est vrai pour toute fonction $f_{1}$ qui est $\tilde{M}$-cuspidale.

(ii) On suppose que le corps de base est p-adique. Alors $$I^{\tilde{G}}_{spec}(\omega,f_{1},f_{2})-\sum_{{\bf G}'}i(\tilde{G},{\bf G}')SI^{\bf G'}_{spec}(f_{1}^{\bf G'},f_{2}^{\bf G'})=0$$ pour toute   fonction $f_{1}\in I(\tilde{G})$  si et seulement si cela est vrai pour toute fonction $f_{1}\in I(\tilde{G})$ o\`u l'on suppose que $f_{1}$ est $\tilde{M}$ cuspidale et a une composante non elliptique nulle pr\`es des \'el\'ements isol\'es.
\end{cor}
Comme on suppose que $f_{2}$ est $\tilde{M}$ cuspidale, dans (1) il suffit de sommer sur les espaces de Levi contenant $\tilde{M}$.
On fixe $f_{1}$ et on note $f'_{1}$ une fonction dans $I(\tilde{G})$ qui a m\^eme composante de Paley-Wiener que $f_{1}$ pour tous les espaces de Levi contenant $\tilde{M}$ et qui vaut 0 ailleurs.  V\'erifions que $f'_{1}$ est  $\tilde{M}$ cuspidal: puisque les traces de cette fonction pour toute induite \`a partir d'un espace de Levi ne contenant pas $\tilde{M}$ sont nulles, on applique  le th\'eor\`eme 0 de \cite{ftlt} 5.5. Pour cela on fixe $\tilde{L}$ un espace de Levi de $\tilde{G}$ dont aucun conjugu\'e ne contient $\tilde{M}$ et une repr\'esentation $\pi_{\tilde{L}}$ de cet espace de Levi. Comme la trace de $f'_{1}$ sur l'induite de $\pi_{\tilde{L}}$ est nulle, la trace du terme constant $f'_{1,\tilde{L}}$ sur $\pi_{\tilde{L}}$ est nulle.  On sait alors que  le terme constant de $f'_{1}$ pour un tel espace de Levi a toutes ses int\'egrales orbitales semi-simples r\'eguli\`eres nulles et donc que $f'_{1}$ est bien $\tilde{M}$-cuspidale. Cela d\'emontre le (i) du corollaire.

Pour le  (ii) du corollaire, on suppose comme dans l'\'enonc\'e que le corps de base est p-adique.
Comme on l'a remarqu\'e dans le corollaire pr\'ec\'edant, dans (1) n'interviennent que des espaces de  Levi propres de $\tilde{G}$  et des repr\'esentations elliptiques de ces espaces de Levi. Ainsi (1) est certainement nul si $f_{1}$   est cuspidale. On peut donc d\`es le d\'epart supposer que $f_{1}$ a sa composante elliptique nulle.

Ainsi avec l'hypoth\`ese on sait  que (1) est nulle si en plus $f_{1}$ est nulle pr\`es des \'el\'ements exceptionnels, on voit que (1) restreint \`a $f_{2}$ fix\'e comme dans l'\'enonc\'e est une forme lin\'eaire sur les \'el\'ements de $I(\tilde{G})$ de composante elliptique nulle qui est n\'ecessairement une combinaison lin\'eaire convenable des int\'egrales orbitales   en les points isol\'es. A la suite de \cite{arthurcrelle}, on a montr\'e qu'une int\'egrale orbitale est    une somme d'int\'egrales sur les espaces $i{\mathcal{A}}^*_{\tilde{L},F}$ (quand $\tilde{L}$ parcourt l'ensemble des sous-espace de Levi de $\tilde{G}$) de traces d'induites de repr\'esentations elliptiques de l'espace de Levi $\tilde{L}$ (cf \cite{appendice}, th\'eor\`eme du paragraphe 5). Par hypoth\`ese sur la composante elliptique de $f_{1}$, seuls les espaces $\tilde{L}$ propres interviennent de fa\c{c}on \'eventuellement non nulle. Mais en revenant \`a la d\'efinition de (1) comme int\'egrale sur $i{\mathcal{A}}^*_{\tilde{G},F}$ l'\'egalit\'e des deux distributions est impossible  sans que chacune des distributions soit nulle. Ainsi (1) est nul sous la seule hypoth\`ese que  $f_{1}$ a sa composante elliptique nulle, ce qui prouve (ii).

\begin{rmq} En jouant sur la sym\'etrie entre $f_{1}$ et $f_{2}$, on vient de montrer que (1) est nulle pour tout couple de fonctions $(f_{1},f_{2})$ dans $I(\tilde{G})$ dont l'une est $\tilde{M}$ cuspidale si et seulement si cela est vrai pour tout couple de fonctions $(f_{1},f_{2})$ dont on suppose que les deux sont $\tilde{M}$ cuspidales et que les deux ont une composante non elliptique nulle pr\`es des \'el\'ements isol\'es (condition que l'on n'impose que pour les places p-adiques).
\end{rmq}
\begin{rmq} On fixe un espace de Levi $\tilde{M}$ de $\tilde{G}$. Supposons que l'on sache que pour toute fonction $f$ et pour tout espace de Levi contenant strictement $\tilde{M}$ on ait tout $\gamma$ \'el\'ement semi-simple r\'egulier de $\tilde{L}$,  $$I^{\tilde{G}}_{\tilde{L}}(\omega,\gamma,f)=I^{\tilde{G},{\mathcal E}}_{\tilde{L}}(\gamma,f)\eqno(2)$$ 
et que l'on sache aussi que 
tout $\gamma$ \'el\'ement semi-simple r\'egulier de $\tilde{M}$,  $$I^{\tilde{G}}_{\tilde{M}}(\omega,\gamma,f)=I^{\tilde{G},{\mathcal E}}_{\tilde{M}}(\gamma,f)\eqno(3)$$ si $f$ est $\tilde{M}$ cuspidale, alors pour une telle fonction $f$, (2) est vrai pour tout espace de Levi $\tilde{L}$ de $\tilde{G}$ sans restriction. En particulier (2) est vrai pour toute fonction cuspidale sur $\tilde{G}$ et pour tout espace de Levi $\tilde{L}$.
\end{rmq}

En effet, d'apr\`es l'hypoth\`ese, la formule des traces locales pour deux fonctions $f_{1},f_{2}$ de $\tilde{G}$ toutes deux $\tilde{M}$ cuspidales est stabilisable, c'est-\`a dire que (1) est nulle sous ces conditions. D'apr\`es la remarque pr\'ec\'edente, (1) est nulle sous la seule hypoth\`ese que $f_{1}$ est $\tilde{M}$ cuspidale, donc sans hypoth\`ese sur $f_{2}$. Cela force la nullit\'e du c\^ot\'e g\'eom\'etrique de \ref{redlocgeo} pour toute fonction $f_{2}$ et pour $f_{1}=f$ comme dans l'\'enonc\'e. En utilisant l'hypoth\`ese (2) et la $\tilde{M}$ cuspidalit\'e de $f=f_{1}$, on obtient qu'une combinaison lin\'eaire des termes 
$$
\int_{\tilde{L}_{ell}/\sim}i'(\gamma) \overline{(I^{\tilde{G}}_{\tilde{L}}(\omega,\gamma,f)-I^{\tilde{G},{\mathcal E}}_{\tilde{L}}(\gamma,f))}I^{\tilde{L}}(\gamma,\omega,f_{2,\tilde{L}})$$
est nulle sans hypoth\`ese sur $f_{2}$. Cela force la nullit\'e des fonctions 
$$
I^{\tilde{G}}_{\tilde{L}}(\omega,\gamma,f)-I^{\tilde{G},{\mathcal E}}_{\tilde{L}}(\gamma,f).
$$

\subsection{Stabilisation g\'eom\'etrique sous hypoth\`eses\label{stabgeohyp}}
\bf Hypoth\`eses de r\'ecurrence locales g\'eom\'etriques propres: \rm
on suppose que pour tout espace de Levi propre $\tilde{L}$ de  $\tilde{G}$ et pour tout espace de Levi $\tilde{R}$ de $\tilde{L}$, on sait que pour tout $f\in I(\tilde{L})$ et pour tout \'el\'ement $\gamma$ semi-simple r\'egulier de $\tilde{R}$, on a $I^{\tilde{L}}_{\tilde{R}}(\omega,\gamma,f)=I^{\tilde{L},{\mathcal E}}_{\tilde{R}}(\gamma,f)$.

\

\bf Hypoth\`ese de r\'ecurrence locales g\'eom\'etriques d\'ependant d'un espace de Levi $\tilde{M}$ de $\tilde{G}$ fix\'e: \rm
on suppose aussi que pour tout espace de Levi, $\tilde{L}$ de  $\tilde{G}$ contenant strictement $\tilde{M}$ et pour tout $f\in I(\tilde{G})$ et tout $\gamma$ \'el\'ement semi-simple r\'egulier de $\tilde{L}$,  $I^{\tilde{G}}_{\tilde{L}}(\omega,\gamma,f)=I^{\tilde{G},{\mathcal E}}_{\tilde{L}}(\gamma,f)$.

\

Ce sont des hypoth\`eses de r\'ecurrence que l'on a \'evidemment le droit de faire puisqu'on les initialise sans probl\`eme pour $\tilde{M}=\tilde{G}$. Par contre les deux hypoth\`eses ci-dessous sont d'une autre nature, on les d\'emontrera via un argument global en \ref{preuvegeoloc}. On fixe $\tilde{M}$ un espace de Levi propre de $\tilde{G}$. On fixe aussi ${\bf M'}$ une donn\'ee endoscopique elliptique de $\tilde{M}$;  on fixe  des donn\'ees auxiliaires $M'_{1}$,..., $\Delta_{{\bf M}'}$; soit $\delta$ un \'el\'ement semi-simple r\'egulier de $\tilde{M}'(F)$ que l'on rel\`eve en un \'el\'ement $\delta_{1}\in \tilde{M}'_{1}(F)$. On pose, pour tout $f\in I(\tilde{G})$ nul pr\`es des \'el\'ements exceptionnels (quand il y en a)
$$
\epsilon^{{\bf M}'}_{\tilde{M}}(f)(\delta_{1}):=\sum_{\gamma\in \tilde{M}} \Delta_{{\bf M}'}(\delta_{1},\gamma)z(\gamma)^{-1}
I^{\tilde{M}}(\gamma,\omega,\epsilon_{\tilde{M}}(f)), \eqno(1)
$$
o\`u $z(\gamma)$ est le nombre d'\'el\'ements de $Cent_{M}(\gamma,F)/Cent^0_{M}(\gamma,F)$. La fonction $\epsilon_{\tilde{M}}(f)$ est celle d\'efinie en \cite{stabilisationVIII} et \cite{stabilisationIX}, c'est-\`a-dire que l'on a $I^{\tilde{M}}(\gamma,\omega,\epsilon_{\tilde{M}}(f))=I^{\tilde{G},{\cal E}}_{\tilde{M}}(\gamma,f)-I^{\tilde{G}}_{\tilde{M}}(\gamma,\omega,f)$. Dans le cas $p$-adique, elle n'est d\'efinie que si $f$ est nulle pr\`es des \'el\'ements exceptionnels. 
\bigskip

{\bf{hypoth\`ese cl\'e 1}} \sl On suppose que pour toute donn\'ee endoscopique elliptique ${\bf M}'$ de $\tilde{M}$ comme ci-dessus il existe une fonction lisse sur l'ensemble des \'el\'ements semi-simples r\'eguliers fortement r\'eguliers  de  $\tilde{M}'(F)$, $\epsilon({\bf M}',\delta)$ telle que pour toute fonction $f$ sur $\tilde{G}$, nulle pr\`es des \'el\'ements exceptionnels  et pour tout \'el\'ement semi-simple fortement r\'egulier $\delta$  de $\bf{M}'$, on a 
$$
\epsilon^{{\bf M}'}_{\tilde{M}}(f)(\delta_{1})=\epsilon({\bf M}',\delta) f_{\tilde{M}}^{{\bf M}'}(\delta_{1}),
$$
o\`u le terme de droite est \'evidemment d\'efini comme dans (1)

Puisque les deux fonctions de $\delta_{1}$ se transforment selon le m\^eme caract\`ere du tore central $C_{1}(F)$, le terme  $\epsilon({\bf M}',\delta)$ ne d\'epend bien que de $\delta$. On voit facilement qu'il ne d\'epend pas du choix des donn\'ees auxiliaires. Dans la suite, on oublie ces donn\'ees qui ne sont pas importantes et on note simplement $\delta$ l'\'el\'ement $\delta_{1}$.

\rm

Dans les cas courants (avec $\omega=1$ ou seulement quadratique) par exemple les cas de GL(n) tordu, cette seule hypoth\`ese va suffire pour la stabilisation locale des int\'egrales orbitales pond\'er\'ees. Mais dans le cas g\'en\'eral et en particulier si $\omega$ n'est pas quadratique on a besoin de plus. On va utiliser l'hypoth\`ese suivante, qui est tr\`es facile \`a d\'emontrer au moment o\`u on d\'emontre l'hypoth\`ese cl\'e 1:

\

{\bf hypoth\`ese cl\'e 2} \sl on suppose que la situation locale est une composante locale d'une situation globale: c'est-\`a-dire que $F$ est un corps global, que $v_{0}$ est une place de $F$ et que ce qui nous int\'eresse est ce qui se passe en la place $v_{0}$. On suppose aussi que ${\bf M}'$ est une donn\'ee endoscopique elliptique de la situation globale et que $\delta$ est un \'el\'ement semi-simple r\'egulier de ${\bf M}'(F)$. On fixe encore un ensemble $V$ de places de $F$ contenant $V_{ram}$ et tel que, pour $v\not\in V$, ${\bf M}'$ soit non ramifi\'ee en $v$ et  $\delta\in \tilde{K}_{v}^{M'}$ et on a donc d\'efini pour tout $v\in V$, $\epsilon({\bf M}',\delta_{v})$ not\'e $\epsilon_{v}({\bf M}',\delta)$  alors:
$$
\sum_{v\in V}\epsilon_{v}({\bf M}',\delta)=0.
$$

\rm

\begin{thm} Sous les hypoth\`eses faites, pour toute fonction $f$ sur $\tilde{G}$ et pour tout \'el\'ement  $\gamma$ fortement r\'egulier de $\tilde{M}$, $I_{\tilde{M}}^{\tilde{G}}(\omega,\gamma,f)=I_{\tilde{M}}^{\tilde{G},{\mathcal E}}(\gamma,f)$.
\end{thm}
On fait d'abord la r\'eduction suivante: il suffit de montrer que $\epsilon({\bf M}',\delta)$ est une fonction identiquement nulle pour toute donn\'ee endoscopique elliptique ${\bf M}'$ de $\tilde{M}$. En effet supposons qu'il en soit ainsi, par inversion des facteurs de transfert, on en d\'eduit que pour toute fonction $f\in I(\tilde{G})$, nulle pr\`es des  des \'el\'ements exceptionnels si $F_{v}$ est p-adique et sans restriction si $F_{v}$ est archim\'edienne, la fonction $\epsilon_{\tilde{M}}(\gamma,f)$ est identiquement nulle. C'est le r\'esultat cherch\'e pour ces fonctions $f$. 
Il faut donc enlever la restriction dans le cas o\`u $F_{v}$ est p-adique. Avec les formules de descente, on se ram\`ene imm\'ediatement au cas o\`u $\gamma$ est elliptique dans $\tilde{M}$.

Comme le th\'eor\`eme est d\'ej\`a montr\'e pour les fonctions cuspidales (cf. la deuxi\`eme remarque de \ref{eliminationsingularite}), on a le th\'eor\`eme pour toute fonction $f$ dont la composante de Paley-Wiener non cuspidale est nulle pr\`es des \'el\'ements exceptionnels. On consid\`ere $f_{1},f_{2}$ un couple d'\'el\'ements dans $I(\tilde{G})$; on suppose que la fonction $f_{2}$ est  $\tilde{M}$ cuspidale et nulle pr\`es des \'el\'ements exceptionnels. 
Pour tout $f_{1}$ $\tilde{M}$-cuspidale et nulle pr\`es des \'el\'ements exceptionnels, on a: 
$$
I_{geo}^{\tilde{G}}(\omega,f_{1},f_{2})-I_{geo}^{\tilde{G},{\mathcal E}}(f_{1},f_{2})=0.
$$
Ainsi les hypoth\`eses du corollaire de \ref{redlocspec} sont satisfaites et on sait que
$$
I_{spec}^{\tilde{G}}(\omega,f_{1},f_{2})-I_{spec}^{\tilde{G},{\mathcal E}}(f_{1},f_{2})=0,
$$
sous la seule hypoth\`ese que $f_{2}$ est $\tilde{M}$ cuspidale et sans hypoth\`ese sur $f_{1}$. 
On r\'ecrit
$$
I_{geo}^{\tilde{G}}(\omega,f_{1},f_{2})-I_{geo}^{\tilde{G},{\mathcal E}}(f_{1},f_{2})=0,
$$avec l'hypoth\`ese que $f_{2}$ est  $\tilde{M}$-cuspidale. En utilisant \ref{redlocgeo} cela se traduit par:$$
\forall f_{1} \in I(\tilde{G}) \qquad 0=
$$
$$ \sum_{\tilde{L}}w(\tilde{L})(-1)^{a_{\tilde{L},\tilde{G}}}\int_{\tilde{L}_{ell}/\sim} \biggl(I_{\tilde{L}}^{\tilde{G}}(\gamma,\omega,f_{2})-I_{\tilde{L}}^{\tilde{G},{\mathcal E}}(\gamma,f_{2})\biggr)\overline{I^{\tilde{L}}(\gamma,\omega,f_{1,\tilde{L}}})\, d\gamma\eqno(1)$$
$$
+ w(\tilde{M})(-1)^{a_{\tilde{M},\tilde{G}}}\int_{\tilde{M}_{ell}/\sim} \overline{\biggl(I_{\tilde{M}}^{\tilde{G}}(\gamma,\omega,f_{1})-I_{\tilde{M}}^{\tilde{G},{\mathcal E}}(\gamma,f_{1})\biggr)}I^{\tilde{M}}(\gamma,\omega,f_{2,\tilde{M}})\, d\gamma \eqno(2)
$$
$$+w(\tilde{M})(-1)^{a_{\tilde{M},\tilde{G}}} \int_{\tilde{M}_{ell}/\sim}\overline{I^{\tilde{M}}(\gamma,\omega,f_{1,\tilde{M}})}
\biggr(I_{\tilde{M}}^{\tilde{G}}(\gamma,\omega,f_{2})-I_{\tilde{M}}^{\tilde{G},{\mathcal E}}(\gamma,f_{2})\biggr)\, d\gamma,\eqno(3)
$$
o\`u la premi\`ere somme ne porte que sur les espaces de Levi $\tilde{L}$ ne contenant pas $\tilde{M}$ \`a conjugaison pr\`es (pour ces espaces de Levi, par d\'efinition, $f_{2,\tilde{L}}=0$).

Rajoutons l'hypoth\`ese que $f_{1}$ est nulle pr\`es des \'el\'ements exceptionnels ce qui permet d'avoir la nullit\'e des termes (2); ainsi on obtient que pour toute fonction $f_{1}$ nulle pr\`es des \'el\'ements exceptionnels, on a
$$
\sum_{\tilde{L}}w(\tilde{L})(-1)^{a_{\tilde{L},\tilde{G}}}\int_{\tilde{L}_{ell}/\sim}\biggl(I_{\tilde{L}}^{\tilde{G}}(\gamma,\omega,f_{2})-I_{\tilde{L}}^{\tilde{G},{\mathcal E}}(\gamma,f_{2})\biggr)\overline{I^{\tilde{L}}(\gamma,\omega,f_{1,\tilde{L}}})\, d\gamma=0,
$$
o\`u la somme porte sur les espaces de Levi de $\tilde{G}$  ne contenant pas $\tilde{M}$ \`a conjugaison pr\`es  et $\tilde{M}$ lui-m\^eme.

En modifiant les coefficients, on peut ne sommer que sur les espaces de Levi pris \`a conjugaison pr\`es. Mais en faisant varier $f_{1}$, les fonctions $I^{\tilde{L}}(\gamma,\omega,f_{1,\tilde{L}})$ d\'efinie sur $\tilde{L}_{ell}$ sont soumises uniquement au fait qu'elles sont invariantes sous le normalisateur dans $G$ de $\tilde{L}$. Les fonctions $\biggl(I_{\tilde{L}}^{\tilde{G}}(\gamma,\omega,f_{2})-I_{\tilde{L}}^{\tilde{G},{\mathcal E}}(\gamma,f_{2})\biggr)$ ont la m\^eme propri\'et\'e d'invariance. Ainsi chaque terme est nul et les fonctions que l'on int\'egre sont nulles elles aussi au moins presque partout. Par continuit\'e elles sont nulles sur tous ces \'el\'ements. D'o\`u pour tout espace de Levi $\tilde{L}$ de $\tilde{G}$ et pour toute fonction $f$ que l'on suppose $\tilde{M}$-cuspidale: 
$$
I_{\tilde{L}}^{\tilde{G}}(\gamma,\omega,f)-I_{\tilde{L}}^{\tilde{G},{\mathcal E}}(\gamma,f)=0.\eqno(4)
$$
On reporte cette nullit\'e dans (1),(2), (3), o\`u on ne fait plus d'hypoth\`ese sur $f_{1}$ et on obtient alors la nullit\'e:
$$
\forall f_{1}\in I(\tilde{G}), \quad I_{\tilde{M}}^{\tilde{G}}(\gamma,\omega,f_{1})=I_{\tilde{M}}^{\tilde{G},{\mathcal E}}(\gamma,f_{1}).
$$
Ceci est la r\'eduction cherch\'ee. Au passage on a d\'emontr\'e que les hypoth\`eses de r\'ecurrence et l'hypoth\`ese cl\'e, entra\^{\i}nent la stabilisation locale g\'eom\'etrique pour toute fonction $f$ suppos\'ee $\tilde{M}$-cuspidale et pour tout espace de Levi $\tilde{L}$; mais cette propri\'et\'e ne sert pas dans la suite.

\subsubsection{D\'ebut de la preuve de la nullit\'e de  $\epsilon({\bf M}',\delta)$.\label{debutpreuveloc}}
On fait les hypoth\`eses locales de r\'ecurrence et l'hypoth\`ese cl\'e 1.
\begin{lem} Avec les notations pr\'ec\'edentes, $\epsilon({\bf M}',\delta)+\overline{\epsilon({\bf M}',\delta)}=0.$
\end{lem}
C'est la m\^eme preuve que \cite{partie3} 6.5.
On r\'ecrit \ref{redlocgeo} pour des fonctions $f_{1}$ et $f_{2}$ que l'on suppose $\tilde{M}$ cuspidales, nulles pr\`es des \'el\'ements isol\'es si  la place est $p$-adique. On suppose aussi que le terme constant $f_{1,\tilde{M}}$ est \`a support dans l'ensemble des \'el\'ements r\'eguliers de $\tilde{M}$. On voit que
$$\int_{\tilde{M}_{ell}/\sim} I^{\tilde{M}} (\gamma,\omega,\epsilon_{\tilde{M}}(f_{1}))\overline{f_{2,\tilde{M}}(\gamma)}+f_{1,\tilde{M}}(\gamma)\overline{I^{\tilde{M}} (\gamma,\omega,\epsilon_{\tilde{M}}(f_{2}))}\, d\gamma\eqno(1)
$$
vaut $I^{\tilde{G}}(\omega,f_{1},f_{2})-\sum_{{\bf G}'}i(\tilde{G},{\bf G}')I^{G'}(f_{1}^{\bf G'},f_{2}^{\bf G'})$. On r\'ecrit cela avec le c\^ot\'e spectral et (1) vaut donc:
$$
\sum_{\tilde{L}, \chi_{L},\tau,\tau'}\int_{i{\mathcal{A}}^*_{{\tilde{G},F}}}d\lambda\, b(\tilde{L},\chi_{L},\tau,\tau')tr\, \tau_{\lambda} f_{1,\tilde{L}} \overline{tr\tau_{\lambda}'\, f_{2,\tilde{L}}}.\eqno(2)
$$
Dans le terme de droite, seuls interviennent les espaces de Levi contenant $\tilde{M}$ puisque les fonctions consid\'er\'ees sont $\tilde{M}$ cuspidales. On montre par r\'ecurrence descendante que le terme de (2) index\'e par un espace de Levi $\tilde{L}$ contenant strictement $\tilde{M}$ est nul: en effet le terme de (2) index\'e par un espace de Levi $\tilde{L}$ est nul si et seulement si il est nul pour les fonctions $\tilde{L}$-cuspidales. Quand $\tilde{L}$ contient strictement $\tilde{M}$, on a d\'ej\`a d\'emontr\'e la stabilisation locale g\'eom\'etrique (hypoth\`ese de r\'ecurrence) et (1) est donc nul; c'est l'argument utilis\'e pour $\tilde{L}=\tilde{G}$ dans le premier corollaire de \ref{redlocspec}. D'o\`u l'assertion.

Ainsi (2) est r\'eduit au terme  index\'e par $\tilde{M}$ que l'on r\'ecrit:
$$
\int_{i{\mathcal{A}}^*_{{\tilde{G},F}}}d\lambda \sum_{\chi_{M}\tau',\tau''}\, b(\chi_{M},\tau',\tau'')tr\, \tau_{\lambda}' f_{1,\tilde{M}} \overline{tr
\tau_{\lambda}''\, f_{2,\tilde{M}}}, \eqno(3)
$$ o\`u la somme porte sur une base des repr\'esentations elliptiques de $\tilde{M}$ modulo l'action par tensorisation de $i{\mathcal{A}}^*_{{\tilde{G},F}}$ identifi\'e \`a des caract\`eres de $\tilde{M}$.

On revient \`a (1) que l'on interpr\`ete comme le terme g\'eom\'etrique de la formule des traces locale pour $\tilde{M}$ pour des fonctions cuspidales. On \'ecrit ce c\^ot\'e g\'eom\'etrique comme transfert, donc une somme sur les groupes endoscopiques elliptiques ${\bf M}'$ de $\tilde{M}$ des termes
$$
\int_{{\bf M}'_{ell}/\sim}d\delta\, \epsilon_{\tilde{M}}^{{\bf M}'}(f_{1})(\delta)\overline{f_{2,\tilde{M}}^{{\bf M}'}(\delta)}
$$
plus un terme o\`u on \'echange les r\^oles de $f_{1}$ et $f_{2}$, \`a la conjugaison complexe pr\`es. Avec l'hypoth\`ese cl\'e 1, cela se r\'ecrit comme somme sur les donn\'ees endoscopiques elliptiques ${\bf M}'$ de
$$
\int_{{\bf M}'_{ell}/\sim}d\delta (\epsilon({\bf M}',\delta)+ \overline{\epsilon({\bf M}',\delta)})f^{\bf M'}_{1,\tilde{M}}(\delta)\overline{f_{2,\tilde{M}}^{{\bf M}'}(\delta)}.
$$

Ici on a la forme stable de la formule des traces pour ${\bf M}'$ et la paire de fonction $((\epsilon({\bf M}',\delta)+ \overline{\epsilon({\bf M}',\delta)})f^{\bf M'}_{1,\tilde{M}}(\delta), f_{2,\tilde{M}}^{{\bf M}'}(\delta))$.
 On r\'ecrit (1) en utilisant la forme spectrale de cette formule des traces donn\'ee en \ref{locspecfin} (b). On se rappelle que les fonctions $f_{2,\tilde{M}}$ sont cuspidales et on obtient donc une int\'egrale sur $i{\mathcal{A}}^*_{\tilde{M},F}$ de traces de repr\'esentations n\'ecessairement elliptiques  de ${\bf M}'$ (on n'a pas besoin de les \'ecrire).  On compare (3) avec ce r\'esultat comme distributions en $f_{2,\tilde{M}}$. Ces deux distributions ne peuvent \^etre \'egales que si elles sont toutes les deux nulles.
\subsubsection{Fin de la preuve du th\'eor\`eme dans certains cas\label{casquadratique}}
En \cite{antipode} en particulier dans le paragraphe 2.4 de cette r\'ef\'erence, il a \'et\'e d\'efini une bijection entre les classes d'isomorphie de donn\'ees endoscopiques elliptiques pour $\tilde{M},\omega$ et celles pour $\tilde{M},\omega^{-1}$, not\'ee ${\bf M}' \mapsto {\bf M}'_{\nabla}$. A une classe de conjugaison stable $\delta$ de ${\bf M}'$ correspond une classe de conjugaison stable de ${\bf M}'_{\nabla}$ mais ces deux classes s'identifient naturellement (cf. \cite{antipode} 2.2 (2)) et on garde la notation $\delta$ pour la classe dans ${\bf M}'_{\nabla}$. Les facteurs de transfert sont modifi\'es par le passage au complexe conjugu\'e.
\begin{lem} Avec les notations du paragraphe pr\'ec\'edent, $\overline{\epsilon({\bf M}',\delta)}=\epsilon({\bf M}'_{\nabla},\delta)$.
\end{lem}
En \cite{antipode} 2.4, il est d\'emontr\'e que pour toute donn\'ee endoscopique elliptique ${\bf G}'$ de $\tilde{G},\omega$ et pour toute fonction $f\in I(\tilde{G})$, on a 
$$
\overline{f^{{\bf G}'}}=(\overline{f})^{{\bf G}'_{\nabla}}.
$$
Avec les propri\'et\'es des facteurs de transfert, on a imm\'ediatement, pour tout \'el\'ement $\gamma$ semi-simple r\'egulier de $\tilde{G}$:
$$
\overline{I_{\tilde{M},\omega}^{\tilde{G},{\mathcal E}}(\gamma,f)}=I_{\tilde{M},\omega^{-1}}^{\tilde{G},{\mathcal E}}(\gamma,\overline{f});
$$on rajoute dans la notation le caract\`ere puisqu'il varie.
$$
\overline{\epsilon_{\tilde{M}}(\gamma, \omega,f)}=I_{\tilde{M}}^{\tilde{G}}(\gamma,\omega^{-1},\overline{f})-I_{\tilde{M}, \omega^{-1}}^{\tilde{G},{\mathcal E}}(\gamma,\overline{f})={\epsilon_{\tilde{M}}(\gamma,\omega^{-1},\overline{f})};
$$
$$
\overline{\epsilon^{{\bf M}'}_{\tilde{M}}(\delta,f)}= \sum_{\gamma\in \tilde{M}}\overline{\Delta_{{\bf M}'}(\gamma,\delta)}I^{\tilde{M}}(\gamma,\omega^{-1},\epsilon_{\tilde{M},\omega^{-1}}(\overline{f}))=\epsilon^{{\bf M}'_{\nabla}}(\delta,\overline{f}).
$$
D'autre part, on a aussi $\overline{f_{\tilde{M}}^{{\bf M}'}}(\delta)=\overline{f_{\tilde{M}}}^{{\bf M}'_{\nabla}}(\delta)$ d'o\`u en r\'ecrivant la d\'efinition de $\epsilon({\bf M}',\delta)$ ou plut\^ot du conjugu\'e de ce nombre:
$$
\epsilon^{{\bf M}'_{\nabla}}(\delta,\overline{f})=\overline{\epsilon({\bf M}',\delta)} \overline{f_{\tilde{M}}}^{{\bf M}'_{\nabla}}(\delta)
$$
et en comparant \`a la  d\'efinition de $\epsilon({\bf M}'_{\nabla},\delta)$ on obtient le lemme.

\begin{cor} On suppose que pour toute donn\'ee endoscopique elliptique ${\bf M}'$ de $\tilde{M}$, on a ${\bf M}'={\bf M}'_{\nabla}$ (cela force $\omega^2=1$), alors   pour tout \'el\'ement semi-simple r\'egulier de $\tilde{M}$, on a $$I_{\tilde{M}}^{\tilde{G}}(\omega,\gamma,f)=I_{\tilde{M}}^{\tilde{G},{\mathcal E}}(\omega,\gamma,f).$$
\end{cor}
Avec les formules d'inversion des facteurs de transfert, et les r\'eductions d\'ej\`a faites,  il suffit de montrer que pour toute donn\'ee endoscopique elliptique de $\tilde{M}$ et pour tout $f\in I(\tilde{G})$  nulle pr\`es des \'el\'ements exceptionnels s'il y en a, la fonction $\epsilon_{\tilde{M}}^{{\bf M}'}(f)$ est identiquement nulle. Il suffit donc de montrer que la fonction $\epsilon({\bf M}',\delta )$ est identiquement nulle. Or on a vu en \ref{debutpreuveloc} que cette fonction est \`a valeurs dans l'ensemble des nombres complexes purement imaginaires. Le lemme pr\'ec\'edent assure que puisque ${\bf M}'={\bf M}'_{\nabla}$ elle est aussi \`a valeurs dans l'ensemble des nombres r\'eels. Elle est donc nulle.
\begin{rmq} L'hypoth\`ese du corollaire pr\'ec\'edent \'equivaut \`a ce que   les facteurs de transfert soient \`a valeurs r\'eelles (ce qui force $\omega ^2=1$). On peut d\'emontrer ce corollaire plus directement.
\end{rmq}
 
En effet, pour $\gamma,\gamma'$ des \'el\'ements semi-simples r\'eguliers de $\tilde{M}$, on pose
$$
a(\gamma,\gamma'):=c^{-1} \sum_{{\bf M}', \delta}\Delta_{{\bf M}'}(\delta,\gamma) (\Delta_{{\bf M}'}(\delta,\gamma'))^{-1} \epsilon({\bf M}',\delta),
$$
o\`u $c$ est le nombre d'\'el\'ements de l'ensemble de sommation. 
On calcule directement, en utilisant les inversions des facteurs de transfert:
$$
I_{\tilde{M}}^{\tilde{G}}(\gamma,\omega,f)-I_{\tilde{M}}^{\tilde{G},\mathcal E}(\gamma,f)=\sum_{\gamma'} a(\gamma,\gamma')I^{\tilde{M}}(\gamma',\omega, f_{\tilde{M}}),$$o\`u la somme sur $\gamma'$ porte sur des repr\'esentants  des classes de conjugaison d'\'el\'ements semi-simples dans $\tilde{M}$. 
Les nombres $a(\gamma,\gamma')$ sont purement imaginaires sous l'hypoth\`ese de la remarque puisque c'est le cas des nombres $\epsilon({\bf M}',\delta)$. On applique l'\'egalit\'e ci-dessus \`a $f$ remplac\'e par $\overline{f}$. Le terme de gauche est remplac\'e par son conjugu\'e puisque $\omega^2=1$ et $I^{\tilde{M}}(\gamma',\omega, f_{\tilde{M}})$ aussi. Donc on obtient, en appliquant d'abord l'\'egalit\'e \`a ${\overline{f}}$
$$
\overline{I_{\tilde{M}}^{\tilde{G}}(\gamma,\omega,f)-I_{\tilde{M}}^{\tilde{G},\mathcal E}(\gamma,f)}=\sum_{\gamma'} a(\gamma,\gamma')\overline{I^{\tilde{M}}(\gamma',\omega, f_{\tilde{M}})}
$$
puis en conjuguant simplement l'\'egalit\'e pour $f$
$$
\overline{I_{\tilde{M}}^{\tilde{G}}(\gamma,\omega,f)-I_{\tilde{M}}^{\tilde{G},\mathcal E}(\gamma,f)}=-\sum_{\gamma'} a(\gamma,\gamma')\overline{I^{\tilde{M}}(\gamma',\omega, f_{\tilde{M}})}.
$$
Cela force la nullit\'e du terme de gauche.

\begin{rmq} Les hypoth\`eses de la remarque pr\'ec\'edente sont satisfaites si $\omega =1$ pour $GL(n)$ tordu que ce soit dans la situation du changement de base ou de la torsion avec l'automorphisme $g\mapsto \, ^tg^{-1}$.
\end{rmq}
\subsection{Une construction uniforme d'extensions de corps de nom\-bres\label{extensiondecorps}}
\begin{lem}Ici $F$ est un corps global  et $S$ est un ensemble fini de places de $F$.   On fixe aussi un entier $d$ et pour tout $v\in S$, on fixe une extension galoisienne $K^v$ de $F_{v}$ de degr\'e inf\'erieure ou \'egale \`a $d$. Alors il existe un corps $K$ qui est une extension galoisienne  de $F$ de degr\'e inf\'erieur ou \'egal \`a $d!$ telle que   tout $v\in S$ et pour toute place $v'$ de $K$ au-dessus de $v$, on ait $K_{v'}=K^v$.
\end{lem}
Pour $v\in S$, on fixe un polyn\^ome de degr\'e $d$ \`a coefficients dans $F_{v}$ dont les racines sont distinctes (mais peuvent \^etre dans $F_{v}$) et engendrent $K^v$ sur $F_{v}$. On note $$
X^d+a_{1,v}X^{d-1}+\cdots +a_{d,v}$$
ces polyn\^omes. On fixe $\epsilon>0$. Par approximation forte, on fixe un polyn\^ome \`a coefficients dans $F$ 
$$
X^d+a_{1}X^{d-1}+ \cdots +a_{d}
$$
tel que pour tout $i\in [1,d]$ et pour tout $v\in S$, on ait $\vert a_{i}-a_{i,v}\vert_{v}\leq \epsilon$. Et on note $K$ l'extension de $F$ engendr\'ee par toutes les racines de ce polyn\^omes. Pour $\epsilon$ suffisamment petit, $K$ r\'epond aux conditions du lemme.
\subsection{Une r\'eduction \'etonnament simple \label{reduc}}
On revient momentan\'ement \`a la situation locale mais on suppose qu'elle vient d'une situation globale; on verra en \ref{globalisationgrossiere} et \ref{globalisation} que l'on peut approximer les situations locales et que ce n'est donc pas restrictif. On note donc $F$ le corps de nombres, $v_{0}$ la place qui nous int\'eresse. On fixe aussi $\tilde{M},\omega$ et ${\bf M}'$ une donn\'ee endoscopique elliptique de $\tilde{M},\omega$. On a aussi un \'el\'ement $\gamma\in \tilde{M}(F)$ un \'el\'ement rationnel que l'on suppose semi-simple et r\'egulier,  et un \'el\'ement $\delta\in \tilde{M}'(F)$ qui correspond \`a $\gamma$. On fixe un ensemble de places $V$ v\'erifiant les conditions de la deuxi\`eme hypoth\`ese cl\'e de 3.5 pour ${\bf M}'$ et $\delta$. 
\begin{lem} On suppose que $\tilde{M}$ n'est pas un tore et que les hypoth\`eses de \ref{stabgeohyp} sont satisfaites. Alors $\epsilon_{v_{0}}({\bf M}',\delta)=0$.
\end{lem}  
On reprend la notation ${\tilde{T}}$ pour un tore tordu maximal contenant $\gamma$. Si $\tilde{T}$ n'est pas elliptique dans $\tilde{M}$, par les formules de descente, en toute place $v$, $\epsilon_{v}({\bf M}',\delta)=0$ pour toute donn\'ee endoscopique elliptique ${\bf M}'$ et pour tout $\delta$ dont la classe stable correspond \`a $\gamma_{v}$ et il n'y a rien de non trivial \`a d\'emontrer. Puisque le lemme porte sur la place $v_{0}$, on suppose plus pr\'ecis\'ement que $\tilde{T}(F_{v_{0}})$ est un tore elliptique de $\tilde{M}(F_{v_{0}})$. On suppose comme dans l'\'enonc\'e que $\tilde{M}$ n'est pas un tore et on va conclure. On note $T$ le tore sous-jacent \`a $\tilde{T}$.

On note $d$ le degr\'e d'une extension galoisienne de $F$ qui d\'eploie ${T}$ et on note $S$ l'ensemble des places de $F$ qui contient $v_{0}$, $V$ et l'ensemble des places $v'$ de $F$ tel que l'on ait $p_{v'}\leq N(G)e_{v'}d!$ o\`u ces termes sont d\'efinis en \cite{stabilisationVI} 1.1 (c'est une condition technique). Pour tout $v\in S$ on fixe $K^v$: si $v=v_{0}$, $K^{v_{0}}=F_{v_{0}}$ et si $v\neq v_{0}$, $K^v$ est une extension galoisienne de $F_{v}$ qui d\'eploie $T_{v}$ et qui est de degr\'e inf\'erieur ou \'egal \`a $d$. Avec le lemme de \ref{extensiondecorps}, on fixe une extension galoisienne, $K$, de $F$ de degr\'e inf\'erieur ou \'egal \`a $d!$ et qui en toute place $v'$ de $K$ au-dessus d'une place $v$  de $F$ qui est dans $S$, $K_{v'}\simeq K^v$. On consid\`ere la situation sur $K$ au lieu de la consid\'erer sur $F$. Au dessus de $v_{0}$, on n'a pas chang\'e grand chose, on a d\'emultipli\'e la place $v_{0}$; en particulier ${\bf M}'$ d\'efinit naturellement une donn\'ee endoscopique de $\tilde{M}$ sur $K$ qui se localise en toute place au dessus de   $v_{0}$ en la donn\'ee endoscopique de d\'epart. Donc elle est elliptique. Le nouveau $V_{ram}$ est inclus dans $S_{K}$ o\`u $S_{K}$ est l'ensemble des places de $K$ situ\'ees au dessus d'une place de $S$:  gr\^ace \`a la d\'efinition de $S$, la condition technique qui doit \^etre satisfaite hors de $V_{ram}$ de \cite{stabilisationVI} 1.1, l'est gr\^ace \`a  la d\'efinition de $S$ et  les conditions de non ramification hors de $V_{ram}$ proviennent de celles qui sont d\'ej\`a v\'erifi\'ees avant le changement de corps. Mais en toute place $v'$ de $S_{K}$ qui n'est pas au-dessus de $v_{0}$, $\gamma_{v'}$ n'est plus elliptique. Donc on sait que  $\epsilon_{\tilde{M}}(f_{v'})$ est nulle sur la classe de conjugaison stable de $\gamma_{v'}$ et donc que  $\epsilon_{v'}({\bf M}',\delta)=0$. La deuxi\`eme hypoth\`ese cl\'e de \ref{stabgeohyp} donne alors simplement:
$$
D\epsilon_{v_{0}}({\bf M}',\delta)=0,\eqno(1)
$$o\`u $D$ est le degr\'e  de $ K$ sur $F$. 
Ceci est la nullit\'e cherch\'ee.
\begin{rmq} Pour \'eviter de faire deux fois la m\^eme preuve, on remarque que la d\'emonstration que l'on vient de faire donne aussi le r\'esultat suivant: supposons que l'on sache que,  pour tout corps local $F_{v'}$, $\epsilon_{\tilde{M}}( f_{v'})$ est nul (avec les notations ci-dessus) si $\tilde{M}$ est un tore d\'eploy\'e sur $F_{v'}$  , alors on a aussi $\epsilon_{v_{0}}({\bf M}',\delta)=0$ sans restriction sur $\tilde{M}$.
\end{rmq}
Gr\^ace au r\'esultat pr\'ec\'edent, on suppose que $\tilde{M}$ est un tore. 
Comme dans la preuve ci-dessus, on d\'eploie $\tilde{M}$ sur une extension galoisienne puis on fixe $K$ qui d\'eploie $\tilde{M}$ sur $S-\{v_{0}\}$ o\`u  $S$ est  comme ci-dessus et tel que, au dessus de $v_{0}$, l'extension de corps est totalement scind\'ee. Et on a encore (1)
qui donne le r\'esultat cherch\'e.
\subsection{Le cas des tores d\'eploy\'es\label{derreduc}}
Dans le cas de l'endoscopie ordinaire, le cas des tores d\'eploy\'es est trivial, il n'en est pas de m\^eme pour l'endoscopie tordue. C'est ici que l'on va utiliser \ref{casquadratique}. On rappelle qu'\`a toute donn\'ee endoscopique ${\bf M}'$ est attach\'e un caract\`ere $\omega_{\sharp}$ du groupe $M_{\sharp}(F):=(M/Z(M)^{\theta})(F)$ (\cite{stabilisationI} 2.7), cela se fait d'abord localement, puis comme ce caract\`ere est une propri\'et\'e de transformation des facteurs de transfert, ces caract\`eres locaux donnent un caract\`ere automorphe du groupe $M_{\sharp}$. Le caract\`ere $\omega_{\sharp}$ d\'epend de la donn\'ee endoscopique, on le note donc plut\^ot $\omega_{\sharp}^{{\bf M}'}$. La restriction de ce caract\`ere \`a l'image de $M(F)$ dans $M_{\sharp}(F)$ est $\omega$. Par la d\'efinition m\^eme, on voit que $$\omega_{\sharp}^{{\bf M}'_{\nabla}}=\overline{\omega_{\sharp}^{{\bf M}'}}=
(\omega_{\sharp}^{{\bf M}'})^{-1}.
$$
\begin{lem} La situation est locale et on suppose que $\tilde{M}$ est un tore d\'eploy\'e. L'application qui a une classe  d'isomorphie, ${\bf M}'$ de  donn\'ees endoscopiques elliptiques de $\tilde{M}$, $\omega$, associe le caract\`ere $\omega_{\sharp} ^{{\bf M}'}$ induit une bijection de cet ensemble de classe d'isomorphie sur l'ensemble des caract\`eres de $M_{\sharp}(F)$ dont la restriction \`a l'image de $M(F)$ est $\omega$. En particulier les donn\'ees endoscopiques elliptiques fixes par l'op\'eration $\nabla$ (cf. \ref{casquadratique}) sont exactement les donn\'ees endoscopiques pour lesquelles $(\omega_{\sharp}^{{\bf M}'})^2=1$ \end{lem}
Comme on suppose que $\tilde{M}$ est un tore d\'eploy\'e, tout se simplifie: une donn\'ee endoscopique est un triplet $( M',{\mathcal M}',\hat{s})$. L'\'el\'ement $\hat{s} \in \hat{M}$ est d\'etermin\'e modulo $Z(\hat{M})=\hat{M}$. Ainsi on prend $\hat{s}=\hat{\theta}$. Alors $\hat{M}'=\hat{M}^{\hat{\theta},0}$. Et ${\mathcal M}'$ est uniquement d\'etermin\'e par un homomorphisme  de $W_{F}$ dans $\hat{M}/\hat{M}^{\hat{\theta},0}$ qui se rel\`eve en  une application $w\in W_{F} \mapsto a_{w} \in \hat{M}$ v\'erifiant que l'application $w\mapsto a_{w}^{-1}\hat{\theta}(a_{w})$ est le cocycle de $W_{F}$ dans $\hat{M}$ donnant le caract\`ere $\omega$ de $M(F)$. Or $M_{\sharp}=M/M^\theta$ est lui aussi un tore; il s'identifie (comme groupe alg\'ebrique) \`a $(1-\theta)M$. Son groupe dual est alors $\hat{M}/\hat{M}^{\hat{\theta},0}$. Et l'application $w\mapsto a_{w}$ apr\`es passage au quotient par $\hat{M}^{\hat{\theta},0}$ est par d\'efinition un morphisme de groupe. Il est facile de voir que ce morphisme correspond exactement au caract\`ere $\omega_{\sharp}$ de la donn\'ee endoscopique. La r\'eciproque est claire: soit $w\mapsto a_{w}$ une application de $W_{F}$ dans $\hat{M}$ qui apr\`es passage  au quotient par $\hat{M}^{\hat{\theta},0}$ est un morphisme et qui est tel que $w \mapsto a_{w}^{-1}\hat{\theta}(a_{w})$ soit le cocycle correspondant \`a $\omega$. On d\'efinit alors ${\mathcal M}'$ comme le sous-groupe de $\hat{M}$ engendr\'e par $\hat{M}^{\hat{\theta},0}$ et les \'el\'ements $(a_{w},w)$ et c'est la donn\'ee endoscopique n\'ecessairement elliptique cherch\'ee.

\begin{cor} La situation est locale; on suppose que $\tilde{M}$ est un tore d\'eploy\'e et que $\omega_{\sharp}^{{\bf M}'}$ est un caract\`ere trivial. Alors $\epsilon({\bf M}',\delta)=0$.
\end{cor}
On v\'erifie que sous les hypoth\`eses faites, ${\bf M}'$ satisfait aux conditions du corollaire de \ref{casquadratique} gr\^ace au lemme pr\'ec\'edent.
\begin{lem} On suppose que $\omega$ est d'ordre fini, alors $\epsilon({\bf M}',\delta)=0$ pour toute donn\'ee endoscopique elliptique de $\tilde{M},\omega$.
\end{lem}
Le fait que $\omega$ soit d'ordre fini, entra\^{\i}ne le fait que pour tout ${\bf M}'$, le caract\`ere $\omega_{\sharp}^{{\bf M}'}$ est lui aussi d'ordre fini. Par exemple si la place est une place complexe, cela force $\omega_{\sharp}^{{\bf M}'}$ \`a \^etre trivial et on applique alors le corollaire pr\'ec\'edent. En g\'en\'eral on se ram\`ene au cas o\`u $\omega_{\sharp}^{{\bf M}'}$ est trivial. On se remet dans une situation globale, on fixe la place $v_{0}$ qui nous int\'eresse et $\omega$ est toujours d'ordre fini. On note $d$ le degr\'e d'une extension galoisienne  de $F$ qui trivialise $\omega$ et on proc\`ede comme dans \ref{reduc}, en construisant une extension galoisienne de $F$ telle qu'on ne change rien au-dessus de $v_{0}$, et en toute place  $v\not=v_{0}$  appartenant \`a l'ensemble  $S$ construit en 3.7, le caract\`ere $\omega_{\sharp}^{{\bf M}'}$ est trivial. On conclut alors comme dans \ref{reduc}.

\begin{rmq}On suppose qu'il existe un caract\`ere $\mu$ de $G(F)$ tel que $\omega(g)=\mu(g \theta(g)^{-1})$ o\`u $\theta$ est $ad_{\gamma_{0}}$ pour $\gamma_{0}\in \tilde{G}(F)$. Alors l'application $ f \mapsto \epsilon_{\tilde{M}}(\, ,\omega,f)$ est nulle et en particulier toutes les fonctions $\epsilon({\bf M}',\,)$ sont identiquement nulles.
\end{rmq}
Le choix de $\theta$ d\'epend du choix d'un \'el\'ement $\gamma_{0}$  de $\tilde{G}(F)$; fixons donc $\gamma_{0}\in \tilde{G}(F)$; on peut changer de $\theta$, cela ne change pas $\mu$. On d\'efinit la fonction $\tilde{\mu}$ sur $\tilde{G}$ en posant $\tilde{\mu}(\gamma)=\mu(g_{0})$ o\`u $g_{0}$ est l'unique \'el\'ement de $G(F)$ tel que $g_{0}\gamma_{0}=\gamma$. Evidemment pour tout $g,g'\in G(F)$, on a: $$\tilde{\mu}(g\gamma g')={\mu}(gg_{0}\theta(g'))=\mu(g) \mu(\theta(g)')\tilde{\mu}(\gamma).$$On associe ainsi \`a tout \'el\'ement, $f$,  de $I(\tilde{G})$ un \'el\'ement $f_{\mu}:=\tilde{\mu}f$ de $I(\tilde{G})$. Les $\omega$-int\'egrales orbitales pond\'er\'ees de $f$ deviennent des int\'egrales orbitales pond\'er\'ees de $f_{\mu}$:
$$
I_{\tilde{M}}^{\tilde{G}}(\omega,\gamma,f)=\tilde{\mu}(\gamma) I_{\tilde{M}}^{\tilde{G}}(\gamma,f_{\mu}).$$
 Le point \`a v\'erifier est donc que $$
I_{\tilde{M}}^{\tilde{G},{\mathcal E}}(\omega,\gamma,f)=\tilde{\mu}(\gamma)I_{\tilde{M}}^{\tilde{G},{\mathcal E}}(\gamma,f_{\mu}).
$$
On v\'erifie ais\'ement, qu'en notant $a_{\mu}()$ un cocycle \`a valeurs dans $Z(\hat{G})$ d\'efinissant $\mu$, il existe une bijection entre les classes d'isomorphie de donn\'ees endoscopiques elliptiques pour $\tilde{G}$ et le cocycle fix\'e d\'efinissant le caract\`ere $\omega$ de $G$ et les donn\'ees endoscopiques elliptiques de $\tilde{G}$ sans caract\`ere et cette bijection est donn\'ee simplement en modifiant ${\mathcal G}'$ en ${\mathcal H}'$ o\`u si $(x(w),w)\in {\mathcal G}'$ pour $w\in W_{F}$ alors $(a_{\mu}(w)x(w),w)\in {\mathcal H}'$. Notons ${\bf G}' \mapsto {\bf H}'$ cette correspondance; il faut montrer que les facteurs de transfert se d\'eduisent par multiplication par ${\mu}(g_{0})$ o\`u $g_{0}$ relie $\gamma$ au point de base choisi pour normaliser les facteurs de transfert. Pr\'ecis\'ement:
$$
\Delta^{\tilde{G},\omega,{\bf G}'}(\gamma,\delta)=\mu(g_{0})\Delta^{\tilde{G},{\bf H}'}(\gamma,\delta).
$$
On revient aux d\'efinitions rappel\'ees en \cite{stabilisationI} 2.2. Ce qui change est le  cocycle $w \mapsto \hat{V}_{1}(w)$ qui est multipli\'e par $a_{\mu}(w)^{-1}$ (cela intervient dans la d\'efinition de $t_{T}(w)$). Donc $\Delta_{imp}$ a bien la propri\'et\'e de transformation \'ecrite.

Avec cela, on en d\'eduit que $f^{{\bf G}'}=(f_{\tilde{\mu}})^{{\bf H}'}$ en regardant le transfert des int\'egrales orbitales des \'el\'ements semi-simples r\'eguliers. D'o\`u l'assertion cherch\'ee.
 Ce raisonnement prouve la remarque car il nous ram\`ene au cas $\omega=1$ qui est trait\'e par le lemme pr\'ec\'edent.
\subsection{Fin des r\'eductions}
\begin{prop} On fait les hypoth\`eses de \ref{stabgeohyp} et alors l'application $$f \mapsto \epsilon_{\tilde{M}}(f)$$ est identiquement nulle.
\end{prop}
Il ne nous reste plus qu'\`a faire le cas o\`u $\tilde{M}$ est un tore d\'eploy\'e et o\`u $\omega$ est un caract\`ere  qui n'est pas de la forme $(1-\theta)\mu$. On suppose donc que $\tilde{M}$ est un tore d\'eploy\'e et il n'y a rien \`a d\'emontrer si $\omega$ n'est pas trivial sur $M(F)^\theta$ o\`u $\theta$ est l'automorphisme obtenu par conjugaison sous n'importe quel \'el\'ement de $\tilde{M}$.

On globalise; si $v_{0}$ est une place complexe on globalise avec une extension totalement imaginaire $F$ de ${\mathbb Q}$.  D'abord par une extension convenable, $K'$ de $F$, on se ram\`ene au cas o\`u $\omega \circ Norm_{K'/F}$ est non ramifi\'e en toutes les places finies. On note $d$ le degr\'e de cette extension. On remarque que pour toute place $v$ finie de $K'$, si un caract\`ere de $M(K')$ est non ramifi\'e, il est uniquement d\'etermin\'e par l'image du Frobenius, $Fr_{v}$ dans $Z(\hat{M})$. On note $a_{v,Fr,\omega}$ l'\'el\'ement de $\hat{M}$ correspondant \`a cet \'el\'ement. Comme $\omega$ est la restriction des caract\`eres $\omega_{\sharp}$ d\'ecrit en  \ref{derreduc}, $a_{v,Fr,\omega}$ est n\'ecessairement dans $(1-\hat{\theta})\hat{M}$. (Cela traduit le fait que $\omega$ est nul sur $M(F)^{\theta}$). On fixe $z_{v,Fr}\in \hat{M}$ tel que $a_{v,Fr,\omega}=z_{v,Fr} \hat{\theta}(z_{v,Fr})^{-1}$. On note   $\mu_{v}$ le caract\`ere non ramifi\'e pour lequel qui $z_{v,Fr}$ est l'image  du Frobenius (via la th\'eorie du corps de classe). Alors  $(\omega\circ Norm_{K'/F})_{v}= \mu_{v}( \mu_{v}\circ \theta)^{-1}$. Par la m\'ethode de \ref{reduc}   et gr\^ace \`a la remarque de 3.8 on se d\'ebarrasse ainsi des places finies pour obtenir le r\'esultat aux places complexes. Ensuite, puisqu'on a le r\'esultat aux places complexes, on peut encore appliquer la m\'ethode de \ref{reduc}  et le m\^eme raisonnement que ci-dessus pour n'avoir que la place $v_{0}$ qui nous int\'eresse et on conclut.

\section{Les caract\`eres pond\'er\'es $\omega$-\'equivariants et leur stabilisation}
Comme pour les int\'egrales orbitales pond\'er\'ees invariantes, on a une d\'efinition locale des caract\`eres pond\'er\'es invariants et une d\'efinition semi-globale de ces caract\`eres. La d\'efinition semi-globale se ramenant \`a la d\'efinition locale par une formule de scindage. Avant de revenir sur ces d\'efinitions, il faut consid\'erer les places non ramifi\'ees hors de $V$, celles o\`u l'on ne rend pas les objets invariants. Donc ici $F$ est un corps de nombres.
\subsection{Caract\`ere pond\'er\'e  aux places non ramifi\'ees et stabilisation}
On fixe $V$ un ensemble fini de places de $F$ contenant $V_{ram}$ (cf. \cite{stabilisationVI} 1.1 pour la d\'efinition de $V_{ram}$). On fixe pour $v\notin V_{ram}$, $\tilde{K}_{v}$ un sous-espace hypersp\'ecial (cf. \cite{stabilisationI} 6.1 pour la d\'efinition de $\tilde{K}_{v}$). On note $1_{\tilde{K}^V}$ la fonction caract\'eristique du sous-espace compact $\tilde{K}^V$ de $\tilde{G}(F^V)$.  La notion d'invariance n'a pas de  sens dans cette situation, par contre on peut calculer $J_{\tilde{M}}^{\tilde{G}}(\pi^V,1_{\tilde{K}^V})$ et stabiliser cette formule. C'est fait en \cite{fonctionLendoscopique} que l'on reprend ici pour adapter les coefficients \`a la situation tordue. On utilise la classification de Langlands des repr\'esentations non ramifi\'ees des groupes p-adiques. Ainsi les caract\`eres pond\'er\'es pour les repr\'esentations non ramifi\'ees et pour  l'unit\'e de l'alg\`ebre de Hecke sph\'erique hors de $V$ s'expriment en termes de fonctions $L$. On va les calculer.

\subsubsection{Rappel\label{nonconnexite}}
On rappelle le lemme suivant qui s'inspire fortement du lemme 2 de  {\cite{fonctionLendoscopique}} pour pouvoir y r\'ef\'erer facilement.
Soit  $\tilde{M}$ un sous-espace de Levi de $\tilde{G}$. On dispose des $L$-groupes, $^LG,^LM$ et de l'automorphisme $\hat{\theta}$ et des composantes neutres $\hat{M}$ et $\hat{G}$ des $L$-groupes. On fixe aussi une donn\'ee endoscopique elliptique ${\bf G}'$ de $\tilde{G},\omega$ et ${\bf M}'$ de $\tilde{M},\omega$. On suppose que ${\bf M}'$ est un espace de Levi de ${\bf G}'$. Alors:

\begin{lem} On a l'\'egalit\'e $Z(\hat{M})^{\hat{\theta}}=Z(\hat{G})^{\theta}(Z(\hat{M})^{\hat{\theta}})^0$. Et aussi
$Z(\hat{M})^{\Gamma,\hat{\theta}}=Z(\hat{G})^{\Gamma,\hat{\theta}}(Z(\hat{M}))^{\Gamma,\hat{\theta},0}$
\end{lem}
\begin{cor} On a $Z(\hat{G})\cap Z(\hat{M})^{\Gamma,\hat{\theta},0}\subset Z(\hat{G}')$ et $$i_{\tilde{M}}(\tilde{G},{\bf G}')=\vert Z(\hat{G}')\cap Z(\hat{M})^{\Gamma,\hat{\theta},0}/Z(\hat{G})\cap Z(\hat{M})^{\Gamma,\hat{\theta},0}\vert^{-1}.
$$
\end{cor}
La premi\`ere assertion est claire car $Z(\hat{M})^{\Gamma,\hat{\theta},0}=Z(\hat{M}')^{\Gamma,0}$ par ellipticit\'e. D'apr\`es le lemme pr\'ec\'edent, l'application naturelle
$$
Z(\hat{G}')^{\Gamma}/Z(\hat{G}')^{\Gamma}\cap Z(\hat{G})\rightarrow Z(\hat{M}')^{\Gamma}/Z(\hat{M}')^{\Gamma}\cap Z(\hat{M})$$
est surjective. Et $i_{\tilde{M}}(\tilde{G},\tilde{G}')$ est l'inverse du  nombre d'\'el\'ements du noyau. Or un \'element $z'\in Z(\hat{G}')^{\Gamma}$ est dans le noyau si son image dans $Z(\hat{G}')/Z(\hat{G}')\cap Z(\hat{M})^{\Gamma,\hat{\theta}}$ est nulle. C'est-\`a-dire $z'=zm$ avec $z\in Z(\hat{G})^{\Gamma,\hat{\theta}}$ et $m\in Z(\hat{M})^{\Gamma,\hat{\theta},0}$. Ainsi $m\in \hat{G}'$ et $z\in Z(\hat{G})\cap \hat{G'}$ et $z\in Z(\hat{G}')$, plus pr\'ecis\'ement $z\in Z(\hat{G}')^{\Gamma} \cap Z(\hat{G})$. Ainsi le noyau est exactement isomorphe \`a
$$Z(\hat{G}')\cap Z(\hat{M})^{\Gamma,\hat{\theta},0}/Z(\hat{G})\cap Z(\hat{M})^{\Gamma,\hat{\theta},0}.$$
Cela termine la preuve.

\subsubsection{Factorisation des facteurs $L$\label{factorisation}}
On fixe un sous-espace de Levi $\tilde{M}$ de $\tilde{G}$ et on fixe $\alpha$ un caract\`ere, non trivial, de $Z(\hat{M})^{\Gamma,\hat{\theta},0}$. On suppose que ce caract\`ere intervient dans l'action de $Z(\hat{M})^{\hat{\theta},0}$ dans l'alg\`ebre de Lie de $\hat{G}$; sinon la fonction attach\'ee \`a un tel $\alpha$ est, par d\'efinition, la fonction identiquement 1. Avec cette hypoth\`ese, il est clair qu'il existe un unique prolongement de $\alpha$ \`a $Z(\hat{M})^{\Gamma,\hat{\theta}}$ trivial sur $Z(\hat{G})^{\Gamma,\hat{\theta}}$. On peut donc voir indif\'eremment $\alpha$ comme un caract\`ere de $Z(\hat{M})^{\Gamma,\hat{\theta}}$ trivial sur  $Z(\hat{G})^{\Gamma,\hat{\theta}}$ ou de la composante neutre de ce groupe.

On note $\hat{{\mathfrak g} }[{\alpha}]$  l'espace propre de valeur propre $\alpha$ pour l'action de $Z(\hat{M})^{\Gamma,\hat{\theta}}$ dans l'alg\`ebre de lie de $\hat{G}$. On remarque que $\hat{M}$ et plus g\'en\'eralement le $L$-groupe de $M$, agit trivialement par conjugaison dans $Z(\hat{M})^{\Gamma,\hat{\theta}}$ et agit donc par conjugaison dans $\hat{\mathfrak g}[\alpha]$.

Soit $v$ une place de $F$ et  $c_{v}$ un homomorphisme de $W_{F_{v}}$ dans le groupe dual de $M$. On a le facteur $L$, $L_{v}(\alpha,c_{v},s)$ qui correspond \`a la repr\'esentation de $W_{F_{v}}$ dans $\hat{\mathfrak{g}}[\alpha]$ via le compos\'e de $c_{v}$ avec la repr\'esentation adjointe de $^L{M}$ dans $\hat{\mathfrak{g}}[\alpha]$. Il est clair que ce facteur $L$ est un produit fini d'objets d\'efinis en \cite{fonctionLendoscopique}, pr\'ecis\'ement le produit des $L_{v}(\beta,c_{v},s)$ o\`u $\beta$ parcourt l'ensemble des caract\`eres de $Z(\hat{M})^{\Gamma}$ ayant $\alpha$ comme restriction \`a $Z(\hat{M})^{\Gamma,\hat{\theta}}$. On fixe encore  $V$ un ensemble fini de places de $V$ contenant $V_{ram}$ et on suppose que $c_{v}$ est donn\'e pour tout $v\notin V$; on note $c^V$ le produit de ces $c_{v}$. On d\'efinit alors $L^V(\alpha,c^V,s)$ en faisant le produit sur toutes les places $v$ en tant que fonction m\'eromorphe de $s$: ceci est bien d\'efini sous des hypoth\`eses assez faibles sur $c^V$, puisqu'il suffit que ce soit vrai pour les fonctions $L^V(\beta,c^V,s)$, on peut donc prendre les hypoth\`eses de \cite{fonctionLendoscopique} qui suivent la formule (1) du paragraphe 2. Pour nous, il suffit par exemple que $c^V$ soit la composante non ramifi\'ee hors de $V$ d'une repr\'esentation automorphe de $\tilde{M}$ puisque cela donne aussi une repr\'esentation automorphe de $M$.
Dans la suite, on aura besoin d'une hypoth\`ese plus faible que $c^V$ soit la composante non ramifi\'ee hors de $V$ d'une repr\'esentation  d'une donn\'ee endoscopique elliptique ${\bf M}'$ de $\tilde{M}, \omega$; en fait ce sont les repr\'esentations des donn\'ees auxiliaires attach\'ees \`a une donn\'ee endoscopique elliptique qui interviennent. Les repr\'esentations consid\'er\'ees ont alors en restriction au   tore induit central celui d\'etermin\'e un caract\`ere fix\'e  par la donn\'ee endoscopique et le choix de la donn\'ee auxiliaire (cf \cite{fonctionLendoscopique} paragraphe 2, le m\^eme ph\'enom\`ene existant dans le cas non tordu).  On dira pour simplifier qu'une telle repr\'esentation  est une repr\'esentation de ${\bf M}'$. C'est  la factorisation que l'on va pr\'ecis\'ement d\'emontrer  ci-dessous en suivant la m\'ethode de \cite{fonctionLendoscopique} qui montre que ces fonctions $L$ partielles sont bien  d\'efinies comme fonctions m\'eromorphes.

\

Soit ${\bf M}'=:(M', {\mathcal M}',\hat{s}_{M'})$ un espace endoscopique elliptique de $\tilde{M}$; on suppose que $c^V$ est \`a valeurs dans ${\mathcal M}'$ et que $c^V$ correspond au caract\`ere des alg\`ebres de Hecke sph\'eriques agissant sur une repr\'esentation automorphe de ${\bf M}'$  non ramifi\'ee hors de $V$.

On \'ecrit $\hat{s}_{M'}=s_{M'}\hat{\theta}$. Soit $s\in s_{M'}Z(\hat{M}^{\Gamma,\hat{\theta}})/Z(\hat{G}^{\Gamma,\hat{\theta}})$ on associe un espace endoscopique ${\bf G}'=(G'(\hat{s}),{\mathcal G}',\hat{s}=s\hat{\theta})$ de $\tilde{G}$ o\`u tout simplement, $\hat{G'}:=\hat{G}'(\hat{s})$ est la composante neutre du centralisateur de $\hat{s}$ dans $\hat{G}$,  ${\mathcal G}'=\hat{G}'{\mathcal M}'$, ce qui donne une action de $\Gamma_{F}$ sur $\hat{G}'$ et un groupe quasid\'eploy\'e $G'$.
Avec l'\'egalit\'e 
$
(Z(\hat{M}')^{\Gamma})^0=(Z(\hat{M})^{\Gamma,\hat{\theta}})^{0},
$
on voit $\alpha$ comme un caract\`ere de $(Z(\hat{M'})^{\Gamma})^0$ et on le note alors $\alpha'$. On d\'efinit alors  $L_{G'(\hat{s})}(\alpha',c^V,s)$ qui est une fonction m\'eromorphe de $s$ comme on l'a rappel\'e ci-dessus (avec un double emploi de la lettre $s$ qu'on esp\`ere compr\'ehensible). On remarque que $\hat{\mathfrak{g}}'[\alpha']=\hat{\mathfrak{g}}[\alpha]\cap \hat{\mathfrak{g}'}$ et que cet espace peut \^etre vide, la fonction est alors \'egale \`a 1. Clairement $\hat{\mathfrak{g}}'$ est l'alg\`ebre de Lie du centralisateur de $\hat{s}$ agissant par conjugaison dans $\hat{\mathfrak{g}}$.
 Soit $z\in Ker\, \alpha$ et consid\'erons la donn\'ee endoscopique associ\'ee \`a $\hat{s}z$. On a:
$\hat{\mathfrak{g}}[\alpha]\cap Cent_{\hat{\mathfrak{g}}}s \subset \hat{\mathfrak{g}}[\alpha]\cap Cent_{\hat{\mathfrak{g}}}\hat{s}z$ puisque $z$ commute \`a $\hat{\mathfrak{g}}[\alpha]$ et par la m\^eme raison l'inclusion inverse est vraie, d'o\`u l'\'egalit\'e 
$L^V_{G'(\hat{s})}(\alpha',c^V,s)= L^V_{G(sz)}(\alpha',c^V,s)$.

\begin{prop} $L^V_{G}(\alpha,c^V,s)=\prod_{\hat{s}\in \hat{s}_{M'}Z(\hat{M})^{\Gamma,\hat{\theta}}/Z(\hat{G})^{\Gamma,\hat{\theta}}Ker\, \alpha }L^V_{G'(\hat{s})}(\alpha',c^V,s)$ et le membre de gauche est une fonction m\'eromorphe de $s$.
\end{prop}
Comme $\hat{s}_{M'}$ commute aux \'el\'ements de $Z(\hat{M})^{\Gamma,\hat{\theta}}$, il agit par conjugaison dans $\hat{\mathfrak{g}}[\alpha]$ et on   d\'ecompose $\hat{\mathfrak{g}}[\alpha]$ sous l'action de $\hat{s}_{M'}$. L'\'el\'ement $\hat{s}_{M'}$  op\`ere donc de fa\c{c}on semi-simple dans cet espace vectoriel. Ainsi, 
$$
\hat{\mathfrak{g}}[\alpha]=\oplus_{x} \hat{\mathfrak{g}}[\alpha,x], \eqno(1)
$$
o\`u  $x$ parcourt l'ensemble des valeurs propres de l'action de  $\hat{s}_{M'}$; en particulier $x\neq 0$ puisque ${\hat{s}}_{M'}$ est un \'el\'ement inversible. Pour toute valeur propre $x$ intervenant, on fixe $t_{x}\in Z(\hat{M})^{\Gamma,\hat{\theta},0}$ tel que $\alpha(t_{x})=x^{-1}$ et on pose $\hat{s}_{x}:=\hat{s}_{M'}t_{x}$; \'evidemment, $t_{x}$ est bien d\'etermin\'e modulo $Ker\, \alpha$. On note ${\bf G}'(\hat{s}_{x})$ l'espace endoscopique correspondant \`a $\hat{s}_{x}$ et ${\bf M}'$. On remarque que $\hat{s}_{x}$ commute aux \'el\'ements de $\hat{\mathfrak{g}}[\alpha,x]$ et  cet espace est donc inclus dans $\hat{\mathfrak{g}}'(\hat{s}_{x})$, d'o\`u
$$
\hat{\mathfrak{g}}[\alpha,x]\subset \hat{\mathfrak{g}}'(\hat{s}_{x})[\alpha'].
$$
L'inclusion inverse est tout aussi claire d'o\`u  $\hat{\mathfrak{g}}[\alpha,x]=\hat{\mathfrak{g}}'(\hat{s}_{x})[\alpha']$.

 R\'eciproquement soit $\hat{s}\in \hat{s}_{M'}Z(\hat{M})^{\Gamma,\hat{\theta}}$, c'est-\`a-dire $\hat{s}=\hat{s}_{M'}z$ avec $z\in Z(\hat{M})^{\Gamma,\hat{\theta}}$.  On suppose que $\hat{\mathfrak{g}}'(\hat{s})[\alpha']\neq 0$; cela veut dire qu'il existe $v\in  \hat{\mathfrak{g}}'(\hat{s})[\alpha']$ non nul commutant \`a $\hat{s}$ ou avec ce qui pr\'ec\`ede appartenant \`a $\hat{\mathfrak{g}}[\alpha,\alpha(z)^{-1}]$. Ainsi $\alpha(z)^{-1}$ est l'un des $x$ consid\'er\'es ci-dessus et $z\in t_{x} Ker\, \alpha$. Ainsi (1) se r\'ecrit
$$
\hat{\mathfrak{g}}[\alpha]=\oplus_{\hat{s}\in \hat{s}_{M'}Z(\hat{M})^{\Gamma,\hat{\theta}}/Ker\, \alpha}\hat{\mathfrak{g}}[\alpha]\cap \hat{\mathfrak{g}}'(\hat{s}).
$$
Et le lemme r\'esulte de cette d\'ecomposition par d\'efinition des facteurs $L$. C'est exactement la d\'emonstration de \cite{fonctionLendoscopique} que l'on a recopi\'ee.

\subsubsection{D\'efinition des caract\`eres pond\'er\'es non ramifi\'es dans le cas d'un espace de Levi maximal\label{caracpondnonramifdef}}
Soit $\lambda\in {\mathcal{A}}_{\tilde{M},{\mathbb C}}^*$. Ici on consid\`ere les \'el\'ements invariants sous l'action de $\tilde{M}$. On identifie, ${\mathcal{A}}_{\tilde{M},{\mathbb C}}^*$ \`a un sous-groupe des homomorphismes non ramifi\'es de $W_{F^V}$ dans $Z(\hat{M})^{\Gamma,\hat{\theta},0}$. Pour tout $\lambda$ on tensorise $c^V$ (notations du paragraphe pr\'ec\'edent) par $\lambda$ et on note $c^V_{\lambda}$ le r\'esultat. On remarque alors    que si  $c^V$ est \`a valeurs dans ${\mathcal M}'$ comme pr\'ec\'edemment, il en est de m\^eme de $c^V_{\lambda}$. Pour $\alpha$ comme dans le paragraphe pr\'ec\'edent, on d\'efinit $L^V(\alpha,c^V,\lambda)$ ou $L(\alpha,c^V,\lambda)$ (le $V$ \'etant d\'ej\`a dans la notation $c^V$). On remarque que pour $s$ fix\'e, cette fonction de $\lambda$ est m\'eromorphe non identiquement nulle: en effet en voyant $\lambda$ comme un homomorphisme de $W_{F^V}$ dans $Z(\hat{M})^{\Gamma,\hat{\theta},0}$, on consid\`ere  $\alpha\circ \lambda$  comme un homomorphisme de $W_{F^V}$ dans ${\mathbb C}^*$, non ramifi\'e. On voit $d\alpha$ comme une forme lin\'eaire sur ${\mathcal{A}}^*_{\tilde{M},{\mathbb C}}$ et $d\alpha(\lambda)\in {\mathbb C}$ est tel que le caract\`ere pr\'ec\'edent soit $\vert\, \vert_{W_{F}}^{d\alpha(\lambda)}$. On a ainsi$L(\alpha,c^V_{\lambda},s)=L(\alpha,c^V,d\alpha({\lambda})+s)$. 
On pose alors
$$
r(\alpha,c^V,\lambda):=L(\alpha,c^V_{\lambda},0)/L(\alpha,c^V_{\lambda},1).
$$
Et sous l'hypoth\`ese que $c^V$ est \`a valeurs dans ${\mathcal M}'$, on a   encore une factorisation:
$$
r(\alpha,c^V,\lambda)=\prod_{\hat{s}\in \hat{s}_{M'}Z(\hat{M})^{\Gamma,\hat{\theta},0}/Ker(\alpha)}
r^{G(s)}(\alpha',c^V,\lambda).
$$
On pose donc, pour tout $\lambda,\Lambda \in {\mathcal{A}}^*_{\tilde{M},{\mathbb C}}$
$$
r(\alpha,c^V,\Lambda,\lambda):=$$
$$ r(\alpha,c^V,\lambda+1/2\Lambda)^{-1}r(\alpha,c^V,\lambda+\Lambda)r(-\alpha,c^V,\lambda+1/2\Lambda)^{-1}r(-\alpha,c^V,\lambda+\Lambda),
$$
o\`u $-\alpha$ est la notation additive pour le caract\`ere inverse de $\alpha$ qui jouit des m\^emes propri\'et\'es que $\alpha$.

La fonction en $\Lambda$ ci-dessus est, comme on l'a vu en \ref{factorisation} une fonction de $d\alpha(\Lambda)$; on \'ecrit alors plus simplement
$\tilde{r}^V(\alpha,c^V_{\lambda},s):=r^V(\alpha,c^V,\Lambda,\lambda)$ o\`u $s=d\alpha(\Lambda)$. On dit que $c^V$ est unitaire si chaque repr\'esentation correspondant \`a $c_{v}$ pour $v\notin V$ est unitaire.
\begin{lem} On suppose que $c^V$ est unitaire. La d\'eriv\'ee en $s=0$ de $\tilde{r}^V(\alpha,c^V_{\lambda},s)$ est holomorphe en $\lambda$ au voisinage de tout point $\lambda_{0}$ imaginaire.
\end{lem}

On \'ecrit $\tilde{r}(\alpha,c^V_{\lambda},s)$ comme fonction m\'eromorphe de $s$ et de $\lambda$. On  calcule sa d\'eriv\'ee logarithmique   en $s$ et on obtient:
$$\tilde{r}'(\alpha,c^V_{\lambda},s)\tilde{r}(\alpha,c^V_{\lambda},s)^{-1}=-1/2 r'(\alpha,c^V_{\lambda},s/2)r(\alpha,c^V_{\lambda},s/2)^{-1} $$
$$ +
r'(\alpha,c^V_{\lambda},s)r(\alpha,c^V_{\lambda},s)^{-1}-1/2r'(-\alpha,c^V_{\lambda},s/2)r(-\alpha,c^V_{\lambda},s/2)^{-1} $$
$$ +
r'(-\alpha,c^V_{\lambda},s)r(-\alpha,c^V_{\lambda},s)^{-1} $$ 
On fait $s=0$ et on obtient $$\tilde{r}'(\alpha,c^V_{\lambda},0)=1/2\biggl(r'(\alpha,c^V_{\lambda},0)/r(\alpha,c^V_{\lambda},0)+r'(-\alpha,c^V_{\lambda},0)/r(-\alpha,c^V_{\lambda},0)\biggr)\eqno(1)$$
ce qui fort heureusement est une formule sym\'etrique en $\alpha$ et $-\alpha$. On fixe $\lambda_{0}$ imaginaire et on pose $\lambda=\lambda_{0}+\mu$. On a alors
$$
r'(\alpha,c^V_{\lambda},0)/r(\alpha,c^V_{\lambda},0)=r'(\alpha,c^V_{\lambda_{0}},d\alpha(\mu))/r(\alpha,c^V_{\lambda_{0}},d\alpha(\mu));$$
$$ r'(-\alpha,c^V_{\lambda},0)/r(-\alpha,c^V_{\lambda},0)=r'(-\alpha,c^V_{\lambda_{0}},-d\alpha(\mu))/r(-\alpha,c^V_{\lambda_{0}},-d\alpha(\mu))$$
Pour \'eviter les confusions, on pose $\sigma:=d\alpha(\mu)$ et on d\'eveloppe la fonction \'ecrite en (1) au voisinage de $\sigma=0$. Pour $\eta=\pm$, on note $d_{\eta}$ l'ordre en $\sigma=0$ de la fonction $r(\eta\alpha,c^V_{\lambda_{0}},s)$ et on voit que (1) est la somme de fonctions holomorphes en $\sigma=0$ avec $(d_{+}-d_{-})\sigma^{-1}$. Le lemme est donc \'equivalent \`a montrer que $d_{+}=d_{-}$.

Comme on suppose que $\phi$ est unitaire et que  $\lambda_{0}$ est imaginaire, on peut remplacer $c^V$ par $c^V_{\lambda_{0}}$ et oublier $\lambda_{0}$.
On note $c^{V*}$ le morphisme correspondant \`a la repr\'esentation contragr\'ediente de la repr\'esentation d\'efinie par $c^V$. On v\'erifie sur les d\'efinitions que $r(-\alpha,c^V,s)=r(\alpha,c^{V*}-s)$; cela vient du fait que l'espace propre pour la valeur propre $\alpha$ pour l'action de $Z(\hat{M})^{\Gamma,\hat{\theta},0}$ est en dualit\'e avec son analogue pour la valeur propre $-\alpha$. En utilisant l'unitarit\'e, on a encore $\overline{r(\alpha,c^{V*}-s)}=r(\alpha,c^V,-\overline{s})$. D'o\`u $r(-\alpha,c^V,s)=\overline{r(\alpha,c^V,-\overline{s})}$. Ainsi l'ordre en $s=0$ de la fonction $r(-\alpha,c^V,s)$ est \'egal \`a l'ordre en $s=0$ de la fonction $r(\alpha,c^V,-s)$, c'est-\`a-dire aussi de la fonction $r(\alpha,c^V,s)$.

\
Pour un caract\`ere $\alpha$ de  $Z(\hat{M})^{\Gamma,\hat{\theta},0}$ qui intervient dans l'alg\`ebre de Lie de $\hat{G}$,  on a interpr\'et\'e ci-dessus  $d\alpha$ comme une forme lin\'eaire sur ${\cal A}_{\tilde{M},{\mathbb C}}^*$, autrement dit comme un \'el\'ement de ${\cal A}_{\tilde{M},{\mathbb C}}$. Il est clair que $d\alpha\in {\cal A}_{\tilde{M}}^{\tilde{G}}$. On se rappelle que cet espace est muni d'une mesure. Supposons que $\tilde{M}$ est un espace de Levi maximal de $\tilde{G}$. On pose alors
$$d^{\tilde{G}}(\alpha)=\frac{1}{2}vol({\cal A}_{\tilde{M}}^{\tilde{G}}/{\mathbb Z}d\alpha).$$ 

\bf D\'efinition: \rm Dans le cas o\`u $\tilde{M}$ est un espace de Levi maximal de $\tilde{G}$,  avec les notations pr\'ec\'edentes, on pose:
$$
r_{\tilde{M}}^{\tilde{G}}(c^V,\lambda):=\sum_{\alpha}r_{\tilde{M}}^{\tilde{G}}(\alpha,c_{\lambda}^V),$$
  o\`u   $\alpha$ parcourt les caract\`eres de $Z(\hat{M})^{\Gamma,\hat{\theta},0}$ qui interviennent dans l'alg\`ebre de Lie de $\hat{G}$ 
et 
$$r_{\tilde{M}}^{\tilde{G}}(\alpha,c_{\lambda}^V)=d^{\tilde{G}}(\alpha)r'(\alpha,c^V_{\lambda},0)/r(\alpha,c^V_{\lambda},0).
$$.  

\subsubsection{Stabilisation dans le cas d'un espace de Levi maximal}
 On d\'efinit aussi $s_{\tilde{M}}^{\tilde{G}}(\alpha,c^V_{\lambda})$ sous l'hypoth\`ese que $\tilde{G}$ est \`a torsion int\'erieure, que $\omega$ est trivial et que $\tilde{M}$ est maximal; cette derni\`ere hypoth\`ese assure que la composante neutre du groupe $Ker(\alpha)$ est incluse dans $Z(\hat{G})^{\Gamma}$ et donc dans $Z(\hat{G})\cap Z(\hat{M})^{\Gamma,0}$. On pose
$$
s_{\tilde{M}}^{\tilde{G}}(\alpha,c^V_{\lambda})=\vert Ker(\alpha)/Z(\hat{G})\cap Z(\hat{M})^{\Gamma,0})\vert^{-1}d^{\tilde{G}}(\alpha) r'(\alpha,c^V_{\lambda},0)/r(\alpha,c^V_{\lambda},0);
$$et dans le cas que l'on consid\`ere ici o\`u $\tilde{M}$ est un espace de Levi maximal de $\tilde{G}$, on pose aussi
$$
s_{\tilde{M}}^{\tilde{G}}(c^V,\lambda)=\sum_{\alpha}s_{\tilde{M}}^{\tilde{G}}(\alpha,c^V_{\lambda}),
$$
o\`u $\alpha$ parcourt l'ensemble des caract\`eres de $Z(\hat{M})^{\Gamma}$  (ici, dans le cas de torsion int\'erieure, $\hat{\theta}$ est trivial) agissant dans l'alg\`ebre de Lie de $\hat{G}$.
Et cette d\'efinition co\"{\i}ncide avec celle du cas non tordu quand $\tilde{G}=G$, ce qui est indispensable, \'evidemment, pr\'ecis\'ement parce que l'on a repris la d\'efinition de \cite{fonctionLendoscopique}.  La d\'efinition pr\'ec\'edente s'adapte au formalisme des donn\'ees endoscopiques, ce qui donne naissance aux fonctions qui interviennent ci-dessous. 
\begin{lem}Avec les notations d\'ej\`a introduites et sans hypoth\`ese sur $\tilde{G},\omega$
$$r_{\tilde{M}}^{\tilde{G}}(\alpha,c^V_{\lambda })=\sum_{\hat{s}\in \hat{s}_{M'}Z(\hat{M})^{\Gamma,\hat{\theta}}/Z(\hat{G})^{\Gamma,\hat{\theta}}}i_{\tilde{M}}(\tilde{G},{\bf G}'(\hat{s}))s_{{\bf M}'}^{{\bf G}'(\hat{s})}(\alpha',c^V_{\lambda}).$$ D'o\`u l'\'egalit\'e de fonctions m\'eromorphes:
$$
r_{\tilde{M}}^{\tilde{G}}(c^V,\lambda)=\sum_{\hat{s}\in \hat{s}_{M'}Z(\hat{M})^{\Gamma,\hat{\theta}}/Z(\hat{G})^{\Gamma,\hat{\theta}}}i_{\tilde{M}}(\tilde{G},{\bf G}'(\hat{s}))s_{{\bf M}'}^{{\bf G}'(\hat{s})}(c^V,\lambda).
$$
\end{lem}
Pour $\hat{s}$ intervenant dans cette somme, tel que $\alpha$ soit aussi une racine dans $\hat{G}'(\hat{s})$ not\'ee $\alpha'$, l'isomorphisme ${\cal A}_{\tilde{M}}^{\tilde{G}}\simeq {\cal A}_{\tilde{M}'}^{\tilde{G}'(\hat{s})}$ pr\'eserve les mesures par d\'efinitions et il envoie $d\alpha$ sur $d\alpha'$. Donc $d^{\tilde{G}}(\alpha)=d^{\tilde{G}'(\hat{s})}(\alpha')$.
De \ref{factorisation}, on a la formule
$$
r_{\tilde{M}}^{\tilde{G}}(\alpha,c^V_{\lambda})=\sum_{\hat{s}\in {\hat s}_{M'}Z(\hat{M})^{\Gamma,\hat{\theta}}/Ker(\alpha)} (r^{G(s)})'(\alpha',c^V_{\lambda},0)/r^{G(s)}(\alpha',c^V_{\lambda},0).
$$
On peut sommer sur $\hat{s}\in {\hat s}_{M'}Z(\hat{M})^{\Gamma,\hat{\theta}}/Z(\hat{G})^{\Gamma}\cap Z(\hat{M})^{\Gamma,\hat{\theta},0}$ \`a condition de mettre devant chaque terme le 
bon coefficient, c'est-\`a-dire $\vert Ker(\alpha)/Z(\hat{G})^{\Gamma}\cap Z(\hat{M})^{\Gamma,\hat{\theta},0}\vert^{-1}$. Il faut donc montrer pour tout $\hat{s}$ fix\'e, l'\'egalit\'e
$
i_{\tilde{M}}(\tilde{G},{\bf G}'(\hat{s}))=$
$$\vert Ker(\alpha)/Z(\hat{G})^{\Gamma}\cap Z(\hat{M})^{\Gamma,\hat{\theta},0}\vert^{-1}\vert Ker(\alpha')/Z(\hat{G}'(\hat{s}))\cap Z(\hat{M}')^{\Gamma,0})\vert.\eqno(1)
$$
Evidemment $Ker(\alpha)=Ker(\alpha')$. Ainsi (1) vaut $$\vert 
Z(\hat{G})^{\Gamma}\cap Z(\hat{M})^{\Gamma,\hat{\theta},0}/Z(\hat{G}'(\hat{s}))\cap Z(\hat{M}')^{\Gamma,0})\vert.$$
Par ellipticit\'e $Z(\hat{M}')^{\Gamma,0}= Z(\hat{M})^{\Gamma,\hat{\theta},0}$ et on trouve
$i_{\tilde{M}'}(\tilde{G},{\bf G}'(\hat{s}))$ d'apr\`es \ref{nonconnexite}.
\subsubsection{Stabilisation dans le cas non ramifi\'e\label{stabnonramif}}
Le caract\`ere pond\'er\'e dans le cas non ramifi\'e  et pour un espace de Levi qui n'est pas maximal se calcule par une formule de descente que l'on va r\'ecrire.

On fixe $V$ un ensemble fini de places contenant $V_{ram}$ et un espace de Levi $\tilde{M}$. 
Soit ${\bf M}'$ une donn\'ee endoscopique  elliptique de $\tilde{M}$ et soit $c^V$ un caract\`ere des alg\`ebres de Hecke sph\'eriques de ${\bf M}'$ hors de $V$.     Via la fonctorialit\'e de Langlands aux places non ramifi\'ees on voit $c^V$ comme un caract\`ere de l'alg\`ebre de Hecke sph\'erique de $\tilde{M}({\mathbb A}_{F}^V)$.  Pour un   caract\`ere $\alpha$ de $Z(\hat{M})^{\hat{\theta},0}$ intervenant dans l'alg\`ebre de Lie de $\hat{G}$, on a d\'ej\`a d\'efini $r_{\tilde{M}}^{\tilde{G}}(\alpha,c^V,\lambda)$ dans le cas o\`u $\tilde{M}$ \'etait un espace de Levi propre maximal de $\tilde{G}$. Si cette condition n'est pas v\'erifi\'ee, on pose $r_{\tilde{M}}^{\tilde{G}}(\alpha,c^V,\lambda)=0$. En rempla\c{c}ant $\tilde{G}$ par $\tilde{L}$, on a ainsi 
 d\'efini $r_{\tilde{M}}^{\tilde{L}}(\alpha,c^V,\lambda)$ comme fonction m\'eromorphe de $\lambda \in {\mathcal{A}}^*_{\tilde{M},{\mathbb C}}$ pour tout espace de Levi $\tilde{L}$ contenant $\tilde{M}$.   
\

\noindent\bf formule de descente.\rm

  Soit $A$ un ensemble de  caract\`eres de $Z(\hat{M})^{\hat{\theta},0}$ intervenant dans l'alg\`ebre de Lie de $\hat{G}$.  Si $A$ est r\'eduit \`a un seul \'el\'ement $\alpha$, on pose
$$ r_{\tilde{M}}^{\tilde{G}}(A,c^V,\lambda):= r_{\tilde{M}}^{\tilde{G}}(\alpha,c^V,\lambda).$$
  En g\'en\'eral, soit $A=A_{1}\cup A_{2}$ une d\'ecomposition de $A$. On pose alors en utilisant les d\'efinitions par r\'ecurrence:
$$
r_{\tilde{M}}^{\tilde{G}}(A,c^V,\lambda):=\sum_{\tilde{L},\tilde{L}'\in {\mathcal L}(\tilde{M})}d_{\tilde{M}}(\tilde{L},\tilde{L}') r_{\tilde{M}}^{\tilde{L}}(A_{1},c^V,\lambda)r_{\tilde{M}}^{\tilde{L}'}(A_{2},c^V,\lambda),\eqno(1)
$$ce qui donne une d\'efinition.
Il faut \'evidemment remarquer que le membre de gauche est aussi celui qui est associ\'e \`a la $\tilde{G}, \tilde{M}$-famille form\'ee avec les fonctions $L$ partielles d\'efinies avec les \'el\'ements de $A$ comme dans le paragraphe 4 de \cite{fonctionLendoscopique} et ainsi la d\'efinition du membre de gauche ne d\'epend pas de la d\'ecomposition $A=A_{1}\cup A_{2}$ choisie. C'est la formule de descente pour les $\tilde{G},\tilde{M}$ familles coupl\'ee au lemme 7.1 de \cite{eisensteinII} dans le cas non tordu et \cite{labwald} lemme 2.10.2.

\noindent\bf version stable des caract\`eres pond\'er\'es.\rm

 Dans le cas o\`u $\tilde{G}$ est \`a torsion int\'erieure avec un groupe sous-jacent quasi-d\'eploy\'e et un caract\`ere trivial, on a une formule de descente pour le caract\`ere pond\'er\'e stable. On pose, avec $A=A_{1}\cup A_{2}$ comme ci-dessus.
$$
s_{\tilde{M}}^{\tilde{G}}(A,c^V,\lambda)=\sum_{\tilde{L},\tilde{L}'\in {\mathcal L}(\tilde{M})}e_{\tilde{M}}(\tilde{L},\tilde{L}') s_{\tilde{M}}^{\tilde{L}}(A_{1}, c^V,\lambda)s_{\tilde{M}}^{\tilde{L}'}(A_{2},c^V,\lambda).\eqno(2)
$$
 Comme on a d\'efini les caract\`eres pond\'er\'es stables dans le cas d'un espace de Levi maximal, cela donne une d\'efinition en toute g\'en\'eralit\'e qui a priori d\'epend de la d\'ecomposition. Cette d\'efinition co\"{\i}ncide avec celle de \cite{fonctionLendoscopique} theorem 5 dans les cas d\'ej\`a connus, c'est \`a dire les cas o\`u il n'y a pas de torsion et quand on aura d\'emontr\'e la proposition ci-dessous, on aura en prime  que la d\'efinition  est ind\'ependante de la d\'ecomposition de $A$. Ici encore, la d\'efinition passe au formalisme des donn\'ees endoscopiques.

\

\noindent\bf version endoscopique des caract\`eres pond\'er\'es.\rm

On revient \`a la situation de $\tilde{G}$ g\'en\'eral consid\'er\'ee ci-dessus et on pose:
$$
r_{{\bf M}'}^{\tilde{G},{\mathcal E}}(A,c^V,\lambda)=\sum_{\hat{s}\in \hat{s}_{M'}Z(\hat{M})^{\Gamma,\hat{\theta}}/Z(\hat{G})^{\Gamma,\hat{\theta}}}i_{\tilde{M}'}(\tilde{G},{\bf G}'(\hat{s}))s_{{\bf M}'}^{{\bf G}(\hat{s})}(A,c^V,\lambda).
$$
\begin{prop}On a l'\'egalit\'e entre fonctions m\'eromorphes $$r_{\tilde{M}}^{\tilde{G}}(A,c^V,\lambda)=r_{{\bf M}'}^{\tilde{G},{\mathcal E}}(A,c^V,\lambda).$$
\end{prop}
On a d\'ej\`a montr\'e la proposition dans le cas des espaces de Levi maximaux.
Et pour la montrer en toute g\'en\'eralit\'e,  il suffit donc  de d\'emontrer que la fonction $r_{{\bf M}'}^{\tilde{G},{\mathcal E}}(c^V,\lambda)$ v\'erifie la m\^eme formule de descente \'ecrite en (1) que $r_{\tilde{M}}^{\tilde{G}}(c^V,\lambda)$. C'est la combinatoire de \cite{ftlt} 1.14 (i) (dans le cas non tordu c'est le paragraphe 6 de \cite{arthurduke99}).

\

\bf D\'efinition: \rm On pose $r_{\tilde{M}}^{\tilde{G}}(c^V,\lambda):= r_{\tilde{M}}^{\tilde{G}}(A,c^V,\lambda)$ et $s_{\tilde{M}}^{\tilde{G}}(c^V,\lambda):= s_{\tilde{M}}^{\tilde{G}}(A,c^V,\lambda)$, o\`u $A$ est exactement la r\'eunion des caract\`eres non triviaux pour l'action de $Z(\hat{M})^{\Gamma,\hat{\theta},0}$ agissant dans l'alg\`ebre de Lie de $\hat{G}$.

\subsubsection{Une propri\'et\'e de croissance\label{croissance}}
On g\'en\'eralise dans ce paragraphe le lemme 3.2 de \cite{partie1}.
\begin{lem} La fonction de $\lambda$, $r_{\tilde{M}}^{\tilde{G}}(c^V,\lambda)$ d\'efinie dans le paragraphe pr\'ec\'edent est holomorphe sur l'axe imaginaire et comme fonction de $\lambda\in i{\mathcal{A}}^*_{\tilde{M}}$ c'est une distribution temp\'er\'ee ce qui veut dire qu'il existe un nombre entier $N$ tel que $r_{\tilde{M}}^{\tilde{G}}(c^V,\lambda)(1+\| \lambda \|)^{-N}$ soit int\'egrable en $\lambda$.
\end{lem}
On commence par consid\'erer le cas o\`u $\tilde{M}$ est un espace de Levi maximal et on va se ramener aux r\'esultats du cas non tordu d\'emontr\'es par J. Arthur. On a donn\'e une formule explicite pour $r_{\tilde{M}}^{\tilde{G}}(c^V,\lambda)$ comme somme de d\'eriv\'ees logarithmiques de fonctions $L$-partielles. De plus chacun de ces termes est lui-m\^eme somme de d\'eriv\'ees logarithmiques $$L_{{\bf G}'(\hat{s})}^{'V}(\alpha',c^V_{\lambda},0)/L_{{\bf G}'(\hat{s})}^{V}(\alpha',c_{\lambda}^V,0)+L_{{\bf G}'(\hat{s})}^{'V}(-\alpha',c^V_{\lambda},0)/L_{{\bf G}'(\hat{s})}^{V}(-\alpha',c_{\lambda}^V,0)$$ o\`u ${\bf G}'(\hat{s})$ est une donn\'ee endoscopique de $\tilde{G}$ et o\`u $\alpha'$ est une racine convenable pour l'action de $Z(\hat{M})^{\Gamma_{F}}$ dans $\hat{G}'(\hat{s})$. Il suffit donc d'avoir le r\'esultat pour une telle fonction. La fonction est maintenant exactement celle qui intervient pour le Levi $\tilde{M}'$ dans l'espace $\tilde{G}'(\hat{s})$  (ou plus exactement, en fixant des donn\'ees auxiliaires,  pour le Levi $\tilde{M}'_{1}$ de l'espace $\tilde{G}'(\hat{s})_{1}$). Le r\'esultat de \cite{partie1} 3.2 s'applique \`a la somme sur $\alpha'$ et non pas \`a une fonction ind\'ependamment. Ce probl\`eme est r\'esolu dans la preuve de la proposition 1 de \cite{fonctionLendoscopique} o\`u il est montr\'e que quitte \`a changer encore de groupe on se ram\`ene facilement au cas o\`u la somme est r\'eduite \`a un seul terme.

On passe maintenant au cas g\'en\'eral: la propri\'et\'e de croissance de l'\'enonc\'e, \`a montrer,  se traduit conform\'ement \`a \cite{partie1} lemma 3.2, par le fait qu'il existe un entier $N$ tel que l'int\'egrale sur $i({\mathcal{A}}^{*}_{\tilde{M}}/{\mathcal{A}}^*_{\tilde{G}})$ de la fonction $r_{\tilde{M}}^{\tilde{G}}(c^V,\lambda)(1+\| \lambda\|)^{-N}$ soit convergente. On revient \`a la formule (1) de \ref{caracpondnonramifdef}; pour tout couple $\tilde{L},\tilde{L}'$ y intervenant, on peut appliquer le lemme par r\'ecurrence puisque l'on a d\'ej\`a d\'emontr\'e le cas des espaces de Levi maximaux; le fait que $d_{\tilde{M}}(\tilde{L},\tilde{L}')$ soit non nul, entra\^{\i}ne la d\'ecomposition en somme directe:
$$
i({\mathcal{A}}^{*}_{\tilde{M}}/{\mathcal{A}}^*_{\tilde{G}})=i({\mathcal{A}}^{*}_{\tilde{M}}/{\mathcal{A}}^*_{\tilde{L}})\oplus i({\mathcal{A}}^{*}_{\tilde{M}}/{\mathcal{A}}^*_{\tilde{L}'}).
$$Donc l'int\'egration en $\lambda$ se d\'ecompose en produit d'int\'egrations et chacune d'elle est convergente.
Cela entra\^{\i}ne le lemme.

\subsection{Caract\`eres pond\'er\'es invariants}
Dans tout ce paragraphe le corps de base est local.
\subsubsection{Rappel des d\'efinitions}
Soit $\tilde{M}$ un espace de Levi de $\tilde{G}$ et soit $\lambda\in {\mathcal{A}}^*_{\tilde{M},{\mathbb C}}$. On voit $\lambda$ comme un caract\`ere de $\tilde{M}$ ce qui n\'ecessite d'avoir choisi une application $H_{\tilde{M}}$ de $\tilde{M}/M^1$ dans ${\mathcal{A}}_{\tilde{M}}$. On fixe donc une telle application comme en \cite{stabilisationII} 1.6.

Soit encore $\pi$ une $\omega$ repr\'esentation de $\tilde{M}$. On a d\'efini pour tout $f\in I(\tilde{G})$, la fonction m\'eromorphe de $\lambda$: $J_{\tilde{M}}^{\tilde{G}}(\pi_{\lambda},f)$. Soit $X\in {\mathcal{A}}_{\tilde{M}}$; pour tout espace de Levi $\tilde{L}$ contenant $\tilde{M}$, on note $X_{\tilde{L}}$ la projection de $X$ sur ${\mathcal{A}}_{\tilde{L}}$. Pour tout $\lambda\in {\mathcal{A}}^*_{\tilde{M},{\mathbb C}}$ tel que $\lambda+i{\mathcal{A}}^*_{\tilde{M}}$ ne coupe pas les hyperplans singuliers de la fonction $J_{\tilde{M}}^{\tilde{G}}(\pi_{\lambda},f)$, on pose:
$$
J_{\tilde{M}}^{\tilde{G}}(\pi,{\lambda},X,f):=\int_{\mu\in \lambda+i{\mathcal{A}}^*_{\tilde{M}}}d\mu\, e^{-\mu(X)}J(\pi_{\mu},f).
$$
Cette fonction de $\lambda$,  est localement constante sur son espace de d\'efinition; elle ne d\'epend que de la restriction de $f$ \`a l'ensemble des \'el\'ements $\gamma\in \tilde{G}$ tel que  $H_{\tilde{G}}(\gamma)=X$. Cette distribution s'\'etend donc aux \'el\'ements de $I_{ac}(\tilde{G})$ et  cela permet de d\'efinir, par induction, les distributions invariantes:
$$f\in I_{ac}(\tilde{G})\mapsto $$
$$
I_{\tilde{M}}^{\tilde{G}}(\pi,{\lambda},X,f)=:J_{\tilde{M}}^{\tilde{G}}(\pi_{\lambda},X,f)-\sum_{\tilde{L}\in {\mathcal E}}I_{\tilde{M}}^{\tilde{L}}(\pi,{\lambda},X,\phi_{\tilde{L}}^{\tilde{G}}(f)),$$
o\`u encore $\lambda$ est tel que $\lambda+i{\mathcal{A}}^*_{\tilde{M}}$ ne coupe pas un  r\'eseau d'hyperplans d\'ependant de $\pi$ et o\`u $\phi_{\tilde{L}}^{\tilde{G}}$  est l'application de $I(\tilde{G})$ dans $I_{ac}(\tilde{L})$ d\'efinie par la proposition de \cite{ftlt} 6.4.

La propri\'et\'e qui r\'esulte de la d\'efinition des fonctions $\phi_{\tilde{L}}^{\tilde{G}}$ est que pour tout $\pi$ temp\'er\'ee, pour tout $X\in {\mathcal{A}}_{\tilde{M}}$ et tout $f\in I(\tilde{G})$,
$I_{\tilde{M}}^{\tilde{G}}(\pi,\lambda,X,f)=0$ si $\lambda$ est imaginaire.
De plus dans tous les cas, $I(\pi,\lambda,X,f)$ ne d\'epend que de la restriction de $f$ \`a l'ensemble  des \'el\'ements $\gamma\in \tilde{G}$ tel que  $H_{\tilde{G}}(\gamma)=X$. Et en tant que fonction de $\lambda$, $I(\pi,\lambda,X,f)$ est  une fonction localement constante mais non partout d\'efinie.

\begin{lem} ({\bf formule de descente}) On fixe $\tilde{R}$ un sous-espace de Levi de $\tilde{M}$ et on suppose qu'il existe une $\omega$ repr\'esentation de $\tilde{R}$ tel que $\pi$ soit l'induite de $\sigma$. Alors, on a pour tout $f\in I(\tilde{G})$ et pour tout $X\in {\mathcal{A}}_{\tilde{M}}$
$$
I_{\tilde{M}}^{\tilde{G}}(\pi,{\lambda},X,f)=\sum_{\tilde{L} \in {\mathcal L}(\tilde{R})} d_{\tilde{R}}^{\tilde{G}}(\tilde{L},\tilde{M})\int_{Y\in {\mathcal{A}}_{\tilde{R}}; Y_{\tilde{M}}=X} dY\, I_{\tilde{R}}^{\tilde{L}}(\sigma,{\lambda},Y,f_{\tilde{L}}),
$$formule qui m\'erite les explications ci-dessous.
 \end{lem}
 D'abord on commence par remarquer que si la place consid\'er\'ee est non archim\'edienne, l'int\'egrale est une  somme finie: en effet, il n'y a qu'un nombre fini d'\'el\'ements dans ${\mathcal L}(\tilde{R})$ et le coefficient $d_{\tilde{R}}^{\tilde{G}}(\tilde{L},\tilde{M})$ est non nul seulement si l'application naturelle de ${\mathcal{A}}_{\tilde{R}}^{\tilde{L}}\oplus {\mathcal{A}}_{\tilde{R}}^{\tilde{M}}$ dans ${\mathcal{A}}_{\tilde{R}}^{\tilde{G}}$ est un isomorphisme. Ainsi on somme sur $Y\in X +{\mathcal{A}}_{\tilde{R}}^{\tilde{M}}$ et sur cet ensemble la projection orthogonale $Y \mapsto Y_{\tilde{L}}$ dans ${\mathcal{A}}_{\tilde{L}}$ est injective: en effet le noyau est ${\mathcal{A}}_{\tilde{R}}^{\tilde{L}}\cap {\mathcal{A}}_{\tilde{R}}^{\tilde{M}}$ qui est nul d'apr\`es ce que l'on vient de voir.   Ainsi dans la somme sur les $Y$ la composante de $Y$ dans ${\mathcal{A}}_{\tilde{M}}$ est fix\'ee et par contre la composante de $Y$ sur ${\mathcal{A}}_{\tilde{L}}^{\tilde{G}}$ est laiss\'ee libre; mais cette derni\`ere composante concerne le centre de $\tilde{L}$ et $f_{\tilde{L}}$ est \`a support compact. D'o\`u la finitude. Dans le cas archim\'edien, le m\^eme argument dit que l'int\'egrale porte sur un compact, elle est donc bien d\'efinie.

Montrons le r\'esultat cherch\'e: on part de la formule de descente pour les caract\`eres pond\'er\'es non invariants donn\'ee en \cite{ftlt} 5.4 (iv) (comme les termes ne sont pas invariants, il faut des donn\'ees auxiliaires pour d\'efinir les termes constants comme ceci est expliqu\'e en loc.cite):
$$
J_{\tilde{M}}^{\tilde{G}}(ind \sigma_{\lambda},f)=\sum_{\tilde{L} \in {\mathcal L}(\tilde{R})} d_{\tilde{R}}^{\tilde{G}}(\tilde{L},\tilde{M}) J_{\tilde{R}}^{\tilde{L}}(\sigma_{\lambda},f_{\tilde{L}}).
$$
On int\`egre cette formule sur $i{\mathcal{A}}_{\tilde{M}}^{*}$:
$$
\int_{\mu\in i{\mathcal{A}}^*_{\tilde{M}}}d\mu\, J_{\tilde{M}}^{\tilde{G}}(ind \sigma_{\lambda+\mu},f)e^{-(\lambda+\mu)(X)}=
J_{\tilde{M}}^{\tilde{G}}(ind\sigma_{\lambda},X,f)$$
car c'est la d\'efinition. Pour traiter le c\^ot\'e droit, on reprend la d\'efinition pour $\tilde{L}\in {\mathcal L} (\tilde{R})$ et pour $Z\in {\mathcal{A}}_{\tilde{R}}$ de $J_{\tilde{R}}^{\tilde{G}}(\sigma_{\lambda},Z,f_{\tilde{L}})$ comme une int\'egrale sur $i{\mathcal{A}}_{\tilde{R}}^*$. Puisque l'on n'int\`egre que sur $i{\mathcal{A}}_{\tilde{M}}^*$, on obtient une formule d'inversion de Fourier, que l'on \'ecrit dans le cas d'un corps de base non archim\'edien:
$$\sum_{\tilde{L} \in {\mathcal L}(\tilde{R})} d_{\tilde{R}}^{\tilde{G}}(\tilde{L},\tilde{M}) \sum_{Z\in {\mathcal{A}}_{\tilde{R}}^{\tilde{M}}} J_{\tilde{R}}^{\tilde{L}}(\sigma_{\lambda},X+Z,f_{\tilde{L}}).
$$
Autrement dit, on a la m\^eme formule que dans l'\'enonc\'e avec les distributions $J$ au lieu de $I$. D'apr\`es les d\'efinitions, il suffit donc d'obtenir une formule analogue pour 
$$
\sum_{\tilde{L}'}I_{\tilde{M}}^{\tilde{L}'}(ind \sigma,\lambda,X,\phi_{\tilde{L}'}^{\tilde{G}}(f)),
$$
o\`u l'on somme sur les espaces de Levi propres $\tilde{L}'$ contenant $\tilde{M}$. 
Soit $\tilde{L}'\in {\mathcal L}(\tilde{M})$ un tel espace de Levi; on consid\`ere $I_{\tilde{M}}^{\tilde{L}'}(ind \sigma,\lambda,X,\phi_{\tilde{L}'}^{\tilde{G}}(f))$. On applique la formule de l'\'enonc\'e par r\'ecurrence \`a  ce terme, ce qui est loisible car on peut restreindre $\phi_{\tilde{L}'}^{\tilde{G}}(f)$ \`a l'ensemble des $\gamma\in \tilde{L}'$ tels que $H_{\tilde{L}'}(\gamma)=X_{\tilde{L}'}$ et alors la fonction ainsi restreinte est \`a support compact. 
On obtient:
$$I_{\tilde{M}}^{\tilde{L}'}(ind \sigma,\lambda,X,\phi_{\tilde{L}'}^{\tilde{G}}(f))=\sum_{\tilde{L}''\in {\mathcal L}(R),\tilde{L}''\subset \tilde{L}'}d^{\tilde{L}'}_{\tilde{R}}(\tilde{L}'',\tilde{M})$$
$$\sum_{Y\in {\mathcal{A}}_{\tilde{R}};Y_{\tilde{M}}=X}I_{\tilde{R}}^{\tilde{L}''}(\sigma,\lambda,Y,(\phi_{\tilde{L}'}^{\tilde{G}})(f)_{\tilde{L}''}).
$$
Il faut encore utiliser la formule de scindage pour les termes $(\phi_{\tilde{L}'}^{\tilde{G}})(f)_{\tilde{L}''}$ donn\'ee en \cite{ftlt} 6.4 (11):
$$
(\phi_{\tilde{L}'}^{\tilde{G}})(f)_{\tilde{L}''}=\sum_{\tilde{L}\in {\mathcal{L}}(\tilde{R}); \tilde{L}''\subset \tilde{L}}d_{\tilde{L}''}^{\tilde{G}}(\tilde{L}',\tilde{L})
(\phi_{\tilde{L}''}^{\tilde{L}})(f_{\tilde{L}}).
$$
D'o\`u $$\sum_{\tilde{L}'\in {\mathcal L}(\tilde{M}); \tilde{L}'\neq \tilde{G}}I_{\tilde{M}}^{\tilde{L}'}(ind \sigma,\lambda,X,\phi_{\tilde{L}'}^{\tilde{G}}(f))=$$
$$ \sum_{\tilde{L}\in {\mathcal L}(\tilde{R})}$$
$$\sum_{\tilde{L}',\tilde{L}'' \in {\mathcal{L}}(\tilde{R}); \tilde{L}'\neq \tilde{G},\tilde{L}''\subset \tilde{L}'\cap \tilde{L}}d_{\tilde{R}}^{\tilde{L}'}(\tilde{M},\tilde{L}'')d_{\tilde{L}''}^{\tilde{G}}(\tilde{L}',\tilde{L})$$
$$
\sum_{Y\in {\mathcal{A}}_{\tilde{R}};Y_{\tilde{M}}=X}I_{\tilde{R}}^{\tilde{L}''}(\sigma,\lambda,Y,\phi_{\tilde{L}''}^{\tilde{L}}(f_{\tilde{L}})).
$$
Dans la somme de droite on fait dispara\^{\i}tre $\tilde{L}'$ en utilisant \cite{stabilisationII} 1.7 (5): l'ensemble $B$ de loc. cite est ici les triplets $(\tilde{L},\tilde{L}'',\tilde{L}')$ (au lieu de $(\tilde{L}',\tilde{R}',\tilde{R})$ in loc.cite), $\tilde{R}$ ici joue le r\^ole de $\tilde{M}$ in loc.cite et $\tilde{M}$ ici joue le r\^ole de $\tilde{L}$ in loc. cite. On a en particulier l'\'egalit\'e 
$$
d_{\tilde{R}}^{\tilde{G}}(\tilde{M},\tilde{L})=d_{\tilde{R}}^{\tilde{L}'}(\tilde{M},\tilde{L}'')d_{\tilde{L}''}^{\tilde{G}}(\tilde{L}',\tilde{L}).
$$ et on obtient que le c\^ot\'e droit vaut:
$$
\sum_{\tilde{L}\in {\mathcal L}(\tilde{R})}d_{\tilde{R}}^{\tilde{G}}(\tilde{M},\tilde{L}) $$
$$\sum_{Y\in {\mathcal{A}}_{\tilde{R}}, Y_{\tilde{M}}=X}\sum_{\tilde{L}'' \in {\mathcal L}(\tilde{R}); \tilde{L}''\subsetneq \tilde{L}}I_{\tilde{R}}^{\tilde{L}''}(\sigma,\lambda,Y,\phi_{\tilde{L}''}^{\tilde{L}}(f_{\tilde{L}})).
$$
Et on obtient le lemme.

\subsubsection{Les caract\`eres pond\'er\'es, variante compacte\label{variantecompacte}}
Il est d\'efini en \cite{stabilisationVIII} 1.4 et \cite{stabilisationIX} 5.8 une variante compacte de l'application $\phi_{\tilde{M}}^{\tilde{G}}$; dans le cas o\`u le corps local est archim\'edien, il y a  deux variantes compacte et  on va utiliser celle construite \`a partir de $\phi_{\tilde{M}}^{rat}$ et non l'autre. On   adopte la notation, $^c\phi_{\tilde{M}}^{\tilde{G}}$ en oubliant l'exposant $rat$ ce qui unifie les notations du cas archim\'edien et du cas $p$-adique et c'est le choix de \cite{stabilisationIX}. C'est une application de $I(\tilde{G})$ dans $I(\tilde{M})$. On pose
$$
^cI_{\tilde{M}}^{\tilde{G}}(\pi,\lambda,f):=J_{\tilde{M}}^{\tilde{G}}(\pi_{\lambda},f)-\sum_{\tilde{L}\in {\mathcal L}(\tilde{M}); \tilde{L}\neq \tilde{G}}\, ^cI_{\tilde{M}}^{\tilde{L}}(\pi,\lambda,\, ^c\phi_{\tilde{L}}^{\tilde{G}}).\eqno(1)
$$
Cette distribution d\'epend m\'eromorphiquement de $\lambda$ et  l'on n'a pas besoin de prendre de coefficients de Fourier pour d\'efinir cette distribution car $^c\phi_{\tilde{L}}^{\tilde{G}}$ pr\'eserve le caract\`ere compact du support.

Cette distribution v\'erifie une formule de descente, pour des repr\'esentations induites, plus simple que sa variante non compacte: pour tout sous-espace de Levi $\tilde{R}$ de $\tilde{M}$ et pour toute $\omega$ repr\'esentation de $\tilde{R}$:
$$^cI_{\tilde{M}}^{\tilde{G}}(ind\, \sigma,\lambda,f)=\sum_{\tilde{L}\in {\mathcal L}(\tilde{R})}d_{\tilde{R}}(\tilde{L},\tilde{M}) \, ^cI_{\tilde{R}}^{\tilde{L}}(\sigma,\lambda,f_{\tilde{L}}).
$$
Pour comparer cette distribution \`a celle d\'ej\`a d\'efinie, on d\'efinit comme dans le paragraphe pr\'ec\'edent les coefficients de Fourier
$$
^cI_{\tilde{M}}^{\tilde{G}}(\pi,\lambda,X,f):=\int_{\mu \in i{\mathcal{A}}^*_{\tilde{M}}}d\mu\, 
^cI_{\tilde{M}}^{\tilde{G}}(\pi,\lambda+\mu,f) e^{-(\lambda+\mu)(X)}.
$$
Dans les r\'ef\'erences d\'ej\`a mentionn\'ees, il a \'et\'e construit une application $^c\theta_{\tilde{M}}^{\tilde{G}}$ de $I(\tilde{G})$ dans $I_{ac}(\tilde{G})$ qui relie les d\'efinitions de $\phi_{\tilde{M}}^{\tilde{G}}$ et $^c\phi_{\tilde{M}}^{\tilde{G}}$ par la formule:
$$\forall f\in I(\tilde{G}), \quad ^c\theta_{\tilde{M}}^{\tilde{G}}(f)=\phi_{\tilde{M}}^{\tilde{G}}(f)-\sum_{\tilde{L}\in {\mathcal L}(\tilde{M}), \tilde{L}\neq \tilde{G}}\, ^c\theta_{\tilde{M}}^{\tilde{L}}(^c\phi_{\tilde{L}}^{\tilde{G}}(f)).\eqno(2)
$$
Cette application sert \`a relier les distributions d\'efinies ici et dans le paragraphe pr\'ec\'edent. 
\begin{lem} Pour tout $\lambda\in {\mathcal{A}}_{\tilde{M}}^*$ o\`u les termes ci-dessous sont d\'efinis, on a l'\'egalit\'e
$$
^cI_{\tilde{M}}^{\tilde{G}}(\pi,\lambda,X,f)=\sum_{\tilde{L}\in {\mathcal L}(\tilde{M})}I_{\tilde{M}}^{\tilde{L}}(\pi,\lambda,X,\, ^c\theta_{\tilde{L}}^{\tilde{G}}(f)).
$$
\end{lem}
On utilise la d\'efinition (1) pour calculer le membre de gauche et on applique tout de suite le lemme par r\'ecurrence pour les termes faisant intervenir un espace de Levi strictement inclus dans $\tilde{G}$ \`a la place de $\tilde{G}$. On obtient$$
^cI_{\tilde{M}}^{\tilde{G}}(\pi,\lambda,X,f)= J_{\tilde{M}}^{\tilde{G}}(\pi,\lambda,X,f)-
$$
$$
\sum_{\tilde{L}_{1}\in {\mathcal L}(\tilde{M}),\tilde{L}_{1}\neq \tilde{G}}\qquad\sum_{\tilde{L}_{2}\in {\mathcal L}(\tilde{M}), \tilde{L}_{2}\subset \tilde{L}_{1}} I_{\tilde{M}}^{\tilde{L}_{2}}(\pi,\lambda,X,\, ^c\theta_{\tilde{L}_{2}}^{\tilde{L}_{1}}(\, ^c\phi_{\tilde{L}_{1}}^{\tilde{G}}(f))).
$$
On \'echange la somme sur $\tilde{L}_{1}$ et la somme sur $\tilde{L}_{2}$ et on a donc une somme sur les espaces de Levi $\tilde{L}_{2}\in {\mathcal L}(\tilde{M})$ v\'erifiant $\tilde{L}_{2}\neq  \tilde{G}$ de 
$
-I_{\tilde{M}}^{\tilde{L}_{2}}(\pi,\lambda,X,F_{2,1})$ o\`u 
$$
F_{2,1}:= \sum_{\tilde{L}_{1}\in {\mathcal L}(\tilde{L}_{2}),\tilde{L}_{1}\neq \tilde{G}}\, ^c\theta_{\tilde{L}_{2}}^{\tilde{L}_{1}}(\phi_{\tilde{L}_{1}}^{\tilde{G}}(f)).$$
Avec (2), on a $F_{2,1}=\phi_{\tilde{L}_{2}}^{\tilde{G}}(f)-\, ^c\theta_{\tilde{L}_{2}}^{\tilde{G}}(f)$. Ainsi 
$
^cI_{\tilde{M}}^{\tilde{G}}(\pi,\lambda,X,f)=$
$$J_{\tilde{M}}^{\tilde{G}}(\pi,\lambda,X,f)-\sum_{\tilde{L}_{2}\in {\mathcal L}(\tilde{M}); \tilde{L}_{2}\neq \tilde{G}} I_{\tilde{M}}^{\tilde{L}_{2}}(\pi,\lambda,X,\phi_{\tilde{L}_{2}}^{\tilde{G}}(f))\eqno(3)$$
$$
+ \sum_{\tilde{L}_{2}\in {\mathcal L}(\tilde{M}); \tilde{L}_{2}\neq \tilde{G}}I_{\tilde{M}}^{\tilde{L}_{2}}(\pi,\lambda,X,\, ^c\theta_{\tilde{L}_{2}}^{\tilde{G}}(f))\eqno(4).
$$
Le  terme (3) est $I_{\tilde{M}}^{\tilde{G}}(\pi,\lambda,X,f)$ par d\'efinition; or $^c\theta_{\tilde{G}}^{\tilde{G}}(f)=f$ et  en ajoutant (3) et (4),  on obtient le lemme.

\subsubsection{Les caract\`eres pond\'er\'es compacts des repr\'esentations temp\'er\'ees\label{carpondtemp}}
On suppose que $\pi$ est une $\omega$-repr\'esentation temp\'er\'ee.
Il a \'et\'e v\'erifi\'e en \cite{stabilisationIX} (les premi\`eres lignes de 5.10) que le coefficient de Fourier $I^{\tilde{M}}(\pi,X,\, ^c\theta_{\tilde{M}}^{\tilde{G}}(f))$  est la transform\'ee de Fourier d'une fonction sur $i{\mathcal{A}}^*_{\tilde{M}}$; dans le cas $p$-adique, il n'y a pas d'exposant $rat$. Ceci n'est autre par d\'efinition que $ tr \pi_{\lambda}(\, ^c\theta_{\tilde{M}}^{\tilde{G}}(f))$.  La fonction de $\lambda$,  $tr \pi_{\lambda}(\, ^c\theta_{\tilde{M}}^{\tilde{G}}(f))$ se prolonge m\'eromorphiquement en $\lambda$. On note $I^{\tilde{M}}(\pi_{\lambda},\, ^c\theta_{\tilde{M}}^{\tilde{G}}(f))$ ce prolongement.
\begin{lem} On suppose que $\pi$ est une $\omega$-repr\'esentation temp\'er\'ee;   on a l'\'egalit\'e, pour tout $f\in I(\tilde{G})$, des fonctions m\'eromorphes en $\lambda$
$$
^cI_{\tilde{M}}^{\tilde{G}}(\pi,\lambda,f)=I^{\tilde{M}} (\pi_{\lambda},\, ^c\theta_{\tilde{M}}^{\tilde{G}}(f)).
$$
\end{lem}
Chaque terme de l'\'egalit\'e est une fonction m\'eromorphe en $\lambda$ et il suffit donc de montrer l'\'egalit\'e pour $\lambda$ imaginaire. Sous cette hypoth\`ese
$$
^cI_{\tilde{M}}^{\tilde{G}}(\pi,{\lambda},f)=J_{\tilde{M}}^{\tilde{G}}(\pi,{\lambda},f)-
\sum_{\tilde{L}\in {\mathcal L}(\tilde{M}), \tilde{L}\neq \tilde{G}} \, ^cI_{\tilde{M}}^{\tilde{L}}(\pi_{\lambda},\, ^c\phi_{\tilde{L}}^{\tilde{G}}(f)).
$$Par construction de $\phi_{\tilde{M}}^{\tilde{L}}(f)$, on a l'\'egalit\'e:
$
J_{\tilde{M}}^{\tilde{G}}(\pi_{\lambda},f)=tr\, \pi_{\lambda}(\phi_{\tilde{M}}^{\tilde{G}}(f))$ et avec une r\'ecurrence  facile, on obtient l'\'egalit\'e:
$$
^cI_{\tilde{M}}^{\tilde{G}}(\pi,{\lambda},f)= tr \pi_{\lambda}\bigl(\phi_{\tilde{M}}^{\tilde{G}}(f)-\sum_{\tilde{L}\in {\mathcal L}(\tilde{M}), \tilde{L}\neq \tilde{G}}\, ^c\theta_{\tilde{M}}^{\tilde{L}}(^c\phi_{\tilde{L}}^{\tilde{G}}(f))\bigr).
$$
C'est exactement la formule (2) ci-dessus qui donne alors le lemme.

\subsection{Le cas de la torsion int\'erieure}
\subsubsection{Les caract\`eres pond\'er\'es invariants stables, premi\`eres d\'efinitions \label{defcarstable}}
On suppose ici que $\tilde{G}$ est \`a torsion int\'erieure, c'est-\`a-dire que $\hat{\theta}=1$, $\omega=1$ et que le groupe $G$ est quasid\'eploy\'e. Dans ce cas, on a besoin de d\'efinir les caract\`eres pond\'er\'es stables. Soit $\tilde{M}$ un espace de Levi de $\tilde{G}$.
On fixe $\pi$ une repr\'esentation de $\tilde{M}$, on d\'efinit 
 pour $X\in {\mathcal{A}}_{\tilde{M}}$ et $\lambda\in {\mathcal{A}}^*_{\tilde{M}}$ g\'en\'eral,  la distribution 
$f\in I(\tilde{G})\mapsto I_{\tilde{M}}^{\tilde{G}}(\pi,X,\lambda,f)$. La condition sur $\lambda$ pour que cette distribution soit d\'efinie est que $\lambda+i{\mathcal{A}}^*_{\tilde{M}}$ ne coupe pas un r\'eseau d'hyperplans. 
On suppose que $\pi$ est stable. 
Pour tout $f\in I(\tilde{G})$, on pose en recopiant les d\'efinitions de J. Arthur: $$SI_{\tilde{M}}^{\tilde{G}}(\pi,\lambda,X,f):= I_{\tilde{M}}^{\tilde{G}}(\pi,\lambda,X,f)-$$
$$\sum_{s\in Z(\hat{M})^{\Gamma}/Z(\hat{G})^{\Gamma},s\neq 1}i_{\tilde{M}}(\tilde{G},{\bf G}'(s))SI_{\tilde{M}}^{\tilde{G}'(s)}(\pi,\lambda,X,f^{{\bf G'}(s)})\eqno(1)
$$Bien s\^ur cette d\'efinition n'a de sens que si les distributions $SI_{\tilde{M}}^{\tilde{G}'(s)}(\pi,\lambda,X,.)$ sont stables. On l'admet par r\'ecurrence en remarquant que puisque si $s\neq 1$ le groupe $G'(s)$ est plus ''petit'' que $G$ et on doit \'evidemment montrer la stabilit\'e de $SI_{\tilde{M}}^{\tilde{G}}(\pi,\lambda,X,.)$ ce qui sera fait dans le paragraphe \ref{stabilite} qui suit.

On a aussi d\'efini les variantes compactes avec les m\^emes propri\'et\'es de stabilit\'e admises par r\'ecurrence pour ${\bf G}'(s)$ si $s\neq 1$ et \`a d\'emontrer pour $\tilde{G}$. On pose donc aussi pour tout $f\in I(\tilde{G})$:
$$
^c(SI)_{\tilde{M}}^{\tilde{G}}(\pi,\lambda,X,f):= \, ^cI_{\tilde{M}}^{\tilde{G}}(\pi,\lambda,X,f)-$$
$$\sum_{s\in Z(\hat{M})^{\Gamma}/Z(\hat{G})^{\Gamma},s\neq 1}i_{\tilde{M}}(\tilde{G},{\bf G'(s)})\, ^c(SI)_{\tilde{M}}^{\tilde{G}'(s)}(\pi,\lambda,X,f^{\bf{G'}(s)}).\eqno(2)
$$
Pour relier ces deux d\'efinitions, il faut utiliser les applications $^c\theta_{\tilde{L}}^{\tilde{G}}$ et on pose pour tout $f\in I(\tilde{G})$:
$$
^c(S\theta)_{\tilde{L}}^{\tilde{G}}(f):= \, ^c\theta_{\tilde{L}}^{\tilde{G}}(f)-\sum_{s\in Z(\hat{L})^{,\Gamma}/Z(\hat{G})^{\Gamma}, s\neq 1}i_{\tilde{L}}(\tilde{G},{\bf G}'(s))\, ^c(S\theta)_{\tilde{L}}^{\tilde{G}'(s)}(f^{{\bf G}'(s)}).
$$
Il a \'et\'e montr\'e en \cite{stabilisationVIII} 2.2 (dont la preuve est en \cite{stabilisationVIII} 4.2) et \cite{stabilisationIX} 6.1 (dont la preuve est en \cite{stabilisationIX} 7.3) que cette application $f\mapsto\, ^c(S\theta)_{\tilde{L}}^{\tilde{G}}(f)$ est stable, c'est-\`a-dire que $^c(S\theta)_{\tilde{L}}^{\tilde{G}}(f)$ ne d\'epend que de l'image de $f$ dans $SI(\tilde{G})$.
\begin{lem} Pour tout $f\in I(\tilde{G})$, on a l'\'egalit\'e
$$
^c(SI)_{\tilde{M}}^{\tilde{G}}(\pi,\lambda,X,f)=\sum_{\tilde{L}\in {\mathcal L}(\tilde{M})}SI_{\tilde{M}}^{\tilde{L}}(\pi,\lambda,X,\, ^c(S\theta)_{\tilde{L}}^{\tilde{G}}(f) ).
$$
\end{lem}
On \'ecrit $$\sum_{s\in Z(\hat{M})^{\Gamma}/Z(\hat{G})^{\Gamma},s\neq 1}i_{\tilde{M}}(\tilde{G},{\bf G'(s)})\, ^c(SI)_{\tilde{M}}^{\tilde{G}'(s)}(\pi,\lambda,X,f^{\bf{G'}(s)})\eqno(3)$$ en utilisant le lemme par r\'ecurrence puisque $s\neq 1$. On obtient:
$$
\sum_{\tilde{L}\in {\mathcal L}(\tilde{M})}
\sum_{s_{\tilde{L}}\in Z(\hat{M})^{\Gamma}/Z(\hat{L})^{\Gamma}, s_{\tilde{L}}\neq 1}$$
$$\sum_{s\in s_{\tilde{L}}Z(\hat{L})^{\Gamma}/Z(\hat{G})^{\Gamma}}i_{\tilde{M}}(\tilde{G},{\bf G}'(s)) SI_{\tilde{M}}^{\tilde{L}'(s_{\tilde{L}})}(\pi,\lambda,X,\, ^c(S\theta)_{\tilde{L}'(s_{\tilde{L}})}^{\tilde{G}'(s)}(f^{{\bf G}'(s)}))
$$
$$
+ \sum_{\tilde{L}\in {\mathcal L}(\tilde{M})}\sum_{s\in Z(\hat{L})^{\Gamma}/Z(\hat{G})^{\Gamma},s\neq 1}i_{\tilde{M}}(\tilde{G},{\bf G}'(s)) SI_{\tilde{M}}^{\tilde{L}}(\pi,\lambda,X,\, ^c(S\theta)_{\tilde{L}}^{\tilde{G}'(s)}(f^{{\bf G}'(s)})),
$$o\`u cette derni\`ere somme correspond aux $s_{\tilde{L}}=1$ et en gardant \'evidemment l'hypoth\`ese que $s\neq 1$.
Soit $s\in Z(\hat{M})^{\Gamma}/Z(\hat{G})^{\Gamma}$ et pour $\tilde{L}\in {\mathcal L}(\tilde{M})$, on note $s_{\tilde{L}}$ l'image de $s$ modulo $Z(\hat{L})^{\Gamma}$. On a l'\'egalit\'e:
$$
i_{\tilde{M}}(\tilde{G},{\bf G}'(s))=i_{\tilde{M}}(\tilde{L},{\bf L}'(s_{\tilde{L}}))i_{\tilde{L}}(\tilde{G},{\bf G}'(s)).$$
On est dans le cas de torsion int\'erieure, donc la stabilisation des int\'egrales orbitales est connue. Cela entra\^{\i}ne avec \cite{stabilisationVIII}, 3.8 et 4.3 et \cite{stabilisationIX} 6.11 (A), qu'avec les notations ci-dessus pour tout $s_{\tilde{L}}\neq 1$:
$$
\sum_{s\in s_{\tilde{L}}Z(\hat{L})^{\Gamma}/Z(\hat{G})^{\Gamma}}i_{\tilde{L}}(\tilde{G},{\bf G}'(s))
 ^c(S\theta)_{\tilde{L}'(s_{\tilde{L}})}^{\tilde{G}'(s)}(f^{{\bf G}'(s)})= \, ^c\theta_{\tilde{L}}^{\tilde{G}}(f)^{{\bf L}'(s_{\tilde{L}})}.\eqno(4)
 $$
 Donc la somme sur les $s_{\tilde{L}}\neq 1$ ci-dessus est la somme des termes
 $$i_{\tilde{M}}(\tilde{L},{\bf L}'(s_{\tilde{L}}))SI^{\tilde{L}'(s_{\tilde{L}})}_{\tilde{M}}(\pi,\lambda,X,\, ^c\theta_{\tilde{L}}^{\tilde{G}}(f)^{{\bf L}'(s_{\tilde{L}})}).$$
 Par d\'efinition ceci n'est autre que $$I_{\tilde{M}}^{\tilde{L}}(\pi,\lambda,X,\, ^c\theta_{\tilde{L}}^{\tilde{G}}(f)) -SI(\pi,\lambda,X,^c\theta_{\tilde{L}}^{\tilde{G}}(f)).$$
Consid\'erons la somme correspondant \`a $s_{\tilde{L}}=1$; l'\'egalit\'e (4) se remplace par la tautologie
 $$\sum_{s\in Z(\hat{L})^{\Gamma}/Z(\hat{G})^{\Gamma}}i_{\tilde{L}}(\tilde{G},{\bf G}'(s))
 ^c(S\theta)_{\tilde{L}}^{\tilde{G}'(s)}(f^{{\bf G}'(s)})= \, ^c\theta_{\tilde{L}}^{\tilde{G}}(f)$$
 et on obtient que les termes correspondant \`a $s_{\tilde{L}}=1$ contribuent simplement par
 $$
 SI_{\tilde{M}}^{\tilde{L}}(\pi,\lambda,X, \, ^c\theta_{\tilde{L}}^{\tilde{G}}(f))-\, ^c(S\theta)_{\tilde{L}}^{\tilde{G}}(f)).
 $$
 En regroupant avec ce que l'on a trouv\'e ci-dessus, on voit que (3) vaut $$
 \sum_{\tilde{L}\in {\mathcal L}(\tilde{M})}I_{\tilde{M}}^{\tilde{L}}(\pi,\lambda,X,\, ^c\theta_{\tilde{L}}^{\tilde{G}}(f))- SI_{\tilde{M}}^{\tilde{L}}(\pi,\lambda,X, \, ^c(S\theta)_{\tilde{L}}^{\tilde{G}}(f)).
 $$
 Avec le lemme de \ref{variantecompacte} ceci n'est autre que $$^cI_{\tilde{M}}^{\tilde{G}}(\pi,\lambda,X,f)-
  \sum_{\tilde{L}\in {\mathcal L}(\tilde{M})} SI_{\tilde{M}}^{\tilde{L}}(\pi,\lambda,X, \, ^c(S\theta)_{\tilde{L}}^{\tilde{G}}(f)).
  $$
  En revenant \`a la d\'efinition de $^c(SI)_{\tilde{M}}^{\tilde{G}}(\pi,\lambda,X,f)$ on obtient le lemme cherch\'e.

\subsubsection{Preuve de la stabilit\'e\label{stabilite}}
On fixe $\tilde{M}$ un espace de Levi de $\tilde{G}$ et une repr\'esentation stable $\pi$ de $\tilde{M}$ (on rappelle que l'on est dans le cas de la torsion intérieure et que $\omega=1$).
\begin{prop} Les distributions $f\in I(\tilde{G})\mapsto SI_{\tilde{M}}^{\tilde{G}}(\pi,\lambda,X,f)$ et $f\in I(\tilde{G})$ $\mapsto \, ^c(SI)_{\tilde{M}}^{\tilde{G}}(\pi,\lambda,X,f)$  sont stables en tout point $\lambda$ o\`u elles sont d\'efinies\end{prop}
On fixe $f\in I(\tilde{G})$ et on suppose que les int\'egrales orbitales stables de $f$ sont toutes nulles, en bref, on dit qu'une telle fonction $f$ est instable. On doit montrer que $SI_{\tilde{M}}^{\tilde{G}}(\pi,\lambda,X,f)=0$ et $^cSI_{\tilde{M}}^{\tilde{G}}(\pi,\lambda,X,f)=0$. 
On reprend la formule du lemme \ref{defcarstable} et son c\^ot\'e droit. Consid\'erons le terme correspondant \`a $\tilde{L}=\tilde{G}$:  d'apr\`es les formules donn\'ees, on a $^c(S\theta)_{\tilde{G}}^{\tilde{G}}=\, ^c\theta_{\tilde{G}}^{\tilde{G}}(f)=f$. Ce terme est donc $SI_{\tilde{M}}^{\tilde{G}}(\pi,\lambda,X,f)$.  Consid\'erons maintenant les termes $SI_{\tilde{M}}^{\tilde{L}}(\pi,\lambda,X,\, ^c(S\theta)_{\tilde{L}}^{\tilde{G}}(f))$, pour $\tilde{L}$  un espace de Levi propre de $\tilde{G}$. Puisque $f$ est instable, on sait que $^c(S\theta)_{\tilde{L}}^{\tilde{G}}(f)=0$ dans $SI(\tilde{L})$. Par r\'ecurrence on sait aussi que $SI_{\tilde{M}}^{\tilde{L}}(\pi,\lambda,X,.)$ est une distribution stable. Ainsi pour $\tilde{L}\neq \tilde{G}$, on a $SI_{\tilde{M}}^{\tilde{L}}(\pi,\lambda,X,\, ^c(S\theta)_{\tilde{L}}^{\tilde{G}}(f))=0$ et l'\'egalit\'e
$$
^c(SI)_{\tilde{M}}^{\tilde{G}}(\pi,\lambda,X,f)= SI_{\tilde{M}}^{\tilde{G}}(\pi,\lambda,X,f).
$$
Il suffit donc de montrer que le membre de gauche est nul. En utilisant les formules de descente, on se ram\`ene au cas o\`u $\pi$ est temp\'er\'ee modulo le centre. On \'ecrit la d\'efinition:
$
^c(SI)_{\tilde{M}}^{\tilde{G}}(\pi,\lambda,X,f)=$
$$\, ^cI_{\tilde{M}}^{\tilde{G}}(\pi,\lambda,X,f)-\sum_{s\in Z(\hat{M})^{\Gamma}/Z(\hat{G})^{\Gamma}, s\neq 1}i_{\tilde{M}}(\tilde{G},{\bf G}'(s))\, ^c(SI)_{\tilde{M}}^{\tilde{G}'(s)}(\pi,\lambda,X,f^{{\bf G}'(s)}).
$$
Ce sont des coefficients de Fourier des distributions qui se d\'eduisent de $^cI(\pi_{\lambda},f)$ de fa\c{c}on formelle. Par m\'eromorphie, il suffit donc encore de consid\'erer le cas o\`u $\pi$ est temp\'er\'ee (et non seulement temp\'er\'ee modulo le centre). On a donc par d\'efinition  l'\'egalit\'e:
$
^c(SI)_{\tilde{M}}^{\tilde{G}}(\pi_{\lambda},f)=$
$$\, ^cI_{\tilde{M}}^{\tilde{G}}(\pi_{\lambda},f)-\sum_{s\in Z(\hat{M})^{\Gamma}/Z(\hat{G})^{\Gamma}, s\neq 1}i_{\tilde{M}}(\tilde{G},{\bf G}'(s))\, ^c(SI)_{\tilde{M}}^{\tilde{G}'(s)}(\pi_{\lambda},f^{{\bf G}'(s)}).
$$
On sait que $^cI_{\tilde{M}}^{\tilde{G}}(\pi_{\lambda},f)= tr\, \pi_{\lambda}(\, ^c\theta_{\tilde{M}}^{\tilde{G}}(f))$ (cf. \ref{carpondtemp}). On admet par r\'ecurrence pour tout $s\neq 1$ comme ci-dessus que 
$$
^c(SI)_{\tilde{M}}^{\tilde{G}'(s)}(\pi_{\lambda},,f^{{\bf G}'(s)})= tr\, \pi_{\lambda}(\, ^c(S\theta)_{\tilde{M}}^{\tilde{G}'(s)}(f^{{\bf G}'(s)}))$$
puisque l'\'egalit\'e ci-dessus donne alors:
$
^c(SI)_{\tilde{M}}^{\tilde{G}}(\pi_{\lambda},f)= tr \pi_{\lambda}(F)$, avec
$$
F=\, ^c\theta_{\tilde{M}}^{\tilde{G}}(f)- \sum_{s\in Z(\hat{M})^{\Gamma}/Z(\hat{G})^{\Gamma}, s\neq 1}i_{\tilde{M}}(\tilde{G},{\bf G}'(s)) \, ^c(S\theta)_{\tilde{M}}^{\tilde{G}'(s)}(f^{{\bf G}'(s)}).
$$
Par d\'efinition $F=\, ^c(S\theta)_{\tilde{M}}^{\tilde{G}}(f)$ et on sait que si $f$ est instable, l'image de $F$ dans $SI(\tilde{M})$ est nulle. Comme $\pi_{\lambda}$ est une repr\'esentation stable par hypoth\`ese, si $f$ est instable, sa trace contre $F$ est nulle. On a donc bien montr\'e que $^c(SI)_{\tilde{M}}^{\tilde{G}}(\pi,\lambda,X,f)=0$ si $f$ est instable ce qui termine la preuve de la proposition.

\subsection{Les caract\`eres pond\'er\'es endoscopiques}
\subsubsection{D\'efinition}\label{carpondend}
On fixe ${\bf M}'$ une donn\'ee endoscopique non n\'ecessairement elliptique de $\tilde{G}$; on a besoin de la notation $\hat{s}_{M'}$ qui fait partie de la donn\'ee. On fixe aussi ${\pi}'_{M'}$ une repr\'esentation  de ${\bf M}'$.  On entend par l\`a que, si l'on fixe des donn\'ees auxiliaires $M'_{1}$, $C_{1}$ etc... $\pi'_{M'}$ s'identifie \`a une repr\'esentation de $\tilde{M}'_{1}$ qui se tranaforme selon le caract\`ere $\lambda_{1}$ de $C_{1}(F)$ d\'etermin\'e par la situation. On suppose que cette repr\'esentation est  stable. Pour toute donn\'ee endoscopique elliptique $\bf{G}'$ de $\tilde{G}$ contenant $\bf{M}'$ comme espace de Levi, on sait d\'efinir le caract\`ere pond\'er\'e invariant en tant que distribution stable  $SI_{{\bf M}'}^{{\bf G}'}(\pi'_{M'},\lambda,X,.)$. On suppose  qu'il existe un espace de Levi $\tilde{M}$ de $\tilde{G}$ tel que ${\bf M}'$ soit une donn\'ee endoscopique elliptique relevante de $\tilde{M}$. On peut alors normaliser les facteurs de transfert, en fixant les facteurs de transfert pour le couple $\tilde{M}$, ${\bf M}'$ de fa\c{c}on \`a poser  pour tout $f\in I(\tilde{G})$:
$$
I_{{\bf M}'}^{\tilde{G}, {\mathcal E}}(\pi'_{M'},\lambda,X,f):= $$
$$\sum_{\hat{s}\in \hat{s}_{M'}Z(\hat{M})^{\Gamma,\hat{\theta}}/Z(\hat{G})^{\Gamma,\hat{\theta}}\cap Z(\hat{M})^{\Gamma,\hat{\theta}}} i_{\tilde{M}'}(\tilde{G},{\bf G}'(\hat{s}))SI_{{\bf M}'}^{{\bf G}'(\hat{s})}(\pi'_{M'},\lambda,X,f^{{\bf G}'(\hat{s})}). \eqno(1)
$$
C'est ce que l'on appelle le caract\`ere pond\'er\'e endoscopique. On a aussi une variante compacte:
$$
^cI_{{\bf M}'}^{\tilde{G}, {\mathcal E}}(\pi'_{M'},\lambda,X,f):= $$
$$\sum_{\hat{s}\in \hat{s}_{M'}Z(\hat{M})^{\Gamma,\hat{\theta}}/Z(\hat{G})^{\Gamma,\hat{\theta}}\cap Z(\hat{M})^{\Gamma,\hat{\theta}}} i_{\tilde{M}'}(\tilde{G},{\bf G}'(\hat{s}))\, ^c(SI)_{{\bf M}'}^{{\bf G}'(\hat{s})}(\pi'_{M'},\lambda,X,f^{{\bf G}'(\hat{s})}). \eqno(2)
$$
\begin{lem} Avec les deux premi\`eres hypoth\`eses de r\'ecurrence g\'eom\'etriques de \ref{stabgeohyp}, 
pour tout $\lambda$ tel que la distribution suivante est d\'efinie et pour tout $X$, on a l'\'egalit\'e:
$$
^cI_{{\bf M}'}^{\tilde{G}, {\mathcal E}}(\pi'_{M'},\lambda,X,f)=\sum_{\tilde{L}\in {\mathcal L}(\tilde{M})}I_{{\bf M}'}^{\tilde{L}, {\mathcal E}}(\pi'_{M'},\lambda,X,\, ^c\theta_{\tilde{L}}^{\tilde{G}}(f))+
$$
$$
I_{{\bf M}'}^{{\bf M}'}(\pi'_{M'},\lambda,X, (\, ^c\theta_{{\bf M}'}^{\tilde{G},{\mathcal E}})(f)-(\, ^c\theta_{\tilde{M}}^{\tilde{G}}(f))^{{\bf M}'}).
$$
\end{lem}
On utilise le lemme de \ref{defcarstable} pour \'ecrire le terme de droite de (2): d'o\`u avec les m\^emes manipulations que dans cette r\'ef\'erence $$^cI_{{\bf M}'}^{\tilde{G}, {\mathcal E}}(\pi'_{M'},\lambda,X,f)=
\sum_{\tilde{L}\in {\mathcal L}(\tilde{M})}$$
$$\sum_{\hat{s}_{\tilde{L}}\in \hat{s}_{M'}Z(\hat{M})^{\Gamma,\hat{\theta}}/Z(\hat{L})^{\Gamma,\hat{\theta}}\cap Z(\hat{M})^{\Gamma,\hat{\theta}}}i_{\tilde{M}'}(\tilde{L},{\bf L}'(\hat{s}_{\tilde{L}}))\, ^c(SI)_{{\bf M}'}^{{\bf L}'(\hat{s}_{\tilde{L}})}(\pi_{M'}',\lambda,X,F(\hat{s}_{\tilde{L}})),
$$
o\`u $F(\hat{s}_{\tilde{L}}):= \sum_{\hat{s}\in \hat{s}_{\tilde{L}}Z(\hat{L})^{\Gamma,\hat{\theta}}/Z(\hat{G})^{\Gamma,\hat{\theta}}}i_{\tilde{L}'(\hat{s}_{\tilde{L}})}(\tilde{G},{\bf G}'(\hat{s}))\, ^c(S\theta)_{{\bf L}'(\hat{s}_{\tilde{L}})}^{{\bf G}'(\hat{s})}(f^{{\bf G}'(s)}).$ Ainsi $F(\hat{s}_{\tilde{L}})=\, ^c \theta^{\tilde{G},{\mathcal E}}_{{\bf L}'(\hat{s}_{\tilde{L}})}(f)$. Tant que $\tilde{L}\neq \tilde{M}$, la deuxi\`eme hypoth\`ese de r\'ecurrence faite assure que l'on conna\^{\i}t la stabilisation g\'eom\'etrique. Ainsi on conna\^{\i}t aussi l'\'egalit\'e: $$
^c \theta^{\tilde{G},{\mathcal E}}_{{\bf L}'(\hat{s}_{\tilde{L}})}(f)= \, (^c\theta_{\tilde{L}}^{\tilde{G}}(f))^{{\bf L}'(\hat{s}_{\tilde{L}})}$$
d'apr\`es \cite{stabilisationVII} 3.8 et 4.3 et \cite{stabilisationIX} 6.11, que l'on est en droit d'appliquer.
On a donc alors:
$$
^cI_{{\bf M}'}^{\tilde{G}, {\mathcal E}}(\pi'_{M'},\lambda,X,f)= \sum_{\tilde{L}\in {\mathcal L}(\tilde{M}); \tilde{L}\neq \tilde{M}}I_{{\bf M}'}^{\tilde{L},{\mathcal E}}(\pi'_{M'},\lambda,X,\, ^c\theta_{\tilde{L}}^{\tilde{G}}(f))
$$
$$+\sum_{\hat{s}\in \hat{s}_{M'}Z(\hat{M})^{\Gamma,\hat{\theta}}/Z(\hat{G})^{\Gamma,\hat{\theta}}}i_{\tilde{M}'}(\tilde{G},{\bf G}'(\hat{s}))
I_{{\bf M}'}^{{\bf M}'}(\pi'_{M'},\lambda,X,\, ^c(S\theta)_{{\bf M}'}^{{\bf G}'(\hat{s})}(f^{{\bf G}'(\hat{s})})).
$$
La premi\`ere somme de l'\'enonc\'e du lemme donne la premi\`ere somme ci-dessus except\'e le terme correspondant \`a $\tilde{L}=\tilde{M}$ qui (dans l'\'enonc\'e du lemme)  est simplement $$I_{{\bf M}'}^{{\bf M}'}(\pi'_{M'},\lambda,X,^c\theta_{\tilde{M}}^{\tilde{G}}(f)^{{\bf M}'});$$ c'est-\`a-dire que ce terme  manque dans la premi\`ere somme ci-dessus. La deuxi\`eme somme ci-dessus est exactement
\'egale \`a$$ I_{{\bf M}'}^{{\bf M}'}(\pi'_{M'},\lambda,X, (\, ^c\theta_{{\bf M}'}^{\tilde{G},{\mathcal E}})(f)).$$ On en d\'eduit le lemme
\subsubsection{Propri\'et\'es de descente des caract\`eres pond\'er\'es endoscopiques \label{descenteendoscopique}}
Il faut un analogue de la proposition 1.14 de \cite{stabilisationII} pour les caract\`eres pond\'er\'es endoscopiques \`a la place des int\'egrales orbitales pond\'er\'ees endoscopiques. On fixe ${\bf M'}$ une donn\'ee endoscopique non elliptique de $\tilde{G}$ comme dans le paragraphe pr\'ec\'edent et on fixe $\pi'$ une repr\'esentation stable de ${\bf M}'$. On fixe aussi un sous-groupe de Levi $R'$ de $M'$. On peut d\'efinir ${\bf R}'$ \`a partir de ${\bf M}'$ mais m\^eme si il existe un sous-espace de Levi $\tilde{R}$ tel que ${\bf R}'$ en soit une donn\'ee endoscopique elliptique rien n'assure que ${\bf R}'$ est relevant. On fixe une repr\'esentation stable $\sigma'$ de ${\bf R}'$.

Pour la proposition ci-dessous, on suppose que $\pi'$ est l'induite de $\sigma'$.
\begin{prop} (i) On suppose que $\tilde{R}$ existe de telle sorte que ${\bf R}'$ en soit une donn\'ee endoscopique elliptique relevante. Alors si $v$ est une place non archim\'edienne:
$$
I_{{\bf M}'}^{\tilde{G}, {\mathcal E}}(\pi',\lambda,X,f)= \sum_{\tilde{L}\in {\mathcal L}(\tilde{R})}d_{\tilde{R}}(\tilde{M},\tilde{L}) \sum_{Y } I_{{\bf R}'}^{\tilde{L},{\mathcal E}}(\sigma',\lambda,Y,f_{\tilde{L}}),
$$
o\`u l'on somme sur les  $Y\in {\mathcal{A}}_{\tilde{R}}$ ayant $X$  pour projection dans ${\mathcal{A}}_{\tilde{M}}$; si $v$ est une place archim\'edienne, la somme est une int\'egrale sur l'ensemble d\'efini de fa\c{c}on analogue.

(ii) On suppose que ${\bf R}'$ n'est pas une donn\'ee endoscopique elliptique relevante pour un sous-espace de Levi de $\tilde{M}$, alors $
I_{{\bf M}'}^{\tilde{G}, {\mathcal E}}(\pi',\lambda,X,f)=0$.
\end{prop}
Le d\'ebut de la preuve est commun \`a (i) et (ii). On \'ecrit la d\'efinition:
$$
I_{{\bf M}'}^{\tilde{G}, {\mathcal E}}(\pi',\lambda,X,f)=  \sum_{\hat{s}\in \hat{s}_{M'}Z(\hat{M})^{\Gamma,\hat{\theta}}/Z(\hat{G})^{\Gamma,\hat{\theta}}}i_{\tilde{M}'}(\tilde{G},{\bf G}'(\hat{s})) SI^{{\bf  G}'(\hat{s})}_{{\bf M}'}(\pi',\lambda,X,f^{{\bf G}'(\hat{s})}).
$$
On \'ecrit la formule de descente pour les termes du membre de droite et on obtient:
$$
\sum_{\hat{s}\in \hat{s}_{M'}Z(\hat{M})^{\Gamma,\hat{\theta}}/Z(\hat{G})^{\Gamma,\hat{\theta}}}i_{\tilde{M}'}(\tilde{G},{\bf G}'(\hat{s})) \sum_{\tilde{L}'\in {\mathcal L}^{\tilde{G}'(\hat{s})}(\tilde{R}')} $$
$$ \sum_{Y\in {\mathcal{A}}_{R'}; Y_{M'}=X}
e_{\tilde{R}'}(\tilde{M}',\tilde{L}')SI^{{\bf L}'}_{{\bf R}'}(\sigma',\lambda,Y,(f^{{\bf G}'(\hat{s})})_{{\bf L}'}).
$$
On remarque  que dans les formules ci-dessus, $(f^{{\bf G}'(\hat{s})})_{{\bf L}'}=0$ sauf si ${\bf L}'$ est une donn\'ee endoscopique elliptique pour un espace de Levi $\tilde{L}$ de $\tilde{G}$. On regroupe donc les termes correspondant \`a un espace de Levi $\tilde{L}$ fix\'e. N\'ecessairement un tel $\tilde{L}$ intervient s'il existe ${\bf L}'$ une de ses donn\'ees endoscopiques elliptiques telle que $d_{\tilde{R}'}(\tilde{M}',\tilde{L}')\neq 0$. 
Comme en loc. cite, m\^eme si ${\bf R}'$ n'est pas une donn\'ee endoscopique elliptique d'un espace de Levi de $\tilde{G}$, on pose $\hat{R}$ l'unique sous-groupe de Levi de $\hat{G}$ tel que $Z(\hat{R})^{\Gamma,\hat{\theta},0}=Z(\hat{R}')^{\Gamma,0}$.  La condition sur $L'$ se traduit alors par le fait que ${\mathcal{A}}_{\hat{R}}^{\hat{G}}$ est la somme directe de ${\mathcal{A}}_{\hat{R}}^{\hat{M}}$ et de ${\mathcal{A}}_{\hat{R}}^{\hat{L}}$; en consid\'erant les orthogonaux de ces espaces dans ${\mathcal{A}}_{\tilde{R}}^{\tilde{G}}$, on obtient le fait que ${\mathcal{A}}_{\tilde{R}}^{\tilde{G}}$ est la somme directe de ${\mathcal{A}}_{\tilde{M}}^{\tilde{G}}$ avec ${\mathcal{A}}_{\tilde{L}}^{\tilde{G}}$, ce qui est \'equivalent \`a la non nullit\'e de $d_{\hat{R}}(\tilde{M},\tilde{L})$. R\'eciproquement pour un $\tilde{L}$ satisfaisant cela et pour ${\bf L'}$ une donn\'ee endoscopique elliptique de $\tilde{L}$, on a aussi $d_{\tilde{R}'}^{\tilde{G}'}(\tilde{M}',\tilde{L}')\neq 0$ pour toute donn\'ee endoscopique elliptique ${\bf G}'$ de $\tilde{G}$ contenant ${\bf M}'$ et ${\bf L}'$. 

 L'application naturelle de $Z(\hat{M})/Z(\hat{G})$ dans $Z(\hat{R})/Z(\hat{L})$ est surjective et donne encore une application surjective (cf \ref{nonconnexite}) de $$Z(\hat{M})^{\Gamma,\hat{\theta},0} \rightarrow Z(\hat{R})^{\Gamma,\hat{\theta},0}/ Z(\hat{L})\cap Z(\hat{R})^{\Gamma,\hat{\theta},0}.\eqno(1)$$  Donc pour ${\bf L}'$ comme ci-dessus, l'ensemble des donn\'ees endoscopiques elliptiques ${\bf G}'$ contenant ${\bf L}'$ et ${\bf M}'$ est non vide.

 On se place maintenant dans la situation de (i); ici ${\bf R}'$ est une donn\'ee endoscopique elliptique relevante d'un espace de Levi $\tilde{R}$ de $\tilde{M}$; les facteurs de transfert ont \'et\'e normalis\'es \`a l'aide du couple ${\tilde{M}}, {\bf M}'$, ils impliquent une normalisation pour le couple ${\tilde{R}},{\bf R}'$ (puisque ce couple est un ''espace'' de Levi de $\tilde{M},{\bf M}'$) et donc une normalisation directe pour le couple ${\tilde{L}},{\bf L}'$ sans passer par l'un des couples $\tilde{G},{\bf G}'(\hat{s})$ et pour tout choix de $\hat{s}$, on a  $(f^{{\bf G}'(\hat{s})})_{{\bf L}'}=(f_{\tilde{L}})^{{\bf L}'(\hat{s}')}$ o\`u $\hat{s}\in \hat{s}_{M'}Z(\hat{M})^{\Gamma,\hat{\theta},0}/Z(\hat{G})\cap Z(\hat{M})^{\Gamma,\hat{\theta},0}$ a pour image $\hat{s}'$ dans $\hat{s}_{M'}Z(\hat{R})^{\Gamma,\hat{\theta},0}/ Z(\hat{L})\cap Z(\hat{R})^{\Gamma,\hat{\theta},0}$. Ainsi dans la formule de descente \'ecrite ci-dessus, ${\bf G}'(\hat{s})$ dispara\^{\i}t au profit de ${\bf L}'(\hat{s}')$ pour $\hat{s}$ et $\hat{s}'$ ayant la propri\'et\'e pr\'ec\'edente. Il y a donc un calcul de coefficients \`a faire mais qui a \'et\'e fait en \cite{stabilisationII} 1.14 pour passer en (5) et (6). Ce calcul donne exactement le r\'esultat cherch\'e comme en loc.cite.
 
 Pour (ii) on se raccroche \'evidemment \`a 1.14 de \cite{stabilisationII} qui construit un sous-groupe ${\mathcal Z}$ du noyau de l'application (1) tel que l'analogue g\'eom\'etrique de 
 $$
 \sum_{z\in {\mathcal Z}}SI^{{\bf L}'(\hat{s}'z)}_{{\bf R}'}(\sigma',\lambda,Y,(f^{{\bf G}'(\hat{s}'z)})_{{\bf L}'})\eqno(2)$$
 soit nul. C'est un subtil probl\`eme de facteurs de transfert qui n\'ecessite l'introduction des donn\'ees auxiliaires mais qui n'est pas \`a refaire puisqu'il a d\'ej\`a \'et\'e fait. Introduisons des donn\'ees auxiliaires $M'_{1}$, $C_{1}$, etc... pour la donn\'ee ${\bf M}'$, qui se restreignent en des donn\'ees $R'_{1}$, $C_{1}$, etc... Ce que montre ce qui suit (10) de 1.14 \cite{stabilisationII}, est que l'introduction de $z\in {\mathcal Z}$ induit un automorphisme de $C^{\infty}_{c,\lambda_{1}}(\tilde{R}'_{1})$   et donne donc dualement un automorphisme sur l'ensemble des repr\'esentations de $\tilde{R}'_{1}$ dont le caract\`ere central se restreint en le caract\`ere $\lambda_{1}$ de $C_{1}(F)$. C'est l'ensemble de  ce qu'on a appel\'e les repr\'esentations  de ${\bf R}'$. Ainsi (2) est en fait une somme sur un ensemble de repr\'esentations $\sigma'(z)$ d\'ependant de $z\in {\mathcal Z}$. L'automorphisme de $C^{\infty}_{c,\lambda_{1}}(\tilde{R}'_{1})$ est la multiplication par une fonction sur $\tilde{R}'$. Cette fonction est  calcul\'ee en loc.cite pr\`es des \'el\'ements elliptiques. D'o\`u l'int\'er\^et pour nous de se ramener au cas o\`u $\sigma'$ est elliptique; c'est tout \`a fait loisible quitte \`a remplacer $\tilde{R}'$ par un de ses sous-groupes de Levi qui ne sera pas plus relevant et on ne perd pas l'hypoth\`ese de stabilit\'e gr\^ace \`a la section  2 de  \cite{prolongementdescaracteres} et \cite{stabilisationIV} 2.8. L'automorphisme se lit alors sur l'ensemble des fonctions cuspidales dont le support est form\'e d'\'el\'ements elliptiques et il est montr\'e en loc. cite qu'en sommant sur $z$, on obtient 0. Dualement on a donc aussi 0 pour $\sigma'$ elliptique et c'est ce qui \'etait cherch\'e.

\subsection{La stabilisation g\'eom\'etrique et la stabilisation spectrale\label{stabgeostabspec}}
Ici on d\'emontre que la stabilisation locale g\'eom\'etrique entra\^{\i}ne la stabilisation locale spectrale. On fixe $\tilde{M}$ un espace de Levi de $\tilde{G}$ et ${\bf M}'$ une donn\'ee endoscopique elliptique de ${\tilde{M}}$. On fixe aussi une repr\'esentation stable $\pi'_{M'}$ de $M'$ et on note $\pi$ la $\omega$-repr\'esentation de $\tilde{M}$ que l'on obtient par transfert.
\begin{cor}(i) Sous les deux premi\`eres hypoth\`eses de 3.5 (avec $\tilde{M}$ fix\'e comme en loc.cite) et sous l'hypoth\`ese que   pour toutes les int\'egrales orbitales pond\'er\'ees  de la forme $I_{\tilde{L}}^{\tilde{G}}$ et $I_{\tilde{L}}^{\tilde{G},{\cal E}}$  sont \'egales, pour $\tilde{L}\in {\mathcal L}(\tilde{M})$ y compris $\tilde{L}=\tilde{M}$. Alors on a la stabilisation des caract\`eres pond\'er\'es:
$$
I_{{\bf M}'}^{\tilde{G},{\mathcal E}}(\pi'_{M'},\lambda,X,f)=I_{\tilde{M}}^{\tilde{G}}(\pi,\lambda,X,f).
$$ 
(ii) On ne suppose que l'hypoth\`ese de stabilisation g\'eom\'etrique pour $\tilde{G}$ remplac\'e par ses espaces de Levi propres. Alors la conclusion est vraie pour les repr\'esentations $\pi'_{M'}$    de caract\`ere central unitaire et pour $\lambda$ dans un petit voisinage de l'axe unitaire.
\end{cor}
Le (ii) est un point cl\'e de la stabilisation et r\'esulte d'une astuce (cf. la preuve) remarqu\'ee par J. Arthur

(i) L'hypoth\`ese de (i) assure que $^c\theta_{{\bf M}'}^{\tilde{G},{\mathcal E}}(f)= (\, ^c\theta_{\tilde{M}}^{\tilde{G}}(f))^{{\bf M}'}$ et donc dans le lemme de \ref{carpondend}, on a simplement l'\'egalit\'e
$$
^cI_{{\bf M}'}^{\tilde{G},{\mathcal E}}(\pi'_{M'},\lambda,X,f)=\sum_{\tilde{L}\in {\mathcal L}(\tilde{M})} I_{{\bf M}'}^{\tilde{L},{\mathcal E}}(\pi'_{M'},\lambda,X,\, ^c\theta_{\tilde{L}}^{\tilde{G}}(f)).
$$
On rappelle aussi la formule:
$$
^cI_{\tilde{M}}^{\tilde{G}}(\pi,\lambda,X,f)=\sum_{\tilde{L}\in {\mathcal L}(\tilde{M})} I_{\tilde{M}}^{\tilde{L}}(\pi,\lambda,X,\, ^c\theta_{\tilde{L}}^{\tilde{G}}(f)).
$$
Ainsi le lemme est  \'equivalent \`a l'\'egalit\'e
$$^cI_{{\bf M}'}^{\tilde{G},{\mathcal E}}(\pi'_{M'},\lambda,X,f)=\, ^cI_{\tilde{M}}^{\tilde{G}}(\pi,\lambda,X,f).
$$
On doit donc montrer l'\'egalit\'e des fonctions m\'eromorphes en $\lambda$: $$
\sum_{\hat{s}\in \hat{s}_{M'}Z(\hat{M})^{\Gamma,\hat{\theta}}/Z(\hat{G}^{\Gamma,\hat{\theta}})}\, i_{\tilde{M}'}(\tilde{G},{\bf {G}}'(\hat{s}))\,^c(SI)_{{\bf M}'}^{{\bf G}'(\hat{s})}(\pi'_{M',\lambda},f^{{\bf G}'(s)})=\, ^cI_{\tilde{M}}^{\tilde{G}}(\pi_{\lambda},f).\eqno(1)
$$
On se ram\`ene par les formules de descente au cas o\`u $\pi'_{M'}$  et donc $\pi$ sont temp\'er\'ees; on a alors le droit de supposer que $\lambda$ est unitaire par m\'eromorphie et on sait alors que le membre de gauche vaut:
$$
\sum_{\hat{s}\in \hat{s}_{M'}Z(\hat{M})^{\Gamma,\hat{\theta}}/Z(\hat{G})^{\Gamma,\hat{\theta}}}\,  i_{\tilde{M}'}(\tilde{G},{\bf {G}}'(\hat{s}))\, tr\, \pi'_{M',\lambda}(\, ^c(S\theta)_{{\bf M}'}^{{\bf G}'(\hat{s})}(f^{{\bf G}'(\hat{s})}))$$ ce qui vaut aussi par d\'efinition $ tr\, \pi'_{M',\lambda}(\, ^c\theta_{{\bf M}'}^{\tilde{G},{\mathcal E}}(f))$.
Et avec l'hypoth\`ese faite cela vaut $tr\, \pi'_{M',\lambda}(\, ^c\theta_{\tilde{M}}^{\tilde{G}}(f)^{{\bf M}'})$. Par d\'efinition du transfert de $\pi'_{M'}$ cela n'est autre que $tr\, \pi_{\lambda}(\, ^c\theta_{\tilde{M}}^{\tilde{G}}(f))$. Ce qui est le membre de droite de (1). Cela termine la preuve de (i)

(ii) On d\'ecompose $\pi'_{M'}$ dans son groupe de Grothendieck; le caract\`ere central de $\pi'_{M'}$ est  unitaire par hypoth\`ese et cela force $\pi'_{M'}$ d'\^etre la somme d'une repr\'esentation temp\'er\'ee avec des induites propres. Pour les repr\'esentations temp\'er\'ees les caract\`eres pond\'er\'es invariants sont nuls et l'\'egalit\'e cherch\'ee est triviale; il reste le cas des induites propres, en acceptant que $\lambda$ soit dans un voisinage de $0$ puisqu'il peut y avoir un probl\`eme de d\'efinition en certains points de l'axe imaginaire. Les formules de descente ram\`enent \`a d\'emontrer l'assertion cherch\'ee pour $\tilde{G}$ remplac\'e par des espaces de Levi propres; pour les espaces de Levi propres l'hypoth\`ese faite en (i) est v\'erifi\'ee et on a donc le r\'esultat cherch\'e.

\begin{rmq} Dans (ii) si on suppose $\pi'_{M'}$ unitaire, alors $I_{{\bf M}'}^{\tilde{G},{\mathcal E}}(\pi'_{M'},\lambda,X,f)=I_{\tilde{M}}^{\tilde{G}}(\pi,\lambda,X,f)$ pour tout $\lambda$ unitaire.
\end{rmq}
Cela vient du fait que les distributions sont alors d\'efinies en tout $\lambda$ unitaire et sont localement constantes en $\lambda$.
\subsection{Caract\`eres pond\'er\'es semi-globaux\label{carponddef}}
Ici le corps de base est global, c'est-\`a-dire un corps de nombres. On fixe $V$ un ensemble fini de places de $F$ contenant au moins  une place archim\'edienne. On fixe $\tilde{M}$ un espace de Levi de $\tilde{G}$ et $\pi_{V}$ une $\omega$-repr\'esentation de $\tilde{M}(F_{V})$. On note ${\mathcal L}(\tilde{M})$ l'ensemble des sous-espaces de Levi de $\tilde{G}$ contenant $\tilde{M}$. Ceci est un objet global. On d\'efinit une application:
$$
\phi_{\tilde{M},V}^{\tilde{G}}: I(\tilde{G}(F_{V})) \mapsto I_{ac}(\tilde{M}(F_{V}))$$
comme en \cite{stabilisationVI} 1.6 et 1.7. ce qui permet de d\'efinir un caract\`ere pond\'er\'e invariant semi-global: $\qquad \forall f\in I(\tilde{G}(F_{V}))$
$$ I_{\tilde{M}}^{\tilde{G}}(\pi_{V},\lambda,X,f):=J_{\tilde{M}}^{\tilde{G}}(\pi_{V},\lambda,X,f)-\sum_{\tilde{L}\in {\mathcal L}(\tilde{M})} I_{\tilde{M}}^{\tilde{L}}(\pi_{V},\lambda,X,\phi_{\tilde{L},V}^{\tilde{G}}(f)).\eqno(1)$$
Dans la d\'efinition de $J_{\tilde{M}}^{\tilde{G}}(\pi_{V},\lambda,X,f)$ il y a une int\'egrale sur $\lambda+i{\mathcal{A}}_{\tilde{M}}^*$ qui est un espace affine car on a suppos\'e que  $V$ contient au moins une place archim\'edienne.

On peut remplacer l'hypoth\`ese que $V$ contient au moins une place archim\'edienne par l'hypoth\`ese que $V$ ne contient que des places ayant m\^eme caract\'eristique r\'esiduelle; on int\`egre alors sur un tore si la place est p-adique. On peut donc d\'efinir de fa\c{c}on semi-globale $I_{\tilde{M}}^{\tilde{G}}(\pi_{v},\lambda,X,f_{v})$ pour toute place $v$ en consid\'erant que $V=\{v\}$. Montrons que  cette d\'efinition s'exprime en fonction de la d\'efinition locale des paragraphes pr\'ec\'edents ainsi:
$$
I_{\tilde{M}}^{\tilde{G}}(\pi_{v},\lambda,X,f_{v})=$$
$$\sum_{\tilde{L}({v})\in {\mathcal L}(\tilde{M}_{v})}d_{\tilde{M}_{v}}(\tilde{M},\tilde{L}(v))\int_{Y\in {\mathcal{A}}_{\tilde{M}_{v}}; Y_{\tilde{M}}=X}dY\,I_{\tilde{M}_{v}}^{\tilde{L}({v})}(\pi_{v},\lambda,X,f_{v,\tilde{L}({v})}),\eqno(2)
$$o\`u $Y_{\tilde{M}}$ est la projection orthogonale de $Y\in {\mathcal{A}}_{\tilde{M}_{v}}$ dans ${\mathcal{A}}_{\tilde{M}}$.

\noindent En effet, pour prouver une telle formule, il faut d'abord la prouver pour $J_{\tilde{M}}^{\tilde{G}}(\pi_{v},\lambda,X,f_v)$. Cette distribution s'obtient par transformation de Fourier \`a partir de la distribution $J_{\tilde{M}}^{\tilde{G}}(\pi_{v,\lambda},f_v)$. Par les propri\'et\'es de $\tilde{G},\tilde{M}$ famille, on a l'\'egalit\'e de fonctions m\'eromorphes en $\lambda\in {\mathcal{A}}^*_{\tilde{M}}$:
$$
J_{\tilde{M}}^{\tilde{G}}(\pi_{v,\lambda},f_{v})= \sum_{\tilde{L}(v)\in {\mathcal L}(\tilde{M}_{v})} d_{\tilde{M}_{v}}(\tilde{M},\tilde{L}(v))J_{\tilde{M}_{v}}^{\tilde{L}(v)}(\pi_{v,\lambda},f_{v,\tilde{L}(v)}).
$$
Le terme de gauche doit \^etre int\'egr\'e sur $\lambda+i{\mathcal{A}}^*_{\tilde{M}}$ contre $e^{-\lambda(X)}$ pour trouver $J_{\tilde{M}}^{\tilde{G}}(\pi_{v},\lambda,X,f_{v})$. Du c\^ot\'e droit, il faut int\'egrer sur $\lambda+i{\mathcal{A}}^*_{\tilde{M}_{v}}$ contre $e^{-\lambda(X)}$ pour trouver $J_{\tilde{M}_{v}}^{\tilde{L}(v)}(\pi_{v},\lambda,X,f_{v,\tilde{L}(v)})$. On passe de l'une des int\'egrales \`a l'autre en faisant pr\'ecis\'ement une transform\'ee de Fourier sur ${\mathcal{A}}_{\tilde{M}_{v}}/{\mathcal{A}}_{\tilde{M}}$. D'o\`u le r\'esultat annonc\'e.

Ensuite, on passe aux termes en $I$ au lieu de $J$ par un calcul standard utilisant les formules de descente. On ne le fait pas.

\

On revient \`a l'ensemble $V$; on ram\`ene la d\'efinition (1) \`a une d\'efinition locale, par r\'ecurrence sur le nombre de places dans $V$; pour simplifier on suppose que $\pi_{V}$ est un produit tensoriel, $\otimes_{v\in V}\pi_{v}$; on se ram\`ene \`a ce cas en faisant des combinaisons lin\'eaires. 

Pour cela  on d\'efinit aussi ${\mathcal L}(\tilde{M}_{V})$ comme l'ensemble des familles $\tilde{L}^V=(\tilde{L}^v;v\in V)$ o\`u pour tout $v\in V$, $\tilde{L}^v\in {\mathcal L}(\tilde{M}(F_{v}))$ et on a les formules de scindage qui relient cette d\'efinition aux d\'efinitions locales: pour toute fonction $f_V\in I(\tilde{G}(F_V))$ produit de ses composantes locales:
$$I_{\tilde{M}}^{\tilde{G}}(\pi_{V},\lambda,X,f_V)=$$
$$\sum_{\{\tilde{L}^v\} \in {\mathcal L}(\tilde{M}_{V})}
d_{\tilde{M}}(\tilde{G},\{\tilde{L}^v\})\int_{\{X_{v}\in {\mathcal{A}}_{\tilde{M}_{v}}; \sum_{v}X_{v,\tilde{M}}=X\}}
\prod_{v\in V}I_{\tilde{M}_{v}}^{\tilde{L}^v}(\pi_{v},\lambda,X_{v},f_{v;\tilde{L}_{v}})dX_{v}.
$$ 
\subsection{Caract\`eres pond\'er\'es semi-globaux et endoscopie, th\'eor\`eme d'annulation\label{carpondendsemiglob}}
On d\'efinit les variantes stables et endoscopiques de caract\`eres pond\'er\'es semi-globaux.

On suppose d'abord que ${\omega}=1$ que $G$ est quasi-d\'eploy\'e et que $\tilde{G}$ est \`a torsion int\'erieure. On pose, pour toute fonction  $f_V\in I(\tilde{G}(F_V))$ produit de ses composantes locales:
$$SI_{\tilde{M}}^{\tilde{G}}(\pi_{V},\lambda,X,f_V)=\sum_{\{\tilde{L}^v\} \in {\mathcal L}(\tilde{M}_{V})}
e_{\tilde{M}_{V}}(\tilde{M},\{\tilde{L}^v\})$$
$$\int_{\{X_{v}\in {\mathcal{A}}_{\tilde{M}_{v}}; \sum_{v}X_{v,\tilde{M}}=X\}}
\prod_{v\in V}SI_{\tilde{M}_{v}}^{\tilde{L}^v}(\pi_{v},\lambda,X_{v},f_{v;\tilde{L}_{v}})dX_{v}. \eqno(1)
$$ Avec une telle d\'efinition, la distribution est certainement stable parce que c'est le cas des distributions locales. Et cela est conforme \`a la d\'efinition usuelle (cf. \cite{stabilisationVI} (i) de la proposition 4.2)
$$
SI_{\tilde{M}}^{\tilde{G}}(\pi_{V},\lambda,X,f_V)=I_{\tilde{M}}^{\tilde{G}}(\pi_{V},\lambda,X,f_V)-\sum_{s\in Z(\hat{M})^{\Gamma}/Z(\hat{G})^{\Gamma}; s\neq 1} SI_{\tilde{M}}^{\tilde{G}'(s)}(\pi_{V},\lambda,X,f^{{\bf G}'(s)}_V).
$$

On ne fait plus d'hypoth\`ese sur $\tilde{G}$; soit   ${\bf M}'$ une donn\'ee endoscopique non n\'ecessairement elliptique ni m\^eme relevante de $\tilde{G}$. On fixe une repr\'esentation $\pi'_{V}$, stable, de ${\bf M}'$ sur $F_{V})$; c'est n\'ecessairement par d\'efinition de la stabilit\'e, une combinaison lin\'eaire de repr\'esentations stables elles m\^emes produit tensoriel de repr\'esentations stables en toute place $v\in V$.  Dans $\hat{G}$, on note $\hat{M}$ le sous-groupe de Levi tel que $Z(\hat{M})^{\Gamma,\hat{\theta},0}=Z(\hat{M}')^{\Gamma,0}$; plus pr\'ecis\'ement $\hat{M}$ est la composante neutre du centralisateur dans $\hat{G}$ de $Z(\hat{M}')^{\Gamma,0}$. 
On pose alors pour toute fonction  $f_V\in I(\tilde{G}(F_V))$ produit de ses composantes locales:
$
I_{{\bf M}'}^{\tilde{G},{\mathcal E}}(\pi'_{V},\lambda,X,f_V):=$
$$
\sum_{\hat{s}\in \hat{s}_{M'}Z(\hat{M})^{\Gamma,\hat{\theta}}/Z(\hat{G})^{\Gamma,\hat{\theta}}}i_{\tilde{M}'}(\tilde{G},{\bf G}'(\hat{s}))(SI)_{{\bf M}'}^{{\bf G}'(\hat{s})}(\pi'_{V},\lambda,X,f^{{\bf G}'(\hat{s})}_V).
$$
Les facteurs de transfert sont normalis\'es globalement.  Dans le cas o\`u ${\bf M}'$ est une donn\'ee endoscopique elliptique d'un espace de Levi $\tilde{M}$ de $\tilde{G}$,  on a la d\'ecompostion, sous l'hypoth\`ese que $\pi'_{V}$ est un produit tensoriel:
$$
I_{\bf{M}'}^{\tilde{G},{\mathcal E}}(\pi'_{V},\lambda,X,f_V)=\sum_{\{\tilde{L}^v\} \in {\mathcal L}(\tilde{M}_{V})}
d_{\tilde{M}_{V}}(\tilde{M},\{\tilde{L}^v\})$$
$$\int_{\{X_{v}\in {\mathcal{A}}_{\tilde{M}_{v}}; \sum_{v}X_{v,\tilde{M}}=X\}}
\prod_{v\in V}I_{{\bf M}'_{v}}^{\tilde{L}^v,{\mathcal E}}(\pi'_{v},\lambda,X_{v},f_{v;\tilde{L}_{v}})dX_{v}. \eqno(2)
$$
On renvoie au (i) de la proposition de  \cite{stabilisationVI} 4.5  pour la preuve de cette formule, c'est juste un calcul de coefficients fait en loc. cite.

\begin{lem} Avec les notations pr\'ec\'edentes, on suppose qu'en toute place $v\in V$ la repr\'esentation $\pi'_{v}$ est une induite de la forme $ind_{{\bf R}'(v)}^{{\bf M}'_{v}}\sigma_{v}$ avec $\sigma_{v}$ une repr\'esentation stable et  elliptique modulo le centre  de ${\bf R}'({v})$ o\`u ${\bf R}'(v)$ est un espace de Levi de ${\bf M}'_{v}$.
On suppose aussi que pour au moins une place $v$, ${\bf R}'(v)$ n'est pas une donn\'ee endoscopique relevante de $\tilde{G}(F_{v})$, alors
$$
I_{{\bf M}'}^{\tilde{G},{\mathcal E}}(\pi'_{V},\lambda,X,f_V)=0$$
\end{lem}
On a utilis\'e l'\'ecriture ${\bf R}'(v)$ pour signifier que ${\bf R}'(v)$ est un espace de Levi de ${\bf M}'_{v}$ qui ne provient pas par localisation d'une donn\'ee globale.

C'est bien l'analogue g\'eom\'etrique de cet \'enonc\'e qui est d\'emontr\'e en 6.10, (31) de \cite{stabilisationVI} mais comme ce n'est pas exactement l'hypoth\`ese faite en \cite{stabilisationVI} 6.1 et suivant on commence par s'y ramener. On commence donc par supposer que ${\bf M}'$ est une donn\'ee endoscopique relevante de $\tilde{G}$. Et on montre la nullit\'e: on consid\`ere (2) et on fixe $v$ tel que ${\bf R'}(v)$ ne soit pas une donn\'ee endoscopique relevante de $\tilde{G}(F_{v})$. On a alors montr\'e en \ref{descenteendoscopique} (ii) que les termes $I_{{\bf M}'_{v}}^{\tilde{L}^v,{\mathcal E}}(\pi'_{v},\lambda,X_{v},f_{v,\tilde{L}_{v}})$ sont tous nuls. D'o\`u l'assertion dans ce cas.

Maintenant on a exactement l'hypoth\`ese que ${\bf M}'$ n'est pas une donn\'ee endoscopique relevante de $\tilde{G}$ comme dans 6.6 de \cite{stabilisationVI} et la m\^eme preuve s'applique avec beaucoup moins (voire pas du tout) de difficult\'es aux places archim\'ediennes.

\subsection{Caract\`eres pond\'er\'es semi-globaux et endoscopie, th\'eor\`eme de transfert\label{transfertcarpond}}
Dans ce paragraphe on fait l'hypoth\`ese de r\'ecurrence g\'eom\'etrique de \ref{redlocgeo}. On fixe $\tilde{M}$ un espace de Levi de $\tilde{G}$ et ${\bf M}'$ une donn\'ee endoscopique elliptique relevante de $\tilde{M}$. On fixe une repr\'esentation $\pi'_{V}$, stable, de ${\bf M}'$ dont on note $\pi_{V}$ le transfert en une $\omega$-repr\'esentation de $\tilde{M}$. On suppose que ces repr\'esentations ont un caract\`ere central.
\begin{prop} Avec les hypoth\`eses et notations pr\'ec\'edentes et en supposant que le caract\`ere central de  $\pi'_{V}$ est unitaire,
on a l'\'egalit\'e des distributions pour tout $\lambda$ tr\`es voisin de $0$ o\`u ces distributions sont d\'efinies
$$\forall f_{V}\in I(\tilde{G}(F_{V})), \qquad I_{\tilde{M}}^{\tilde{G}}(\pi_{V},\lambda,X,f_{V})=I_{{\bf M}'}^{\tilde{G},{\mathcal E}}(\pi'_{V},\lambda,X,f_{V}).
$$
\end{prop}
On se ram\`ene ais\'ement au cas o\`u $\pi_{V}$ est un produit tensoriel. On compare les d\'efinitions de \ref{carponddef} et \ref{carpondendsemiglob} (2). Il suffit donc de montrer avec les notations de ces r\'ef\'erences que pour tout $\tilde{L}^v\in {\mathcal L}(\tilde{M}_{v})$, on a l'\'egalit\'e $$\forall f\in I(\tilde{G}(F_{v}))\quad I_{{\bf M}'_{v}}^{\tilde{L}^v,{\mathcal E}}(\pi'_{v},\lambda,X_{v},f_{v})= I_{\tilde{M}_{v}}^{\tilde{L}^v}(\pi_{v},\lambda,X_{v},f_{v}).$$ Ceci a \'et\'e montr\'e en \ref{stabgeostabspec} (ii) parce que l'on a suppos\'e que  $\pi'_{V}$ a un caract\`ere central unitaire, pour tout $\lambda$ dans un voisinage de $0$ comme dans l'\'enonc\'e. 

\begin{rmq}
Les fonctions, en $\lambda$,  de l'\'enonc\'e sont localement constantes. Si l'une est d\'efinie en $\lambda=0$ l'autre l'est aussi au moins par continuit\'e. Donc avec la seule hypoth\`ese que soit $\pi'_{V}$ soit $\pi_{V}$ est somme de repr\'esentations unitaires, on aura une \'egalit\'e de distribution en $\lambda=0$. On supprime alors $\lambda$ de la notation.
\end{rmq}
Une difficult\'e mineure de cette th\'eorie est que l'on ne sait pas que le transfert pr\'eserve l'unitarit\'e, il n'y a m\^eme sans doute pas de raison que ce soit vrai, d'o\`u l'int\'er\^et de la remarque. 
\subsection{Caract\`eres pond\'er\'es globaux}
\subsubsection{D\'efinition des caract\`eres pond\'er\'es globaux\label{defcarglob}}
On fixe $\tilde{M}$ un espace de Levi de $\tilde{G}$.  On fixe un ensemble fini de places $V$ contenant $V_{ram}$. Pour tout $v\notin V$, on fixe un morphisme $c_{v}: W_{F_{v}} \rightarrow \, ^LM(F_{v})$  (en notant ainsi le $L$-groupe de $M$ vu sur $F_{v}$).

On a d\'efini en \ref{caracpondnonramifdef} la fonction m\'eromorphe de $\lambda\in {\mathcal{A}}^*_{\tilde{M},{\mathbb C}}$, $r_{\tilde{M}}^{\tilde{G}}(c_{\lambda}^V)$; cette fonction est invariante sous ${\mathcal{A}}^*_{\tilde{G},{\mathbb C}}$. On g\'en\'eralise cette d\'efinition en rempla\c{c}ant $\tilde{G}$ par un espace de Levi $\tilde{L}\in {\mathcal L}(\tilde{M})$. Et la fonction m\'eromorphe de $\lambda\in {\mathcal{A}}^*_{\tilde{M},{\mathbb C}}$, not\'ee $r_{\tilde{M}}^{\tilde{L}}(c^V_{\lambda})$, est invariante sous ${\mathcal{A}}^*_{\tilde{L},{\mathbb C}}$ et elle ne d\'epend donc que de la projection $\lambda^{\tilde{L}}$ de $\lambda$ dans ${\mathcal{A}}^{*\tilde{L}}_{\tilde{M},{\mathbb C}}$. On suppose que les repr\'esentations param\'etr\'ees par $c_{v}$ pour tout $v\notin V$ sont unitaires. La fonction m\'eromorphe ainsi d\'efinie est donc holomorphe en tout $\lambda$ tel que $\lambda^{\tilde{L}}$ est  unitaire.

On fixe $\lambda\in {\mathcal{A}}^*_{\tilde{M},{\mathbb C}}$; pour $\tilde{L}\in {\mathcal L}(\tilde{M})$, on d\'ecompose $\lambda=\lambda^{\tilde{L}}+\lambda_{\tilde{L}}$ suivant la d\'ecomposition orthogonale ${\mathcal{A}}^*_{\tilde{M},{\mathbb C}}= {\mathcal{A}}^*_{\tilde{L},{\mathbb C}}\oplus {\mathcal{A}}^{*\tilde{L}}_{\tilde{M},{\mathbb C}}$.
Dans les paragraphes pr\'ec\'edents, pour toute $\omega$ repr\'esentation $\pi_{V}$ de $\tilde{G}(F_{V})$, on a d\'efini la distribution $$f_{V}\in I(\tilde{G}(F_{V})) \mapsto I_{\tilde{L}}^{\tilde{G}}(ind_{\tilde{M}}^{\tilde{L}}(\pi_{V}\otimes \lambda^{\tilde{L}}),\lambda_{\tilde{L}},X_{\tilde{L}},f_{V}),$$ o\`u $X_{\tilde{L}}\in {\mathcal{A}}_{\tilde{L}}$.
Cette distribution n'est pas partout d\'efinie, elle est localement constante en $\lambda_{\tilde{L}}$ (l\`a o\`u elle est d\'efinie) et elle d\'epend m\'eromorphiquement de $\lambda^{\tilde{L}}$: cela r\'esulte des formules de descente. Si $\pi_{V}$ est unitaire la distribution est d\'efinie pour $\lambda$ unitaire. On peut donc la calculer en $\lambda_{\tilde{L}}=0$.

On suppose que $\pi_{V}$ est unitaire et on revient \`a $c^V$ que l'on suppose aussi unitaire; pour $\lambda\in {\mathcal{A}}^*_{\tilde{M},{\mathbb C}}$ on d\'efinit $r_{\tilde{M}}^{\tilde{L}}(c_{\lambda}^V)$. Cette fonction m\'eromorphe de $\lambda$ ne d\'epend que de $\lambda^{\tilde{L}}$. On d\'efinit alors pour tout
$ f_{V}\in I(\tilde{G}(F_{V}))$
$$
I_{\tilde{M}}^{\tilde{G}}(\pi_{V}\otimes c^V,f)=\sum_{\tilde{L}\in {\mathcal L}(\tilde{M})}\int_{i{\mathcal{A}}^*_{\tilde{M}}/i{\mathcal{A}}^*_{\tilde{L}}}d\lambda^{\tilde{L}}\, r_{\tilde{M}}^{\tilde{L}}(c^V_{\lambda}) I_{\tilde{L}}^{\tilde{G}}(ind_{\tilde{M}}^{\tilde{L}} (\pi_{V}\otimes \lambda^{\tilde{L}}),0,0,f_{V}). \eqno(1)
$$Cette int\'egrale converge gr\^ace \`a \ref{croissance}. C'est le caract\`ere pond\'er\'e global.

\subsubsection{Caract\`eres pond\'er\'es globaux stables (cas de la torsion int\'erieure\label{carglobstab})}
Dans ce paragraphe, on suppose que $\tilde{G}$ est \`a torsion int\'erieure avec ${\omega}=1$ et $G$ quasid\'eploy\'e. On d\'efinit une version stable des caract\`eres pond\'er\'es globaux. On fixe encore $\tilde{M}$,   $c^V$  et $\pi_{V}$ comme dans le paragraphe pr\'ec\'edent. 
Pour $\lambda\in {\mathcal{A}}^*_{\tilde{M},{\mathbb C}}$ et pour $\tilde{L}\in {\mathcal L}(\tilde{M})$, on consid\`ere encore la d\'ecomposition $\lambda=\lambda_{\tilde{L}}+\lambda^{\tilde{L}}$. Fixons $\tilde{L}\in {\mathcal L}(\tilde{M})$. On a d\'efini en \ref{stabnonramif} les fonctions  m\'eromorphes de $\lambda$, $s_{\tilde{M}}^{\tilde{G}}(c^V_{\lambda})$ que l'on consid\`ere pour $\tilde{G}$ remplac\'e par $\tilde{L}$, ce sont donc des fonctions $s_{\tilde{M}}^{\tilde{L}}(c^V_{\lambda})$ qui ne d\'ependent  que de $\lambda^{\tilde{L}}$. Et on pose, avec les notations pr\'ec\'edentes
$$
(SI)_{\tilde{M}}^{\tilde{G}}(\pi_{V}\otimes c^V,f)=$$
$$\sum_{\tilde{L}\in {\mathcal L}(\tilde{M})}\int_{i{\mathcal{A}}^*_{\tilde{M}}/i{\mathcal{A}}^*_{\tilde{L}}}d\lambda^{\tilde{L}}\, s_{\tilde{M}}^{\tilde{L}}(c^V_{\lambda^{\tilde{L}}}) (SI)_{\tilde{L}}^{\tilde{G}}(ind_{\tilde{M}}^{\tilde{L}}( \pi_{V}\otimes \lambda^{\tilde{L}}),0,0,f_{V}).\eqno(2)
$$
Il est clair d'apr\`es \ref{carpondendsemiglob} qu'une telle  distribution est stable. 
\subsubsection{Caract\`eres pond\'er\'es globaux endoscopiques, transfert \label{carglobtransfert}}
On ne fait plus d'hypoth\`eses sur $\tilde{G}$ et on d\'efinit aussi les  variantes endoscopiques des caract\`eres pond\'er\'es globaux. On fixe $\tilde{M}$ un espace de Levi de $\tilde{G}$. Soit ${\bf M}'$ une donn\'ee endoscopique elliptique de $\tilde{M}$ et $\pi'_{V}$ une repr\'esentation unitaire de ${\bf M}'(F_{V})$ que l'on suppose stable. On suppose aussi donn\'e $c^{'V}$ produit tensoriel pour tout $v\notin V$ de morphisme de $W_{F_{v}}$ dans ${\mathcal M}'$.  Et on suppose aussi que $c^{'V}$ est unitaire.
Par inclusion de ${\mathcal M}'$ dans le $L$-groupe de $M$, on obtient un morphisme du type de ceux consid\'er\'es pr\'ec\'edemment et on le note encore $c^{'V}$. Si ${\bf M}'$ n'est pas relevant, on pose $\pi_{V}=0$ et sinon on note $\pi_{V}$ la repr\'esentation de $\tilde{M}(F_{V})$ obtenue \`a partir de $\pi'_{V}$ par transfert. On pose pour tout
$f_{V}\in I(\tilde{G}(F_{V}))$
$$
I_{{\bf M}'}^{\tilde{G},{\mathcal E}}(\pi'_{V}\otimes c^{'V},f_{V})=$$
$$ \sum_{\hat{s}\in \hat{s}_{M'}Z(\hat{M})^{\Gamma,\hat{\theta}}/Z(\hat{G})^{\Gamma,\hat{\theta}}} i_{\tilde{M}'}(\tilde{G},{\bf G}'(\hat{s})) (SI)_{{\bf M}'}^{{\bf G}'(\hat{s})}(\pi'_{V}\otimes c^{'V},f_{V}^{{\bf G}'(\hat{s})}).\eqno(3)
$$
\begin{prop} Avec les hypoth\`eses et notations pr\'ec\'edentes, pour tout $f_V\in I(\tilde{G}(F_{V}))$ on a l'\'egalit\'e:
$$
I_{{\bf M}'}^{\tilde{G},{\mathcal E}}(\pi'_{V}\otimes c^{'V},f_{V})=I_{\tilde{M}}^{\tilde{G}}(\pi_{V}\otimes c^{'V},f_{V}).
$$ En particulier si ${\bf M}'$ n'est pas relevant le terme de gauche est nul.
\end{prop}
Il faut montrer que le terme de gauche v\'erifie une propri\'et\'e de descente analogue \`a \ref{defcarglob} (1). Pour  ${\bf L}'$ une donn\'ee endoscopique de $\tilde{G}$ on dit que ${\bf L}'\in {\mathcal L}({\bf M}')$ si ${\bf L}'=(L',{\mathcal L}',\hat{s}_{L'})$ et ${\bf M}'=(M',{\mathcal M}',\hat{s}_{M'})$ sont telles que $M'$ est un sous-groupe de Levi de $L'$, $\hat{s}_{L'}\in \hat{s}_{M'}Z(\hat{M})^{\Gamma,\hat{\theta}}$ et ${\mathcal L}'=\hat{L}'{\mathcal M}'$. A une telle donn\'ee on associe le sous-groupe de Levi ${\hat{L}}$ de ${\hat{G}}$ d\'efini par le fait que $\hat{L}$ est le commutant dans $\hat{G}$ de $Z(\hat{L}')^{\Gamma,0}$. A chaque terme $(SI)_{{\bf M}'}^{{\bf G}'(\hat{s})}(\pi'_{V}\otimes c^{'V},f_{V}^{{\bf G}'(\hat{s})})$ on applique la d\'efinition \ref{carglobstab} (2), ce qui fait intervenir une somme sur les ${\bf L}'$ contenant ${\bf M}'$ et inclus dans ${\bf G}'(\hat{s})$. On peut modifier l'ordre des sommes (cf \cite{stabilisationVI}, 6.6) pour avoir une somme sur les sous-groupes de Levi, ${\hat{L}}$ de ${\hat{G}}$ contenant ${\hat{M}}$ et sur les donn\'ees endoscopiques  ${\bf L}'$ donnant ${\hat{L}}$ puis sur les donn\'ees ${\bf G}'(\hat{s})$ contenant ${\bf L}'$: comme toujours la somme en $\hat{s}$ se d\'ecompose en une somme sur $\hat{s}_{L'}\in \hat{s}_{M'}Z(\hat{M})^{\Gamma,\hat{\theta}}/Z(\hat{L})^{\Gamma,\hat{\theta}}$ et $\hat{s}\in \hat{s}_{L'}Z(\hat{L})^{\Gamma,\hat{\theta}}/Z(\hat{G})^{\Gamma,\hat{\theta}}$. On a d\'ej\`a vu que quand on fixe ${\bf L}'$ la deuxi\`eme somme donne le terme (5) ci-dessous et on obtient donc  que le terme correspondant \`a un ${\bf L}'$ fix\'e est $$i_{\tilde{M}'}(\hat{L},\hat{L}')
\int_{\lambda^{L'}\in i{\mathcal{A}}^*_{M'}/i{\mathcal{A}}^*_{L'}}d\lambda^{L'}\, s_{M'}^{L'}(c^{'V}_{\lambda^{{L}'}})\eqno(4)
$$
$$
I_{{\bf L}'}^{\tilde{G},{\mathcal E}}(ind_{{\bf M}'}^{{\bf L}'}(\pi'_{V}\otimes \lambda^{L'}), 0,0,f_{V}).\eqno(5)
$$
Si ${\bf M}'$ n'est pas relevant tous les termes (5) sont nuls d'apr\`es 4.7. On est donc ramen\'e au cas o\`u ${\bf M}'$ est relevant. Dans ce cas l\`a, les ${\bf L}'$ sont aussi relevants et on note $\tilde{L}$ l'espace de Levi de $\tilde{G}$ correspondant et, d'apr\`es \ref{stabnonramif}, on a:
$$
\sum_{\hat{s}\in \hat{s}_{M'}Z(\hat{M})^{\Gamma,\hat{\theta}}/Z(\hat{L})^{ {\hat{\Gamma}},\hat{\theta}}}i_{\tilde{M}'}(\hat{L},\hat{L}'(\hat{s}))s_{{\bf M}'}^{{\bf L}'(\hat{s})}(c^{'V}_{\lambda})=r_{\tilde{M}}^{\tilde{L}}(c^{V}_{\lambda}),
$$
o\`u $c^{V}_{\lambda}$ est l'image de $c^{'V}_{\lambda}$ dans le $L$-groupe de $M$ via l'inclusion de ${\mathcal M}'$ dans ce groupe. On utlise encore la proposition de \ref{transfertcarpond} pour transformer (5) en $I_{\tilde{ L}}^{\tilde{G}}(ind_{\tilde{M}}^{\tilde{L}}\pi_{V}, 0,0,f_{V})$. Et on obtient exactement (1), ce qui prouve la proposition.

\
 Oublions l'espace de Levi $\tilde{M}$ et consid\'erons une donn\'ee endoscopique ${\bf M}'$ de $\tilde{G}$, $\omega$ qui ne correspond \`a aucun espace de Levi de $\tilde{G}$. Soit $\pi'_{V}$  et ${c'}^V$ comme ci-dessus. On peut associer \`a ${\bf M}'$ un Levi $\hat{M}$ de $\hat{G}$. On d\'efinit encore $I_{{\bf M}'}^{\tilde{G},{\mathcal E}}(\pi'_{V}\otimes c^{'V},f_{V})$ par la formule (3) qui conserve un sens. La m\^eme preuve que ci-dessus montre que:
\begin{prop} Sous ces hypoth\`eses,  $I_{{\bf M}'}^{\tilde{G},{\mathcal E}}(\pi'_{V}\otimes c^{'V},f_{V})=0$. 
\end{prop}

\section{Le c\^ot\'e spectral de la formule des traces}
\subsection{Rappel des termes discrets\label{rappeldestermesdiscrets}}
Il y a deux notions pour ''discret'' en ce qui concerne les repr\'esentations automorphes; la notion habituelle est que ce sont les repr\'esentations intervenant discr\`etement dans la d\'ecomposition spectrale des fonctions de carr\'e int\'egrable. L'autre notion qui nous concerne ici est le fait que ces repr\'esentations interviennent dans la partie discr\`ete de la formule des traces, qui va \^etre red\'ecrite ci-dessous. Il serait plus logique mais beaucoup trop lourd de distinguer en parlant de repr\'esentations discr\`etes pour les premi\`eres et de repr\'esentation $t$-discr\`ete pour les secondes. On va quand m\^eme le faire dans ce paragraphe, en utilisant la notation g\'en\'erique $\rho$ pour les repr\'esentations puis on posera $\pi_{disc}=\rho_{t-disc}$ et donc le $disc$ utilis\'e avec $\pi$ sera pour les repr\'esentations $t$-discr\`etes et c'est la notation qui sera ensuite utilis\'ee.

On reprend \cite{labwald} 6.1.

On ne peut pas parler de repr\'esentations discr\`etes sans caract\`ere unitaire du centre du groupe. Pour faire simple (et pas le plus g\'en\'eral possible) on reprend les d\'efinitions ${\mathfrak {A}}_{G}$ et ${\mathfrak A}_{\tilde{G}}$ introduites en \cite{stabilisationVI} fin de 1.3: pour ${\mathfrak A}_{G}$ c'est la composante neutre topologique des points sur ${\mathbb R}$ d'un tore d\'eploy\'e maximal dans $G$ vu comme groupe sur ${\mathbb Q}$ et ${\mathfrak A}_{\tilde{G}}$ en sont les \'el\'ements invariants sous $\theta$. Ces espaces s'identifie naturellement (via l'application log) \`a ${\mathcal A}_{G}$ et ${\mathcal A}_{\tilde{G}}$. On utilise alors  l'inclusion diagonale de ${\mathfrak A}_{{G}}$ dans ${G}(F_{\infty})$ pour identifier ${G}({\mathbb A}_{F})=\mathfrak{A}_{\tilde{G}}G({\mathbb A}_{F})^1$.   Les repr\'esentations que nous allons consid\'erer sont invariantes ont un caract\`ere sous ${\mathfrak A}_{\tilde{G}}$ unitaire et un caract\`ere pour l'action de  ${\mathfrak A}_{G}$ qui est $\theta$- semi-invariant au sens $\theta(\chi)=\chi \omega$ si $\chi$ est ce caract\`ere; on a suppos\'e que l'action de $\tilde{G}(F)$ sur le centre de $G(F)$ est  d'ordre fini et que $\omega$ est unitaire. Donc un tel caract\`ere $\chi$ est n\'ecessairement unitaire et donc toute $\omega$-repr\'esentation irr\'eductible
de $\tilde{G}({\mathbb A}_{F})$ ayant un caract\`ere  unitaire pour l'action de $\mathfrak{A}_{\tilde{G}}$ a un caract\`ere  unitaire pour l'action de ${\mathfrak A}_{{G}}$.

Le groupe des points rationnels $\tilde{G}(F)$ agit naturellement par conjugaison dans l'ensemble des formes automorphes de carr\'e int\'egrale et on pose pour toute forme automorphe de carr\'e int\'egrable, $\phi$, invariante sous $\mathfrak{A}_{\tilde{G}}$ 
 pour tout   $g\in G({\mathbb A}_{F})$, $\gamma\in \tilde{G}({\mathbb A}_{F})$
$$
\rho_{disc}^{\tilde{G}}(\gamma,\omega)\phi(g)=\omega(\delta^{-1}g\gamma)\phi(\delta^{-1}g\gamma),
$$
o\`u $\delta$ est n'importe quel \'el\'ement de $\tilde{G}(F)$ et o\`u $\delta^{-1}g \gamma \in G({\mathbb A}_{F})$ est d\'efini par l'\'egalit\'e $g \gamma= \delta (\delta^{-1}g \gamma)$ dans $\tilde{G}({\mathbb A}_{F})$.
On remarque que cette formule est ind\'ependante de $\delta$.

Soit $x,y\in G({\mathbb A}_{F})$ est appliquons cette formule avec $\gamma$ remplac\'e par $x\gamma y$: 
$$
\rho_{disc}^{\tilde{G}}(x\gamma y,\omega)\phi(g)= \omega (\delta^{-1} gx\gamma y)\phi (\delta^{-1} gx \gamma y)=\omega(y) \omega (\delta^{-1} gx\gamma ) \rho(y)(\phi)(\delta^{-1}gx\gamma).
$$Or
$$
\omega (\delta^{-1} gx\gamma ) \rho(y)(\phi)(\delta^{-1}gx\gamma)=\rho_{disc}^{\tilde{G}}(\gamma,\omega) (\rho(y)\phi)(gx).
$$
D'o\`u
$$\rho_{disc}^{\tilde{G}}(x\gamma y,\omega)\phi(g)= \omega(y)\rho(x) \rho_{disc}^{\tilde{G}}(\gamma,\omega)(\rho(y)\phi)(g).
$$
ce qui est bien une $\omega$ repr\'esentation.

Les autres termes \`a ajouter pour avoir la partie discr\`ete de la formule spectrale viennent des sous groupes de Levi de $G$, $M$ tel que $Norm_{\tilde{G}(F)}M$ non seulement est non vide mais en plus contient un \'el\'ement r\'egulier, c'est-\`a-dire un \'el\'ement $u$ qui agit sans point fixe 
 sur ${\mathcal A}^*_{M}/{\mathcal A}^*_{\tilde{G}}$. Soit $u$ un tel \'el\'ement.
 La partie discr\`ete associ\'ee \`a $M$ et $u$  est \'ecrite en \cite{labwald} dans la section 14.3 (avant l'\'enonc\'e de 14.3.2) dans le cas non invariant. Rappelons la construction. On fixe un sous-groupe parabolique $P$ de $G$ de sous-groupe de Levi $M$; on note $u.P$ le sous-groupe parabolique de $G$ obtenu en conjuguant $P$ par $u$. 
On note, pour tout sous-groupe parabolique $Q$ de $G$ de sous-groupe de Levi $M$,  $U_{Q}$ son radical unipotent.

Soit $\phi$ une fonction sur $G({\mathbb A}_{F})$ invariante \`a gauche sous $M(F)U_{P}({\mathbb A}_{F})$. Pour tout $\gamma\in \tilde{G}({\mathbb A}_{F})$, on d\'efinit la fonction $\gamma_{\omega}.\phi$ en posant pour tout $g\in G({\mathbb A}_{F})$:
$$
(\gamma_{\omega}.\phi)(g)=\omega(u^{-1}g\gamma) \phi(u^{-1}g \gamma).
$$
C'est une fonction sur $G({\mathbb A}_{F})$ invariante \`a gauche sous $M(F)U_{u.P}({\mathbb A}_{F})$.

Pour revenir en une fonction invariante \`a gauche sous $M(F)U_{P}({\mathbb A}_{F})$ on utilise l'op\'erateur d'entrelacement standard. Pour que cette op\'erateur soit holomorphe il suffit que $\phi$ se transforme \`a gauche sous ${\mathfrak{A}}_{M}$ par un caract\`ere unitaire.

Pour cela, on note ${\mathcal{A}}^2(U_{P}({\mathbb A}_{F})M(F){\mathfrak{A}}_{\tilde{G}}(1-u)_{\omega}{\mathfrak{A}}_{M}\backslash G({\mathbb A}_{F}))$ les formes automorphes de carr\'e int\'egrable au sens de \cite{mw} I.2.17,  invariantes pour l'action \`a gauche de $U_{P}({\mathbb A}_{F})M(F){\mathfrak{A}}_{\tilde{G}}$ et se transformant sous l'action \`a gauche par un \'el\'ement de  ${\mathfrak {A}}_{M}$ de la forme $a^{-1}(u.a)$ (avec $a\in {\mathfrak{A}}_{M})$ par la multiplication par $\omega(a)$ et qui sont ${\mathfrak{A}}_{M}$ finies pour l'action de ce groupe. 
Comme $u$ n'a pas de point fixe dans ${\mathfrak{A}}^*_{M}/{{\mathfrak{A}}}^*_{G}$ il suffit de v\'erifier que  si l'on d\'ecompose les fonctions consid\'er\'ees suivant l'action de ${\mathfrak{A}}_{M}$ seuls des caract\`eres unitaires interviennent. Par hypoth\`ese ${\mathfrak{A}}_{\tilde{G}}$ agit trivialement. Il faut donc regarder l'action de ${\mathfrak{A}}_{M}/{\mathfrak{A}}_{\tilde{G}}$.  Dans ${\mathcal{A}}_{{G}}/{\mathcal{A}}_{\tilde{G}}$, l'action de $u$ est \'evidemment sans point fixe donc avec l'hypoth\`ese faite sur $u$, l'action de $u$ dans ${\mathcal{A}}_{M}/{\mathcal{A}}_{\tilde{G}}$ est sans point fixe et on a demand\'e que les \'el\'ements $a^{-1} (u.a)$ agissent via le caract\`ere unitaire $\omega$. Donc la partie r\'eelle des caract\`eres de $\mathfrak{A}_{M}$ intervenant est invariante sous l'action adjointe de $u$, cette partie r\'eelle est donc nulle.

On sait alors d\'efinir l'op\'erateur d'entrelacement par prolongement m\'eromorphe
$$
M_{P\vert u.P}(u): {\mathcal{A}}^2(U_{u.P}({\mathbb A}_{F})M(F)\mathfrak{A}_{\tilde{G}}(1-u)_{\omega}\mathfrak{A}_{M}\backslash G({\mathbb A}_{F})) \rightarrow $$
$$
{\mathcal{A}}^2(U_{P}({\mathbb A}_{F})M(F)\mathfrak{A}_{\tilde{G}}(1-u)_{\omega}\mathfrak{A}_{M}\backslash G({\mathbb A}_{F})) .
$$ 
Et l'action de $\gamma$ sur $\phi$ est alors $$\rho_{t-disc}(\gamma)(\phi):=
M_{P\vert u.P}(u) \gamma_{\omega} \phi.
$$On v\'erifie comme ci-dessus que cela donne une $\omega$-repr\'esentation de $\tilde{G}({\mathbb A}_{F})$ dans ${\mathcal{A}}^2(U_{P}({\mathbb A}_{F})M(F)\mathfrak{A}_{\tilde{G}}(1-u)_{\omega}\mathfrak{A}_{M}\backslash G({\mathbb A}_{F}))$. On la note $\rho_{disc,M,u,\omega}$.

On pose $w^M:=\vert Norm_{G(F)}(M)/M(F)\vert$ et:
 $$\rho_{t-disc,\omega}:=\sum_{M,u}w(M) \frac{1}{\vert {\mathrm det}({u}-1)_{\vert {\mathcal{A}}_{M}/{\mathcal{A}}_{G }\vert }} \rho_{disc,M,u,\omega},$$  o\`u $M$ parcourt l'ensemble des sous-groupes de Levi de $G$, pris \`a conjugaison pr\`es et o\`u $u$ parcourt l'ensemble des \'el\'ements de $G(F)$ normalisant $M$ r\'eguliers au sens ci-dessus et pris \`a translation pr\`es par l'action de $M(F)$ op\'erant \`a gauche ou \`a droite (cela n'a pas d'importance).

 \bf D\'efinition et Notations: \rm on fixe un ensemble fini de places, $V$, de $F$ contenant $V_{ram}$. Pour tout $v$ non dans $V$, on note $c_{v}$ un morphisme non ramifi\'e  de $W_{F_{v}}$ dans $^LG(F_{v})$, o\`u encore $^LG(F_{v})$ d\'esigne le $L$-groupe de $G$ sur $F_{v}$. On note $c^V$ le produit tensoriel de ces morphismes; ainsi $c^V$ donne une repr\'esentation non ramifi\'ee de $G({\mathbb A}_{F}^{V})$ ou suivant le point de vue un caract\`ere du produit des alg\`ebres de Hecke sph\'eriques en toute place hors de $V$. On suppose que $c^V$ est $\omega, \hat{\theta}$ invariant. En faisant op\'erer trivialement $\tilde{K}^V$ sur les vecteurs $K^V$-invariants de la repr\'esentation, hors de $V$ la repr\'esentation s'\'etend  en une repr\'esentation de $\tilde{G}({\mathbb A}_{F}^V)$; c'est la normalisation, canonique, que nous utiliserons syst\'ematiquement hors de $V$.

Aux places archim\'ediennes, on fixe $\nu$ un caract\`ere infinit\'esimal tel que $\theta(\nu)=\nu+d\omega$; on r\'ealise ${\mathfrak{A}}_{\tilde{G}}$ comme un sous-groupe du centre de $G(F_{\infty})$ form\'e d'\'el\'ements invariants sous $\tilde{G}(F_{\infty})$. On suppose que $\nu$ est trivial sur l'alg\`ebre de Lie de ce groupe (vu comme sous espace du centre de l'alg\`ebre enveloppante).

On note alors $\pi_{\nu,disc}(c^V)$ la $\omega$ repr\'esentation de $\tilde{G}({\mathbb A}_{F})$ qui est le produit du scalaire $\vert det (1-\theta)_{\vert {\mathcal A_{G}}/{\mathcal A}_{\tilde{G}}}\vert^{-1}$ et de  la somme des sous-repr\'esentations de $\rho_{t-disc,\omega}$ ayant des vecteurs invariants sous $\tilde{K}^V$ et sur lesquelles l'alg\`ebre enveloppante de $G(F_{\infty})$ op\`ere par le caract\`ere central $\nu$ et o\`u toute fonction sur $K^V\backslash G({\mathbb A}_{F}^V)/K^V$ op\`ere, sur les vecteurs $\tilde{K}^V$ invariants par le caract\`ere $c^V$. On voit cette repr\'esentation comme une $\omega_{V}$-repr\'esentation de $\tilde{G}(F_{V})$ en se limitant aux \'el\'ements $\tilde{K}^V$ invariants.  Le scalaire  $\vert det (1-\theta)_{\vert {\mathcal A_{G}}/{\mathcal A}_{\tilde{G}}}\vert^{-1}$ a \'et\'e mis un peu formellement pour que les formules soient coh\'erentes entre l'espace ${\tilde{G}}$ et ses espaces de Levi, (cf. \cite{stabilisationVI})

D'apr\`es les constructions, cet espace de repr\'esentations se r\'ealise dans un espace de fonctions invariantes \`a gauche sous ${\mathfrak{A}}_{\tilde{G}}$.
On fixe $\lambda\in i{\mathcal{A}}^*_{\tilde{G}}$ cela donne un caract\`ere unitaire de ${\mathfrak{A}}_{\tilde{G}}$ qui s'\'etend en un caract\`ere unitaire de $\tilde{G}({\mathbb A}_{F})$ trivial \'evidemment sur ${G}({\mathbb A}_{F})^1$. On note alors $\pi_{\nu,disc,\lambda}(c^V)$ le produit tensoriel de ce caract\`ere $\lambda$ avec $\pi_{\nu,disc}(c^V)$ et on voit cette repr\'esentation comme une $\omega_{V}$-repr\'esentation de $\tilde{G}(F_{V})$. Il faut faire attention aux places hors de $V$, ce n'est plus $c^V$ le caract\`ere de l'alg\`ebre de Hecke sph\'erique mais $c^V\otimes \lambda$.
C'est la distribution suivante sur laquelle on a prise via la formule des traces:
 $$
 f_{V}\in I(\tilde{G}(F_{V}))\mapsto tr\, \pi_{\nu,disc,\lambda}(c^V)(f_{V}). \eqno(1)
 $$
On a \'evidemment $\pi_{\nu,disc,\lambda+\mu}(c^V)=\pi_{\nu,disc,\lambda}(c^V)\otimes \mu$ pour tout $\lambda,\mu \in i{\mathcal{A}}^*_{\tilde{G}}$.

Mais dans la formule des traces, c'est le coefficient de Fourier de cette distribution qui intervient: pour $X\in {\mathcal{A}}_{\tilde{G}}$, on d\'efinit la distribution, $ \forall f_{V}\in I(\tilde{G}(F_{V}))$
 $$
 I_{\nu}(c^V,X,f_{V}):=\int_{\lambda\in i{\mathcal{A}}^*_{\tilde{G}}}d\lambda\,  tr\, \pi_{\nu,disc,\lambda}(c^V)(f_{V}) e^{-\lambda(X)}.
 $$
 Et la partie discr\`ete de la formule des traces invariante dans $V$ est la distribution $ \oplus_{\nu,c^V}I_{\nu}(c^V,0,f_{V})$, ce qui veut dire que l'on ne consid\`ere que les traces des repr\'esentations restreintes \`a $\tilde{G}(F_{V})\tilde{K}^V\cap {G}(\mathbb A)^1\tilde{G}(F)$.
 
\subsection{Rappel des termes continus}
On fixe un espace de Levi $\tilde{M}$ de $\tilde{G}$; on g\'en\'eralise la d\'efinition $\pi_{\nu,disc,\lambda}(c^V)$ en supposant ici que $\nu$ est un caract\`ere infinit\'esimal de $M(F_{\infty})$ et $c^V$ est un syst\`eme de caract\`ere des alg\`ebres de Hecke sph\'eriques hors de $V$ pour $M(F^V)$; on note alors plut\^ot $\pi^{\tilde{M}}_{\nu,disc,\lambda}(c^V)$. On a donc d\'efini en \ref{defcarglob} la distribution sur $I(\tilde{G}(F_{V}))$, $$f_{V}\mapsto I_{\tilde{M}}^{\tilde{G}}(\pi^{\tilde{M}}_{\nu,disc}(c^V),f_{V}).
$$
La partie continue de la formule des traces invariantes est d'apr\`es \cite{invglob}  Theorem 4.4 la distribution sur $I(\tilde{G}(F_{V}))$ obtenue en sommant (on reviendra ci-dessous sur les probl\`emes de convergence)
$$
\sum_{\tilde{M}} w(\tilde{M}) \sum_{\nu,c^V} I_{\tilde{M}}^{\tilde{G}}(\pi^{\tilde{M}}_{\nu, disc}(c^V),f_{V});
$$ la pr\'esentation ici est un peu diff\'erente qu'en loc. cite car, d'une part, on a choisi de faire entrer ce qui est not\'e $a^M(\pi)$ dans \cite{invglob} (4.5)  pour la  partie ''multiplicit\'e'' dans $\tilde{M}$ (not\'e $a_{disc}^{M_{1}}$, le $M_{1}$ est notre $\tilde{M}$)  dans $\pi^{\tilde{M}}_{\nu, disc}(c^V)$ et  d'autre part le $r_{M_{1}}^{M}(\pi_{1,\lambda})$ a \'et\'e pris en compte dans la d\'efinition de $I_{\tilde{M}}^{\tilde{G}}(\pi^{\tilde{M}}_{\nu, disc}(c^V),f_{V})$ dans \ref{defcarglob}. Ensuite on remplace dans l'\'enonc\'e de \cite{invglob} theorem 4.4, l'int\'egrale sur $\Pi(M,t)$ par sa valeur donn\'ee avant (4.5) de loc. cite; cela devient une int\'egrale sur $\Pi_{disc}(M_{1},t)$, ce qui est essentiellement notre formulation.
La diff\'erence avec \cite{invglob} est le fait que l'on n'a pas mis de somme sur la partie imaginaire du caract\`ere infinit\'esimal (le $t$); ici on utilise les r\'esultats de Finis, Lapid et M\"uller qui assurent la convergence comme pour la formule des traces non invariantes dans le paragraphe 14.3 de \cite{labwald}.

Pour tout $\tilde{M}$ sous espace de Levi propre de $\tilde{G}$, on a une int\'egrale sur $i{\mathcal{A}}^*_{\tilde{M}}$ ce qui explique que l'on a une partie continue. 

\subsection{Repr\'esentations semi-finies}
\subsubsection{D\'efinition\label{definitionsemifinie}}
On a besoin de cette d\'efinition: soit $\pi$ une $\omega$ repr\'esentation de $\tilde{G}({\mathbb A}_{F})$ ou plut\^ot un \'el\'ement du groupe de Grothendieck  de ces $\omega$-repr\'esentations; on dit que $\pi$ est semi-finie si les conditions suivantes sont v\'erifi\'ees:

$\overline{\pi} \simeq \check{\pi}$

pour tout ensemble fini de $K$-types, ${\mathcal K}$, il existe $R\in {\mathbb R}_{>0}$ tel que les composantes de $\pi$ admettant des vecteurs se transformant sous au moins un de ces $K$-types ont toutes un caract\`ere infinit\'esimal dont la partie r\'eelle est de norme inf\'erieure ou \'egale \`a $R$;
si en plus,  en fixant $R'\in {\mathbb R}_{\geq 0}$, parmi ces repr\'esentations il n'y en a plus qu'un nombre fini (en comptant les multiplicit\'es) dont la partie imaginaire du  caract\`ere infinit\'esimal est de norme born\'ee par $R'$;

pour $\nu$ un caract\`ere infinit\'esimal fix\'e, on note $\pi_{\nu}^{\mathcal K}$ la somme des composantes irr\'eductibles de $\pi$ ayant des vecteurs se transformant sous au moins un $K$-type dans ${\mathcal K}$ et ayant pour caract\`ere infinit\'esimal $\nu$. Cette repr\'esentation est de longueur finie d'apr\`es les propri\'et\'es demand\'ees et on   demande en plus  que la somme sur $\nu$ des traces ces repr\'esentations \'evalu\'ees sur des fonctions $f\in C^{\infty}_{c}(\tilde{G}({\mathbb A}_{F}))$, se transformant selon des $K$-types inclus dans l'ensemble fix\'e ${\mathcal K}$,  converge absolument.

\

Fixons $\pi$ une repr\'esentation semi-finie et $V$ un ensemble fini de places. On suppose que toutes les composantes de $\pi$ sont non ramifi\'ees hors de $V$. On d\'ecompose alors   $\pi=\sum_{\nu}\pi_{\nu}$  o\`u $\pi_{\nu}$ est la somme des termes de caract\`ere infinit\'esimal $\nu$ et pour $\nu$ fix\'e on d\'ecompose $\pi_{\nu}=\sum_{c^V}\pi_{\nu}(c^V)$, o\`u   $c^V$ parcourt l'ensemble  des caract\`eres des alg\`ebres de Hecke sph\'eriques hors de $V$.

Soit ${\bf M}'$ une donn\'ee endoscopique de $\tilde{G}$ que l'on ne suppose pas elliptique mais que l'on suppose relevante; plus exactement on note $M'_{1}$ une donn\'ee auxiliaire relative \`a cette donn\'ee endoscopique de sorte que    le groupe ${\mathcal M}'$ (de la donn\'ee endoscopique) soit plong\'e dans le groupe dual, $^LM'_{1}$ de $M'_{1}$; ainsi $M'_{1}$ est une extension par un tore induit du groupe $M'$ de la donn\'ee endoscopique ${\bf M}'$ et il existe un caract\`ere de ce tore (donn\'e par ${\bf M}'$ et  le plongement de ${\mathcal M}'$ dans $^LM'_{1}$) tel que tous les objets attach\'es \`a ${M}'_{1}$ intervenant dans l'endoscopie d\'efinie par ${\bf M}'$ se transforment par ce caract\`ere sous l'action du tore induit. Surtout pour ce que l'on fait ici, on peut simplifier les notations en ''oubliant'' que $\tilde{M}'\neq \tilde{M}'_{1}$:  au lieu d'objets attach\'es \`a $\tilde{M}'_{1}$ se transformant selon ce caract\`ere, on parlera symboliquement d'objets attach\'es \`a ${\bf M}'$.  A tout caract\`ere infinit\'esimal $\nu'$  de ${\bf M}'$ est associ\'e un caract\`ere infinit\'esimal $\nu$ de $G$ tel que le transfert d'un paquet stable de repr\'esentations de ${\bf M}'$ de caract\`ere infinit\'esimal $\nu'$ \`a un espace de Levi de $\tilde{G}$ ait $\nu$ pour caract\`ere infinit\'esimal. On \'ecrira $\nu' \mapsto \nu$ pour signifier cette relation qui n'est ni injective ni surjective. 
Soit aussi $c^{'V}$ un caract\`ere de l'alg\`ebre de Hecke sph\'erique de ${\bf M}'$ hors de $V$. Par la correspondance de Langlands, cela d\'efinit un caract\`ere pour l'alg\`ebre de Hecke sph\'erique hors de $V$ de $\tilde{G}$. On note encore $c^{'V}\mapsto c^V$ cette relation qui n'est l\`a aussi ni injective ni surjective.

\subsubsection{Les repr\'esentations semi-finies et la partie discr\`ete de la formule des traces}
On fixe $V$ un ensemble fini de places contenant les places archim\'ediennes. On reprend la notation $I_{disc}^{\tilde{G}}(\omega,f_{V}1_{\tilde{K}^V})$ pour la partie discr\`ete de la formule des traces. 

\begin{rmq} Cette distribution s'exprime comme la trace d'une repr\'esentation semi-finie\end{rmq}

Ceci est loin d'\^etre \'evident mais les arguments sont d\'ej\`a dans la litt\'erature comme on va l'expliquer. Le premier point qu'il nous faut est la propri\'et\'e suivante. On fixe un ensemble fini de $K_{\infty}$-types. Soit $\pi$ une des repr\'esentations irr\'eductibles intervenant dans $I_{disc}$; on suppose que $\pi_{\infty}$ a au moins un $K_{\infty}$-type dans l'ensemble fix\'e. Alors la partie r\'eelle du caract\`ere infinit\'esimal de $\pi$ est born\'ee ind\'ependamment de $\pi$. C'est utilis\'e en \cite{partie3} preuve de 4.1 qu'on peut remonter \`a \cite{invglob} fin de la preuve de 6.5. Faute de r\'ef\'erence dans ces citations (le r\'esultat doit \^etre bien connu des experts d\`es la fin des ann\'ees 80, cf. \cite{hermitien} (4) du corollaire 7.2), on explique le r\'esultat ainsi: une repr\'esentation unitaire, $\pi$, est unitairement induite \`a partir d'une repr\'esentation unitaire, $\pi'$ d'un de ses sous-groupes de Levi ayant un caract\`ere infinit\'esimal r\'eel (voir \cite{voganlanglands} paragraphe 3). Pour les repr\'esentations unitaires, M\"uller a montr\'e en \cite{muller} (8.1) qui renvoie \`a (3.4) (et qui n'est pas difficile) que pour $\pi'$ une repr\'esentation unitaire d'un groupe de Lie r\'eel et pour $\sigma$ l'un de ses $K_{\infty}$-types l'action de l'op\'erateur de Casimir agissant sur $\pi'$ est born\'ee par l'action du Laplacien de $K_{\infty}$ agissant sur $\sigma$. Avec l'hypoth\`ese que le caract\`ere infinit\'esimal de $\pi'$ est r\'eel cela borne ce caract\`ere infinit\'esimal. Comme le fait que l'on impose \`a $\pi$ de contenir au moins certains $K_{\infty}$-types se propagent \`a $\pi'$ (pour des types qui d\'ependent du sous-groupe de Levi mais pas des repr\'esentations), on obtient l'assertion.

Ensuite la propri\'et\'e de convergence absolue r\'esulte des travaux de M\"uller sur la tra\c{c}abilit\'e du spectre discret (\cite{muller1}) et les propri\'et\'es restantes sont des propri\'et\'es g\'en\'erales de finitude du nombre de repr\'esentations automorphes quand on a fix\'e le caract\`ere infinit\'esimal et la ramification. Cela prouve la remarque.

\subsubsection{Utilisation des multiplicateurs sur les repr\'esentations semi-finies\label{multiplicateurs}}
Cette th\'eorie des multiplicateurs est due \`a Arthur et permet de s\'eparer les caract\`eres infinit\'esimaux (\cite{multiplicateur} 4.2). On l'utilise comme dans \cite{partie3}, en montrant ici la propri\'et\'e qui servira.

Soit $V$ un ensemble fini de places contenant  $V_{ram}$ et soit ${\mathcal F}$ un ensemble de fonctions dans $I(\tilde{G}(F_{V}))$ (donc en particulier $K_{V}$ finies). On suppose que ${\mathcal F}$ est stable par convolution pour tout multiplicateur $\alpha$. On fixe aussi $\pi$ une repr\'esentation semi-finie de $\tilde{G}({\mathbb A}_{F})$ dont toutes les composantes sont non ramifi\'ees hors de $V$. On reprend la notation $\pi=\sum_{\nu,c^V}\pi_{\nu}(c^V)$ de \ref{definitionsemifinie}. 

\begin{lem} (i) On suppose que $tr\, \pi(f_{V} 1_{\tilde{K}}^V)=0$ pour tout $f_{V}\in {\mathcal F}$.
Alors $tr\, \pi_{\nu}(f_{V} 1_{\tilde{K}^V})=0$ pour tout $f_{V}\in {\mathcal F}$ et tout caract\`ere infinit\'esimal $\nu$.

(ii) On suppose en plus que pour tout ensemble fini de place $V'$ de $F$ contenant $V$ et pour toute fonction $f_{V'-V}$ dans $I(\tilde{G}(F_{V'-V}))$ non ramifi\'ee, on a $tr\, \pi(f_{V}f_{V'-V}1_{\tilde{K}^{V'}})=0$ alors $tr\, \pi_{\nu}(c^V)(f_{V}1_{\tilde{K}^V})=0$ pour tout couple $\nu,c^V$ comme pr\'ec\'edemment.
\end{lem}

On fixe $f\in {\mathcal F}$; il suffit de supposer que ${\mathcal F}$ est exactement l'ensemble des \'el\'ements $f_{\alpha}$ pour $\alpha$ parcourant l'ensemble des mutliplicateurs; on rappelle que $\alpha$ est une fonction \`a support compact sur ${\mathfrak h}$  o\`u ${\mathfrak h}$ est l'alg\`ebre de Lie (complexifi\'ee) d'un tore de $G(F_{\infty})$, invariante sous $W_{\infty}$, le groupe de Weyl de ce tore complexifi\'e. On note $\check{\alpha}$ la transform\'ee de Fourier de $\alpha$ et $\check{\alpha}$ parcourt donc l'ensemble des fonctions de Paley Wiener sur ${\mathfrak h}^{*}$ invariante sous $W_{\infty}$. Comme $f_V$ est $K_{V}$-fini et que $\pi$ est non ramifi\'e hors de $V$, l'hypoth\`ese de semi-finitude assure qu'il existe $R\in {\mathbb R}_{>0}$ tel que si $\pi_{\nu}$ a des vecteurs invariants sous l'un des $K$-types fix\'es, alors $\vert \vert Re\, \nu\vert \vert \leq R$. On fixe $\nu_{0}$ avec cette propri\'et\'e sur les invariants. On fixe aussi, gr\^ace \`a \cite{AC}, 2e partie, lemme 15.2,  $\alpha$ tel que $\check{\alpha}(\nu)\in [0, 1]$ pour tout  caract\`ere infinit\'esimal $\nu$  de partie r\'eelle   born\'ee en norme par $R$ et tel que $\nu$ soit le caract\`ere infinit\'esimal d'une repr\'esentation $\pi$ v\'erifiant $\overline{\pi}\simeq \check{\pi}$;  et on impose en plus (avec la m\^eme r\'ef\'erence)  que  $\check{\alpha}(\nu)= 1$ pour un tel $\nu$ uniquement si $\nu=\nu_{0}$.

 Pour tout $m\in {\mathbb N}$ on note $*^m\alpha$ le convol\'e de $\alpha$ $m$-fois. Alors $tr\, \pi(f_{*^m\alpha})=\sum_{\nu}\check{\alpha}(\nu)^m tr\, \pi_{\nu}(f_{V}1_{\tilde{K}^V})$. Pour $\nu$ intervenant dans cette somme, avec $\nu\not=\nu_{0}$, $\check{\alpha}(\nu)^m$ tend vers $0$ quand $m$ tend vers l'infini. On utilise alors l'hypoth\`ese de convergence absolue pour montrer que quand $m$ tend vers l'infini cette somme converge vers $tr\, \pi_{\nu_{0}}(f_{V}1_{\tilde{K}}^V)$. On a ainsi montr\'e (i).

Montrons (ii). On fixe $\nu_{0}$ et $f_{V}\in {\mathcal F}$. On fixe aussi un ensemble de $K$-types qui sont les $K^V$ types triviaux hors de $V$ et  tel que $f_{V}$ se transforme suivant ces $K_{V}$-types. On sait alors que $\pi_{\nu_{0}}$ a un nombre fini de composantes irr\'eductibles admettant ces $K$-types et on les d\'ecompose en $\sum_{c^V}\pi_{\nu_{0}}(c^V)$. On fixe $V'$ un ensemble fini de places contenant $V$ et tel que les caract\`eres des repr\'esentations constituant $\pi_{\nu_{0}}(c^V)$ restreints en des distributions sur $\otimes_{v\in V'}I(\tilde{G}(F_{v}))$ soient lin\'eairement ind\'ependants, ou encore que $\otimes_{v\in V'}I(\tilde{G}(F_{v}))$  s\'eparent les diff\'erentes repr\'esentations $\pi_{\nu_{0}}(c^V)$. Par hypoth\`ese $tr\, \pi(f_{V}f_{V'-V}1_{\tilde{K}^{V'}})=0$ donc d'apr\`es (i) appliqu\'e \`a $V'$ ce qui est loisible, $tr\, \pi_{\nu_{0}}(f_{V}f_{V'-V}1_{\tilde{K}^{V'}})=0$ et en faisant varier $f_{V'-V}$ on obtient $tr\, \pi_{\nu_{0}}(c^V)(f_{V}f_{V'-V}1_{\tilde{K}^{V'}})=0$ pour tout choix de $f_{V'-V}$ et pour tout $c^V$. Cela entra\^{\i}ne a fortiori (ii).
\subsection{Repr\'esentation semi-finie et stabilit\'e\label{decompositionstable}}
Soit $\pi$ une repr\'esentation semi-finie de $\tilde{G}$, on la suppose stable; cela veut dire que pour toute place $v_{0}$, pour tout ensemble fini de places $V'$ contenant $V_{ram}$ et $v_{0}$ et pour toute fonction $f_{V'}\in I(\tilde{G}(F_{V'}))$ d\'ecompos\'ee, dont la composante en la place $v_{0}$ annule toutes les int\'egrales orbitales stables de $\tilde{G}(F_{v_{0}})$, $tr\, \pi(f_{V'}1_{\tilde{K}^{V'}})=0$. On fixe $V$ contenant $V_{ram}$ comme pr\'ec\'edemment et on  d\'ecompose $\pi$ en $\sum_{\nu, c^V} \pi_{\nu}(c^{V})$
\begin{cor} Pour tout $\nu,c^{V}$, les distributions $$f_{V}\in I(\tilde{G}(f_{V}))\mapsto tr\, \pi_{\nu}(c^V)(f_{V}1_{\tilde{K}^V})$$ sont stables.
\end{cor}
On fixe $\nu,c^V$ et $v_{0}\in V$. et on doit montrer que pour toute fonction $f_{v_{0}}\in I(\tilde{G}(F_{v_{0}}))$ qui annule les int\'egrales orbitales stables en la place $v_{0}$ et pour toute fonction $f_{V}\in I(\tilde{G}(F_{v}))$ d\'ecompos\'ee et ayant $f_{v_{0}}$ comme composante en la place $v_{0}$, on a $tr\, \pi_{\nu}(c^V)(f_{V}1_{\tilde{K}^V})=0$.
On applique \ref{multiplicateurs} \`a l'ensemble ${\mathcal F}$ qui est pr\'ecis\'ement l'ensemble des $f_{V,\alpha}$ o\`u $\alpha$ parcourt l'ensemble des multiplicateurs. Mais pour cela il faut v\'erifier que pour tout $f'_{V}$ dans cet ensemble $f'_{v_{0}}$ annule encore toutes les int\'egrales orbitales stables. C'est \'evident d'apr\`es les d\'efinitions si $v_{0}$ n'est pas une place archim\'edienne. Supposons donc que $v_{0}$ est une place archim\'edienne. Pour toute repr\'esentation stable, $\pi'_{v_{0}}$ de $\tilde{G}(F_{v_{0}})$, on a $$tr\, \pi'_{v_{0}}(f_{v_{0},\alpha})=\check{\alpha}(\nu_{\pi'_{v_{0}}})tr\, \pi'_{v_{0}}(f_{v_{0}})=0,
$$
par hypoth\`ese sur $f_{v_{0}}$. D'apr\`es \cite{stabilisationIV} 2.3, la fonction $f_{v_{0},\alpha}$ annule donc toutes les int\'egrales orbitales stables en la place $v_{0}$. Cela permet d'appliquer \ref{multiplicateurs} (i) et (ii) pour obtenir le corollaire.

\subsection{Enonc\'e du lemme fondamental tordu\label{lemmefondamental}}
La situation ici est locale et on note $v$ la place consid\'er\'ee. On suppose que $v\notin V_{ram}$.
On fixe ${\bf G}'$ une donn\'ee endoscopique elliptique non ramifi\'ee de $\tilde{G}(F_{v})$. A tout morphisme non ramifi\'e de $W_{F_{v}}$ dans ${\mathcal G}'$,   on associe un morphisme non ramifi\'e de $W_{F_{v}}$ dans le $L$-groupe de $G$.  Dans notre cas non ramifi\'e, on peut fixer un isomorphisme ${\cal G}'\simeq {^LG}'$.  Cela donne une correspondance entre les repr\'esentations non ramifi\'ees de $G'(F_{v})$ dans l'ensemble des repr\'esentations non ramifi\'ees   $\pi$ de $G(F_{v}) $ telles que $\pi^{\theta}\simeq \pi\otimes \omega$. Comme on a suppos\'e que $v\notin V_{ram}$, une telle repr\'esentation non ramifi\'ee de $G(F_{v})$ s'\'etend canoniquement en une repr\'esentation de $\tilde{G}(F_{v})$ en demandant que $\tilde{K}_{v}$ agisse trivialement sur l'espace des invariants sous $K_{v}$ (cet espace est de dimension 1, bien entendu). Si la donn\'ee endoscopique a elle aussi une torsion (n\'ecessairement int\'erieure) on \'etend aussi les repr\'esentations non ramifi\'ees de $G'(F_{v})$ \`a l'espace tordu. Le lemme fondamental tordu pour toute l'alg\`ebre de Hecke sph\'erique dit que si $f_{v}\in I(\tilde{G}(F_{v}))$ est bi-invariante sous $K_{v}$ et alors on peut choisir une fonction $f_{v}^{{\bf G}'}$ dans $I(\tilde{G}'(F_{v}))$ bi-invariante sous $K_{G',v}$ qui soit un transfert de $f_{v}$ et qui v\'erifie $tr \pi(f_{v})= tr \pi'(f^{{\bf G}'}_{v})$ si $\pi$ et $\pi'$ se correspondent dans la correspondance que l'on vient de d\'ecrire. C'est la m\^eme formulation que \cite{stabilisationI} 6.4.

Dans le cas non tordu et pour l'\'el\'ement unit\'e de l'alg\`ebre de Hecke ceci est d\'emontr\'e par Ngo et la d\'emonstration de Ngo a fait sauter le verrou bloquant la stabilisation de la formule des traces. Hales a montr\'e en \cite{hales} que l'assertion pour l'\'el\'ement unit\'e de l'alg\`ebre de Hecke entra\^{\i}ne le lemme fondamental pour tous les \'el\'ements de l'alg\`ebre de Hecke sph\'erique, dans le cas non tordu. Dans le cas tordu,  il est  d\'emontr\'e que le r\'esultat de Ngo entra\^{\i}ne le lemme fondamental pour l'\'el\'ement unit\'e de l'alg\`ebre de Hecke dans le cas tordu (\cite{waldslemmefond}). Il reste \`a \'etendre le r\'esultat de Hales, c'est un travail en cours  mais d'ores et d\'ej\`a Arthur a remarqu\'e que le cas de GL(n) tordu r\'esulte de \cite{selecta} et de son extension au cas tordu \cite{prolongementdescaracteres}.

On admet ici ce lemme fondamental.

\subsection{La variante stable de la partie discr\`ete de la formule des traces}
On suppose ici que $\tilde{G}$ est \`a torsion int\'erieure et que $G$ est quasi-d\'eploy\'e avec $\omega=1$. On fixe un ensemble fini de places $V$ de $F$ contenant $V_{ram}$.
On consid\`ere  la distribution  qui \`a $f_{V}\in I(\tilde{G}(F_{V}))$ associe
$$SI_{disc}^{\tilde{G}}(f_{V}1_{\tilde{K}^V}):=I_{disc}^{\tilde{G}}(f_{V}1_{\tilde{K}^V})-\sum_{{\bf G}'; G'\neq G}i(\tilde{G},{\bf G}')(SI)_{disc}^{{\bf G}'}(f_{V}^{\bf G'}1_{{K_{G'}^V}});\eqno(2)$$
seuls les groupes endoscopiques elliptiques non ramifi\'es hors de $V$ interviennent non trivialement. Ainsi on v\'erifie  facilement que cette distribution
est la trace d'une repr\'esentation semi-finie (parce que le transfert d'une repr\'esentation semi-finie est semi-finie).

On a une d\'ecomposition d'apr\`es \ref{definitionsemifinie}, $$(SI)_{disc}^{\tilde{G}}(f_{V}1_{\tilde{K}^V})=\sum_{\nu,c^V}tr\, \pi_{\nu,st}(c^V)(f_{V}1_{\tilde{K}^V}),$$
la notation anticipe la proposition suivante o\`u on va montrer les propri\'et\'es de stabilit\'e

\begin{prop} On fixe $\nu$ et $c^V$. Alors
$
\pi_{\nu,st}(c^V)=$
$$\pi_{\nu,disc}(c^V)-\sum_{{\bf G}'; {G}'\neq G}i(\tilde{G},{\bf G}')\, \sum_{\nu'\mapsto \nu; c^{'V}\mapsto c^V}transfert(\pi_{\nu',st}^{{\bf G}'}(c^{'V})) \eqno(2)
$$ et 
les repr\'esentations $\pi_{\nu,st}(c^{V})$ sont stables.\end{prop} 
L'\'egalit\'e (2) est facile: on part de la d\'efinition (1) que l'on \'ecrit, pour tout $f_{V}\in I(\tilde{G})(F_{V})$
$$
SI_{disc}^{\tilde{G}}(f_{V}1_{\tilde{K}^V})-\biggl(I_{disc}^{\tilde{G}}(f_{V}1_{\tilde{K}^V})-\sum_{{\bf G}'; G'\neq G}i(\tilde{G},{\bf G}')(SI)_{disc}^{{\bf G}'}(f_{V}^{\bf G'}1_{{K_{G'}^V}})\biggr)=0.
$$
Il faut remarquer que si $V'$ est un ensemble fini de places contenant $V$, on a une \'egalit\'e analogue en rempla\c{c}ant $V$ par $V'$ avec une compatibilit\'e \'evidente si on prend pour $f_{V'}=f_{V}1_{\tilde{K}^{V'-V}}$. Comme on a d\'emontr\'e que les repr\'esentations sous-jacentes aux deux membres de (2) sont semi-finies, on applique \ref{multiplicateurs} pour obtenir (2) et il faut savoir \'evidemment comment se comporte le transfert vis \`a vis des caract\`eres infinit\'esimaux et des caract\`eres des alg\`ebres de Hecke sph\'erique. Pour le caract\`ere infinit\'esimal on renvoie \`a \cite{stabilisationI} corollaire de 2.8 et pour les alg\`ebres de Hecke sph\'erique, c'est le lemme fondamental tel que rappel\'e en \ref{lemmefondamental}.

Montrons la stabilit\'e.
On sait que la distribution $$f_{V}\in I(\tilde{G}(F_{V}))\mapsto I^{\tilde{G}}(f_{V}1_{\tilde{K}^V})-\sum_{{\bf G}'; G'\neq G}i(\tilde{G},{\bf G}')(SI)^{{\bf G}'}(f_{V}^{\bf G'}1_{K_{G'}^V}),$$
\'ecrite sous forme g\'eom\'etrique
est stable (\cite{stabilisationVII} 3.4 sachant que le th\'eor\`eme 3.3 de \cite{stabilisationVII} est d\'emontr\'e de \cite{stabilisationVII} 3.5 \`a 3.8). Sous forme spectrale, elle est une somme de termes index\'es par les espaces de Levi $\tilde{M}$ de $\tilde{G}$. Et on va isoler  le terme qui nous int\'eresse c'est-\`a-dire celui correspondant \`a $\tilde{M}=\tilde{G}$. Fixons $\tilde{M}$ un sous-espace de Levi propre de $\tilde{G}$ et \'ecrivons le terme lui correspondant; on utilise la combinatoire de \cite{stabilisationVI} 6.5. C'est une somme sur les couples $\nu_{M},c_{M}^V$ o\`u $\nu_{M}$ parcourt l'ensemble des caract\`eres infinit\'esimaux de $M(F_{\infty})$ triviaux sur $\mathfrak{A}_{\tilde{M}}$ et $c_{M}^V$ l'ensemble des caract\`eres des alg\`ebres de Hecke sph\'eriques hors de $V$, du terme ci-dessous multipli\'e par le coefficient $w(\tilde{M})$:
$$
I_{\tilde{M}}^{\tilde{G}}(\pi_{\nu_{M}}(c_{M}^V)f_{V}1_{\tilde{K}^V})\eqno(1)$$
$$-\sum_{{\bf M}'\neq M}i(\tilde{M},{\bf M}')\sum_{\nu_{M'},c_{M'}^V\mapsto \nu_{M},c_{M}^V}I_{{\bf M}'}^{\tilde{G},{\mathcal E}}(\pi_{\nu_{M'},st},0,f_{V}1_{\tilde{K}^V})\eqno(2)
$$
$$
- \sum_{s\in Z(\hat{M})^{\Gamma}/Z(\hat{G})^{\Gamma}, s\neq 1)}i_{\tilde{M}}(\tilde{G},{\bf G}'(s))(SI)_{\tilde{M}}^{{\bf G}'(s)}(\pi_{\nu,st},0,f^{{\bf G}'(s)_{V}}1_{\tilde{K}^{'V}}).\eqno(3)
$$
On a montr\'e que dans les termes de (1), on peut enlever l'exposant ${\mathcal E}$ en rempla\c{c}ant les repr\'esentations stables de ${\bf M}'$ par leur transfert \`a $\tilde{M}$. Ainsi (1) moins (2) est exactement $I_{\tilde{M}}^{\tilde{G}}(\pi_{\nu_{M},st}(c_{M}^V),f_{V}1_{\tilde{K}^V})$. Quand on enl\`eve \`a ce terme les termes \'ecrits en (3), on a exactement la d\'efinition de $(SI)_{\tilde{M}}^{\tilde{G}}(\pi_{\nu_{M},st}(c_{M}^V),f_{V}1_{\tilde{K}^V})$. 
Et on a montr\'e que cette distribution est stable en \ref{carglobstab}. Ainsi la somme des trois termes (1), (2), (3) est une distribution stable. Et on en d\'eduit  que la distribution $f_{V}\mapsto SI_{disc}^{\tilde{G}}
(f_{V}1_{\tilde{K}^V})$ est stable.
Il reste \`a s\'eparer les termes en fonction des couples $\nu,c^V$ et montrer que chaque terme est stable, ce qui a \'et\'e fait en \ref{decompositionstable}

\subsection{Enonc\'e de la stabilisation spectrale \label{enoncestabspec}}
On a d\'efini les repr\'esentations semi-finies $\pi_{\nu',st}^{{\bf G}'}(c^{'V})$
pour toute donn\'ee endoscopique de $\tilde{G}$ o\`u $\nu'$ est un caract\`ere infinit\'esimal de ${\bf G}'$ trivial sur $\mathfrak{A}_{G'}$ et $c^{'V}$ est un caract\`ere du produit des alg\`ebres de Hecke sph\'eriques hors de $V$. On pose
$$
\pi_{\nu}^{\tilde{G},{\mathcal E}}(c^V)=\sum_{{\bf G}'}i(\tilde{G},{\bf G}') \sum_{\nu'\mapsto \nu; c^{'V}\mapsto c^V}transfert \, \pi_{\nu',st}^{{\bf G}'}(c^{'V}).
$$
On rappelle que $V$ contient $V_{ram}$ et que ${\bf G}'$ ne parcourt que les groupes endoscopiques elliptiques non ramifi\'es hors de $V$ et relevants.

La stabilisation de la partie discr\`ete spectrale de la formule des traces consiste \`a montrer que $$\pi^{\tilde{G}}_{\nu,disc}(c^V)=\pi_{\nu}^{\tilde{G},{\mathcal E}}(c^V)\eqno(1)$$
pour tout couple $\nu,c^V$. En termes plus explicites, on veut montrer que pour toute fonction $f_{V}\in   I(\tilde{G}(F_{V}))$, on a
$$
tr\, \pi^{\tilde{G}}_{\nu, disc}(c^V)(f_{V}1_{K}^V)=\sum_{{\bf G}'}i(\tilde{G},{\bf G}') \sum_{\nu'\mapsto \nu; c^{'V}\mapsto c^V}tr\, \pi_{\nu',st}^{{\bf G}'}(c^{'V})(f^{{\bf G}'}_{V}1_{K^{'V}}).
$$
Comme la d\'ependance en $V$ n'est qu'un probl\`eme de ramification et que l'on a admis le lemme fondamental (\ref{lemmefondamental}) cela veut dire que pour toute fonction $f\in I(\tilde{G}({\mathbb A}_{F}))$, qui hors d'un nombre fini de places, $V'$ contenant $V$ est la fonction caract\'eristique du compact $K_{v}$ (pour $v\notin V'$) et qui est bi-invariante par $K_{v}$ pour tout $v\in V'-V$, on a
$$
tr \pi^{\tilde{G}}_{\nu,disc}(c^V)(f)=\sum_{{\bf G}'}i(\tilde{G},{\bf G}') \sum_{\nu'\mapsto \nu; c^{'V}\mapsto c^V}\pi_{\nu',st}^{{\bf G}'}(c^{'V})(f^{{\bf G}'}).
$$
\subsection{L'hypoth\`ese spectrale de r\'ecurrence\label{hypspecrec}}
On peut \'evidemment g\'en\'eraliser les d\'efinitions en rempla\c{c}ant $\tilde{G}$ par l'un de ses espaces de Levi et l'hypoth\`ese spectrale de r\'ecurrence est exactement que pour tout espace de Levi propre $\tilde{M}$ de $\tilde{G}$, on a l'\'egalit\'e:
$$
\pi^{\tilde{M}}_{\nu_{M},disc}(c_{M}^V)=\pi_{\nu_{M}}^{\tilde{M},{\mathcal E}}(c_{M}^V).$$

\subsection{R\'eduction de la stabilisation spectrale\label{redstabspecglob}}
\begin{prop} Sous l'hypoth\`ese spectrale de r\'ecurrence et sous l'hypoth\`ese de r\'ecurrence locale g\'eom\'etrique de \ref{stabgeohyp},  pour toute fonction $f$ de la forme $f_{V}1_{\tilde{K}^V}$, on a:
$$
I^{\tilde{G}}(V,f)-\sum_{{\bf G}'}i(\tilde{G},{\bf G}')SI^{{\bf G}'}(V,f^{{\bf G}'})=
I_{disc}^{\tilde{G}}(f)-\sum_{{\bf G}'}i(\tilde{G},{\bf G}')SI_{disc}^{{\bf G}'}(f^{{\bf G}'}).
$$
En particulier le membre de gauche ne d\'epend pas de $V$. La somme sur ${\bf G}'$ ne porte que sur les donn\'ees endoscopiques de $\tilde{G},\omega$ non ramifi\'ees hors de $V$ et relevantes.
\end{prop}
La pr\'esence de $V$ m\'erite explications: les distributions $I$ et $SI$ d\'ependent du choix de $V$ car elles ne sont invariantes qu'en tant que distributions sur $I(\tilde{G}(F_{V}))$. Dans le membre de droite  $I_{disc}^{\tilde{G}}$ qui est spontan\'ement une distribution invariante, n'en d\'epend pas.

En tenant compte de \ref{multiplicateurs} la stabilisation de la partie discr\`ete de la formule des traces, c'est-\`a-dire la nullit\'e du membre de droite de l'\'enonc\'e, est exactement \'equivalente \`a la preuve de l'\'egalit\'e (1) de \ref{enoncestabspec} pour tout $\nu,c^V$.

Montrons la proposition. On fixe $V$ et on  d\'ecompose:
$$\sum_{{\bf G}'}
i(\tilde{G},{\bf G}')(SI^{{\bf G}'}(V,f^{{\bf G}'})-SI_{disc}^{{\bf G}'}(f^{{\bf G}'})).
$$
Comme en \cite{stabilisationVI} 6.5, on transforme cette expression en une somme sur les sous-groupes de Levi des donn\'ees endoscopiques ${\bf G}'$; si ${\bf M}'$ est une donn\'ee endoscopique elliptique pour un espace de Levi $\tilde{M}$ de $\tilde{G}$, le terme correspondant est (\`a un scalaire pr\`es)
$$
i(\tilde{M},{\bf M}')\sum_{\nu_{M'},c_{M'}^V}I_{{\bf M}'}^{\tilde{G},{\mathcal E}}(\pi_{\nu_{M'},st}(c_{M'}^V),f_{V}1_{\tilde{K}^V}).\eqno(1)
$$
Si ${\bf M}'$ n'est relevant pour aucun espace de Levi de $\tilde{G}$, on a une expression similaire avec une d\'efinition formelle de $i(\tilde{M},{\bf M}')$ et on a montr\'e en \ref{carglobtransfert} que ce terme correspondant (1) est nul.  Revenant au cas o\`u ${\bf M}'$ est une donn\'ee relevante et elliptique pour un espace de Levi $\tilde{M}$,
on a d\'ej\`a montr\'e en \ref{carglobtransfert} que sous l'hypoth\`ese g\'eom\'etrique de r\'ecurrence,  l'expression (1) vaut 
$$
i(\tilde{M},{\bf M}')\sum_{\nu_{M'},c_{M'}^V}I_{\tilde{M}}^{\tilde{G}}(transfert\, \pi_{\nu_{M'},st}(c_{M'}^V),f_{V}1_{\tilde{K}^V}).\eqno(2)
$$
Et ceci vaut encore par les d\'efinitions
$\sum_{\nu_{M},c^V_{M}}
I_{\tilde{M}}^{\tilde{G}}(\pi^{\tilde{M},{\mathcal E}}_{\nu_{M}}(c_{M}^V),f_{V}1_{\tilde{K}^V}).
$
Avec l'hypoth\`ese spectrale de r\'ecurrence, on obtient le terme index\'e par $\tilde{M}$ dans l'\'ecriture spectrale de la distribution $f_{V}\mapsto I^{\tilde{G}}(f_{V}1_{\tilde{K}^V})$. Ainsi en sommant sur les espaces de Levi $\tilde{M}$, la somme des contributions des donn\'ees non elliptiques relevantes de $\tilde{G}$ est exactement $I^{\tilde{G}}(V,f)-I_{disc}^{\tilde{G}}(f)$.
\section{Digression, automorphismes de la situation}
\subsection{Action du groupe adjoint ou de son analogue dans le cas tordu\label{groupediese}}
On se place dans une situation locale en fixant une place $v$. On note $G_{\sharp}$ le groupe $G/Z(G)^{\theta}$. Il a \'et\'e montr\'e en \cite{stabilisationI} 2.6 que ce groupe op\`ere sur les facteurs de transfert: soit ${\bf G}'$ une donn\'ee endoscopique elliptique de $\tilde{G},\omega$ alors il existe un caract\`ere $\omega_{\sharp}^{{\bf G}'}$ de $G_{\sharp}(F)$ qui vaut $\omega$ sur l'image de $G(F)$ dans $G_{\sharp}(F)$ et tel que, pour tout $g_{\sharp}\in G_{\sharp}(F)$, pour tout $\gamma\in \tilde{G}_{reg}(F)$ et $\delta$ dans une donn\'ee auxiliaire attach\'ee \`a ${\bf G}'$:
$$
\Delta(\delta,ad(g_{\sharp})^{-1}\gamma)=\omega^{\bf G'}_{\sharp}(g_{\sharp})\Delta(\delta,\gamma).
$$
D'autre part $G_{\sharp}(F)$ op\`ere sur $I(\tilde{G})$ et op\`ere donc sur les $\omega$-repr\'esentations de $\tilde{G}(F)$. Par d\'efinition, le transfert commute \`a ces actions. 

Pour ce qui suit, on suppose que $\omega=1$.

Soit $\xi$ un caract\`ere de $G_{\sharp}(F)/G(F)$; on fixe un syst\`eme de repr\'esentants du groupe fini $G_{\sharp}(F)/G(F)$, not\'e ${\mathcal R}$. Pour tout $f\in I(\tilde{G})$ et tout $\gamma\in \tilde{G}(F)$ on pose $$
f^{\xi} (\gamma):=\vert G_{\sharp}(F)/G(F)\vert^{-1}\sum_{g_{\sharp}\in {\mathcal R}} \xi(g_{\sharp})f(ad(g_{\sharp})^{-1}\gamma).
$$
L'image de $f^\xi$ dans $I(\tilde{G})$ ne d\'epend pas du choix de ${\mathcal R}$.
\begin{lem} Avec les notations pr\'ec\'edentes:

(i) l'\'el\'ement $f^\xi$ de $I(\tilde{G})$ est nul  s'il n'existe pas  une donn\'ee endoscopique elliptique ${\bf G}'$ telle que $\xi=\omega^{{\bf G}'}_{\sharp}$. Et pour toute donn\'ee endoscopique elliptique ${\bf G}'$ de $\tilde{G}$, le transfert $(f^\xi)^{{\bf G}'}$ v\'erifie:
$
(f^\xi)^{{\bf G}'}=0$ si $\omega_{\sharp}^{{\bf G}'}\neq \xi$ et $(f^\xi)^{{\bf G}'}=f^{{\bf G}'}$ sinon.

(ii) On suppose que $\tilde{G}=G$, que la place locale n'est pas dans $V_{ram}$ et que $f$ est bi-invariante par le compact maximal $K$. Alors la composante elliptique de $f^{1_{G_{\sharp}}}$ est nulle sauf \'eventuellement si $G$ est un tore (ici $1_{G_{\sharp}}$ est la notation pour le caract\`ere trivial de $G_{\sharp}$).

(iii) On revient \`a $\tilde{G}$ mais on  suppose que $G$ n'est pas un tore. On suppose aussi que $G$ est d\'eploy\'e et que les donn\'ees endoscopiques elliptiques de $\tilde{G}$ sont \`a torsion int\'erieure triviale. On suppose encore  que $f$ est bi-invariante par le compact maximal $K$. Alors la composante elliptique de $f^{1_{G_{\sharp}}}$ est nulle
\end{lem}

Soit $\pi$ une repr\'esentation temp\'er\'ee de $\tilde{G}(F)$. On \'ecrit $\pi$ comme transfert \`a partir de repr\'esentations temp\'er\'ees stables de donn\'ees endoscopiques elliptiques de ${\tilde{G}(F)}$, c'est-\`a-dire que pour tout ${\bf G}'$ une donn\'ee endoscopique elliptique de $\tilde{G}(F)$, il existe $\pi_{st}^{G'}$ tel que pour tout $f\in I(\tilde{G})$
$$
tr \pi(f)=\sum_{{\bf G}'} tr\pi^{G'}_{st}(f^{{\bf G}'}).
$$
On fait agir $G_{\sharp}(F)$. On calcule le transfert $(f^\xi)^{{\bf G}'}$ en partant des d\'efinitions: pour $g_{\sharp}\in G_{\sharp}(F)$, le transfert $(^{g_\sharp}f)^{{\bf G}'}$ vaut $\omega_{\sharp}^{{\bf G}'}(g_\sharp)^{-1}f^{{\bf G}'}$. On obtient alors les formules de la deuxi\`eme partie de (i) et  $$
tr \pi(f^{\xi})=\sum_{{\bf G}'; \omega_{\sharp}^{{\bf G}'}=\xi} tr\pi^{G'}_{st}(f^{{\bf G}'}). \eqno(1)
$$
 Avec (1) on d\'emontre la totalit\'e de (i).

 On utilise cette \'egalit\'e (1) pour la preuve de (ii) et (iii) et dans (1) on se limite aux repr\'esentations $\pi$ qui sont elliptiques; les repr\'esentations stables du c\^ot\'e droit de l'\'egalit\'e sont alors elles aussi elliptiques. On se place d'abord sous l'hypoth\`ese de (ii). Ici $\tilde{G}=G$ et la situation est non ramifi\'ee. En particulier $f$ est non ramifi\'ee. Pour que le membre de gauche de (1) soit non nul il faut que l'une des repr\'esentations $g_{\sharp}.\pi$ pour $g_{\sharp} \in G_{\sharp}(F)=G_{ad}(F)$ soit non ramifi\'ee. On peut \'evidemment supposer que $\pi$ est ''irr\'eductible'' au sens qu'elle correspond \`a un triplet $(M,\sigma,r)$ de \cite{ftlt} 2.10, 2.11. Cela veut dire que $g_{\sharp}.\pi$ est une repr\'esentation elliptique incluse dans une s\'erie principale non ramifi\'ee  donc $M$ est un tore. On sait que $G_{ad}(F)/G(F)$ agit transitivement par permutation  dans l'ensemble des sous-quotients irr\'eductibles d'une telle induite. Ainsi $\sum_{g_{\sharp}\in G_{ad}(F)/G(F)}g_{\sharp}.\pi$ vaut un multiple de:
$$
(\sum_{\mu\in \hat{R}}\mu(r)) (\sum_{\pi_{\tau}}\pi_{\tau}),
$$
o\`u $\mu$ parcourt les caract\`eres du R-groupe \'evalu\'e en l'\'el\'ement $r\in R$ d\'eterminant $\pi$ et o\`u $\pi_{\tau}$ parcourt l'ensemble des sous-quotients irr\'eductibles de la s\'erie principale consid\'er\'ee. Ceci vaut 0 puisque $r$ est un \'el\'ement r\'egulier de $R$,  sauf si  $M=G$. Dans ce dernier cas, $\pi$ est la s\'erie principale n\'ecessairement irr\'eductible et  $G$ est un tore. Cela montre (ii).

Montrons (iii), on a mis comme hypoth\`ese que les donn\'ees endoscopiques de $\tilde{G}$ sont des groupes non tordus. Le c\^ot\'e droit de (1) est nul sauf s'il existe ${\bf G}'$ une donn\'ee endoscopique de caract\`ere $\omega_{\sharp}^{\bf G'}$ trivial et ayant des repr\'esentations elliptiques stables n'annulant pas une fonction non ramifi\'ee; ici on utlise encore le lemme fondamental tordu. Fixons une telle donn\'ee ${\bf G}'$. Le groupe adjoint de $G'$ agit trivialement dans $SI_{cusp}(G')$. Donc on peut reprendre l'argument donn\'e pour d\'emontrer (ii) qui entra\^{\i}ne que ${\bf G}'$ ne contribue au c\^ot\'e droit de (1) que si $G'$ est un tore. On exploite maintenant le fait que le caract\`ere $\omega_{\sharp}^{{\bf G}'}$ est trivial. On revient \`a la d\'efinition de ce caract\`ere donn\'ee en \cite{stabilisationI} 2.7. Soit $w\in W_{F}$ et $g_{w}=(g(w),w)\in {\mathcal G}'$ donnant l'action de $\Gamma_{F}$ sur $\hat{G}'$; on utilise tout de suite le fait que $\Gamma_{F}$ agit trivialement dans $\hat{G}$. Donc par hypoth\`ese $g({w})$ respecte un \'epinglage de $\hat{G}'$ et v\'erifie, en \'ecrivant $\hat{s}=s\hat{\theta}$ l'\'element de la donn\'ee endoscopique dont le centralisateur est $\hat{G}'$
$$
s \hat{\theta} (g(w))s^{-1}= g(w).
$$
On rel\`eve $s$ et $g(w)$ dans le rev\^etement simplement connexe, $\hat{G}_{SC}$, de $\hat{G}$ en $s_{sc}\in \hat{G}_{SC}$ et $g(w)_{sc}$  et il existe donc $a_{sc}(w)\in Z(\hat{G}_{SC}) $ tel que ($\hat{\theta}$ se rel\`eve canoniquement en un automorphisme de $\hat{G}_{SC}$)
$$
s_{sc} \hat{\theta}(g(w)_{sc})s_{sc}^{-1}=a_{sc}(w)g(w)_{sc}.
$$
On \'ecrit $g(w)=z(w)g(w)_{sc}$, avec $z(w)\in Z(\hat{G})$, en identifiant $g(w)_{sc}\in \hat{G}_{SC}$ \`a son image dans $\hat{G}$. 
Le caract\`ere $\omega_{\sharp}^{{\bf G}'}$ est donn\'e par le cocycle  qui, \`a $w\in W_{F}$, associe l'image de $(z(w),a_{sc}(w))$ dans $Z(\hat{G}_{\sharp})$ (cf. ce qui pr\'ec\`ede l'\'enonc\'e du lemme en \cite{stabilisationI} 2.7). Par hypoth\`ese ce cocycle est trivial et comme $\Gamma_{F}$ agit trivialement sur $\hat{G}$,  cela signifie qu'il existe $z(w)_{sc}\in Z(\hat{G}_{SC})$ de sorte que $z(w)\in z(w)_{sc}\hat{T}^{\hat{\theta},0}$ et $a_{sc}(w)=(1-\hat{\theta})(z(w)_{sc})$. On modifie la d\'ecomposition $g(w)=z(w)g(w)_{sc}$ en rempla\c{c}ant $z(w)$ par $z(w)z(w)_{sc}^{-1}$ et $g(w)_{sc}$ par $z(w)_{sc}g(w)_{sc}$. Les relations pr\'ec\'edentes deviennent
  $z(w)\in \hat{T}^{\hat{\theta},0}$ et  $a_{sc}(w)=1$. Alors $g(w)_{sc}$ est dans le centralisateur de l'automorphisme $ad(s_{sc})\circ \hat{\theta}$. Comme $\hat{G}_{SC}$ est simplement connexe le centralisateur d'un automorphisme semi-simple est connexe. Cette connexit\'e entra\^{\i}ne que l'image du centralisateur de $ad(s_{sc})\circ \hat{\theta}$ est exactement $\hat{G}'$ et donc que l'image de $g(w)_{sc}$ est dans $\hat{G}'$. Puisque $z(w)\in \hat{T}^{\hat{\theta},0}\subset \hat{G}'$, $g(w)$ est aussi dans $\hat{G}'$. Comme $g(w)$ stabilise un \'epinglage de $\hat{G}'$, $g(w)$ est central dans $\hat{G}'$. Ainsi $g_{w}$ agit trivialement sur $\hat{G}'$. Comme $\hat{G}'$ fait partie d'une donn\'ee endoscopique elliptique de $\tilde{G}$ qui est un tore, on a:
$$
\hat{G}'=\hat{G}^{',\Gamma_{F},0}=Z(\hat{G})^{\Gamma,\hat{\theta},0}.
$$
Puisque $\hat{G}'$ contient $\hat{T}^{\hat{\theta},0}$, ce tore est inclus dans $Z(\hat{G})$, ce qui ne peut se produire que si $\hat{G}$ est un tore. Cela  termine la preuve.

\subsection{Fonction caract\'eristique du compact et action du groupe adjoint\label{groupediesespherique}}
On garde la situation locale de \ref{groupediese}. On suppose que la place fix\'ee n'est pas dans $V_{ram}$, en particulier, on est dans une situation locale p-adique non ramifi\'ee. On note encore $1_{\tilde{K}_{v}}$ la fonction caract\'eristique du compact hyperspecial fix\'e et pour tout caract\`ere $\xi$ de $G_{\sharp}(F)/G(F)$, on a d\'efini $1_{\tilde{K}_{v}}^\xi$. On note $K_{\sharp,v}$ le compact maximal hyperspecial de $G_{\sharp}(F)$ contenant  l'image de $K_{v}$ dans $G_{\sharp}(F)$.
\begin{lem} (i) La fonction $1_{\tilde{K}_{v}}^\xi$ est non nulle seulement si $\xi$ est non ramifi\'e, c'est-\`a-dire trivial sur $K_{\sharp,v}$.

(ii) On suppose que $\xi=1$; la composante elliptique de $1_{\tilde{K}_{v}}^\xi$ est nulle dans le cas non tordu, pourvu que $G$ ne soit pas un tore. C'est aussi vrai dans le cas tordu au moins si $\tilde{G}(F_{v})$ a des classes de conjugaison exceptionnelles.
\end{lem}
(i) est une reformulation du (ii) de la proposition de  \cite{stabilisationVII} 2.1 appliqu\'e \`a $\tilde{M}=\tilde{G}$. En effet en appliquant cette proposition \`a $\tilde{M}=\tilde{G}$ on voit que le transfert de la fonction caract\'eristique de $\tilde{K}_{v}$ est nulle pour tout espace endoscopique qui n'est pas non ramifi\'e. D'o\`u avec le (i) du lemme de \ref{groupediese} pour que $1_{\tilde{K}_{v}}^\xi$ soit non nul il faut que $\xi=\omega_{\sharp}^{{\bf G}'}$ pour une donn\'ee endoscopique non ramifi\'ee ${\bf G}'$. Pour une telle donn\'ee $\omega_{\sharp}^{{\bf G}'}$ est non ramifi\'e d'apr\`es \cite{stabilisationVII} 2.1 (3).

 L'assertion (ii) dans le cas non tordu est cons\'equence du (ii) du lemme 6.1. Dans le cas tordu, on suppose que $\tilde{G}(F_{v})$ a des classes de conjugaison exceptionnelles. D'apr\`es \cite{stabilisationIII} 6.3, cela entra\^{\i}ne que $\tilde{G}$ provient par changement de base (ce qui est innofensif) d'un espace pour lequel les hypoth\`eses du (iii) du lemme 6.1 sont satisfaites. Cela entra\^{\i}ne de plus que $G$ est simplement connexe. Alors l'assertion r\'esulte de cette assertion  (iii) du lemme 6.1.

\subsection{ Action globale du groupe adjoint et de son analogue dans le cas tordu\label{groupedieseglobal}}
On revient \`a une situation globale. On d\'efinit sur $F$, le groupe $G_{\sharp}:= G/Z(G)^{\theta}$. Et les caract\`eres automorphes de $G_{\sharp}({\mathbb A}_{F})$ sont d\'efinis de fa\c{c}on usuelle. On fixe un ensemble $V$ de places de $F$ contenant $V_{ram}$. Clairement les \'el\'ements de $\tilde{G}(F)$ agisse par conjugaison sur $G_{\sharp}$
\begin{rmq} Il n'existe qu'un nombre fini de caract\`eres automorphes de $G_{\sharp}$ non ramifi\'es hors de $V$, triviaux sur l'image de $G({\mathbb A}_{F})$ ou plus g\'en\'eralement prolongeant un caract\`ere fix\'e de $G({\mathbb A}_{F})$.
\end{rmq}
On fixe un ensemble $V'$ de places de $F$, contenant $V$ et tel que $$G_{\sharp}({\mathbb A}_{F})=G_{\sharp}(F) G_{\sharp}(F_{V'})K_{\sharp}^{V'}\eqno(1).$$ La remarque r\'esulte alors facilement du fait que $G_{\sharp}(F_{V'})/G(F_{V'})$ est un groupe fini.

Pour tout $v$ on fixe un caract\`ere $\xi(v)$ de $G_{\sharp}(F_{v})$ prolongeant $\omega_v$.  On  a donc d\'efini $f_{v}^{\xi({v})}$ pour tout $f_{v}\in I(\tilde{G}(F_{v}))$. On note $\xi(V):=\prod_{v}\xi(v)$ et on d\'efinit $f_{V}^{\xi(V)}$ pour tout $f_{V}\in I(\tilde{G}(F_{V}))$, d\'ecompos\'ee, en faisant le produit des $f_{v}^{\xi(v)}$.

On note ${\mathcal C}(V)$ l'ensemble des caract\`eres  de la remarque pr\'ec\'edente prolongeant $\omega$ et pour $\nu$ un caract\`ere infinit\'esimal de $G(F_{\infty})$, on note $\pi_{disc,\nu}^{\tilde{G},V}$ la somme des $\pi^{\tilde{G}}_{\nu,disc}(c^V)$ de \ref{rappeldestermesdiscrets} quand $c^V$ varie. C'est la partie discr\`ete de la formule des traces pour un caract\`ere infinit\'esimal fix\'e et non ramifi\'e hors de $V$.

\begin{lem} Soit $f_{V}\in I(\tilde{G}(F_{V}))$. 
Alors $tr\, \pi_{disc,\nu}^{\tilde{G},V}(f^{\xi(V)}_{V}1_{K^V})=0$ s'il n'existe par $\xi\in {\mathcal C}(V)$ tel que $\xi(V)=\xi_{V}$ c'est-\`a-dire que $\xi(V)$ soit la restriction de $\xi$ \`a $G_{\sharp}(F_{V})$.
\end{lem}
Pour d\'emontrer ce lemme, on est en droit d'augmenter $V$: en effet soit $V'\supset V$. On a clairement $tr\, \pi_{disc,\nu}^{\tilde{G},V}(f^{\xi(V)}_{V}1_{K^V})=tr\, \pi_{disc\nu}^{\tilde{G},V'}(f^{\xi(V)}_{V}1_{K^V})$.

On d\'ecompose $1_{K_{V'-V}}$ suivant les caract\`eres n\'ecessairement non ramifi\'es du groupe $G_{\sharp}(F_{V'-V})$ et il suffit \'evidemment de d\'emontrer la proposition pour $V'$ et le produit de $\xi(V)$ avec l'un de ces caract\`eres. On suppose donc que $V$ est tel que (1) ci-dessus soit v\'erifi\'e pour $V'=V$. Pour d\'emontrer la proposition, on doit donc montrer que $\xi_{V}$ est trivial sur $G_{\sharp}(F)K_{\sharp}^V\cap G(F_{V})$. Soit $g_{\sharp,V}$ dans cet intersection. On calcule $^{g_{\sharp}}f^{\xi(V)}=\xi(V)(g_{\sharp})^{-1}f^{\xi(V)}$. On \'ecrit aussi $g_{\sharp,V} g^V=g_{F}$ o\`u $g_{F}\in G_{\sharp}(F)$ et $g^V\in K_{\sharp}^V$. Le groupe $G_{\sharp}(F)$ op\`ere directement sur l'espace de la repr\'esentation $\pi_{disc,\nu}^{\tilde{G},V}$ \`a une double difficult\'e pr\`es: la partie $t$-discr\`ete mais non discr\`ete d\'ecrite en \ref{rappeldestermesdiscrets} fait intervenir des induites, un \'el\'ement de $G_{\sharp}(F)$ transforme cette induite en une autre induite mais les semi-simplifi\'es sont les m\^emes. L'autre difficult\'e vient du fait que les espaces de formes automorphes consid\'er\'es sont $K_{\infty}$-finis et $G_{\sharp}(F)$ ne conserve pas cette finitude puisque les \'el\'ements de $G_{\sharp}(F)$ remplacent  $K_{\infty}$ en un de ses conjugu\'es sous $G(F_{\infty})$ (les compacts maximaux sont conjugu\'es en les places archim\'ediennes); mais on peut \'elargir l'espace des repr\'esentations automorphes consid\'er\'ees pour qu'il soit conserv\'e par l'action de $G_{\sharp}(F)$ et ensuite on ne calcule que la trace que sur les \'el\'ements de $I(\tilde{G}(F_{\infty}))$, or cet espace est en fait ind\'ependant du choix de $K_{\infty}$ puisque les diff\'erents choix sont conjugu\'es.

On remarque aussi que l'action de $G_{\sharp}(F)$ pr\'eserve le caract\`ere infinit\'esimal. On fait agir un \'el\'ement $g_{F}=g_{\sharp,V}g^V$ avec $g^V$ dans $\tilde{K}^V$ ainsi $g_{F}$ pr\'eserve la non ramification hors de $V_{ram}$.  Et pour toute fonction $f\in I(\tilde{G}({\mathbb A}_{F}))$ et pour tout $g_{F}\in G_{\sharp}(F)$ comme pr\'ec\'edemment, on a
$$
tr\, \pi_{disc,\nu}^{\tilde{G},V} (^{g_{F}}f)=tr\, \pi_{disc,\nu}^{\tilde{G},V}(f).
$$
Les \'el\'ements de $K_{\sharp}^V$ agisse trivialement sur $f_{V}1_{K^V}$ et on en d\'eduit donc que l'action de $$\xi(V)(g_{\sharp,V})^{-1}tr\, \pi_{disc,\nu}^{\tilde{G},V}(f^{\xi(V)}_{V}1_{K^V})=tr\, \pi_{disc,\nu}^{\tilde{G},V}(\, ^{g_{\sharp,V}}f^{\xi(V)}_{V}1_{K^V})=tr\, \pi_{disc,\nu}^{\tilde{G},V}(f^{\xi(V)}_{V}1_{K^V})$$ce qui force $\xi(V)(g_{\sharp,V})=1$ comme voulu.

\

On fixe $V$ suffisamment grand pour que (1) soit satisfait. Soit $f_{V}\in I(\tilde{G}(F_{V}))$ et $\xi \in {\mathcal C}$. On note $\xi(V)$ la restriction de $\xi$ \`a $G_{\sharp}(F_{V})$.
\begin{prop}  Avec ces hypoth\`eses et notations, on a pour toute place $v'$ de $F$ non dans $V$:
$$
tr\, \pi_{disc,\nu}^{\tilde{G},V}(f^{\xi(V)}_{V}1_{K^V})=tr\, \pi_{disc,\nu}^{\tilde{G},V\cup \{v'\}}(f^{\xi(V)}_{V}1^{\xi_{v'}}_{K_{v'}}1_{K^{V- \{v'\}}})
$$
\end{prop}
On a \'evidemment $tr\, \pi_{disc,\nu}^{\tilde{G},V}(f^{\xi(V)}_{V}1_{K^V})=tr\, \pi_{disc,\nu}^{\tilde{G},V\cup \{v'\}}(f^{\xi(V)}_{V}1_{K^V})=$
$$
\sum_{\xi(v')}tr\, \pi_{disc,\nu}^{\tilde{G},V\cup \{v'\}}(f^{\xi(V)}_{V}1^{\xi({v'})}_{K_{v'}}1_{K^{V- \{v'\}}}),
$$
o\`u la somme porte sur les caract\`eres de $G_{\sharp}(F_{v'})/
G(F_{v'})$. On a vu en \ref{groupediesespherique} que $\xi(v')$ est n\'ecessairement non ramifi\'e pour r\'eelement intervenir. Par le lemme pr\'ec\'edent on sait que $\xi(V)\xi(v')$ est la restriction \`a $G_{\sharp}(F_{V\cup\{v'\}})$ d'un caract\`ere dans ${\mathcal C}$ et comme $V$ est suffisamment grand pour que (1) soit satisfait, ce caract\`ere est uniquement d\'etermin\'e par $\xi(V)$ et vaut donc $\xi$. D'o\`u la proposition.

\section{Fin de la stabilisation locale g\'eom\'etrique}

\subsection{Mise en place des objets\label{miseenplace}}
\bf{Hypoth\`ese de r\'ecurrence} \sl  On suppose que l'hypoth\`ese locale g\'eom\'etrique de r\'ecurrence de 3.5  est satisfaite pour tout espace de Levi $\tilde{L}$ contenant strictement $\tilde{M}$, c'est-\`a-dire que pour tout $\tilde{L}$ contenant strictement $\tilde{M}$, pour toute place $v\in V$, et pour tout $f_{v}\in I(\tilde{G}(F_{v}))$, $I_{\tilde{L}}^{\tilde{G}}(\omega,f_{v})=I_{\tilde{L}}^{\tilde{G},{\mathcal E}}(f_{v})$.

On suppose aussi que pour toute place $v\in V$ 
 et pour toute paire $(\tilde{L},\tilde{L}')$ o\`u $\tilde{L}'$ est un espace de Levi propre de $\tilde{G}$ contenant $\tilde{L}$,  on a la stabilisation locale g\'eom\'etrique c'est-\`a-dire que pour tout $f_{v}\in I(\tilde{L}'(F_{v}))$, $I_{\tilde{L}}^{\tilde{L}'}(\omega,f_{v})=I_{\tilde{L}}^{\tilde{L}',{\mathcal E}}(f_{v})$.

 On suppose aussi que l'hypoth\`ese spectrale globale de r\'ecurrence est satisfaites, c'est-\`a-dire qu'avec les notations de \ref{hypspecrec}, pour tout sous-espace de Levi propre, $\tilde{L}$ de $\tilde{G}$, et pour tout $f_{V}\in I(\tilde{L}(F_{V}))$, on a $tr\, \pi_{disc,\nu}^{\tilde{L}}(c^V)(f_{V}1_{\tilde{K}^{V,\tilde{L}}})=tr\, \pi_{disc,\nu}^{\tilde{L},{\mathcal E}}(c^V)(f_{V}1_{\tilde{K}^{V,\tilde{L}}})$. 
\rm

Soit $f_{V}\in I(\tilde{G}(F_{V}))$ et fixons $v_{0}\in V$. On reprend les constructions de l'application $\epsilon_{\tilde{M},v_{0}}$ donn\'ee en \cite{stabilisationVIII} 4.4 pour $v_{0}$, $p$-adique et dans le paragraphe 8 de  \cite{stabilisationIX}  si $v_{0}$ est une place archim\'edienne; cette application envoie l'ensemble des fonctions sur $\tilde{G}(F_{v})$ nulles pr\`es des \'el\'ements exceptionnels de $\tilde{G}(F_{v_{0}})$, s'il y en a et si $v_{0}$ est $p$-adique,  dans $I_{ac}(\tilde{M}(F_{v_{0}}))$. On pose, pour de telles fonctions $f_{V}$: $$h^f(v_{0}):= \epsilon_{\tilde{M},v_{0}}(f_{v_{0}})\prod_{v\in V, v\neq v_{0}}f_{v,\tilde{M}}.$$ 
C'est un \'el\'ement de $I_{ac}(\tilde{M}(F_{V}))$. Soit $\pi_{V}$ une repr\'esentation de longueur finie de $\tilde{M}(F_{V})$ triviale sur $\mathfrak{A}_{\tilde{M}}$ et ayant un caract\`ere infinit\'esimal, not\'e $\mu_{\pi}$. On suppose en plus que $\pi_{V}$ est unitaire. 

Pour $X\in {\mathcal{A}}_{\tilde{M}}$, on a d\'efini $I^{\tilde{M}}(\pi_{V},X,h^f({v_{0}}))$, ce qui utilise l'unitarit\'e de $\pi_{V}$. On va montrer que c'est un coefficient de Fourier d'une fonction que l'on va un peu expliciter pour pouvoir contr\^oler l'action des multiplicateurs de $\tilde{G}$ agissant sur $f_{V}$.

\begin{lem} Il existe une fonction m\'eromorphe sur ${\mathcal{A}}^*_{\tilde{M}}$, not\'ee $I^{\tilde{M}}(\lambda,h^f(v_{0}))$ \`a d\'ecroissance rapide quand on restreint cette fonction \`a l'axe des  $\lambda$ unitaires et telle que
$$
I^{\tilde{M}}(\pi_{V},X,h^f({v_{0}}))=\int_{i{\mathcal{A}}^*_{\tilde{M}}}d\lambda\, I^{\tilde{M}}(\lambda,h^f({v_{0}}))e^{-\lambda(X)}.
$$
De plus pour tout \'el\'ement $z$ du centre de l'alg\`ebre enveloppante de $G(F_{\infty})$, on a:
$$
I^{\tilde{M}}(\lambda,h^{z.f}(v_{0}))=z(\mu_{\pi}+\lambda)I^{\tilde{M}}(\lambda,h^f(v_{0})).
$$
\end{lem}
On note $\pi_{V,v_{0}}$ l'unique repr\'esentation de $\tilde{M}(F_{V})$, elliptique (modulo le centre) en la place $v_{0}$ et qui a m\^eme trace que $\pi_{V}$ sur toutes les fonctions dans $I(\tilde{M}(F_{V}))$ cuspidales en la place $v_{0}$. Cela revient \`a d\'ecomposer dans le groupe de Grothendieck en la place $v_{0}$ suivant la base form\'ee par les induites de repr\'esentations elliptiques et en enlevant toutes les induites propres.  Montrons que
$$
I^{\tilde{M}}(\pi_{V},X,h^f({v_{0}}))=I^{\tilde{M}}(\pi_{V,v_{0}},X,h^f({v_{0}}));$$
on sait que les int\'egrales orbitales de $h^f(v_{0})$ sont nulles en les points de $\tilde{M}(F_{V})$ dont la composante en la place $v_{0}$ est conjugu\'ee d'un \'el\'ement non elliptique de $\tilde{M}(F_{v_{0}})$ (cf. \cite{stabilisationVIII} 4.4 et  \cite{stabilisationIX} 8.9 (ii)). Donc le seul probl\`eme est que $h^f(v_{0})$ n'a pas dans ses propri\'et\'es d'\^etre \`a support compact. Soit $b$ une fonction \`a support compact sur ${\mathcal{A}}_{\tilde{M}}$ qui vaut $1$ en $X$. On voit $b$ comme une fonction sur $\tilde{M}(F_{V})$ et le produit  $bh^f(v_{0})$ est maintenant \`a support compact. Ce produit  v\'erifie toujours les propri\'et\'es de nullit\'e de certaines  int\'egrales orbitales au m\^eme titre que  $h^f(v_{0})$. D'o\`u
$$
I^{\tilde{M}}(\pi_{V},X,bh^f({v_{0}}))=I^{\tilde{M}}(\pi_{V,v_{0}},X,bh^f({v_{0}}));$$
mais aussi $$I^{\tilde{M}}(\pi_{V},X,h^{f}(v_{0}))=I^{\tilde{M}}(\pi_{V},X,bh^{f}(v_{0}))$$
$$\hbox{ et }I^{\tilde{M}}(\pi_{V,v_{0}},X,h^f({v_{0}}))=I^{\tilde{M}}(\pi_{V,v_{0}},X,bh^f({v_{0}})),$$ d'o\`u le r\'esultat annonc\'e.

Fixons $v\in V$; si $v\neq v_{0}$, pour toute repr\'esentation $\rho_{v}$ de $\tilde{M}(F_{v})$ et pour tout $X_{v}\in {\mathcal{A}}_{\tilde{M}(F_{v})}$, on a
$$
I^{\tilde{M}}(\rho_{v},X_{v},f_{\tilde{M},v})=\int_{i{\mathcal{A}}^*_{\tilde{M}(F_{v})}}d\lambda_{v}I^{\tilde{M}}(\rho_{v},\lambda_{v},f_{v}),\eqno(2)
$$
o\`u $I^{\tilde{M}}(\rho_{v},\lambda_{v},f_{v})= \, tr\, ind (\rho_{v\otimes \lambda_{v}})(f_{v})e^{-\lambda_{v}(X_{v})}.$

En $v_{0}$, il existe encore, (cf. le d\'ebut de la preuve de \cite{stabilisationIX}, 8.9) une fonction m\'eromorphe de $\lambda_{v_{0}} \in {\mathcal{A}}^*_{\tilde{M}(F_{v_{0}})}$, not\'ee $I^{\tilde{M}}(\rho_{v_{0}},\lambda_{v_{0}},\epsilon_{\tilde{M},v_{0}}(f_{v_{0}}))$ telle que l'on ait l'analogue de (2) pour $\rho_{v_{0}}$ une repr\'esentation temp\'er\'ee de $\tilde{M}(F_{v_{0}})$, c'est-\`a-dire
$$
I^{\tilde{M}}(\rho_{v_{0}},X_{v_{0}},\epsilon_{\tilde{M},v_{0}}(f_{v_{0}}))=\int_{i{\mathcal{A}}^*_{\tilde{M}(F_{v_{0}})}}d\lambda_{v_{0}}I^{\tilde{M}}(\rho_{v_{0}},\lambda_{v_{0}},\epsilon_{\tilde{M},v_{0}}(f_{v_{0}})).$$
 Et si $v_{0}$ est une place archim\'edienne, remplacer $f$ par $zf$ avec $z$ dans le centre de l'alg\`ebre enveloppante de $\tilde{G}(F_{v_{0}})$ multiplie par la fonction m\'eromorphe de $\lambda_{v_{0}}$ par $z(\mu_{\rho_{v_{0}}}+\lambda_{v_{0}})$ (si $\rho_{v_{0}}$ a un caract\`ere infinit\'esimal, $\mu_{\rho_{v_{0}}}$) d'apr\`es la m\^eme r\'ef\'erence dans le cas o\`u $v_{0}$ est une place archim\'edienne.

Dans le cas o\`u $v$ est une place archim\'edienne diff\'erente de $v_{0}$ la transformation sous le centre de l'alg\`ebre enveloppante est \'evidemment compl\`etement claire.

Puisque l'on peut remplacer $\pi_{V}$ par $\pi_{V,v_{0}} $ pour calculer $I^{\tilde{M}}(\pi,X,h^f(v_{0}))$, on peut r\'ealiser cette fonction de $X$ comme transform\'ee de Fourier de la fonction de $\lambda\in {\mathcal{A}}^*_{\tilde{M}}$ qui vaut en $\lambda$ (cf. \ref{carponddef} appliqu\'e \`a $\tilde{G}=\tilde{M}$),  $$ I(\rho_{v_{0}},\lambda,\epsilon_{\tilde{M},v_{0}}(f_{v_{0}}))\prod_{v\not=v_{0}}I(\rho_{v},\lambda,f_{v}),
$$ o\`u on a \'ecrit $\pi_{V,v_{0}}=\otimes_{v\in V}\rho_{v}$.
Cela donne le lemme.
\begin{cor} Pour tout multiplicateur $\alpha$ de $\tilde{G}$, on a l'\'egalit\'e
$$I^{\tilde{M}}(\pi_{V},X,h_{v_{0}}^{f_{\alpha}})=\int_{i{\mathcal{A}}^*_{\tilde{M}}}d\lambda\,\hat{\alpha}(\nu_{\pi_{V}}+\lambda) I^{\tilde{M}}(\lambda,h^f_{v_{0}})e^{-\lambda(X)}.
$$
\end{cor}
La famille de formes lin\'eaires (d\'ependant m\'eromorphiquement de $\lambda\in {\mathcal{A}}^*_{\tilde{M}})$
$$
\lambda \in {\mathcal{A}}^*_{\tilde{M}}\mapsto \ell(\lambda,f_{V}):=I^{\tilde{M}}(\lambda,h^f(v_{0}))$$ se voit sur l'ensemble des fonctions de Paley-Wiener sur l'espace des $\omega$-repr\'esentations temp\'er\'ees de $\tilde{G}(F_{V})$; ces applications lin\'eaires, $f_{V}\mapsto \ell(\lambda,f_{V})$, v\'erifient $\ell(\lambda,z.f_{V})=z(\nu_{\pi_{V}}+\lambda)\ell(\lambda,f_{V})$. 

Fixons $\lambda$ telle que la forme lin\'eaire $\ell(\lambda,f_{V})$ soit d\'efinie. C'est une forme lin\'eaire sur $I(\tilde{G}(F_{V}))/{\mathcal M}_{\chi_{\lambda}}I(\tilde{G}(F_{V}))$, o\`u ${\mathcal M}_{\chi_{\lambda}}$ est l'id\'eal maximal du centre de l'alg\`ebre enveloppante correspondant au caract\`ere infinit\'esimal $\nu_{\pi_{V}}+\lambda$. Soit $q$ l'application naturelle de $I(\tilde{G}(F_{V}))$ dans ce quotient et on va montrer l'assertion:
\sl pour tout $\alpha,f_{V}$ comme ci-dessus, $q(f_{V,\alpha})=\check{\alpha}(\nu_{\pi_{V}}+\lambda)q(f_{V})$.\rm

Pour d\'emontrer cette assertion, il faut le faire pour $f_{v_{\infty}}$ dans une composante de Paley-Wiener. On fixe donc une telle composante que l'on note ${\mathcal F}$ et on utilise le lemme de \cite{stabilisationIV} 2.5 de la fa\c{c}on suivante.
On note $\Pi_{\chi_{\lambda}}$ l'ensemble des repr\'esentations induites de repr\'esentations elliptiques de $G(F_{v_{\infty}})$ ayant m\^eme caract\`ere infinit\'esimal que l'induite de $\pi_{\infty}\otimes \lambda$ de $\tilde{M}$ \`a $\tilde{G}$. C'est un ensemble fini. Si la composante de Paley-Wiener ne coupe pas $\Pi_{\chi_{\lambda}}$ alors ${\mathcal M}_{\chi_{\lambda}}{\mathcal F}={\mathcal F}$ d'apr\`es le (i) de ce lemme. Si l'intersection est non vide, alors ${\mathcal M}_{\chi_{\lambda}}{\mathcal F}$ contient les \'el\'ements de ${\mathcal F}$ nuls \`a un ordre suffisant, notons le $N$, en les points de l'intersection; c'est le (ii) du lemme. Soit $\phi$ un \'el\'ement de ${\mathcal F}$ et $\alpha$ un multiplicateur; si $f_{v_{\infty}}\in I(G(F_{v_{\infty}}))$ correspond \`a $\phi$ alors $f_{v_{\infty},\alpha}$ correspond au produit $\check{\alpha}\phi$ o\`u $\check{\alpha}\phi$ en une repr\'esentation $\pi$ est le produit de $\phi(\pi)$ avec la valeur de $\check{\alpha}$ sur le caract\`ere infinit\'esimal de $\pi$. On peut approximer la restriction de  $\check{\alpha}$ \`a $\Pi_{\chi_{\sigma}} \cap {\mathcal F}$ en utilisant le centre de l'alg\`ebre enveloppante, c'est \`a dire qu'il existe $z_{\alpha}$ dans ce centre tel que $\check{\alpha}-z_{\alpha}$ s'annule en l'ordre au moins $N$ en le caract\`ere infinit\'esimal $\chi_{\lambda}$. Et alors $$q(f_{v_{\infty},\alpha})=q(z_{\alpha}f)=z_{\alpha}(\nu_{\pi_{V}}+\lambda)q(f)=\check{\alpha}(\nu_{\pi_{V}}+\lambda)q(f),$$ ce qui est l'assertion cherch\'ee.
Et le corollaire s'en d\'eduit.
\subsection{Stabilisation de la formule des traces pour certaines fonctions\label{casparticulier}}
Soit $f_{V}\in I(\tilde{G}(F_{V}))$; on suppose que pour tout $v\in V$,  $f_{v}$ est nulle pr\`es des \'el\'ements exceptionnels, s'il y en a et si $v$ est $p$-adique, et on suppose que si $v$ est une place archim\'edienne, $f_{v}$ est \`a support dans les \'el\'ements r\'eguliers (c'est une condition qui vient du probl\`eme technique suivant: tant que l'on ne sait pas stabiliser la formule des traces locales, les termes $I_{\tilde{M}}^{\tilde{G},{\mathcal E}}(\gamma_{v},v_{v})$ ne sont pas d\'efinis pour $\gamma_{v}$ un \'el\'ement singulier). On pose pour une telle fonction $f_{V}$:
$$
h_{V}^f:=\sum_{v\in V}\epsilon^{\tilde{G}}_{\tilde{M},v}(f_{v})\prod_{v''\in V-\{ v\}}f_{\tilde{M},v''}.
$$
L'ensemble $V$ contient au moins deux places (il contient $V_{ram}$ qui contient les places divisant $2$, $3$ et $5$) et on fixe $v_{1},v_{2}\in V$. On suppose que les fonctions $f_{V}$ consid\'er\'ees sont $\tilde{M}$-cuspidales en $v_{1}$ et $v_{2}$ (ce qui veut dire que $f_{\tilde{M},v_{i}}$ est cuspidale pour $i=1,2$). Ainsi $h_{V}^f$ est somme de fonctions  cuspidales en au moins deux places; on rappelle que $\epsilon_{\tilde{M},v}^{\tilde{G}}$ est une application de $I(\tilde{G}(F_{v}))$ dans $I_{cusp,ac}(\tilde{M}(F_{v}))$. On rappelle aussi que la premi\`ere somme sur $v$ ne porte en fait que sur les \'el\'ements de $v\in V$ tel que ${\mathcal A}_{\tilde{M}}={\mathcal A}_{\tilde{M},v}$ car sinon on sait d\'ej\`a avec les formules de descente que $\epsilon_{\tilde{M},v}^{\tilde{G}}$ est une application identiquement nulle. Notons $V'$ l'ensemble des places de $F$ appartenant \`a  $V$ o\`u l'\'egalit\'e ${\mathcal A}_{\tilde{M}}={\mathcal A}_{\tilde{M},v}$ est satisfaite.

\begin{lem} Sous l'hypoth\`ese de r\'ecurrence de \ref{miseenplace}, on a pour tout $f_{v\in V}$ comme ce qui pr\'ec\`ede l'\'enonc\'e, avec les notations des paragraphes pr\'ec\'edents:
$$
I^{\tilde{G}}_{geo}(\omega,f_{V}1_{\tilde{K}^V})-I^{\tilde{G},{\mathcal E}}_{geo}(\omega,f_{V}1_{\tilde{K}^V})=x I^{\tilde{M}}_{geo}(\omega,h^f_{V}1_{\tilde{K}^{M,V}}),
$$o\`u $x$ est un scalaire non nul explicite (ind\'ependant de $f_{V}$)
\end{lem}
On reprend les notations de \cite{stabilisationVI} 6.11. Le terme de gauche de l'\'enonc\'e est une somme sur les espaces de Levi $\tilde{L}$ de $\tilde{G}$ des termes:
$$\sum_{{\mathcal O}^{\tilde{L}}}
I_{\tilde{L}}^{\tilde{G}}(A^{\tilde{L}}({\mathcal O}^{\tilde{L}},V,\omega),f)$$
$$-\sum_{{\bf L}' \in {\mathcal E}(\tilde{L}, \omega,V)}i(\tilde{L},{\bf L}')\sum_{{\mathcal O}^{{\bf L}'}\in \tilde{L}'_{ss}(F_{V})/st-conj; {\mathcal O}^{L'}\mapsto {\mathcal O}^{\tilde{L}}}I_{\tilde{L}}^{\tilde{G},{\mathcal E}}({\bf L}',SA^{L'}({\mathcal O}^{{\bf L}'}),V,f).
$$C'est la diff\'erence entre les deux valeurs de $X$ de \cite{stabilisationVI} 6.11.
On sait que les termes correspondants \`a $\tilde{L}=\tilde{G}$ sont nuls, d'apr\`es \cite{stabilisationVII} 3.3  sauf (\'eventuellement) pour les \'el\'ements exceptionnels o\`u le r\'esultat est encore \`a montrer. Mais l'hypoth\`ese que $f_{V}$ est nulle pr\`es de ces \'el\'ements  en les places non-archim\'ediennes  assure aussi la nullit\'e de ces termes.

Supposons que $\tilde{L}$ soit un espace de Levi propre de $\tilde{G}$; par r\'ecurrence on sait que $A^{\tilde{L}}({\mathcal O}^{\tilde{L}},V,\omega)$ est le transfert de $\sum_{{\bf L}'}i(\tilde{L},{\bf L}')SA^{{\bf L}'}({\mathcal O}^{{\bf L}'},V)$; le terme \'ecrit est donc une combinaison lin\'eaire de distributions de la forme $I_{\tilde{L}}^{\tilde{G}}(\omega, \gamma,f_{V})-I_{\tilde{L}}^{\tilde{G},\mathcal E}(\omega,\gamma,f_{V})$ pour $\gamma$ des \'el\'ements rationnels de $\tilde{L}$ pris \`a conjugaison pr\`es. Par hypoth\`ese de r\'ecurrence et les formules de descente, il ne reste qu'une combinaison lin\'eaire de termes, de la forme \`a un scalaire pr\`es
$$\sum_{v_{0}\in V}
\biggl(I_{\tilde{L}}^{\tilde{G}}(\omega_{v_{0}},\gamma_{v_{0}},f_{v_{0}})-I_{\tilde{L}}^{\tilde{G},{\mathcal E}}(\omega_{v_{0}},\gamma_{v_{0}},f_{v_{0}})\biggr)
\prod_{v\in V-\{v_{0}\}} I_{\tilde{L}}^{\tilde{L}}(\omega_{v},\gamma_{v},(f_{v})_{L}).\eqno(2)
$$
En utilisant les propri\'et\'es de cuspidalit\'e de la fonction $f_{V}$ en les places $v_{1}$ et $v_{2}$ un tel terme est certainement nul si $\tilde{L}$ ne contient pas (\`a conjugaison pr\`es) $\tilde{M}$; il est donc ici important de sommer sur les classes de conjugaison d'espaces de Levi et non pas sur les espaces de Levi semi-standard, d'o\`u des coefficients   qu'il n'est pas important de calculer. Si $\tilde{L}$ contient strictement $\tilde{M}$, les hypoth\`eses de r\'ecurrence donnent la stabilisation locale g\'eom\'etrique pour les \'el\'ements semi-simples r\'eguliers et on a donc cette stabilisation g\'eom\'etrique pour tout \'el\'ement gr\^ace \`a \cite{stabilisationV} 1.11.

Finalement il ne reste que les termes avec $\tilde{L}=\tilde{M}$; en utilisant les propri\'et\'es de descente, la premi\`ere somme ne porte que sur $v_{0}\in V'$ (o\`u $V'$ a \'et\'e d\'efini avant l'\'enonc\'e). Par cuspidalit\'e, il reste encore uniquement les termes correspondant \`a des \'el\'ements semi-simples r\'eguliers elliptiques de $\tilde{M}$. Or (2) n'est autre que $I_{\tilde{M}}^{\tilde{M}}(\omega,\gamma_{V},h^{f}_{V})$. En revenant \`a $A^{\tilde{M}}({\mathcal O}^{\tilde{M}})$, on voit que le terme de gauche de l'\'enonc\'e est exactement
$$
\sum_{{\mathcal O}^{\tilde{M}}}I^{\tilde{M}}(A^{\tilde{M}}({\mathcal O}^{\tilde{M}},V,\omega),h^{f}_{V}).
$$
Maintenant on utilise les propri\'et\'es de cuspidalit\'e de $h^f_{V}$ pour v\'erifier que ceci est exactement $I_{geo}^{\tilde{M}}(\omega,V,h^f_{V})$.
\begin{prop} Avec les hypoth\`eses et notations du lemme pr\'ec\'edent et celles de \ref{miseenplace}, pour tout $\nu,c^V$
$$
tr\, \pi_{\nu,disc}^{\tilde{G}}(c^V)(f_{V}1_{\tilde{K}^V})=tr\, \pi_{\nu,disc}^{\tilde{G},{\mathcal E}}(c^V)(f_{V}1_{\tilde{K}^V})$$
\end{prop}
Le lemme pr\'ec\'edent donne l'\'egalit\'e 
$$
I_{geo}^{\tilde{G}}(\omega,f_{V})-I_{geo}^{\tilde{G},{\mathcal E}}(\omega,f_{V})=x I^{\tilde{M}}_{geo}(\omega,h^f_{V}). \eqno(1)
$$
On remplace ces distributions, \'ecrites \`a l'aide du c\^ot\'e g\'eom\'etrique, par les distributions qui leur sont \'egales mais \'ecrites \`a l'aide du c\^ot\'e spectral. On a vu en \ref{redstabspecglob} que le terme de gauche vaut 
$$I_{disc}^{\tilde{G}}(\omega,f_{V})-I_{disc}^{\tilde{G},{\mathcal E}}(\omega,f_{V}).\eqno(2)
$$
Le terme de droite vaut $x I_{disc}^{\tilde{M}}(\omega,h^f_{V}1_{\tilde{K}^{\tilde{M},V}})$ car $h^f_{V}$ est une somme finie de fonctions cuspidales en deux places.

On fixe $f_{V}$ et $\nu$ et on a d\'ej\`a vu qu'en fixant convenablement un multiplicateur $\alpha$, en rempla\c{c}ant $f_{V}$ par $f_{\alpha_{m},V}$ avec $m\in {\mathbb N}$ o\`u $\alpha_{m}$ est le convol\'e $m$ fois de $\alpha$, (2) a une limite qui est la somme sur les caract\`eres $c^V$ de 
$$\int_{i{\mathcal{A}}^*_{\tilde{G}}}d\lambda\,\left( tr\, \pi_{\nu,disc,\lambda}^{\tilde{G}}(c^V)(f_{V}1_{\tilde{K}^V})-tr\, \pi_{\nu,disc,\lambda}^{\tilde{G},{\mathcal E}}(c^V)(f_{V}1_{\tilde{K}^V})\right)\eqno(2)$$
o\`u l'indice $\lambda$ indique la tensorisation par le caract\`ere $\lambda$ de $\mathfrak{A}_{\tilde{G}}$.

Quant \`a  $I_{disc}^{\tilde{M}}(\omega,h^{f_{\alpha^m}}_{V}1_{\tilde{K}^{\tilde{M},V}})$ c'est la somme absolument convergente  des coefficients de Fourier $I^{\tilde{M}}(\pi_{disc,\nu_{{M}}}^{\tilde{M}}(c^{V}_{{M}}), 0,h^{f_{\alpha^m}}_{V})$ calcul\'es en \ref{miseenplace}. On a vu en \ref{miseenplace} qu'un tel terme valait 
$$
\int_{i{\mathcal{A}}^*_{\tilde{M}}}d\lambda\, \hat{\alpha}^m(\lambda+\nu_{M})I(\pi_{disc,\nu_{M}}^{\tilde{M}}(c^V),\lambda,h^f_{V}).\eqno(3)
$$
 On sait, d'apr\`es \cite{stabilisationIX} 8.9,  que la fonction $h^f_{V}$ est $K\cap M$ finie. Ainsi dans le calcul de  l'expression (3),  gr\^ace aux propri\'et\'es de semi-finitude de $\pi_{disc}^{\tilde{M}}$, seuls  interviennent les caract\`eres infinit\'esimaux de partie r\'eelle born\'ee, la borne d\'ependant de  $K\cap M$. Ainsi, pour un bon choix de $\alpha$ (qui ne d\'epend que des $K$-types d\'etermin\'es par $f$) la limite quand $m$ tend vers l'infini de $\hat{\alpha}^m(\lambda+\nu_{M})$  est $0$ presque partout.  On en d\'eduit que la limite de (3) quand $m$ tend vers l'infini est nulle gr\^ace \`a la propri\'et\'e de convergence absolue d\'emontr\'ee dans la proprosition ci-dessous (cf. \ref{convergenceabsolue}).

 On trouve que la limite quand $m$ tend vers l'infini de $I_{disc}^{\tilde{M}}(\omega,h^f_{\alpha_{m},V}1_{\tilde{K}^{\tilde{M},V}})$ est nulle. Ainsi, pour $\nu$ fix\'e,
la somme sur $c^V$ des termes \'ecrits en (2) est nulle pour toute fonction $f_{V}$. Par inversion de Fourier sur $i{\mathcal{A}}^*_{\tilde{G}}$ qui est dual de $\mathfrak{A}_{\tilde{G}}$, on en d\'eduit:
$$
\sum_{c^V} tr\,\bigl(\pi^{\tilde{G}}_{\nu,disc}(c^V)-\pi_{\nu}^{\tilde{G},{\mathcal E}}(c^V)\bigr)(f_{V}1_{\tilde{K}^V})=0.
$$
On en d\'eduit l'\'egalit\'e cherch\'ee pour $c^V$ fix\'e en utilisant \ref{multiplicateurs} (ii).

\subsection{Propri\'et\'e de convergence absolue pour la formule des traces \label{convergenceabsolue}}
On veut d\'emontrer que (3) de \ref{casparticulier} (ou une variante) converge absolument. On veut donc montrer 
\begin{prop} $\sum_{\pi^{\tilde{M}}} \int_{i{\mathcal A}_{\tilde{M}}^*}d\lambda\,  \vert I(\lambda,\pi^{\tilde{M}},h_{V}^f)\vert $ converge o\`u la somme porte sur les repr\'esentations discr\`etes (au sens de la formule des traces) de $\tilde{M}({\mathbb A}_{F})$, non ramifi\'ees hors de $V$ et ayant des $K_{V}$ types dans un ensemble fix\'e a priori par le choix de $f$. 
\end{prop}
Comme on se ram\`ene aux r\'esultats de M\"uller sur le spectre discret de $M({\mathbb A}_{F})$, on modifie l'\'enonc\'e. Pour toute repr\'esentation automorphe irr\'eductible unitaire, $\pi_{M}$, de $M({\mathbb A}_{F})$ v\'erifiant:
$$
\forall \gamma\in \tilde{G}({\mathbb A}_{F}),\,\, ^\gamma\pi_{M}\simeq \pi_{M} \otimes \omega
$$
on fixe une $\omega$-repr\'esentation de $\tilde{M}({\mathbb A}_{F})$ unitaire. On note ${\mathcal A}^{2}[\pi]$ la somme des sous-espaces des formes automorphes de carr\'e int\'egrable sur $M({\mathbb A}_{F})$ isomorphes \`a $\pi$; on note $m(\pi):=dim \, Hom_{M({\mathbb A}_{F})}(\pi, {\mathcal A}^2[\pi]$). Cet espace ${\mathcal A}^{2}[\pi]$ est une $\omega$-repr\'esentation de $\tilde{M}({\mathbb A}_{F})$ (cf. \ref{rappeldestermesdiscrets}) et on v\'erifie que pour tout $f\in I(\tilde{M}(A))$
$$
\vert tr _{{\mathcal A}^2[\pi]}(f)\vert= \vert z(\tilde{\pi})\vert \vert tr\, \tilde{\pi}(f)\vert,\eqno(1)
$$avec $z(\tilde{\pi})$ un nombre complexe v\'erifiant $\vert z(\tilde{\pi})\vert \leq m(\pi)$: le $z(\tilde{\pi})$ est la trace d'un op\'erateur unitaire d\'efini par un \'el\'ement de $\tilde{M}(F)$ dans l'espace de dimension $m(\pi)$, $Hom_{M({\mathbb A}_{F})}(\pi,{\mathcal A}^2(M(F)\backslash M({\mathbb A}_{F})))$ et la valeur absolue de ce nombre ne d\'epend pas du choix de $\tilde{\pi}$. 
On rappelle le r\'esultat cl\'e de \cite{muller}: il existe un entier $k$ tel que
$$
\sum_{\pi} m(\pi) (1+\vert\nu(\pi_{\infty})\vert)^{-k} \leq \infty\eqno(2)
$$
converge, o\`u la somme porte sur les repr\'esentations $\pi$ incluses dans le spectre discret de $M({\mathbb A}_{F})$, non ramifi\'ees hors de $V$ et ayant des $K_{V}$-types dans un ensemble fix\'e a priori. Ce r\'esultat s'\'etend au spectre $t$-discret en suivant toujours  \cite{muller}, comme  on va l'expliquer ici. Le spectre $t$-discret de $\tilde{M}({\mathbb A}_{F})$ fait intervenir, plus g\'en\'eralement (cf. \ref{rappeldestermesdiscrets}), les repr\'esentations dans $${\mathcal A}^{2}(U_{P}({\mathbb A}_{F})L(F) {\mathfrak A}_{L}\backslash M({\mathbb A}_{F}))$$ o\`u $P$ est un sous-groupe parabolique de $M$ de radical unipotent $U_{P}$ et de sous-groupe de Levi $L$. Comme on ne regarde que les repr\'esentations non ramifi\'ees hors de $V$, cela fixe dans ces induites uniquement la composante hors de $V$ et l'induite se d\'ecompose donc en au plus $N_{V}$-repr\'esentations irr\'eductibles o\`u $N_{V}$ est le produit du cardinal de $V$ avec $$sup_{v\in V} \vert Norm_{M(F_{v})}(L(F_{v}))/L(F_{v})\vert.$$ Ainsi le r\'esultat (2) s'\'etend en sommant sur les repr\'esentations $t$-discr\`etes de $M({\mathbb A}_{F})$, non ramifi\'ees hors de $V$ et ayant des $K_{V}$-types dans un ensemble fix\'e a priori. 
Il n'y a \'evidemment aucune difficult\'e \`a v\'erifier que (1) s'\'etend aussi pour toute repr\'esentation $\pi$ de $M({\mathbb A}_{F})$ donnant lieu \`a une repr\'esentation $t$-discr\`ete de $\tilde{M}({\mathbb A}_{F})$. On a donc
$$
\sum_{\pi} \vert z(\tilde{\pi}) (1+\vert \nu(\pi_{\infty})\vert)^{-k} \leq \infty, \eqno(3)
$$
pour $k$ convenable et o\`u la somme porte sur les repr\'esentations $\pi$ de $M({\mathbb A}_{F})$ intervenant dans la construction du spectre $t$-discret de $\tilde{M}({\mathbb A}_{F})$ comme expliqu\'e en \ref{rappeldestermesdiscrets}. En tenant compte de (3), pour d\'emontrer la proposition, il suffit de prouver que pour tout entier $N$, il existe une constante positive $c_{N}$ tel que pour tout $\lambda\in i{\mathcal A}^*_{\tilde{M}}$, on ait, pour tout $\pi$ comme ci-dessus:
$$
\vert I(\tilde{\pi},\lambda,h_{V}^f)\vert \leq c_{N}(1+\vert \nu(\pi_{\infty})\vert+\vert \lambda\vert)^{-N}. \eqno(4)
$$
Comme $h_{V}^f$ est une somme sur les places de $v$, il suffit de montrer cette assertion pour tout place $v$ de $V$ en rempla\c{c}ant $h_{V}^f$ par $\epsilon_{\tilde{M}}(f_{v})f_{\tilde{M}}^v$. De plus avec les formules de descente, il suffit de prouver que les in\'egalit\'es suivantes:

si $v$ est une place archim\'edienne:
$$
\vert I(\tilde{\pi}_{v},\lambda,f_{\tilde{M},v})\vert \leq c_{N}(1+\vert \nu(\pi_{\infty})\vert+\vert \lambda\vert)^{-N}
\eqno(5)
$$

si $v$ est une place $p$-adique
$$
\vert I(\tilde{\pi}_{v},\lambda,f_{\tilde{M},v})\vert\leq c_{N}\eqno(6)
$$
donc une in\'egalit\'e ind\'ependante de $\lambda$

si $v$ est une place archim\'edienne
$$
\vert I(\tilde{\pi}_{v},\lambda,\epsilon_{\tilde{M}}(f_{v}))\vert \leq c_{N}(1+\vert \nu(\pi_{\infty})\vert+\vert \lambda\vert)^{-N}
\eqno(7)
$$

si $v$ est une place $p$-adique
$$
\vert I(\tilde{\pi}_{v},\lambda,\epsilon_{\tilde{M}}(f_{v}))\vert\leq c_{N}
\eqno(8)
$$

Pour les assertions (5) et (6), on commence par remarquer que puisque le terme constant $f_{\tilde{M},v}$ est \`a support compact, $I(\lambda,\tilde{\pi}_{v},f_{\tilde{M},v})$ n'est autre que la trace de $\tilde{\pi}_{\lambda,v}$ sur $f_{\tilde{M},v}$. 
L'assertion (5) est alors d\'emontr\'ee en 2.1 de \cite{clozeldelorme}; l'exponentiel en la partie r\'eelle de $\nu_{\infty}(\pi)$ est born\'ee ici puisque l'on a fix\'e les $K_{\infty}$-types possibles.

L'assertion (6); la trace est une somme de coefficients matriciels o\`u on ne consid\`ere que les $K_{v}$-types sous-lesquels se transforment $f_{\tilde{M},v}$. Or ces $K_{v}$-types sont fix\'es et leur multiplicit\'e dans $\tilde{\pi}_{\lambda,v}$ est born\'ee ind\'ependamment de cette repr\'esentation. Donc le terme de gauche de (6) est born\'e par une constante d\'ependant uniquement des $K_{v}$ type fois $\int_{\tilde{M}(F_{v})}\vert f_{\tilde{M},v}(\gamma)\vert \, d\gamma$ qui est fini puisque $f_{\tilde{M},v}$ est \`a support compact. D'o\`u (6).

Montrons (8). On va d'abord  montrer qu'il n'y a qu'un nombre fini de repr\'esentations $\tilde{\pi}_{v}$ pour lesquelles $I(\lambda,\tilde{\pi}_{v},\epsilon_{\tilde{M}}(f_{v}))$ est non nul.
Les repr\'esentations globales qui interviennent ont un caract\`ere central trivial sur ${\mathfrak A}_{\tilde{M}}$. Elles sont non ramifi\'ees hors de $V$ et leur caract\`ere sous les points ad\'eliques du tore $A_{\tilde{M}}({\mathbb A}_{F})$ ne parcourt donc qu'un ensemble fini de caract\`eres automorphes unitaires. De plus la fonction $\epsilon_{\tilde{M},v}(f_{v})$ n'intervient que si ${\mathcal A}^*_{\tilde{M}}={\mathcal A}^*_{\tilde{M},v}$ ce que l'on suppose. Donc l'action de $\tilde{\pi}_{v}$ restreinte \`a $A_{\tilde{M},v}(F_{v})$ se fait via un caract\`ere unitaire dans un ensemble fini. On utilise le fait que $\epsilon_{\tilde{M},v}(f_{v})$ est une fonction cuspidale donc $I(\tilde{\pi}_{v},\lambda,\epsilon_{\tilde{M}}(f_{v}))=I(\rho_{v},\lambda,\epsilon_{\tilde{M}}(f_{v}))$ si $\tilde{\pi}_{v}-\rho_{v}$ est une induite propre. On peut trouver $\rho_{v}$ ayant cette propri\'et\'e, \'etant elliptique,  avec m\^eme caract\`ere que $\pi_{v}$ sur $A_{\tilde{M},v}(F_{v})$ et pour que $I(\rho_{v},\lambda,\epsilon_{\tilde{M}}(f_{v}))$ soit non nul, $\rho_{v}$ doit contenir certains $K_{v}$-types d\'etermin\'es uniquement par $f_{v}$. Il n'y a donc qu'un nombre fini de possibilit\'es; en particulier le support cuspidal de $\rho_{v}$ varie dans un ensemble fini et il en est donc de m\^eme de $\pi_{v}$, en tout cas pour que $I(\pi_{v},\lambda,\epsilon_{\tilde{M}}(f_{v}))$ soit non nul. D'o\`u l'assertion de finitude; cette assertion entra\^{\i}ne (8) car $\lambda$ parcourt un ensemble compact, la place \'etant finie, et la d\'ependance en $\lambda$ est continue.

Montrons (9). L'assertion de finitude est analogue, on remplace l'usage du support cuspidal par l'action du centre de l'alg\`ebre enveloppante. Ici $\lambda$ ne varie pas dans un tore compact. La fonction $\lambda\mapsto I(\tilde{\pi}_{v},\lambda,\epsilon_{\tilde{M}}(f_{v}))$ est \`a d\'ecroissance rapide, il existe donc, pour tout entier $N$,  une  constante $c_{N}$ ind\'ependante de $\tilde{\pi}_{v}$ tel que $$
\vert I(\tilde{\pi}_{v},\epsilon_{\tilde{M}}(f_{v}))\vert \leq c_{N} (1+\vert \lambda\vert)^{-N}.
$$
Comme 
$$
(1+\vert \lambda\vert)^{-N}\leq (1+\vert \nu_{\infty}(\pi_{\infty})\vert )^{N}(1+\vert \nu_{\infty}(\pi_{\infty})\vert +\vert \lambda\vert)^{-N},
$$
en modifiant la constante, on obtient (9).
\subsection{Globalisation\label{globalisationgrossiere}}
On note avec des indices $_{0}$ les termes  locaux que l'on veut globaliser. On a un corps local $F_{0}$, un espace tordu $\tilde{G}_{0}$ avec un caract\`ere $\omega_{0}$ et un groupe $G_{0}$, en fait c'est le cocycle $\underline{a}_{0}$ qui d\'etermine $\omega_{0}$ qui compte et on a aussi un espace de Levi $\tilde{M}_{0}$ de $\tilde{G}_{0}$. D'apr\`es \cite{kottwitzrogawski} paragraphe 2, il existe un corps de nombres $F$, des donn\'ees $G,\tilde{G},\tilde{M}$ d\'efinies sur $f$ et une place $v_{0}$ de $F$ qui se localisent en les donn\'ees locales pr\'ec\'edentes.  Et on peut en plus supposer comme nous le ferons que $\tilde{M}(F)$ est dense dans $\tilde{M}(F_{v_{0}})$ et que $A_{\tilde{M}}=A_{\tilde{M},v}$. 

On fixe ${\bf M}'_{0}$ une donn\'ee endoscopique elliptique relevante de $\tilde{M}_{0},\underline{a}_{0}$. Puisque $\tilde{M}(F)$ est dense dans $\tilde{M}(F_{v_{0}})$, l'ensemble des \'el\'ements $\delta_{0}$ dans ${\bf M}'_{0}$ qui appartiennent \`a une classe de conjugaison stable form\'ees d'\'el\'ements semi-simples r\'eguliers  correspondant \`a une classe de conjugaison stable d'\'el\'ements semi-simples r\'eguliers de $\tilde{M}(F)$ est dense dans ${\bf M}'_{0}$. On fixe $\delta_{0}$ dans cet ensemble et $\gamma\in \tilde{M}(F)$ lui correspondant. On note $\tilde{T}$ le tore tordu maximal de $\tilde{M}$ contenant $\gamma$. Le seul cas qui nous int\'eresse est le cas elliptique, donc on suppose que $\tilde{T}$ est elliptique dans $\tilde{M}$.

\begin{prop} Il existe un cocycle $\underline{a}$ de $W_{F}$, une donn\'ee endoscopique de $\tilde{M}, {\underline{a}}$, not\'ee ${\bf M}'$, et un \'el\'ement $\delta$ de cette donn\'ee, d\'efini sur $F$, tel que la localisation de ${\underline{a}}$ en la place $v_{0}$ soit le cocycle ${\underline{a}_{0}}$ et la localisation de ${\bf M}',\delta$ soit ${\bf M'_{0}},\delta_{0}$. 
\end{prop}

La d\'emonstration de cette proposition est report\'ee en \ref{globalisationdonneendo}.

Cette globalisation n'est pas suffisante pour s\'eparer les \'el\'ements de la classe de conjugaison stable de $\gamma$ qui interviennent dans les formules globales, uniquement par des consid\'erations locales. A ${\bf M}'$ est associ\'e un caract\`ere de $M_{\sharp}({\mathbb A}_{F})$; comme $M_{\sharp}=M/Z(M)^{\theta}$, il y a une application naturelle de $M_{\flat}:= M/Z(G)^{\theta}$ dans $M_{\sharp}$. On a besoin des propri\'et\'es suivantes,  pour toute place $v$ de $F$:

l'application naturelle  $M_{\flat}(F_{v})/M(F_{v}) \rightarrow G_{\sharp}(F_{v})/G(F_{v})$ est bijective et l'application naturelle $M_{\flat}(F_{v}) \rightarrow M_{\sharp}(F_{v})$ est surjective.

En effet, $G_{\sharp}(F_{v})$ est engendr\'e par un tore maximalement d\'eploy\'e et par $G_{SC}(F_{v})$, d'o\`u la surjectivit\'e de la premi\`ere application et son injectivit\'e est claire. Pour la deuxi\`eme application, on remarque que l'application de groupes alg\'ebriques,
$$
M_{\flat}= M/Z(G)^{\theta} \rightarrow M_{\sharp}=M/Z(M)^{\theta}$$
a pour noyau le tore $Z(M)^\theta/Z(G)^{\theta}=Z(M )^{\theta,0}/Z(G)^{\theta}\cap Z(M)^{\theta,0}$. Ce tore est induit et la surjectivit\'e annonc\'ee s'en d\'eduit.

On peut donc relever $\omega_{\sharp}$ en $\omega_{\flat}$ un caract\`ere de $M_{\flat}(F_{v})$ (pour toute place $v$ de $F$) et cette application de rel\`evement est injective. Et en utilisant le premier isomorphisme ci-dessus, on voit $\omega_{\sharp}$ comme un caract\`ere de $G_{\sharp}(F_{v})$. C'est le caract\`ere de $G_{\sharp}(F_{v})$ associ\'e \`a ${\bf M}'$ vu comme groupe endoscopique non elliptique de $\tilde{G}$. On note ce caract\`ere $\omega_{\sharp,v}^{\bf M'}$.

\subsection{Propri\'et\'es de finitude du nombre de certaines donn\'ees endoscopiques\label{finitudeendo}}
On se place dans une situation globale o\`u $\tilde{G}$ est d\'efini sur un corps de nombres $F$, que le caract\`ere $\omega$ est donn\'e par un cocycle fix\'e, $\underline{a}$ et o\`u on a un espace de Levi $\tilde{M}$ de $\tilde{G}$. On fixe aussi un nombre fini de places $V$ de $F$ contenant $V_{ram}$. Il n'y a qu'un nombre fini de donn\'ees endoscopiques elliptiques de $\tilde{M},\underline{a}$ non ramifi\'ees hors de $V$. Mais comme on va bouger $V$, il faut une condition de finitude plus ind\'ependante de $V$; ce que l'on contr\^ole mieux est le caract\`ere associ\'e \`a une donn\'ee endoscopique elliptique. On voudrait cette m\^eme propri\'et\'e de finitude en imposant seulement que le caract\`ere associ\'e \`a la donn\'ee soit non ramifi\'e hors de $V$. Ce n'est sans doute pas correct en g\'en\'eral mais cela le devient quand on fixe en plus un tore tordu maximal de $\tilde{M}$ et que l'on ne consid\`ere que les donn\'ees endoscopiques elliptiques ayant elles-m\^emes un tore se transf\'erant en le tore donn\'e et encore quelques hypoth\`eses. On suppose en plus qu'il existe une extension galoisienne $E$ de $F$ qui d\'eploie le tore $T$ et une place $u$ de $F$ tel que $E_{u}$ soit un corps. L'application de localisation de $Gal(E/F)$ dans $Gal(E_{u}/F_{u})$ est alors bijective.

On note ${\mathcal K}^V_{F}(T,\omega)$ l'ensemble des donn\'ees endoscopiques elliptiques de $\tilde{M}$ ayant ces propri\'et\'es.
\begin{lem} Avec les hypoth\`eses faites l'ensemble ${\mathcal K}^V_{F}(T,\omega)$ est fini.
\end{lem}
On montre d'abord l'assertion suivante qui servira aussi dans \ref{globalisation}; on fixe $\gamma\in \tilde{T}(F)$ et on consid\`ere les couples $({\bf M}',\delta)$ form\'es d'un \'el\'ement de ${\mathcal K}^V_{F}(T,\omega)$ et d'une classe de conjugaison stable $\delta$ dans cette donn\'ee qui se transf\`ere en la classe de conjugaison stable de $\gamma$.

Soit ${\bf M}'',\delta''$ satisfaisant \`a ces conditions. On commence par consid\'erer la situation apr\`es extension de $F$ \`a $E$; sur $E$ le tore $T$ se d\'eploie par hypoth\`ese et les donn\'ees endoscopiques relatives \`a des classes de conjugaison stable de $\tilde{T}(E)$ donnent lieu \`a des caract\`eres automorphes de $M_{\sharp}({\mathbb A}_{E})$ qui valent $\omega \circ Norm_{E/F}$ sur l'image de $M({\mathbb A}_{F})$; d'o\`u un homomorphisme de $W_{E}$ dans $Z(\hat{M}_{\sharp})$. Ce dernier groupe est d\'ecrit en \cite{stabilisationI} 2.7. On reprend la suite exacte qui suit (2) de cette r\'ef\'erence:
$$
1\rightarrow Z(\hat{M}_{SC})/Z(\hat{M}_{SC})^{\hat{\theta}} \stackrel{(\pi,1-\hat{\theta})}{\rightarrow}
Z(\hat{M})/(Z(\hat{M})\cap \hat{T}^{\hat{\theta},0})\times Z(\hat{M}_{SC})$$
$$ \rightarrow Z(\hat{M}_{\sharp})\rightarrow 1. \eqno(*)
$$
On pose $\hat{U}:= Z(\hat{M})\hat{T}^{\hat{\theta},0}/\hat{T}^{\hat{\theta},0}$. Le couple $({\bf M}'',\delta'')$ est associ\'e \`a un cocycle dans $$Z^{1,0}(W_{F}; \hat{T}/\hat{T}^{\hat{\theta},0} \stackrel{(1-\hat{\theta})}{\rightarrow} \hat{T}/Z(\hat{M})).$$
On note simplement $\delta''$ ce cocycle.
En restreignant \`a $W_{E}$, on obtient un morphisme de $W_{E}$ dans $\hat{U}$ (qui ne d\'epend pas des choix). On note $\chi_{\delta''}$ ce morphisme. Montrons:

(1) L'application  qui, \`a $({\bf M}'' ,\delta'')$ comme ci-dessus, associe le localis\'e en la place $u$, $({\bf M}''_{u},\delta''_{u})$ et le caract\`ere $\chi_{\delta''}$ est injective:

en effet soit pour $i=1,2$, des donn\'ees $({\bf M}'_{i} ,\delta'_{i})$ ayant m\^eme image. Alors $\delta'_{1}( \delta'_{2})^{-1}$ est un cocycle trivial sur $W_{E}$ puisque $\chi_{\delta_{1}}=\chi_{\delta_{2}}$. Il se factorise donc par $Gal (E/F)$ mais en la place $u$ les localis\'es sont les m\^emes. Comme on a suppos\'e que la localisation de $Gal(E/F)$ dans $Gal(E_{u}/F_{u})$ est bijective, cela veut dire que l'on a l'\'egalit\'e de cocycles $\delta'_{1}=\delta'_{2}$. D'o\`u l'injectivit\'e.

Il faut encore montrer qu'il n'y a qu'un nombre fini de caract\`eres $\chi_{\delta''}$ possibles. Et pour cela on montre d'abord  que $\chi_{\delta''}$ est non ramifi\'e hors des places de $E$ au-dessus de places de $F$ dans  $V$: soit $v$ une place de $F$ hors de  $V$ et $v'$ une place de $E$ au dessus. On sait que $\omega_{\sharp,v'}^{\bf M''}$ est non ramifi\'e puisque ${\bf M}''\in {\mathcal K}^V_{F}(T,\omega)$ (cf. 7.4 dont on a repris la notation). On v\'erifie que $\omega_{\sharp,v'}^{{\bf M}''}$ est l'image de $\chi_{\delta'',v'}$ dans (*): ces cocycles doivent \^etre cohomologues et comme la situation est d\'eploy\'ee sur $E$, ils co\"{\i}ncident. Soit $w\in I_{E_{v'}}$ le groupe d'inertie de $W_{E_{v'}}$. On vient de v\'erifier que $(\chi_{\delta'',v'}(w),1)$ est d'image triviale dans $Z(\hat{M}_{\sharp})$. D'apr\`es l'exactitude de (*) cela entra\^{\i}ne l'existence de $z\in Z(\hat{M}_{SC})/Z(\hat{M}_{SC}^{\hat{\theta}})$ tel que $1=(1-\hat{\theta})z$ et $\chi_{\delta'',v'}(w)=\pi(z)$ o\`u $\pi$ est l'application naturelle. Comme $(1-\hat{\theta})$ est injective sur $Z(\hat{M}_{SC})/Z(\hat{M}_{SC})^{\hat{\theta}}$ cela force $z=1$ et la non ramification cherch\'ee de $\chi_{\delta'',v'}$.

Pour avoir la finitude du nombre de caract\`eres, il ne reste plus qu'\`a remarquer que si l'on fixe un point de base ${\bf M}',\delta$, les cocycles $\chi_{\delta''}\chi_{\delta}^{-1}$ sont \`a valeurs dans $Z(\hat{M})^{\hat{\theta}}/Z(\hat{M})^{\hat{\theta}} \cap \hat{T}^{\hat{\theta},0}$ qui est un groupe fini: en effet $(1-\hat{\theta})\chi_{\delta''}$ est \`a valeurs dans $Z(\hat{M})$ et doit co\"{\i}ncider avec le cocyle $\underline{a}$ d\'eterminant $\omega\circ Norm_{E/F}$. D'o\`u l'assertion qui termine la preuve de la finitude.

\subsection{Globalisation fine \label{globalisation}}
On fixe ${\bf M}_{0}'$ une donn\'ee endoscopique elliptique en la place $v_{0}$. On globalise une premi\`ere fois \`a l'aide de \ref{globalisationgrossiere}, d'o\`u un corps de nombres $F$. On fixe aussi un tore tore tordu elliptique $\tilde{T}$ de $\tilde{M}$ qui en la place $v_{0}$ se localise en un tore tordu elliptique.  On fixe aussi une donn\'ee endoscopique elliptique ${\bf M}'$ de $\tilde{M}$ et un \'el\'ement $\delta'$ de cette donn\'ee dont on suppose qu'il appartient \`a un tore se transf\'erant en $\tilde{T}$.

On dit que $E$, une extension finie de $F,$ est une bonne extension pour la place $v_{0}$ si  les conditions suivantes sont satisfaites. D'une part, il existe une place $v'_{0}$ de $E$ au-dessus de $v_{0}$ telle que $E_{v'_{0}}=F_{v_{0}}$. D'autre part, il existe  une extension galoisienne $E'$ de $E$ qui  d\'eploie $\tilde{T}$ et deux places $u_{1}$ et $u_{2}$ de $E$  distinctes de $v'_{0}$ telles que $E'_{u_{i}}$ soit un corps pour $i=1,2$, donc  telles que la localisation de $Gal(E'/E)$ dans $Gal(E'_{u_{i}}/E_{u_{i}})$ soit un isomorphisme. Il existe de bonnes extensions pour la place $v_{0}$.

 On fixe de telles donn\'ees.  On note $V_{ram}^{E}$ l'analogue de $V_{ram}$ sur le corps $E$. Pour $V$ un ensemble de places de $E$ contenant $V^{E}_{ram}$ et tel que $\tilde{T}$ soit non ramifi\'e hors de $V$, on reprend la notation  ${\mathcal K}_{E}^V(T,\omega)$ de \ref{globalisationgrossiere}, ce sont les groupes endoscopiques elliptiques pour $\tilde{M}$ (vu sur $E$) dont le caract\`ere associ\'e est non ramifi\'e hors de $V$ et qui poss\`edent un tore se transf\'erant en $\tilde{T}$. Pour la suite, on aura besoin de deux places, pour le lemme qui suit il suffit d'une place $u$ ayant les propri\'et\'es de $u_{1}$ et $u_{2}$. On fixe donc une telle place.
 
\begin{lem}
Il existe un ensemble fini de places $V$ de $E$ (contenant $V^{E}_{ram}$) tel que pour toute donn\'ee   endoscopique elliptique ${\bf M}'' \in {\mathcal K}^V_{E}(T,\omega)$   contenant un \'el\'ement $\delta''$ dans un tore se transf\'erant en $\tilde{T}$ telles  que 
 le localis\'e de ${\bf M}'',\delta''$ co\"{\i}ncide avec le localis\'e de ${\bf M}',\delta'$ en la place $u$, on a:

soit $({\bf M}'', \delta'')$ $=$ $({\bf M}',\delta')$ 

soit le caract\`ere automorphe associ\'e \`a ${\bf M}''$ est diff\'erent du caract\`ere automorphe associ\'e \`a {\bf M}'.
\end{lem}
Comme annonc\'e on reprend l'assertion (1) de la preuve de \ref{finitudeendo} en rempla\c{c}ant $F$ et $E$ par $E$ et $E'$; l'hypoth\`ese faite sur $E$ et chacune des places $u_{i}$ pour $i=1,2$ assure que l'on peut l'utiliser. On fixe ${\bf M}'', \delta''$ satisfaisant \`a la condition du lemme. D'apr\`es l'assertion (1) de \ref{finitudeendo}, soit 
 $({\bf M}'', \delta'')$ $=$ $({\bf M}',\delta')$ soit, avec les notations introduites \`a cette place, $\chi_{\delta''}\neq \chi_{\delta'}$. On va montrer que cela force $$\omega_{\sharp}^{{\bf M}'}\circ Norm_{E'/E}\neq \omega_{\sharp}^{{\bf M}''}\circ Norm_{E'/E}.$$ En effet soit $v'$ une place de $E'$ tel que $\chi_{\delta'',v'}\neq \chi_{\delta',v'}$; on note $\omega_{\sharp,v'}^{{\bf M}'}$ la composante en la place $v'$ du caract\`ere $\omega_{\sharp}^{{\bf M}'}\circ Norm_{E'/E}$ et  on utilise une notation analogue avec ${\bf M}'$ remplac\'e par ${\bf M}''$. Et comme dans la preuve de \ref{finitudeendo}, on voit que $\omega_{\sharp,v'}^{{\bf M}'}(\omega_{\sharp,v'}^{{\bf M}''})^{-1}$ est l'image du cocycle $(\chi_{\delta',v'} \chi_{\delta'',v'}^{-1},1)$ dans la suite exacte (*) de cette preuve. Et on conclut comme en loc. cite.

\subsection{Preuve de la stabilisation g\'eom\'etrique locale\label{preuvegeoloc}}
\subsubsection{La premi\`ere hypoth\`ese cl\'e}
On d\'emontre ici la premi\`ere hypoth\`ese cl\'e faite en \ref{stabgeohyp}.

Dans \ref{stabgeohyp} la situation est locale; on la globalise gr\^ace  \`a \ref{globalisationgrossiere} et \ref{globalisation}

On a ainsi construit une situation globale o\`u en la place $v=v_{0}$ on retrouve la situation locale de \ref{stabgeohyp} et tel que ${M}(F)$ est dense dans $M(F_{v_{0}})$; comme $\tilde{M}(F)$ est non vide et que $M(F)$ est dense dans $M(F_{v_{0}})$, $\tilde{M}(F)$ est dense dans $\tilde{M}(F_{v_0})$. D'apr\`es la construction, on a aussi ${\mathcal A}_{M}={\mathcal A}_{M_{0}}$ et ${\mathcal A}_{\tilde{M}}={\mathcal A}_{\tilde{M}_{0}}$. On a construit la situation globale de fa\c{c}on \`a globaliser le caract\`ere de $M(F_{v_{0}})$ et les donn\'ees endoscopiques elliptiques ${\bf M}'_{0}$ en la place locale.

\begin{lem} Il existe une fonction lisse $\epsilon({\bf M}',\delta_{0})$ sur l'ensemble des classes de conjugaison stable d'\'el\'ements fortement r\'eguliers de ${\bf{M}}'(F_{v_{0}})$ telle que pour toute fonction $f_{v_{0}}\in I(\tilde{G}(F_{v_{0}}))$ nulle pr\`es des \'el\'ements exceptionnels s'il y en a  et si $v_{0}$ est $p$-adique:$$
\epsilon^{{\bf M}'}_{\tilde{M}}(f)(\delta_0)=\epsilon({\bf M}',\delta_{0}) f_{\tilde{M}}^{{\bf M}'}(\delta_{0}),
\eqno(1)
$$
\end{lem}
On commence par remarquer que si $\epsilon({\bf M}',\delta_{0})$ existe en tant que fonction elle est n\'ecesssairement lisse. Le probl\`eme est donc de d\'emontrer l'existence de cette fonction. On sait que le membre de gauche est l'int\'egrale orbitale en $\delta_{0}$ d'une fonction cuspidale. Il suffit donc  de  d\'emontrer l'existence de cette fonction pour les \'el\'ements $\delta_{0}$ qui sont elliptiques dans ${\bf{M}}'(F_{v_{0}})$.

On remarque aussi qu'il suffit de prouver l'existence d'une telle fonction en tout point de ${\bf{M}}'(F_{v_{0}})$ qui correspond \`a un \'el\'ement $\gamma\in \tilde{M}(F_{v_{0}})$ qui provient, par localisation, de $\tilde{M}(F)$, puisque l'on a suppos\'e que $\tilde{M}(F)$ est dense dans $\tilde{M}(F_{v_{0}})$. On fixe un tel $\gamma$, on note $\tilde{T}$ le tore tordu maximal le contenant et on peut supposer que $\tilde{T}$ est elliptique en $v_{0}$.

On reprend \`a peu pr\`es telle quelle la d\'emonstration de \cite{partie3} 7.3 avec les modifications n\'ecessaires puisque la globalisation est plus compliqu\'ee.

On remplace \'eventuellement $F$ par une extension finie $E$ de $F$, bonne pour la place $v_{0}$ comme en \ref{globalisation}, munie de ses places $v'_{0}$, $u_{1}$ et $u_{2}$. On fixe $V$ un ensemble de places de $E$ contenant $v'_{0}$, contenant $V^{E}_{ram}$ ainsi que $u_{i}$ pour $i=1,2$ et satisfaisant pour $u=u_{1}$ et $u=u_{2}$ au lemme de \ref{globalisation}. Ici on peut remplacer $F$ par $E$, $v_{0}$ par $v'_{0}$, $E$ dispara\^{\i}t et on a donc ${\mathcal K}^V_{F}(T,\omega)$ l'ensemble des donn\'ees endoscopiques elliptiques de $\tilde{M}$, dont le caract\`ere associ\'e est non ramifi\'e hors de $V$ et ayant un tore qui se transf\`ere en $\tilde{T}$. On sait que cet ensemble est fini. On suppose bien \'evidemment (quitte \`a \'elargir $V$) que  ${\bf M}'$ est dans cet ensemble et on a  $\delta\in {\bf M}'$ une globalisation de $\delta_{v_{0}}$. En utilisant les lemmes de \ref{finitudeendo} et \ref{globalisation}, on fixe $S$ un ensemble de places de $F$ hors de $V$ tel que pour tout couple ${\bf M}'',\delta''$ comme dans \ref{globalisation} (en particulier en la place $u_{1}$ ou (c'est bien ''ou'' et pas ''et'') la place $u_{2}$, le localis\'e de ce couple est le m\^eme que le localis\'e du couple ${\bf M}',\delta$, soit $({\bf M}'',\delta'')=({\bf M}',\delta)$, soit il existe  au moins deux places distinctes  dans $ S$ o\`u le caract\`ere associ\'e \`a ${\bf M}''$ diff\`ere du caract\`ere associ\'e \`a ${\bf M}'$.

On pose $V'=V\cup S$.

On construit maintenant des fonctions auxquelles on va pouvoir appliquer \ref{casparticulier}
On commence par fixer en toute place de $V'$ un voisinnage compact de $\gamma_{v}$. Ainsi  le produit de ces compacts avec hors de $V$ le compact $\tilde{K}^V$ coupe $\tilde{M}(F)$ en un ensemble fini de points (contenant $\gamma$), ensemble not\'e $C$. Toutes les fonctions consid\'er\'ees qui sont des \'el\'ements de $I(\tilde{G})$ sont telles que leur terme constant pour $\tilde{M}$ est \`a support dans le voisinage compact fix\'e. 

On affine les constructions de la fa\c{c}on suivante:
en les places $u_{i}$ pour $i=1,2$, on suppose que $f_{u_{i}}$ v\'erifie en plus, $f_{u_{i},\tilde{M}}$ est cuspidale et
$(f_{u_{i},\tilde{M}})^{{\bf M}''}=0$ si ${\bf M}''$ n'a pas la m\^eme localisation que ${\bf M}'$ en $u_{i}$ et $f^{{\bf M}'}_{\tilde{M},u_{i}}(\delta'_{i})=0$ pour tout $\delta'_{i}$ la composante locale d'un \'el\'ement de $C$ qui n'est pas dans la classe de conjugaison stable de  $\delta_{u_{i}}$ ni dans un de ses conjugu\'es sous $Norm_{G}(\tilde{M})$. 

En toutes les places de $V'$ hormis celles au dessus de $v_{0}$, on impose en plus que $f_{v,\tilde{M}}^{\bf{M}'}(\delta_{v})=1$.

Pour $v\in S$, on impose en plus que la fonction consid\'er\'ee soit de la forme $f_{v}^{\omega_{\sharp}^{\bf M'}}$, et en particulier   se transforme sous l'action de $M_{\flat}(F_{v})$ suivant la composante locale du caract\`ere associ\'e \`a ${\bf M}'$.
 
On applique la proposition de \ref{casparticulier} \`a $f_{V'}$, c'est-\`a-dire que l'on conna\^{\i}t pour une telle fonction la stabilisation spectrale. On revient alors au lemme de \ref{casparticulier} et  le membre de gauche de ce lemme    est nul; il en est donc de m\^eme du membre de droite de ce lemme c'est-\`a-dire :
$$
I_{geo}^{\tilde{M}}(\omega,h^f_{V'}1_{\tilde{K}^{V'}}) =0. \eqno(2)
$$
Le membre de gauche de (2) est une somme d'int\'egrales orbitales sur les \'el\'ements rationnels de $\tilde{M}(F)$. On l'\'ecrit en termes endoscopiques en consid\'erant s\'epar\'ement chaque terme composant $h^f_{V'}$ (qui, on le rappelle, est une somme sur $v\in V'$). A cause des conditions en $u_{i}$ pour $i=1$ et $2$, on est s\^ur que les couples ${\bf M}'', \delta''$ intervenant dans la somme se localisent en au moins une de ces places en le localis\'e de ${\bf M}',\delta$. 
Ainsi soit $({\bf M}'',\delta'')=({\bf M}',\delta)$ soit il existe une place dans $S$ o\`u le caract\`ere $\omega_{\sharp}^{\bf M''}$ diff\'ere du caract\`ere $\omega_{\sharp}^{\bf M'}$. 
Or pour contribuer  la donn\'ee endoscopique doit \'evidemment avoir son caract\`ere qui co\"{\i}ncide avec celui de ${\bf M}'$ en toute place de $S$ (sauf \'eventuellement une) puisque qu'en toute place de $S$, le terme constant de $f_{V',\tilde{M}}$ a cette propri\'et\'e de transformation.

Alors le terme de gauche de (2) n'est autre que $SI^{\bf{M}'}(\delta,h^f_{V'}\otimes 1_{\tilde{K}}^{{\bf M}',V'})$. 
Cela donne l'\'egalit\'e
$$0=
xSI^{\bf{M}'}(\delta_{v_{0}},\epsilon_{\tilde{M}}^{\bf{M}'}(f_{v_{0}})) \prod_{v\in V'; v\neq v_{0}} SI^{{\bf M}'}(\delta, f_{v,\tilde{M}}) + y SI^{{\bf M}'}(\delta_{v_{0}}, f_{v_{0},\tilde{M}}),
$$ o\`u $x$ est un scalaire non nul et $y$ un autre scalaire (\'eventuellement nul) qui est la contribution des places  $v\neq v_{0}$ faisant intervenir $\epsilon_{\tilde{M}}^{\bf M'}(f_{v})$
D'o\`u encore 
$$
SI^{\bf{M}'}(\delta_{v_{0}},\epsilon_{\tilde{M}}^{\bf{M}'}(f_{v_{0}}))=0\eqno(3)
$$
si $SI^{{\bf M}'}(\delta_{v_{0}}, f_{v_{0},\tilde{M}})=0$. Ainsi la premi\`ere forme lin\'eaire est proportionnelle \`a la seconde et   cela termine la preuve de la premi\`ere hypoth\`ese cl\'e.

\subsubsection{Preuve de la deuxi\`eme hypoth\`ese cl\'e}
La situation est encore globale. On dispose de $\tilde{G},\tilde{M},\underline{a}$. Pour toute donn\'ee endoscopique elliptique ${\bf M}'$ de $\tilde{M}$ et en toute place $v$, on a d\'efini la fonction sur les \'el\'ements tr\`es r\'eguliers de ${\bf M}'(F_{v})$, $\epsilon_{v}({\bf M}', \, )$. On sait que cette fonction est nulle sauf \'eventuellement sur les \'el\'ements elliptiques.
\begin{lem} Pour tout ensemble fini $V$ de places contenant $V_{ram}$, pour  toute donn\'ee endoscopique ${\bf M}'$ de $\tilde{M}$ non ramifi\'ee hors de $V$ et pour tout \'el\'ement $\delta\in {\bf M}'(F)$ non ramifi\'e hors de $V$ c'est-\`a-dire $\delta\in \tilde{K}_{v}^{M'}$ pour $v\not\in V$:
$$
\sum_{v\in V}\epsilon_{v}({\bf M}',\delta_{v})=0.
$$
\end{lem}
 On \'ecrit \ref{casparticulier}: on consid\`ere des fonctions $f_{V}=\prod_{v\in V}f_{v}\in I(\tilde{G}(F_{V}))$, qui sont $\tilde{M}$ cuspidales en au moins deux places et qui sont nulles pr\`es des \'el\'ements isol\'es pour tout $v\in V$ s'il y en a  et si $v$ est $p$-adique. On conna\^{\i}t donc la stabilisation spectrale pour ces fonctions et donc $I^{\tilde{M}}_{geo}(\omega,h_{V}^{f}1_{\tilde{K}^{M,V}})=0$ (avec les notations de loc.cite). Par d\'efinition $h_{V}^f$ est une somme sur les places $v_{0}\in V$ des fonctions $h_{v_{0}}:=\epsilon^{\tilde{G}}_{\tilde{M},v_{0}}(f_{v_{0}})\prod_{v''\in V-\{v_{0}\}}f_{\tilde{M},v''}$. On \'ecrit $I^{\tilde{M}}_{geo}$ sous forme endoscopique pour chacune de ces fonctions et on obtient une somme sur $v_{0}\in V$ des termes
 $$
 \sum_{{\bf M}'} i(\tilde{M},{\bf M}')SI_{geo}^{{\bf M}'}(h_{v_{0}}^{{\bf M}'}1_{\tilde{K}^{{\bf M}',V}}).
 $$
 On calcule $$h_{v_{0}}^{{\bf M}'}=(\epsilon_{\tilde{M},v_{0}}(f_{v_{0}}))^{{\bf M'}(F_{v_{0}})}\prod_{v\in V-\{v_{0}\}}(f_{\tilde{M},v})^{{\bf M}'(F_{v})}=\epsilon_{v_{0}}({\bf M}',\, )\prod_{v\in V}(f_{\tilde{M},v})^{{\bf M}'(F_{v})}.
 $$
 Donc finalement on obtient:
 $$
 \sum_{{\bf M}'} i(\tilde{M},{\bf M}')SI_{geo}^{{\bf M}'} ((\sum_{v\in V}\epsilon_{v}({\bf M}',\, ))\prod_{v\in V}(f_{\tilde{M},v})^{{\bf M}'(F_{v})}1_{\tilde{K}^{{\bf M}',V}})=0. \eqno(1)
 $$
Un  terme index\'e par une donn\'ee ${\bf M}'$ est une somme sur les classes de conjugaisons stables coupant les points rationnels de ${\bf M}'$, et quand $f_{V}$ est fix\'ee la somme est n\'ecessairement finie puisque la fonction \`a laquelle on l'applique est \`a support compact. On fixe en toute place de $V$ un sous ensemble compact et on s'autorise encore \`a restreindre le support des fonctions en supposant qu'il se trouve \`a l'int\'erieur de ce compact. On fixe ${\bf M}'_{0}$ et $\delta \in {\bf M}'_{0}(F)$. On note $\delta_{V}$ la composante dans $V$ de $\delta$; on consid\`ere en fait la classe de conjugaison stable de $\delta_{V}$ sous l'action du normalisateur de $\tilde{M}$ dans $G$. En jouant sur les fonctions $f_{V}$, on peut restreindre la somme aux donn\'ees endoscopiques elliptiques qui se localisent en toute place de $V$ en le localis\'e de ${\bf M}'_{0}$; on peut encore imposer  que la somme sur les points rationnels d'une telle donn\'ee ne porte que sur les classes de conjugaison stable qui en toute place de $V$ se localisent en la m\^eme classe que $\delta_{V}$. On obtient alors:
$$
\left(\sum_{v\in V}\epsilon_{v}({\bf M}',\delta_{v})\right) \sum_{{\bf M}',\delta'}1_{\tilde{K}^{{\bf M}',V}}((\delta')^{V})=0,
$$
o\`u la somme porte sur les couples ${\bf M}',\delta'$ d\'ecrits. Mais les fonctions $1_{\tilde{K}^{{\bf M}',V}}((\delta')^{V})$ valent $1$ ou $0$ et au moins l'une d'elles vaut $1$ d'apr\`es l'hypoth\`ese de non-ramification de $\delta$ hors de $V$. Et on obtient le lemme.
\section{Stabilisation de la formule des traces\label{thm}}

\subsection{Stabilisation spectrale \label{stabspecglob}}
\begin{thm} Pour toute fonction $f\in \tilde{G}(\mathbb A)$,
$$
I_{disc}(\omega,f)=\sum_{{\bf G}'}i(\tilde{G},{\bf G}')SI^{{\bf G}'}_{disc}(f^{{\bf G}'}).\eqno(1)
$$ o\`u ${\bf G}'$ parcourt l'ensemble des donn\'ees endoscopiques elliptiques de $\tilde{G}$. 
Ou encore, pour tout  caract\`ere infinit\'esimal $\nu$ de ${G}(F_{\infty})$ et pour tout ensemble fini de places $V$ contenant les places archim\'ediennes et les places o\`u $\tilde{G}$ est ramifi\'e (c'est-\`a-dire $V_{ram}$) et pour tout caract\`ere de l'alg\`ebre de Hecke sph\'erique $c^V$ hors de $V$, on a l'\'egalit\'e (dans le groupe de Grothendieck des $\omega$ repr\'esentations de  $\tilde{G}(\mathbb A)$)
$$
\pi_{disc,\nu}(c^V)=\sum_{{\bf G}'}i(\tilde{G},{\bf G}') \sum_{\nu'\mapsto \nu; c^{'V}\mapsto c^V}transfert(\pi^{{\bf G}'}_{disc,st,\nu'}(c^{'V})), \eqno(2)
$$
o\`u ici ${\bf G}'$ parcourt l'ensemble des donn\'ees endoscopiques elliptiques de $\tilde{G}$ non ramifi\'ees hors de $V$ et relevantes
\end{thm}
On fixe un ensemble fini $V'$ de places contenant $V_{ram}$. On suppose en plus que $V'-V_{ram}$ contient au moins une place et on en fixe une $v_{0}$. 
On a d\'emontr\'e la stabilisation locale des int\'egrales orbitales pond\'er\'ees. On consid\`ere  
$$I_{geo}^{\tilde{G}}(\omega,f_{V'}1_{\tilde{K}^{V'}})-\sum_{{\bf G}'}i(\tilde{G},{\bf G}')SI_{geo}^{{\bf G}'}(f_{V'}^{{\bf G}'}1_{K^{G',V'}}).\eqno(3)$$
Il ne reste que les termes correspondant aux \'el\'ements exceptionnels (cf. la preuve de \ref{casparticulier}). Donc (3) est nulle pour toute fonction $f$ qui en  $v_{0}$,  est nulle pr\`es des \'el\'ements exceptionnels (s'il y en a). Ainsi pour une telle fonction on a l'\'egalit\'e (1) du th\'eor\`eme. On a v\'erifi\'e en \ref{multiplicateurs} que cela entraine l'\'egalit\'e des traces en (2) \`a condition de se limiter aux fonctions que l'on vient de pr\'eciser. On pose $V=V'-\{v_{0}\}$. 

On fixe $f_{V}\in I(\tilde{G}(F_{V}))$, on note alors $f^{v_{0}}$ le produit de cette fonction avec $1_{\tilde{K}^{V'}}$. La forme lin\'eaire
$$f_{v_{0}}\in I(\tilde{G}(F_{v_{0}}))\mapsto$$
$$tr\, \pi_{\nu}(c^{V'})(f_{v_{0}}f^{v_{0}})-\sum_{{\bf G}'}i(\tilde{G},{\bf G}') \sum_{\nu'\mapsto \nu; c^{'V'}\mapsto c^{V'}}\pi^{{\bf G}'}_{st,\nu'}(c^{'V'})(f_{v_{0}}^{{\bf G}'}f^{v_{0},{\bf G}'})
$$
est nulle sur l'intersection des noyaux des formes lin\'eaires d\'efinies par les $\omega$-int\'egrales orbitales en les \'el\'ements isol\'es de $\tilde{G}(F_{v_{0}})$. On sait comme dans \ref{eliminationsingularite} que cette nullit\'e se propage \`a toutes les fonctions $f_{v_{0}}$ dont la composante elliptique est nulle; pour nous il suffit de l'avoir pour les fonctions $f_{v_{0}}$ dans l'alg\`ebre de Hecke sph\'erique.

Malheureusement, il n'est pas vrai en toute g\'en\'eralit\'e que  de telles fonctions ont  leur composante elliptique nulle. Comme ce bug appara\^{\i}t d\'ej\`a dans la stabilisation de la formule des traces non tordue, on explique ici comment faire fonctionner cette preuve dans le cas tordu et dans le cas non tordu. Dans le cas non tordu les \'el\'ements exceptionnels sont les  produits d'\'el\'ements centraux et d'\'el\'ements unipotents et il y en a dans tous les cas (mais ils ne posent pas de probl\`eme si $G$ est un tore) et dans le cas tordu, il n'y en a que si $\omega=1$ et $G$ v\'erifie les hypoth\`eses du (iii) de \ref{groupediese} (pas tout \`a fait, le groupe n'est d\'eploy\'e qu'\`a une restriciton des scalaires pr\`es, op\'eration insignifiante).

On reprend les notations de \ref{groupedieseglobal}, en particulier $\pi_{disc,\nu}^{V}$ et on suppose que $V$ est suffisamment grand pour que le (1) de loc.cite soit v\'erifi\'e. On veut montrer que $$
tr \pi_{disc,\nu}^{V}(f_{V}1_{\tilde{K}^V})=\sum_{{\bf G}'}i(\tilde{G},{\bf G}') tr \pi_{\nu,st}^{{\bf G}',V}(f_{V}^{{\bf G}'}1_{K^V_{{\bf G}'}}).\eqno(4)
$$
On reprend la notation ${\mathcal C}(V)$ de \ref{groupedieseglobal} et pour $\xi\in {\mathcal C}(V)$ on note $\xi(V)$ la restriction de $\xi$ \`a $G_{\sharp}(F_{V})$ et pour tout $v_{0}\notin V$ on note $\xi(v_{0})$ la restriction de $\xi$ \`a $G_{\sharp}(F_{v_{0}})$. On garde aussi les notations $f_{V}^{\xi(V)}$ de loc.cite. Le c\^ot\'e gauche de (4) est, d'apr\`es la proposition de \ref{groupedieseglobal}
$$
tr \pi_{disc,\nu}^{V}(f_{V}1_{\tilde{K}^V})=\sum_{\xi\in {\mathcal C}(V)}tr \pi_{disc,\nu}^{V}(f^{\xi(V)}_{V}1_{\tilde{K}^V}).\eqno(5)
$$
Fixons $\xi\in {\mathcal C}(V)$ et avec \ref{groupedieseglobal} fixons $v_{0}$ tel que $\xi(v_{0})=1$. On sait alors (cf. \ref{groupediesespherique}) que $1_{\tilde{K}_{v_{0}}}^{\xi(v_{0})}$ a sa composante elliptique nulle.  D'ap\`es la proposition 6.3, 
 on a:
$$
tr \pi_{disc,\nu}^{V}(f^{\xi(V)}_{V}1_{\tilde{K}^V})=tr \pi_{disc,\nu}^{V\cup\{v_{0}\}}(f^{\xi(V)}_{V}1_{\tilde{K}_{v_{0}}}^{\xi(v_{0})}1_{\tilde{K}^{V\cup\{v_{0}\}}})\eqno(6)$$
et on va pouvoir appliquer la stabilisation gr\^ace \`a la propri\'et\'e de $1_{\tilde{K}_{v_{0}}}^{\xi(v_{0})}$. Il faut d'abord remarquer que 
$$
\sum_{{\bf G}'}i(\tilde{G},{\bf G}') tr \pi_{\nu,st}^{{\bf G}',V}((f_{V}^{\xi(V)})^{{\bf G}'}1_{K^V_{{\bf G}'}})=$$
$$
\sum_{{\bf G}'}i(\tilde{G},{\bf G}') tr \pi_{\nu,st}^{{\bf G}',V\cup \{v_{0}\}}((f_{V}^{\xi(V)})^{{\bf G}'}(1^{\xi(v_{0})}_{\tilde{K}_{v_{0}}})^{{\bf G}'}1_{K^{V\cup\{v_{0}\}}_{{\bf G}'}}).\eqno(7)
$$
En effet le membre de gauche de (7) ne voit que les ${\bf G}'$ tel que le caract\`ere $\omega_{\sharp}^{{\bf G}'}$ soit dans ${\mathcal C}(V)$: le caract\`ere doit \^etre non ramifi\'e hors de $V$ car les donn\'ees endoscopiques le sont et il doit \^etre trivial sur l'image de $G({\mathbb A}_{F})$. De plus comme la fonction qui intervient est le transfert de $f_{V}^{\xi(V)}$, ce terme ne voit que les \'el\'ements de ${\mathcal C}(V)$ ayant pour  \`a $G_{\sharp}(F_{V})$ le caract\`ere $\xi(V)$ (cf. \ref{groupediese}). Ainsi pour de telles donn\'ees endoscopiques en la place $v_{0}$ on a aussi le caract\`ere $\xi(v_{0})$ et donc, avec \ref{groupediese} $(1^{\xi(v_{0})}_{\tilde{K}_{v_{0}}})^{{\bf G}'}=(1_{\tilde{K}_{v_{0}}})^{{\bf G}'}$ c'est-\`a-dire la fonction caract\'eristique de $K^{G'}_{v_{0}}$. D'o\`u l'\'egalit\'e.
On a l'\'egalit\'e du terme de droite de (6) avec le terme de droite de  (7) puisque la composante elliptique de $1_{\tilde{K}_{v_{0}}}^{\xi(v_{0})}$ est nulle. D'o\`u 
$$
tr \pi_{disc,\nu}^{V}(f^{\xi(V)}_{V}1_{\tilde{K}^V})=\sum_{{\bf G}'}i(\tilde{G},{\bf G}') tr \pi_{\nu,st}^{{\bf G}',V}((f_{V}^{\xi(V)})^{{\bf G}'}1_{K^V_{{\bf G}'}}).
$$
On a aussi
$$
\sum_{\xi\in {\mathcal C}(V)}\sum_{{\bf G}'}i(\tilde{G},{\bf G}') tr \pi_{\nu,st}^{{\bf G}',V}((f_{V}^{\xi(V)})^{{\bf G}'}1_{K^V_{{\bf G}'}})
=\sum_{{\bf G}'}i(\tilde{G},{\bf G}') tr \pi_{\nu,st}^{{\bf G}',V}(f_{V}^{{\bf G}'}1_{K^V_{{\bf G}'}}).\eqno(8)
$$
En effet, $f_{V}=\sum_{\xi'(V)}f_{V}^{\xi'(V)}$ o\`u cette fois les $\xi'(V)$ parcourt tous les caract\`eres de $G_{\sharp}(F_{V})/G(F_{V})$. Mais pour une donn\'ee endoscopique elliptique ${\bf G}'$ fix\'ee, $tr \pi_{\nu,st}^{{\bf G}',V}((f_{V}^{\xi'(V)})^{{\bf G}'}1_{K^V_{{\bf G}'}})=0$ si $\xi'(V)$ n'est pas la restriction de $\omega_{\sharp}^{{\bf G}'}$ \`a $G_{\sharp}(F_{V})$ dont si $\xi'(V)$ n'est pas de la forme $\xi(V)$ avec $\xi\in {\mathcal C}(V)$. D'o\`u (8). Puisque le terme de droite de (6) vaut le terme de gauche de (8), on obtient l'\'egalit\'e du terme de gauche de (6) avec le terme de droite de (8) ce qui est (4) et termine la preuve dans le cas tordu comme dans le cas non tordu.
\subsection{Une d\'ecomposition parfois plus fine de l'\'egalit\'e de stabilisation\label{stabilisationfine}}
On g\'en\'eralise un peu les constructions de \ref{groupediese} pour inclure le cas o\`u $\omega$ n'est pas trivial. On ne s'int\'eresse qu'au cas global. On fixe $\xi$ un caract\`ere automorphe de $G_{\sharp}({\mathbb A}_{F})$. On suppose que la restriction de $\xi$ \`a l'image de $G({\mathbb A}_{F})$ est $\omega$. On fixe $V$ un ensemble de places de $F$ (contenant $V_{ram}$) tel que $$
G_{\sharp}({\mathbb A}_{F})=G_{\sharp}(F)G_{\sharp}(F_{V})K_{\sharp}^V,
$$
o\`u
$K_{\sharp}^V$ est un compact maximal hyperspecial hors de $V$. On reprend la notation $\pi_{disc,\nu}^{\tilde{G}}(c^V)$ des paragraphes pr\'ec\'edent.
On peut d\'efinir la distribution sur $I(\tilde{G}(F_{V}))$,
$$
f_{V}\mapsto tr\, \pi_{\nu,disc}^{\tilde{G}}(c^V)(f_{V}1_{\tilde{K}^V}).
$$
Soit $g_{\sharp,V}\in G_{\sharp}(F_{V})$, cet \'el\'ement agit par conjugaison sur $I(\tilde{G}(F_{V}))$ et donc sur la distribution. On suppose en plus que $g_{\sharp,V}\in (G_{\sharp}(F)K_{\sharp}^V\cap G(F_{V}))$. Un tel \'el\'ement laisse invariante cette distribution. Et si $g_{\sharp,V}$ est dans l'image de $G(F_{V})$, l'action se fait via le caract\`ere $\omega^{-1}$. On suppose que le caract\`ere $\xi$ est non ramifi\'e hors de $V$. La distribution
$$f_{V}\mapsto 
\xi(g_{\sharp,V})tr\, \pi_{\nu,disc}^{\tilde{G}}(c^V)(^{g_{\sharp,V}}f_{V}1_{\tilde{K}^V})$$
ne d\'epend que de l'image de $g_{\sharp,V}$ dans le groupe fini $$H:=G_{\sharp}(F_{V})/G(F_{V})(G_{\sharp}(F)K_{\sharp}^V\cap G(F_{V})).$$ On la note $\xi(g_{\sharp,V})tr\, ^{g_{\sharp,V}}\pi_{\nu,disc}^{\tilde{G}}(c^V)$. Et on d\'efinit la distribution $$tr\, \pi_{disc,\nu}^{\tilde{G},\xi}(c^V):=\vert H\vert^{-1}\sum_{g_{\sharp,V}\in H}\xi(g_{\sharp,V})tr\, ^{g_{\sharp,V}}\pi_{\nu,disc}^{\tilde{G}}(c^V).
$$
\begin{thm} Avec $\xi$ et $V$ satisfaisant les hypoth\`eses pr\'ec\'edentes, c'est-\`a-dire que $V$ est suffisamment grand et $\xi$ est non ramifi\'e hors de $V$, on a
pour tout \'el\'ement $f_{V}\in I(\tilde{G}(F_{V}))$
$$
tr\, \pi_{disc,\nu}^{\tilde{G},\xi}(c^V)(f_{V}1_{K^V})=\sum_{{\bf G}';\omega_{\sharp}^{{\bf G}'}=\xi}i(\tilde{G},{\bf G}') tr\, \pi^{{\bf G}'}_{\nu,st}(c^V)(f_{V}^{{\bf G}'}1_{K_{G'}^V}),
$$o\`u l'on a regroup\'e dans le terme de droite la somme sur les $\nu'$ et les $c^{'V}$ de \ref{thm}.
\end{thm}
La preuve est \'el\'ementaire: le terme de gauche est la valeur de $tr\, \pi_{disc,\nu}^{\tilde{G}}(f^{\xi}_{V}1_{\tilde{K}^V})$ o\`u la fonction $f^\xi$ vaut, apr\`es un choix de repr\'esentant de $H$ dans $G_{\sharp}(F_{V})$
$$
f_{V}^{\xi}=\vert H\vert^{-1}\sum_{g_{\sharp,V}}\xi(g_{\sharp,V})\, ^{g_{\sharp,V}}f_{V}.
$$
On applique la stabilisation spectrale \`a cette fonction, on obtient un analogue du terme de droite en g\'en\'eralisant simplement le (i) de  \ref{groupediese} o\`u la somme porte sur les donn\'ees endoscopiques elliptiques non ramifi\'ees hors de $V$ et dont le caract\`ere associ\'e de $G_{\sharp}({\mathbb A}_{F})$ co\"{\i}ncide avec $\xi$ sur $G_{\sharp}(F_{V})$. Mais comme $V$ est suffisamment grand cela force l'\'egalit\'e de ce caract\`ere avec $\xi$.

\subsection{Un exemple, le cas de $GL(n)$ tordu \label{exemple}} 
On va donner l'exemple de $GL(n)$ tordu par l'automorphisme $g\mapsto \, ^tg^{-1}$. Dans ce cas, $\pi_{disc,\nu}^{\tilde{G}}(c^V)$ est une repr\'esentation irr\'eductible. Le groupe $G_{\sharp}$ vaut $GL(n)/\{\pm 1\}$. Soit $v$ une place de $F$ et  $z$ un \'el\'ement du centre de $GL(n,F_{v})$. On note $E$ l'extension au plus quadratique de $F_{v}$ tel que $z$ y admette une racine carr\'e, not\'ee$z^{1/2}$. L'image de $z^{1/2}$ dans $G_{\sharp}(E)$ est en fait dans $G_{\sharp}(F_{v})$. Et on v\'erifie ais\'ement que $G_{\sharp}(F_{v})$ est engendr\'e par l'image de $GL(n,F_{v})$ et par ces \'el\'ements. Pour un caract\`ere quadratique $\mu$ du centre de $GL(n,F_{v})$ et un caract\`ere $\xi$ de $G_{\sharp}(F_{v})$ trivial sur $GL(n,F_{v})$, on dit que $\xi$ et $\mu$ se correspondent  si et seulement si  $\xi(z^{1/2})=\mu(z)$ pour tout \'el\'ement $z$ du centre de $GL(n,F_{v})$.  Pour des caract\`eres  automorphes $\mu$ du centre de $GL(n,{\mathbb A}_{F})$ et $\xi$ de $G_{\sharp}({\mathbb A}_{F})$ v\'erifiant les m\^emes conditions, on dit qu'ils se correspondent si leurs composantes locales se correspondent en toute place. On note $\xi\simeq \mu$  cette correspondance.
 
\begin{prop} Avec les notations du paragraphe pr\'ec\'edent, la distribution $tr\, \pi_{disc,\nu}^{\tilde{G},\xi}(c^V)$ est nulle sauf si   $\xi$ correspond au caract\`ere central de  la repr\'esentation irr\'eductible  de $\pi_{disc,\nu}^{\tilde{G},\xi}(c^V)$. Et dans ce dernier cas, $tr\, \pi_{disc,\nu}^{\tilde{G},\xi}(c^V)=tr\, \pi_{disc,\nu}^{\tilde{G}}(c^V)$. 
\end{prop}
\begin{cor} On note  $\omega_{(c^V)}$ le caract\`ere central de $\pi_{disc,\nu}^{\tilde{G}}(c^V)$ et la repr\'esentation $
\pi_{disc,\nu}^{\tilde{G}}(c^V)$ est un transfert de $$\sum_{{\bf G}'; \omega_{\sharp}^{{\bf G}'}\simeq\omega_{(c^V)}}i(\tilde{G},{\bf G}')\pi_{\nu,st}^{{\bf G}'}(c^V). \eqno(1)
$$
\end{cor}
Cela n'est  qu'un exemple car \cite{book} donne beaucoup plus de pr\'ecisions. Toutefois, par exemple si $n$ est impair et si $\pi_{disc,\nu}^{\tilde{G}}(c^V)$ est de carr\'e int\'egrable (par exemple cuspidale) cela permet de d\'emontrer assez vite (uniquement avec le paragraphe 3 de \cite{book}) que la somme \`a droite n'a qu'un terme. Si $n$ est pair, sous les m\^emes hypoth\`eses pour $\pi_{disc,\nu}^{\tilde{G}}(c_{V})$, si le caract\`ere central de cette repr\'esentation est non trivial, ce caract\`ere central d\'etermine une extension quadratique de $F$ et dans la somme on ne trouve encore qu'une donn\'ee endoscopique elliptique, celle correspondant au groupe orthogonal pair quasid\'eploy\'e non d\'eploy\'e  qui se d\'eploie dans l'extension de $F$ \`a $E$. Si le caract\`ere central est trivial,
on ne peut pas distinguer par cette m\'ethode une donn\'ee endoscopique avec le groupe orthogonal impair $SO(n+1)$ de la donn\'ee endoscopique avec le groupe orthogonal pair d\'eploy\'e  et on n'\'evite pas l'utilisation des fonctions $L$ partielles et les derniers chapitres de \cite{book}.

\subsection{Une remarque sur la finitude de $\pi_{disc,\nu}(c^V)$ et son calcul pour les groupes classiques\label{finitude}}

Dans le cas o\`u $G=GL(n,F)$ et $\tilde{G}$ est l'espace tordu consid\'er\'e en \ref{exemple}, Arthur a d\'emontr\'e en \cite{book} en toute g\'en\'eralit\'e que  pour ${\bf G}'$ n'importe laquelle des donn\'ees endoscopiques elliptiques de $\tilde{G}$, les repr\'esentations  $\pi_{\nu',disc}^{G'}(c^{G',V})$ sont de longueur finie (ici l'indice $st$ est remplac\'e par $disc$). On peut renforcer cette propri\'et\'e de finitude en 
\begin{rmq} Soit ${\bf G}'$ une donn\'ee endoscopique elliptique de $\tilde{G}$. Alors pour tout $\nu, c^V$ comme ci-dessus, la repr\'esentation automorphe de $G'({\mathbb A}_{F})$, $$\oplus_{\nu'\mapsto \nu, c^{G',V}\mapsto c^V}\pi_{disc,\nu'}^{G'}(c^{G',V})$$ est de longueur finie. Ceci est m\^eme vrai pour $c^V$ fix\'e et en n'imposant rien au caract\`ere infinit\'esimal. C'est-\`a-dire avec des notations \'evidentes:
la repr\'esentation $\oplus_{c^{G',V}\mapsto c^{V}}\pi_{disc}^{G'}(c^{G',V})$ est de longueur finie. 
\end{rmq}
Cela r\'esulte \'evidemment des r\'esultats de \cite{book}. Arthur  montre entre autre que toute composante locale d'une des composantes irr\'eductibles d'une des repr\'esentations $\pi_{\nu',disc}^{G'}(c^{G',V})$ en une place $v\in V$ est n\'ecessairement dans un paquet de repr\'esentations qui ne d\'epend que de la composante locale en la place $v$ de $\pi_{disc,\nu}^{\tilde{G}}(c^V)$ (qui est une repr\'esentation irr\'eductible) et qui est caract\'eris\'e par des relations de transfert (cf. \cite{book}). En particulier il n'y a qu'un nombre fini de possibilit\'es. L'ensemble des repr\'esentations automorphes de carr\'e int\'egrable de $G'({\mathbb A}_{F})$ qui   dans $V$ sont l'une des repr\'esentations d\'etermin\'ees et qui sont non ramifi\'ees hors de $V$ sont en nombre fini (on rappelle que $V$ contient les places archim\'ediennes et que le caract\`ere infinit\'esimal est donc fix\'e dans un ensemble fini). D'o\`u la finitude; la toute petite difficult\'e qu'il a fallu contourner est qu'il y a en g\'en\'eral  un nombre infini de $c^{G',V}$ qui se transf\`erent en $c^{V}$ si $G'$ est un produit de groupes classiques l'un \'etant un groupe special orthogonal pair. Et on vient de voir que seul un nombre fini de $c^{G',V}$ intervient vraiment dans l'ensemble des formes automorphes de carr\'e int\'egrable, ceci est bien s\^ur dans \cite{book}.

\begin{rmq} La repr\'esentation   $\oplus_{c^{G',V}\mapsto c^{V}}\pi_{st}^{G'}(c^{G',V})$ est aussi de longueur finie.
\end{rmq}
C'est la m\^eme d\'emonstration.

\

Le m\^eme r\'esultat de finitude est vrai d'apr\`es \cite{mok} si $E/F$ est une extension quadratique, $G=GL(n,E)$ et $\tilde{G}$ est d\'efini par l'automorphisme $g\mapsto \, ^t\overline{g}^{-1}$ o\`u $\overline{g}$ est la conjugaison induite par l'extension $E$ de $F$. 

\

Soit $G$ une forme int\'erieure d'un groupe classique, c'est-\`a-dire que le groupe quasi-d\'eploy\'e correspondant est  un groupe sp\'ecial orthogonal, symplectique ou unitaire. Les r\'esultats d'Arthur seront tr\`es certainement g\'en\'eralis\'es \`a terme \`a un tel groupe mais on peut d\'ej\`a avoir une description qualitative des repr\'esentations $\pi^G_{\nu,disc}(c^V)$.

\begin{rmq} Les repr\'esentations $\pi_{disc}^{G}(c^V)$ sont de longueur finie
\end{rmq}
On applique (1) de \ref{stabspecglob} avec $\tilde{G}=G$. Il n'y a qu'un nombre fini de termes dans le membre de droite de (1): en effet il n'y a qu'un nombre fini de donn\'ees endoscopiques elliptiques non ramifi\'ees hors de $V$. Fixons l'un de ces groupes, ${\bf G}'$ et $c^{G',V}$ un syst\`eme de caract\`eres de Hecke qui par transfert non ramifi\'e s'envoie sur $c^V$ et tel que pour au moins un caract\`ere infinit\'esimal $\nu'$, $\pi_{disc,\nu'}^{G'}(c^{G',V})$ ne soit pas nul. Alors $G'$ est un produit de groupes classiques et il existe un entier $m=m_{1}+m_{2}$ (d\'ependant de $G'$) et une repr\'esentation automorphe irr\'eductible  de $GL(m_{1},{\mathbb A}_{F})\times GL(m_{2},{\mathbb A}_{F})$ qui est une induite de repr\'esentations automorphes de carr\'e int\'egrable et qui en toutes places hors de $V$ a pour caract\`eres de Hecke le transfert de $c^{G',V}$. Ainsi la repr\'esentation de $GL(m_{1}+m_{2},{\mathbb A}_{F})$ obtenue par induction (ou s\'erie d'Eisenstein) est une induite de repr\'esentations de carr\'e int\'egrable avec hors de $V$ un syst\`eme de caract\`ere de Hecke qui s'obtient directement par transfert non ramifi\'e \`a partir de $c^V$.  Ainsi il n'y a qu'un nombre fini de possibilit\'es pour cette repr\'esentation de $GL(m_{1}+m_{2},{\mathbb A}_{F})$ et donc aussi pour les repr\'esentations de $GL(m_{i},{\mathbb A}_{F})$ pour $i=1,2$. Avec la remarque pr\'ec\'edente le membre de droite de (1) est de longueur finie, son transfert est donc certainement de longueur finie.

\

\bf Question : \sl Soit $G$ un groupe r\'eductif d\'efini sur $F$, et $V, c^V$ comme pr\'ec\'edemment. Soit $\nu$ un caract\`ere infinit\'esimal. Il est naturel de se demander si $\pi_{disc,\nu}^{G}(c^V)$ est de longueur finie. Il est m\^eme vraisemblable que pour $c^V$ fix\'e, la somme $\oplus_{\nu} \pi^G_{disc,\nu}(c^V)$ est de longueur finie\rm

\subsection{V\'erification de toutes les hypoth\`eses de r\'ecurrence, r\'ecapitulatif\label{recapitulatif}}
Pour arriver \`a la stabilisation de la formule des traces, nous avons fait un certain nombre d'hypoth\`eses de r\'ecurrence et il faut donc les v\'erifier pour $\tilde{G}$.

On avait r\'ecapitul\'e les hypoth\`eses en \ref{miseenplace}; il y a des hypoth\`eses locales g\'eom\'etriques. Pour $\tilde{G}$ leurs v\'erification est faite en \ref{preuvegeoloc}. Disons tout de suite que nous n'avions pas d'hypoth\`eses locales spectrales car on les avait r\'esolues en \ref{stabgeostabspec}.

Et il y a les hypoth\`eses de r\'ecurrence globales; l'hypoth\`ese spectrale de r\'ecurrence est le th\'eor\`eme \ref{thm} que l'on vient de d\'emontrer. Il reste l'hypoth\`ese g\'eom\'etrique  c'est-\`a-dire que pour tout ${\mathcal O} $ une classe de conjugaison stable semi-simple de $\tilde{G}(F_{V})$ (notations de \cite{stabilisationVI} 5.3 et 5.4) $$A^{\tilde{G}}({\mathcal O},V)= A^{\tilde{G},{\mathcal E}}({\mathcal{O}},V).\eqno(1)
$$
Des r\'eductions ont \'et\'e faites en \cite{stabilisationVII} 3.5 ramenant la preuve de cette assertion aux \'el\'ements exceptionnels (cf. \ref{eliminationsingularite}), cette r\'eduction a servi pour la preuve de la stabilisation spectrale.

La distribution $f_{V}\mapsto I^{\tilde{G}}(A^{\tilde{G}}({\mathcal O},V),f_{V}1_{\tilde{K}^V})$ est une combinaison lin\'eaire d'int\'egrales orbitales d'apr\`es le (ii) de la proposition de 2.3 dans \cite{stabilisationVI} associ\'ees \`a des \'el\'ements $u\gamma$ (cf. loc. cite) o\`u $u$ est un \'el\'ement unipotent du centralisateur de $\gamma$ et o\`u $\gamma$ est un \'el\'ement de ${\mathcal O}$. Il en est de m\^eme de la distribution construite avec $A^{\tilde{G},{\mathcal E}}({\mathcal O},V)$ par la d\'efinition de \cite{stabilisationVI} 5.4 (cf. le d\'ebut de ce paragraphe). Rappelons que quand $\gamma,u$ varient avec les propri\'et\'es pr\'ec\'edentes, ces int\'egrales orbitales sont des distributions lin\'eairement ind\'ependantes.
Or on sait d\'ej\`a que $I^{\tilde{G}}(\omega,f_{V}1_{\tilde{K}^V})-I^{\tilde{G},{\mathcal E}}(f_{V}1_{\tilde{K}^V})=0$ pour tout $f_{V}\in I(\tilde{G}(F_{V}))$ puisque l'on a d\'emontr\'e la stabilisation spectrale (en \ref{redstabspecglob} et \ref{stabspecglob}). Avec les r\'eductions d\'ej\`a connues, ceci est la somme des distributions$$
I^{\tilde{G}}(A^{\tilde{G}}({\mathcal O},V),f_{V}1_{\tilde{K}^V})-I^{\tilde{G}}(A^{\tilde{G}, {\mathcal E}}({\mathcal O},V),f_{V}1_{\tilde{K}^V})
$$
quand ${\mathcal O}$ varie (cf. \cite{stabilisationVI} 6.11). Ainsi chacune de ces distributions doit \^etre nulle par ind\'ependance lin\'eaire. Cela termine la preuve.
\subsection{Stabilisation g\'eom\'etrique}
Le th\'eor\`eme de \cite{stabilisationVI} 5.9 est d\'emontr\'e: $$I_{geo}^{\tilde{G}}(\omega,f)=\sum_{{\bf G}'}i(\tilde{G},{\bf G}')SI_{geo}^{{\bf G}'}(f^{{\bf G}'}).$$
La d\'emonstration a \'et\'e donn\'ee en \cite{stabilisationVI} 6.11 sous les hypoth\`eses que l'on a v\'erifi\'ees en \ref{recapitulatif} mais bien \'evidemment c'est aussi un corollaire de \ref{stabspecglob}.

\subsection{Stabilisation de la formule des traces locales\label{stabloc}}
\subsubsection{La partie elliptique}
Ici $F$ est un corps local; pour simplifier les formules on suppose que $A_{\tilde{G}}=1$, sinon il faut int\'egrer sur ${\mathcal{A}}^*_{\tilde{G}}$. Pour toute donn\'ee endoscopique ${\bf G}'$ de $\tilde{G},\omega$, on fixe une base orthonormale des caract\`eres elliptiques stables de $G'$, not\'ee ${\mathcal B}_{st}^{G'}$. Pour $\tau$ une repr\'esentation elliptique, on note $i'(\tau)$ son produit scalaire elliptique, il est calcul\'e dans le (ii) du th\'eor\`eme de  \cite{ftlt} 7.3.
\begin{thm}Pour tout couple de fonctions $f_{1},f_{2}\in I(\tilde{G})$, on a l'\'egalit\'e
$$
\sum_{\tau \in Rep_{ell}(\tilde{G})} i'(\tau) \overline{tr\, \tau(f_{1})} tr\, \tau(f_{2})= \sum_{{\bf G}'} i(\tilde{G},{\bf G}') \sum_{\phi^{G'}\in {\mathcal B}_{st}^{G'}} \overline{tr\, \phi^{G'}(f^{{\bf G}'}_{1})} tr\, \phi^{G'}(f^{{\bf G}'}_{2}).
$$
\end{thm}
Ce th\'eor\`eme est un traduction spectrale de la proposition de \cite{stabilisationI} 4.17: en effet supposons d'abord que $f_{1}$ et $f_{2}$ soient des fonctions cuspidales. Alors le c\^ot\'e gauche s'interpr\`ete comme le produit scalaire de $f_{1}$ et $f_{2}$ pour le produit scalaire elliptique gr\^ace au fait que la norme d'une repr\'esentation elliptique $\tau$ est pr\'ecis\'ement $i'(\tau)$ par d\'efinition (cf. \cite{ftlt} 7.3). Cette \'egalit\'e est donc une \'egalit\'e de transfert entre repr\'esentations elliptiques qui est vraie pour les fonctions cuspidales. Elle est donc vraie pour toute fonction d'apr\`es le r\'esultat principal de \cite{prolongementdescaracteres} pour les corps p-adiques et \cite{stabilisationIV} 3.2 pour les corps archim\'ediens. 
\subsubsection{La partie discr\`ete non elliptique}
On rappelle que l'on a montr\'e en \ref{locspecfin} qui s'appuie sur les paragraphes pr\'ec\'edents, que l'on pouvait d\'ecomposer $I^{\tilde{G}}(\omega,f_{1},f_{2})$ et les analogues stables pour les groupes endoscopiques en une somme sur les espaces de Levi:
$$I^{\tilde{G}}(\omega,f_{1},f_{2})=\sum_{\tilde{M}}  I^{\tilde{G}}(\tilde{M},\omega,f_{1},f_{2}),
$$
o\`u la somme porte sur les classes de conjugaison d'espace de Levi de $\tilde{G}$ et o\`u 
$$
I^{\tilde{G}}(\tilde{M},\omega,f_{1},f_{2})=\sum_{\tau} i'(\tau) \overline{tr\, \tau(f_{1,\tilde{M}})} tr\, \tau(f_{2,\tilde{M}}),
$$
o\`u la somme porte sur les repr\'esentations elliptiques $\tau$ de $\tilde{M}$ dont l'induite est une repr\'esentation discr\`ete de $\tilde{G}$ modulo conjugaison sous le normalisateur de $\tilde{M}$ dans $G$ et o\`u $\iota'(\tau)$ est le coefficient intervenant dans loc.cite.

On a aussi montr\'e en \ref{redlocspec} que l'on avait aussi une d\'ecomposition pour la variante stable de cette distribution bien que l'\'ecriture soit moins jolie et surtout moins explicite. Notons $SI^{{\bf G}'}({\bf M}',f'_{1},f'_{2})$ la distribution correspondant.
\begin{prop}Pour tout $\tilde{M}$ espace de Levi de $\tilde{G}$ et pour tout couple de fonctions $f_{1},f_{2}\in I(\tilde{G})$,  on a l'\'egalit\'e:
$$
I^{\tilde{G}}(\tilde{M},\omega,f_{1},f_{2})=\sum_{{\bf M}'}i(\tilde{M},{\bf M}') \sum_{{\bf G}'}i_{\tilde{M}'}(\tilde{G},{\bf G}') SI^{{\bf G}'}({\bf M}',f_{1}^{\bf G'},f_{2}^{{\bf G}'}),
$$  o\`u la somme sur ${\bf M}'$ est la somme sur les donn\'ees endoscopiques elliptiques relevantes de $\tilde{M}$ et o\`u la somme sur ${\bf G}'$, quand ${\bf M}'$ est fix\'e est la somme sur les groupes endoscopiques elliptiques de $\tilde{G}$ contenant $M'$ comme sous groupe de Levi.
\end{prop}
Les distributions sont des combinaisons lin\'eaires de traces de repr\'esentations induites \`a partir de repr\'esentations elliptiques de $\tilde{M}$ pour le c\^ot\'e gauche et ${\bf M}'$ pour le c\^ot\'e droit. Elles s'expriment donc en fonction des termes constants  $f_{1,\tilde{M}}, f_{2,\tilde{M}}$ et de leurs transferts \`a ${\bf M}'$. Il suffit donc de montrer cette \'egalit\'e pour $f_{1}$ et $f_{2}$ des fonctions $\tilde{M}$-cuspidales. Puisqu'on a d\'emontr\'e la stabilisation des int\'egrales orbitales pond\'er\'ees, on sait stabiliser la formule des traces locale; on a donc l'\'egalit\'e cherch\'ee en sommant sur les espaces de Levi mais avec l'hypoth\`ese sur $f_{1}$ et $f_{2}$ il suffit de sommer sur les espaces de Levi contenant $\tilde{M}$ et on obtient la proposition avec une r\'ecurrence facile qui est initialis\'ee par le cas elliptique du paragraphe pr\'ec\'edent.
\section{Preuve de \ref{globalisationgrossiere}\label{globalisationdonneendo}}
Suivant Kottwitz-Rogawski, on peut choisir un corps de nombres $F$, des donn\'ees, $G,\tilde{G}, \tilde{M}$ sur $F$ et une place $v_{0}$ de $F$ de sorte que $F_{v_{0}}\simeq F_{0}$ et  que, modulo cet isomorphisme, les donn\'ees $G,\tilde{G},\tilde{M}$ localis\'ees en $v_{0}$ soient isomorphes \`a $G_{0},\tilde{G}_{0},\tilde{M}_{0}$. On peut de plus supposer que $\tilde{M}(F)$ soit dense dans $\tilde{M}(F_{{0}})$ et que $A_{\tilde{M}_{0}}$ soit le localis\'e de $A_{\tilde{M}}$. L'ensemble des $\delta_{0}\in {\bf M}'_{0}(F_{0})_{ell}$ qui correspondent \`a la classe de conjugaison stable d'un \'el\'ement de $\tilde{M}(F)$ est dense dans ${\bf M}'_{0}(F_{0})_{ell}^{rel}$ (on d\'esigne ainsi l'ensemble des \'el\'ements qui se transf\`erent \`a $\tilde{M}_{0}(F_{0})$). On note $D$ cet ensemble. On fixe $\delta_{0}\in D$ et $\gamma\in \tilde{M}(F)$ dont la classe de conjugaison stable correspond \`a $\delta_{0}$. On note $\tilde{T}$ le tore tordu maximal de $\tilde{M}$ contenant $\gamma$. Il est n\'ecessairement elliptique.

\

Fixons une extension galoisienne $E_{0}$ de $F_{0}$ finie que l'on pr\'ecisera plus loin telle que le tore localis\'e $T_{v_{0}}$ se d\'eploie sur $E_{0}$. On peut choisir une extension galoisienne finie $E'/F$, une place $v'_{0}$ de $E'$ au-dessus de $v_{0}$ de sorte que $T$ se d\'eploie sur $E'$ et $E'_{v'_{0}}$ contienne $E_{0}$ via l'isomorphisme $F_{0}\simeq F_{v_{0}}$. Le groupe $Gal(E'_{v'_{0}}/F_{0})$ s'identifie au fixateur de $v'_{0}$ dans $Gal(E'/F)$. On remplace $F$ par le corps des points fixes de ce fixateur et $v_{0}$ par la place $v'_{0}$ restreinte \`a ce sous-corps. Alors $Gal(E'_{v'_{0}}/F_{0})=Gal(E'/F)$. On note $E$ l'extension interm\'ediaire $F \subset E \subset E'$ telle que $Gal(E'/E)=Gal(E'_{v'_{0}}/E_{0})$. Alors $E/F$ est galoisienne, $Gal(E/F)=Gal(E_{0}/F_{0})$ et les actions galoisiennes globale et locale en $v_{0}$ co\"{\i}ncident sur $X_{*}(T)$. Donc $T$ est d\'eploy\'e sur $E$ ce qui entra\^{\i}ne que $G$ et $M$ le sont aussi. Puisque $\tilde{T}_{v_{0}}$ est elliptique dans $\tilde{M}_{v_{0}}$  et que le plus grand tore d\'eploy\'e $A_{\tilde{M}}$ central dans $\tilde{M}$ a pour localis\'e en $v_{0}$ le plus grand tore d\'eploy\'e $A_{\tilde{M}_{v_{0}}}$ dans $\tilde{M}_{v_{0}}$, $\tilde{T}$ est elliptique dans $\tilde{M}$.

\begin{rmq} On peut supposer que ces derni\`eres propri\'et\'es sont aussi v\'erifi\'ees pour deux autres places $u_{1},u_{2}$ de $F$.
\end{rmq} 
Il suffit pour cela de fixer une extension galoisienne finie $K/F$ telle que $v_{0}$ soit totalement d\'ecompos\'ee dans $K$ avec $[K:F]\geq 3$. On remplace $F$ par $K$ et $E$ par $EK$ qui est n\'ecessairement un corps. On remplace $v_{0}$ par une place au-dessus du $v_{0}$ initial  et on choisit pour $u_{1}, u_{2}$ deux autres places au-dessus du $v_{0}$ initial.

\

\bf Globalisation du caract\`ere $\omega_{0}$ et de l'\'el\'ement $\delta_{0}$. \rm

\

On commence par quelques rappels sur la classe de conjugaison stable de $\gamma$ dans $\tilde{M}_{0}(F_{0})$. Pour ces rappels on suppose que $F_{0}\neq {\mathbb R}$, le cas $F_{0}={\mathbb R}$ se traite de la m\^eme fa{\c c}on mais avec des $K$-espaces.

La classe de conjugaison stable de $\gamma$ est ${\mathcal C}^{st}(\gamma):=\{\gamma' \in \tilde{M}_{0}(F_{0});  \exists m\in M_{0}(\overline{F}_{0}); \gamma'=m^{-1}\gamma m\}$.  Si $\gamma' =m^{-1}\gamma m$ la condition d'\^etre dans $\tilde{M}_{0}(F_{0})$ se traduit exactement par le fait que pour tout $\sigma\in Gal(\overline{F}_{0}/F_{0})$, $m\sigma(m)^{-1}\in T^{\theta}$. On pose:
$$
{\mathcal Y}:=\{m \in M_{0}(\overline{F}_{0}); \forall \sigma\in Gal(\overline{F}_{0}/F_{0}),m\sigma(m)^{-1}\in T^{\theta}\}.
$$
On note $\pi$ la projection naturelle de $M_{0,SC}$ dans $M_{0}$. On a une bijection:
$$
T^{\theta}\backslash {\mathcal Y} / \pi(M_{0,SC}(F_{0})) \simeq H^{1,0}(\Gamma_{F_{0}}; T_{sc}\stackrel{1-\theta}{\rightarrow} (1-\theta)T),
$$
o\`u $T_{sc}$ est l'image r\'eciproque de $T$ dans $M_{0,SC}$.

Rappel de la construction: pour $m\in {\mathcal Y}$, on \'ecrit $m=z\pi(m_{sc})$ avec $z\in Z(M_{0})(\overline{F}_{0})$, $m_{sc}\in M_{0,SC}(\overline{F}_{0})$. Pour $\sigma\in \Gamma_{F_{0}}$, on pose $u(\sigma)=m_{sc}\sigma(m_{sc})^{-1}$. On voit que $u(\sigma) \in T_{sc}$ et que le couple $(u, (1-\theta)z)$ d\'efinit un \'el\'ement de $H^{1,0}(\Gamma_{F_{0}}; T_{sc} \stackrel{1-\theta}{\rightarrow} (1-\theta)T)$.

Notons ${\mathcal K}_{F_{0}}(T,\omega_{0})$ l'ensemble des caract\`eres de $H^{1,0}(\Gamma_{F_{0}}; T_{sc} \stackrel{1-\theta}{\rightarrow} (1-\theta)T)$ qui via la bijection ci-dessus, deviennent des fonctions sur $T^\theta \backslash {\mathcal Y}$ qui se transforment selon le caract\`ere $\omega_{0}$ par multiplication \`a droite par $M_{0}(F_{0})$.

D'apr\`es la th\'eorie g\'en\'erale, l'ensemble des couples $({\bf M}'_{1},\delta_{1})$ o\`u ${\bf M}'_{1}$ est une donn\'ee endoscopique elliptique de $(M_{0},\tilde{M}_{0},\underline{a}_{0})$ et $\delta_{1}$ est un \'el\'ement de ${\bf{M}'_{1}}(F_{0})$ correspondant \`a $\gamma$, couples pris \`a \'equivalence pr\`es, est naturellement en bijection avec ${\mathcal K}_{F_{0}}(T,\omega_{0})$. En particulier, le couple fix\'e $({\bf M}'_{0},\delta_{0})$  correspond \`a un \'el\'ement de ${\mathcal K}_{F_{0}}(T,\omega_{0})$, que l'on note simplement $\delta_{0}$.

On a d\'efini le groupe $M_{0,ab}(F_{0}):=\pi(M_{0,SC}(F_{0}))\backslash M_{0}(F_{0})\simeq H^{1,0}(\Gamma_{F_{0}}; T_{sc} \rightarrow T)$. On a le diagramme commutatif (cf \cite{stabilisationI}, 1.12):

$$
 \xymatrix {\relax &T_{sc} \ar[d]^{=} \ar[r]& T \ar[d]^{(1-\theta)}\\&
 T_{sc}\ar[r]^{(1-\theta)\,} &\, (1-\theta)T.}
 $$
 D'o\`u un homomorphisme
 $$
 H^{1,0}(\Gamma_{F_{0}};T_{sc}\rightarrow T) \rightarrow  H^{1,0}(\Gamma_{F_{0}}; T_{sc} \stackrel {1-\theta}{\rightarrow} (1-\theta)T). \eqno(*)
 $$
 On voit que l'action par multiplication \`a droite de $M_{0,ab}(F_{0})$ sur $T^\theta \backslash {\mathcal Y}/\pi(M_{0,SC}(F_{0}))$ correspond \`a la multiplication par  $H^{1,0}(\Gamma_{F_{0}};T_{sc}\rightarrow T)$ via l'homomorphisme ci-dessus. Donc ${\mathcal K}_{F_{0}}(T,\omega_{0})$ est l'ensemble des caract\`eres de $H^{1,0}(\Gamma_{F_{0}}; T_{sc} \stackrel {1-\theta}{\rightarrow} (1-\theta)T)$ qui pouss\'es par l'homomorphisme ci-dessus deviennent le caract\`ere $\omega_{0}$ de $ H^{1,0}(\Gamma_{F_{0}};T_{sc}\rightarrow T) $. Les groupes de caract\`eres des deux groupes de (*) sont les deux derniers groupes de la suite exacte:
 $$
 H^{1,0}(W_{F_{0}}; \hat{T}/\hat{T}^{\hat{\theta},0} \stackrel {1-\hat{\theta}}{\rightarrow}\hat{T})/Im\, \hat{T}^{\Gamma_{F_{0}},0} \rightarrow H^{1,0}(W_{F_{0}}; \hat{T}/\hat{T}^{\hat{\theta},0} \stackrel {1-\hat{\theta}}{\rightarrow} \hat{T}_{ad})/Im\, \hat{T}_{ad}^{\Gamma_{F_{0}},0}$$
 $$ \rightarrow
 H^{1,0}(W_{F_{0}},\hat{T}\rightarrow \hat{T}_{ad})\simeq H^1(W_{F_{0}}; Z(\hat{M})),
 $$
o\`u ici $\hat{T}_{ad}=\hat{T}/Z(\hat{M})$. Donc $\delta_{0}$ s'interpr\`ete comme un \'el\'ement du groupe central qui s'envoie sur l'\'el\'ement $\underline{a}_{0}$ du dernier groupe.

\

On s'est plac\'e dans l'espace ambiant $\tilde{M}_{0}$. On peut faire la m\^eme construction en rempla\c{c}ant $\tilde{M}_{0}$ par $\tilde{G}_{0}$. Pour pr\'eciser les notations avec $\tilde{M}_{0}$, on a l'ensemble ${\mathcal K}^{\tilde{M}_{0}}_{F_{0}}(T,\omega_{0})$ et avec $\tilde{G}_{0}$, l'ensemble ${\mathcal K}^{\tilde{G}_{0}}_{F_{0}}(T,\omega_{0})$. En fait ces ensembles sont identiques: ceci est bien connu mais nous allons le rev\'erifier. Du c\^ot\'e dual, ce qui change est le tore $\hat{T}_{ad}$ qui dans un cas est $\hat{T}/Z(\hat{M})$ et dans l'autre est $\hat{T}/Z(\hat{G})$. Pour v\'erifier l'\'egalit\'e annonc\'ee, on montre:
\begin{lem}
(i) L'homorphisme 
$$
H^{1,0}(W_{F_{0}}; \hat{T}/\hat{T}^{\hat{\theta},0} \stackrel {1-\hat{\theta}}{\rightarrow} \hat{T}/Z(\hat{G}))/Im\, (\hat{T}/Z(\hat{G}))^{\Gamma_{F_{0}},0} \rightarrow 
$$
$$H^{1,0}(W_{F_{0}}; \hat{T}/\hat{T}^{\hat{\theta},0} \stackrel {1-\hat{\theta}}{\rightarrow} \hat{T}/Z(\hat{M}))/Im\,( \hat{T}/Z(\hat{M}))^{\Gamma_{F_{0}},0} 
$$
est injectif.

(ii) Un \'el\'ement de $H^{1,0}(W_{F_{0}}; \hat{T}/\hat{T}^{\hat{\theta},0} \stackrel {1-\hat{\theta}}{\rightarrow} \hat{T}/Z(\hat{M}))/Im\,( \hat{T}/Z(\hat{M}))^{\Gamma_{F_{0}},0} $ qui s'envoie sur ${\underline{a}_{0}}$ est dans l'image de l'homorphisme de (i).
\end{lem}
Remarquons que $(\hat{T}/Z(\hat{G}))^{\Gamma_{F_{0}},0} $ et $(\hat{T}/Z(\hat{M}))^{\Gamma_{F_{0}},0} $ sont les images naturelles du m\^eme groupe $\hat{T}^{\Gamma_{F_{0}},0}$. Donc le noyau de l'homomorphisme \'ecrit dans l'\'enonc\'e de (i) est la projection du noyau de l'homomorphisme
$$
H^{1,0}(W_{F_{0}}; \hat{T}/\hat{T}^{\hat{\theta},0} \stackrel {1-\hat{\theta}}{\rightarrow} \hat{T}/Z(\hat{G})) \rightarrow 
H^{1,0}(W_{F_{0}}; \hat{T}/\hat{T}^{\hat{\theta},0} \stackrel {1-\hat{\theta}}{\rightarrow} \hat{T}/Z(\hat{M})).\eqno(1)
$$
On a la suite exacte de tores complexes:
$$
\xymatrix{\relax 
&1\ar[r]& Z(\hat{M})/Z({\hat{G}})\ar[r]&\hat{T}/Z(\hat{G}) \ar[r]&\hat{T}/Z(\hat{M})\ar[r]
&1
\\
&&1\ar[r] \ar[u]&\hat{T}/\hat{T}^{\hat{\theta},0} \ar[r]^{=} \ar[u]&\hat{T}/\hat{T}^{\hat{\theta},0}\ar[r] \ar[u]&1}.
$$
Le noyau de (1) est donc l'image de 
$$
(Z(\hat{M})/Z(\hat{G}))^{\Gamma_{F_{0}}}=H^0(W_{F_{0}}; Z(\hat{M})/Z(\hat{G}))\rightarrow $$
$$
H^{1,0}(W_{F_{0}}; \hat{T}/\hat{T}^{\hat{\theta},0}\stackrel{1-\hat{\theta}}{\rightarrow} \hat{T}/Z(\hat{G})).
$$Or $Z(\hat{M})/Z(\hat{G})$ est un tore induit ce qui entra\^{\i}ne que $(Z(\hat{M})/Z(\hat{G}))^{\Gamma_{F_{0}}}$ est connexe et son image tombe dans l'image de $(\hat{T}/Z(\hat{G}))^{\Gamma_{F_{0}},0}$. Cela prouve (i).

Consid\'erons un cocycle $(u,t)\in Z^{1,0}(W_{F_{0}}; \hat{T}/\hat{T}^{\hat{\theta},0} \stackrel{1-\hat{\theta}}{\rightarrow} \hat{T}/Z(\hat{M}))$ qui s'envoie sur $\underline{a}_{0}$ dans 
$H^{1,0}(W_{F_{0}}; \hat{T}/\hat{T}^{\hat{\theta},0}\stackrel{1-\hat{\theta}}{\rightarrow} \hat{T}/Z(\hat{M}))$.

Cette image est $((1-\hat{\theta})u,t)$. Repr\'esentons $\underline{a}_{0}$ par un cocycle \`a valeurs dans $Z(\hat{G})$ en revenant \`a $\underline{a}_{0} \in H^1(W_{F_{0}},Z(\hat{G}))$. Alors $((1-\hat{\theta})u,t)$ est cohomologue \`a $(\underline{a}_{0},1)$. Il existe donc $x\in \hat{T}$ tel que $(1-\hat{\theta})u(\sigma)=a_{0}(\sigma)x^{-1}\sigma(x)$ pour tout $\sigma\in W_{F_{0}}$ et $t=x$ modulo $Z(\hat{M})$. Le couple $(u,x_{ad})$, o\`u $x_{ad}$ est l'image de $x$ dans $\hat{T}/Z(\hat{G})$ est un cocycle dans $Z^{1,0}(W_{F_{0}}; \hat{T}/\hat{T}^{\hat{\theta},0} \stackrel{1-\hat{\theta}}{\rightarrow} \hat{T}/Z(\hat{G}))$. On voit que son image par l'homomorphisme du (i) est le cocycle de d\'epart $(u,t)$. Cela prouve (ii).

\

On peut donc identifier $\delta_{0}$ \`a un \'el\'ement de $$H^{1,0}(W_{F_{0}}; \hat{T}/\hat{T}^{\hat{\theta},0} \stackrel {1-\hat{\theta}}{\rightarrow} \hat{T}/Z(\hat{G}))/ Im\, (\hat{T}/Z(\hat{G}))^{\Gamma_{F_{0},0}} 
$$
qui s'envoie sur $\underline{a}_{0}$ dans $H^{1,0}(W_{F_{0}},\hat{T}\rightarrow \hat{T}/Z(\hat{G}))$. On le repr\'esente par un cocycle $(u_{0},t)$. Soit $K_{0}/F_{0}$ l'extension galoisienne qui d\'eploie $T$. La restriction de $u_{0}$ \`a $W_{K_{0}}$ est un homomorphisme  $W_{K_{0}}\rightarrow \hat{T}/\hat{T}^{\hat{\theta},0}$ qui v\'erifie $(1-\hat{\theta})u_{0}(\sigma)_{ad}=\sigma(t)t^{-1}=1$ pour tout $\sigma\in W_{K_{0}}$ (l'indice $ad$ d\'esigne la projection dans $\hat{T}_{ad}=\hat{T}/Z(\hat{G})$). Or l'ensemble des \'el\'ements $x\in \hat{T}/\hat{T}^{\hat{\theta},0}$ tels que $(1-\hat{\theta})(x)_{ad}=1$ est $\hat{U}:=Z(\hat{G})\hat{T}^{\hat{\theta},0}/\hat{T}^{\hat{\theta},0}$ qui n'est pas connexe en g\'en\'eral. Si $F_{0}={\mathbb R}$, on pose $E_{0}={\mathbb C}$ et comme $W_{E_{0}}={\mathbb C}^*$, la restriction de $u_{0}$ \`a $W_{E_{0}}$ est \`a valeurs dans $\hat{U}^0$. Si $F_{0}$ est $p$-adique, on interpr\`ete la restriction de $u_{0}$ \`a $W_{K_{0}}$ comme un caract\`ere $\chi_{K_{0}}$ de $K_{0}^*$ dans $\hat{U}$. Soit $C$ le noyau de la restriction de ce caract\`ere au groupe des unit\'es $O_{K_{0}}^*$. Parce que $u_{0}$ est un cocycle, $\chi_{K_{0}}(\sigma(x))=\sigma\chi_{K_{0}}(x)$ pour tout $\sigma \in Gal(K_{0}/F_{0})$ et tout $x\in K_{0}^*$. En particulier $\chi_{K_{0}}$ envoie $F_{0}^*$ dans $\hat{U}^{\Gamma_{F_{0}}}$. En notant $n:=[\hat{U}^{\Gamma_{F_{0}}}:\hat{U}^{\Gamma_{F_{0}},0}]$, $\chi_{K_{0}}$ envoie $F_{0}^{*,n}:=\{x^n; x\in F_{0}^*\}$ dans $\hat{U}^{\Gamma_{F_{0}},0}$. Le groupe $CF_{0}^{*,n}$ est ouvert d'indice fini dans $K_{0}^*$. Par le corps de classes, c'est le groupe des normes d'une extension ab\'elienne finie $E_{0}$ et parce que $u_{0}$ est un cocycle, on v\'erifie que l'extension $E_{0}$ de $F_{0}$ est elle aussi galoisienne. On a maintenant d\'efini le corps $E_{0}$ que l'on utilise dans la construction pr\'ec\'edente de globalisation.

\

On va construire l'\'el\'ement $\delta\in H^{1,0}(W_{F_{0}}; \hat{T}/\hat{T}^{\hat{\theta},0}\stackrel{1-\hat{\theta}}{\rightarrow} \hat{T}/Z(\hat{G}))$ qui par localisation s'envoie sur $\delta_{0}$. On dispose des quotients $W_{E/F}$ et $W_{E_{0}/F_{0}}$. On a fait ce qu'il fallait pour que $u_{0}$ se quotiente par $W_{E_{0}/F_{0}}$. Et parce que $Gal(E/F)=Gal(E_{0}/F_{0})$, on voit que $W_{E/F}$ s'identifie au produit semi-direct de ${\mathbb A}_{E}^*/E^* \rtimes W_{E_{0}/F_{0}}$ quotient\'e par $E_{0}^*$ plong\'e antidiagonalement. Pour prolonger $u_{0}$ en un cocycle $u:W_{E/F} \rightarrow \hat{T}/\hat{T}^{\hat{\theta},0}$, il suffit de prolonger le caract\`ere $\chi_{0}:= u_{0 \vert E_{0}^*}: E_{0}^* \rightarrow \hat{U}$ en un caract\`ere $\chi: {\mathbb A}_{E}^*/E^* \rightarrow \hat{U}$ qui soit \'equivariant pour $Gal(E/F)$. Si $F_{0}={\mathbb C}$, on a aussi $E_{0}=F_{0}$ donc $E=F$. La condition d'\'equivariance est triviale. D'autre part, on a vu que $\chi_{0}$ prend ses valeurs dans $\hat{U}^0$. Il est trivial que l'on peut prolonger $\chi_{0}$ en $\chi$. Si $F_{0}$ est p-adique, on s'est arrang\'e pour que $\chi_{0}$ soit non ramifi\'e et prenne ses valeurs dans $\hat{U}^{\Gamma_{F_{0}},0}$. Ecrivons $\hat{U}^{\Gamma_{F_{0}},0}=({\mathbb C}^*
)^N$. Alors $\chi_{0}$ est de la forme $(\vert\, \vert^{s_{1}}_{E_{0}}, \cdots, \vert\, \vert^{s_{N}}_{E_{0}})$ que l'on prolonge par $(\vert\, \vert^{s_{1}}_{{\mathbb A}_{E}}, \cdots, \vert\, \vert^{s_{N}}_{{\mathbb A}_{E}})$.

Reste le cas, $F_{0}={\mathbb R}$, $E_{0}={\mathbb C}$, $E/F$ quadratique. Le caract\`ere $\chi_{0}$ prend ses valeurs dans $\hat{U}^0$. On sait que l'on peut d\'ecomposer $\hat{U}^0$ muni de l'action de $Gal(E/F)=\{1,\sigma\}$ en produits de trois types de tores: ${\mathbb C}^*$ avec $\sigma(x)=x$, ${\mathbb C}^* \times {\mathbb C}^*$ l'action du groupe de Galois \'echangeant les deux copies et ${\mathbb C}^*$ avec $\sigma(x)=x^{-1}$. Il suffit de traiter chaque cas. On traite le dernier qui est le plus difficile.

La condition d'\'equivariance est que $\chi$ doit \^etre trivial sur les normes $N_{E/F}(E^*)$. On v\'erifie que $E^*N_{E/F}({\mathbb A}^*_{E})$ est un sous-groupe ferm\'e de ${\mathbb A}^*_{E}$, on munit ${\mathbb A}^*_{E}/E^*N_{E/F}({\mathbb A}^*_{E})$ de la topologie quotient. C'est un groupe compact. L'homomorphisme $E^*_{0}/N_{E_{0}/F_{0}}(E_{0}^*) \rightarrow {\mathbb A}^*_{E}/E^*N_{E/F}({\mathbb A}^*_{E})$ est continu. On va voir qu'il est injectif. Alors son image est un sous-groupe compact de ${\mathbb A}^*_{E}/E^*N_{E/F}({\mathbb A}^*_{E})$ et on peut consid\'erer que $\chi_{0}$ est d\'efini sur ce sous-groupe. Il est connu qu'un caract\`ere continu d'un sous-groupe compact d'un groupe compact \`a valeurs dans ${\mathbb C}^*$ se prolonge au groupe tout entier en un caract\`ere continu. Cela conclut modulo l'assertion d'injectivit\'e admise. Montrons-la: soit $e_{0}\in E_{0}^*$ et supposons que $e_{0}\in E^*N_{E/F}({\mathbb A}^*_{E})$. Ainsi il existe $e\in E^*$ tel que  $e\in e_{0}N_{E_{0}/F_{0}}(E_{0}^*)$ et  $e\in N_{E_{v}/F_{v}}(E_{v}^*)$ pour toute place $v\neq v_{0}$. Cela entra\^{\i}ne que $e\in F^*$ et, en notant $\kappa$ le caract\`ere quadratique associ\'e \`a $E/F$, on a $\kappa_{v}(e)=1$ pour tout $v\neq v_{0}$. La formule de produit entra\^{\i}ne que $\kappa_{v_{0}}(e)=1$ donc $e_{0}\in N_{E_{0}/F_{0}}(E_{0}^*)$ ce qu'il fallait d\'emontrer.

\

On a maintenant d\'efini $\delta \in H^{1,0}(W_{F_{0}}; \hat{T}/\hat{T}^{\hat{\theta},0}\stackrel{(1-\hat{\theta})}{\rightarrow} \hat{T}/Z(\hat{G}))$ qui par localisation s'envoie sur $\delta_{0}$. On note $\underline{a}$ l'image de $\delta$ dans $H^{1,0}(W_{F}; \hat{T} \rightarrow \hat{T}/Z(\hat{G}))=H^1(W_{F}; Z(\hat{G}))$, ou plut\^ot son image modulo $ker^1(W_{F}; Z(\hat{G}))$. Il est clair que $\underline{a}$ s'envoie sur $\underline{a}_{0}$ par localisation. Il d\'efinit un caract\`ere $\omega$ de $G({\mathbb A}_{F})$, automorphe, qui se restreint en $\omega_{0}$ en la place $v_{0}$.

\begin{rmq} Le caract\`ere $\omega$ est trivial sur $T^\theta({\mathbb A}_{F})$. 
\end{rmq}
En effet la restriction de $\omega $ \`a $T({\mathbb A}_{F})$ est associ\'ee \`a l'\'el\'ement de $H^1(W_{F}; \hat{T})$ qui est \'egal \`a $(1-\hat{\theta})u$, si $\delta$ est repr\'esent\'e par $(u,t)$. Donc cette restriction est compos\'ee de $T({\mathbb A}_{F}) \stackrel{1-\theta}{\rightarrow} (1-\theta)(T)({\mathbb A}_{F})$ et du caract\`ere de ce dernier groupe d\'etermin\'e par $u$. En tout cas, c'est trivial sur $T^\theta({\mathbb A}_{F})$.

colette.moeglin@imj-prg.fr

jean-loup.waldspurger@imj-prg.fr

\end{document}